%% file: chapter.tex
\documentclass{amsbook}

\usepackage{graphicx}

\newtheorem{theorem}{Theorem}[chapter]
\newtheorem{lemma}[theorem]{Lemma}
\newtheorem{corollary}[theorem]{Corollary}

\theoremstyle{definition}
\newtheorem{definition}[theorem]{Definition}

\theoremstyle{remark}
\newtheorem{remark}[theorem]{Remark}

\numberwithin{section}{chapter}
\numberwithin{theorem}{section}
\numberwithin{equation}{section}

\DeclareMathOperator{\Li}{Li}

\newcommand{\lrf}[1]{\left\lfloor #1\right\rfloor}

\begin{document}
\frontmatter
\title[On Prime Values of Bunyakovsky's Polynomials and Prime Products of Schinzel's Polynomials]{On Prime Values of Bunyakovsky's Polynomials and Prime Products of Schinzel's Polynomials}

\author{Shan-guang Tan}
%\address{Zhejiang University, Hangzhou, 310027, CHINA}
\curraddr{Zhejiang University, Hangzhou, 310027, CHINA}
\email{tanshanguang@163.com}

%\date{June 6, 2026}

\begin{abstract}
First, we prove a theorem of the number of primes in short intervals with a much more accurate estimation of prime numbers, the error range of which is less than $x^{1/2-0.0327283}$ for $x\geq10^{41}$. The theorem states as follows: Given a positive real number $\beta$ that determines a real number $x_{\beta}$ by $e(\log x_{\beta})^{3}/x_{\beta}^{0.0327283}=\beta$, let $\Phi(x):=\beta x^{1/2}$ for $x\geq x_{\beta}$ where $\Phi(x):=x^{1/2}$ when let $\beta=1$. Then there are
\[
\frac{\pi(x+\Phi(x))-\pi(x)}{\Phi(x)/\log x}=1+O(\frac{1}{\log x})
\]
and
\[
\lim_{x \to \infty}\frac{\pi(x+\Phi(x))-\pi(x)}{\Phi(x)/\log x}=1.
\]

Next, based on the theorem, Bunyakovsky's polynomials, Schinzel's polynomials, and Schinzel's polynomials with a negative leading coefficient are researched. Results show that

(1) given a Bunyakovsky's polynomial $f(x)$, $f(m)$ is prime for infinitely many positive integers $m$,

(2) given a Schinzel's polynomial $f(x)=\prod_{i=1}^{k}f_{i}(x)$, all numbers $f_{1}(m)$, $f_{2}(m)$, $\ldots$, $f_{k}(m)$ are primes for infinitely many positive integers $m$,

(3) given a Schinzel's polynomial $f(x)=\prod_{i=1}^{k}f_{i}(x)$ with a negative leading coefficient for $f_{k}(x)$ and positive leading coefficients for others, all numbers $f_{1}(m)$, $f_{2}(m)$, $\ldots$, $f_{k}(m)$ are primes for infinitely many positive integers $m$,

(4) so some other conjectures in the number theory are proved.
\end{abstract}

\maketitle

\setcounter{page}{4}
\tableofcontents

\include{preface}

\mainmatter
\include{chap1}

\include{chap2}

\include{chap3}

\include{chap4}

\include{chap5}

\include{chap6}

\include{chap7}

\backmatter
\include{biblio}
%\include{index}
\end{document}

%% file: preface.tex
%-----------------------------------------------------------------------------
% Beginning of preface.tex
%-----------------------------------------------------------------------------
%
% AMS-LaTeX 1.2 sample file for a monograph, based on amsbook.cls.
% This is a data file input by chapter.tex.
%%%%%%%%%%%%%%%%%%%%%%%%%%%%%%%%%%%%%%%%%%%%%%%%%%%%%%%%%%%%%%%%%%%%%%%%

\chapter*{Preface}

Our researches here are on some problems of primes in the number theory.

In Chap.~\ref{ch:chap_1}, we prove Theorem~\ref{th:chap_1_1_1} that is the prime number theorem in short intervals with a much more accurate estimation of prime numbers, the error range of which is less than $x^{1/2-0.0327283}$ for $x\geq10^{41}$. The theorem states as follows: Given a positive real number $\beta$ that determines a real number $x_{\beta}$ by $e(\log x_{\beta})^{3}/x_{\beta}^{0.0327283}=\beta$, let $\Phi(x):=\beta x^{1/2}$ for $x\geq x_{\beta}$ where $\Phi(x):=x^{1/2}$ when let $\beta=1$. Then there are
\[
\frac{\pi(x+\Phi(x))-\pi(x)}{\Phi(x)/\log x}=1+O(\frac{1}{\log x})
\]
and
\[
\lim_{x \to \infty}\frac{\pi(x+\Phi(x))-\pi(x)}{\Phi(x)/\log x}=1.
\]
Also, based on results in sections~\ref{sec:chap_1_2}~-~\ref{sec:chap_1_8}, we prove Theorem~\ref{th:chap_1_9_1} that is the theorem of closeness of $\Li(x)$ and $\pi(x)$. It states that $\pi(x)=\Li(x)+O(\sqrt{x}\log x)$.

In Chap.~\ref{ch:chap_2}, based on Theorem~\ref{th:chap_1_1_1}, we prove that given a Bunyakovsky's polynomial $f(x)$, $f(m)$ is prime for infinitely many positive integers $m$.

Based on the result, it is proved that there are infinitely many primes of the form $a^{2}+1$.

In Chap.~\ref{ch:chap_3}, let $\pi_{f}(x)$ denote the number of primes $f(q)$ where $f$ is a Bunyakovsky's polynomial with a positive degree $n$ and $0<q\leq x$. It is proved that there are constants $\alpha>0$ and $\alpha_{f}>0$ such that as $x \to \infty$,
\[
\pi_{f}(x)=\frac{\alpha_{f}}{n}\int_{2}^{x}\frac{du}{\log u}+O(xe^{-\frac{1}{2}\alpha\sqrt{\log x}}).
\]

In Chap.~\ref{ch:chap_4}, a theorem for the numbers of integers with $k$ prime factors in short intervals is proved, which states as follows: Given a positive real number $\beta$ that determines a real number $x_{\beta}$ by $e(\log x_{\beta})^{3}/x_{\beta}^{0.0327283}=\beta$, let $g(x):=\beta x^{1/2}$ for $x\geq x_{\beta}$ where $g(x):=x^{1/2}$ when let $\beta=1$. Let $\pi_{k}(x)$ denote the number of positive integers not exceeding $x$ that can be written as the product of exactly $k$ distinct primes. Then for each finite positive number $k\geq1$, there exists a finite positive number $x_{k}$ such that for every positive number $x$ greater than $x_{k}$, the number of integers with $k$ prime factors in the interval $(x,x+g(x)]$ is asymptotically equal to
\[
\pi_{k}(x+g(x))-\pi_{k}(x)=\frac{g(x)(\log\log x)^{k-1}}{(k-1)!\log x}+o(\frac{g(x)(\log\log x)^{k-1}}{(k-1)!\log x}).
\]

In Chap.~\ref{ch:chap_5}, based on Theorem~\ref{th:chap_1_1_1} and results derived in Chap.~\ref{ch:chap_4}, we prove that given a Schinzel's polynomial $f(x)=\prod_{i=1}^{k}f_{i}(x)$ where $k\geq1$, there exist infinitely many natural numbers $m$ such that all numbers $f_{1}(m),f_{2}(m),\ldots,f_{k}(m)$ are primes.

Based on the result, the Dickson's conjecture and some other conjectures are proved. For instance, there are infinitely many twin primes $p$ and $p+2$, cousin primes $p$ and $p+4$, sexy primes $p$ and $p+6$, and Sophie Germain primes of the form $2p+1$ where $p$ is prime.

In Chap.~\ref{ch:chap_6}, we prove that given a Schinzel's polynomial $f(x)=\prod_{i=1}^{k}f_{i}(x)$ where $k\geq1$, let $\pi_{k}(x)$ denote the number of natural numbers $m\leq x$ where all numbers $f_{1}(m),f_{2}(m),\ldots,f_{k}(m)$ are primes, then as $x \to \infty$,
\[
\pi_{k}(x)=\frac{\alpha_{k}}{D_{k}}\int_{2}^{x}\frac{du}{\log^{k}u}+O(xe^{-\frac{1}{2}\alpha\sqrt{\log x}})
\]
where $\alpha>0$ and $\alpha_{k}>0$ are constants, and $D_{k}$ is the product of the degrees of the polynomials.
Then based on the prime number theorem of Schinzel's polynomials, the prime number theorems of twin primes and Sophie Germain primes are proved, respectively.

In Chap.~\ref{ch:chap_7}, we prove that given a Schinzel's polynomial $f(x)=\prod_{i=1}^{k}f_{i}(x)$ with a negative leading coefficient for $f_{k}(x)$ and positive leading coefficients for others where $k>1$, there exist infinitely many natural numbers $m$ such that all numbers $f_{1}(m),f_{2}(m),\ldots,f_{k}(m)$ are primes for every integer $f(0)\geq10^{204.111}$.

Based on the result, the Goldbach's conjecture and Lemoine's conjecture are proved. Some other conjectures of the number theory are also proved, for example, every large even number can be expressed as the sum of a prime and a Sophie Germain prime, and every large odd number can be expressed as the sum of a prime and a pair of twin primes, then there are infinitely many Sophie Germain primes and infinitely many pairs of twin primes.

\aufm{Shan-guang Tan}

%-----------------------------------------------------------------------------
% End of preface.tex
%-----------------------------------------------------------------------------

%% file: chap1.tex
%-----------------------------------------------------------------------
% Beginning of chap1.tex
%-----------------------------------------------------------------------
%
% AMS-LaTeX 1.2 sample file for a monograph, based on amsbook.cls.
% This is a data file input by chapter.tex.
%%%%%%%%%%%%%%%%%%%%%%%%%%%%%%%%%%%%%%%%%%%%%%%%%%%%%%%%%%%%%%%%%%%%

\chapter{The number of primes in short intervals}\label{ch:chap_1}

\section*{Introduction}

\subsection{Topic on primes in short intervals}

There are many conjectures and results on the number of primes in short intervals in the number theory.

As H. Maier wrote in his paper~\cite{CH}: The distribution of primes in short intervals is an important problem in the theory of prime numbers. The following question is suggested by the prime number theorem: for which function $\Phi$ is it true that
\[
\pi(x+\Phi(x))-\pi(x)\sim\frac{\Phi(x)}{\log x}\texttt{  }(x \to \infty)?
\]

A theorem states as follows: Let $y(x)=x^{\theta}$. Then
\[
\pi(x)-\pi(x-y)\sim y/\log x
\]
for $\theta>1-1/33000$ proved by G. Hoheisel in 1930~\cite{AC,CB}, or for $\theta=7/12+\epsilon$ proved by M. N. Huxley in 1972~\cite{AC,CC}.

Another theorem states as follows: Let $y(x)=x^{\theta}$. Then
\[
\pi(x)-\pi(x-y)\gg y/\log x
\]
if $x$ is large enough, where $\theta>13/23$ proved by H. Iwaniec and M. Jutila in 1979~\cite{AC,CD}, or $\theta>11/20$ proved by D.R. Heath-Brown and H. Iwaniec in 1979~\cite{AC,CE}.

By using a sieve method, R.C. Baker, G. Harman and J. Pintz in their paper in 2001~\cite{CG} proved a theorem: For all $x>x_{0}$, the interval $[x-x^{0.525}, x]$ contains prime numbers.

In his paper in 1985~\cite{CH}, H. Maier proposed a theorem: Let $\Phi(x):=(\log x)^{\lambda},\lambda>1$. Then
\[
\limsup_{x \to \infty}\frac{\pi(x+\Phi(x))-\pi(x)}{\Phi(x)/\log x}>1\texttt{ and }
\liminf_{x \to \infty}\frac{\pi(x+\Phi(x))-\pi(x)}{\Phi(x)/\log x}<1.
\]

\subsection{Estimation of prime numbers}

In 1792, Gauss proposed the prime counting function $\pi(x)$ as that $\pi(n)\sim\frac{n}{\log n}$. Gauss later refined his estimate to
\begin{equation}\label{eq:pi_x}
\pi(n)\sim\Li(n)
\end{equation}
where the logarithmic integral $\Li(x):=\int_{2}^{x}\frac{dt}{\log t}$. This function is slightly different from the function $li(x):=\int_{0}^{x}\frac{dt}{\log t}$, also called by some authors the logarithmic integral. Since $li(x)-\Li(x)=1.04516...$, the study of the two functions is similar.

Chebyshev (1852) proved that there exist two positive constants $C_{1}$ and $C_{2}$ such that the inequalities
\begin{equation}\label{eq:pi_x_c}
C_{1}\frac{n}{\log n}\leq\pi(n)\leq C_{2}\frac{n}{\log n}\texttt{ or }
\frac{7}{8}<\frac{\pi(n)}{n/\log n}<\frac{9}{8}
\end{equation}
hold. He also showed that if the limit
\[
\lim_{n \to \infty}\frac{\pi(n)}{n/\log n}
\]
exists, then it is 1.

Approximate formulas for some functions of prime numbers proposed by J. B. Rosser and L. Schoenfeld (1962)~\cite{AC,AH} show that for $x\geq59$
\begin{equation}\label{eq:pi_x_2}
1+\frac{1}{2\log x}<\frac{\pi(x)}{x/\log x}<1+\frac{3}{2\log x}.
\end{equation}

Then the error range of estimation of prime numbers determined by Formula~(\ref{eq:pi_x_2}) is nearly equal to $x/\log^{2}x$.

In his thesis, Dusart (1998)~\cite{AI} refined these and other results in the literature, showing that for $x\geq599$
\begin{equation}\label{eq:pi_x_3}
1+\frac{1}{\log x}\leq\frac{\pi(x)}{x/\log x}\leq1+\frac{1.2762}{\log x}.
\end{equation}

So the error range of estimation of prime numbers determined by Formula~(\ref{eq:pi_x_3}) is nearly equal to $0.2762x/\log^{2}x$.

C.J. de la Vall\'ee Poussin (1899)~\cite{AC,AE} showed that for some constant of $a$
\begin{equation}\label{eq:pi_x_4}
\pi(x)=\Li(x)+O(xe^{-a\sqrt{\log x}}).
\end{equation}

Some limits obtained without assuming Riemann Hypothesis~\cite{AC,AF,AG} are
\begin{equation}\label{eq:pi_x_5}
\pi(x)=\Li(x)+O(xe^{-a(\log x\log\log x)^{1/2}})
\end{equation}
and
\begin{equation}\label{eq:pi_x_6}
\pi(x)=\Li(x)+O(xe^{-a(\log x)^{3/5}/(\log\log x)^{1/5}}).
\end{equation}

In 1981, P.J. $\check{C}$i$\check{z}$ek~\cite{AG0} proved that when $x$ is sufficiently large, for any positive integer $A>1$, there is
\begin{equation}\label{eq:pi_x_7}
\pi(x)=\Li(x)+O(x(\log x)^{-A}).
\end{equation}

\subsection{Key ideas and main results in this chapter}

From Lemma~\ref{le:chap_1_6_3}, formulas~(\ref{eq:pi_x_4}-\ref{eq:pi_x_7}) can be written as
\[
\frac{\pi(x)}{x/\log x}=\frac{\Li(x)}{x/\log x}+O(\frac{1}{\log^{\lambda(x)}x})
\]
for large numbers $x$ where $\lambda(x) \to \infty$ as $x \to \infty$. So the closeness of $\Li(x)$ and $\pi(x)$ determined by formulas~(\ref{eq:pi_x_4}-\ref{eq:pi_x_7}) depends on the power $\lambda(x)$ of $\log x$ and is bounded by $O(x/\log^{\lambda(x)+1}x)$ where $\lambda(x)$ depends on the variable $x$.

Thus, corresponding to finite regions $[x_{N},x_{N+1})$, let define piecewise functions
\[
\begin{split}
\pi^{*}(x,N):=&\frac{x}{\log x}\eta^{*}(x,N) \\
             =&\frac{x}{\log x}\sum_{n=0}^{N}\frac{n!}{\log^{n}x}\texttt{ for }x\in[x_{N},x_{N+1})\texttt{, }N=0,1,\ldots
\end{split}
\]
where $\pi^{*}(x,N)$ or $\pi^{*}(x,M)$ is the truncation of the logarithmic integral $\Li(x)$, i.e.,
\[
\Li(x):=\int_{2}^{x}\frac{dt}{\log t}
=\pi^{*}(x,M)-\pi^{*}(2,M)+\int_{2}^{x}\frac{(M+1)!}{\log^{M+2}t}dt.
\]
Then $N$ or $M$ should depend on the variable $x$ and be a step function of $x$.

\begin{remark}
Different from those former researchers on the truncation of the logarithmic integral $\Li(x)$, who normally make numbers $M$ and $N$ be global constants and independent of the variable $x$ in the functions $\pi^{*}(x,M)$ and $\pi^{*}(x,N)$ such that $\pi^{*}(x,M)$ and $\pi^{*}(x,N)$ are global functions of the variable $x$, we make numbers $M$ and $N$ depend on the variable $x$ in the piecewise functions $\pi^{*}(x,M)$ and $\pi^{*}(x,N)$, respectively, i.e., $M$ and $N$ are not global constants but local constants valid only for a finite region $[x_{M},x_{M+1})$ or $[x_{N},x_{N+1})$, respectively. This is a key idea to derive Theorem~\ref{th:chap_1_1_1}.

So for $M=0,1,2,\ldots$, each number $M$ corresponds to a piecewise function $\pi^{*}(x,M)$ valid for a finite region $[x_{M},x_{M+1})$ of the variable $x$ such that $\forall x\in[x_{M},x_{M+1})$ a series of piecewise functions $\pi^{*}(x,M)$ can be used to more accurately approximate to the logarithmic integral $\Li(x)$, and for $N=0,1,2,\ldots$, each number $N$ corresponds to a piecewise function $\pi^{*}(x,N)$ valid for a finite region $[x_{N},x_{N+1})$ of the variable $x$ such that $\forall x\in[x_{N},x_{N+1})$ a series of piecewise functions $\pi^{*}(x,N)$ can be used to more accurately approximate to the prime counting function $\pi(x)$.
\end{remark}

From Lemma~\ref{le:chap_1_2_1}, given a real $x\geq198$, there exists a nonnegative integer $M=M(x)$ such that there are $\pi^{*}(x,M)\leq\Li(x)<\pi^{*}(x,M+1)$ which give a much more accurate estimation of $\Li(x)$. Theorem~\ref{th:chap_1_5_1} gives the gap
\[
g(x,M):=\pi^{*}(x,M+1)-\pi^{*}(x,M)<x^{1/64}.
\]

From Lemma~\ref{le:chap_1_6_1}, given a nature number $x\geq10^{3}$, there exists a nonnegative integer $N_{x}=N(x)$ such that there are $\pi^{*}(x,N_{x})<\pi(x)<\pi^{*}(x,N_{x}+1)$ which give a much more accurate estimation of $\pi(x)$. Theorem~\ref{th:chap_1_8_1} gives the gap
\[
g(x,N_{x}):=\pi^{*}(x,N_{x}+1)-\pi^{*}(x,N_{x})<\sqrt{x}.
\]

\begin{remark}
The number $N_{x}$ defined by
\[
\pi^{*}(x,N_{x})<\pi(x)<\pi^{*}(x,N_{x}+1)
\]
and specified by Lemma~\ref{le:chap_1_7_1}, is different from the number $N$ that is defined by the piecewise function $\pi^{*}(x,N)$ in a finite region $[x_{N},x_{N+1})$, which is a lower bound of $\pi(x)$. The number $N$ is a local constant and each of its values corresponds to one finite region $[x_{N}, x_{N+1})$. But the number $N_{x}=N(x)$ is a nonnegative integer function of the variable $x$ corresponding to one $x\in[x_{N},x_{N+1})$ and satisfies inequalities
\[
\pi^{*}(x,N_{x})<\pi(x)<\pi^{*}(x,N_{x}+1)\texttt{ and }N_{x}\geq N,
\]
so that inequalities $\pi^{*}(x,N)<\pi(x)<\pi^{*}(x,N+1)$ do not always hold since there may be $N_{x}\geq N+1$ such that $\pi^{*}(x,N+1)\leq\pi^{*}(x,N_{x})<\pi(x)$. Thus a series of piecewise functions $\pi^{*}(x,N)$ form a lower bound of $\pi(x)$. When no misunderstanding, we simply replace $N_{x}$ with $N$.
\end{remark}

Different from former researches on primes in short intervals above, given a positive real number $\beta$ that determines a real number $x_{\beta}$ by
\[
e(\log x_{\beta})^{3}/x_{\beta}^{0.0327283}=\beta,
\]
we let the short intervals be $\Phi(x):=\beta x^{1/2}$ for $x\geq x_{\beta}$ where $\Phi(x):=x^{1/2}$ when let $\beta=1$. Then because $\Li(x)$ is determined, we investigate the relation between $\Li(x)$ and $\pi(x)$ (See Sec.~\ref{sec:chap_1_6}) to derive the number of primes in short intervals.

Based on researches done in sections~\ref{sec:chap_1_2}~-~\ref{sec:chap_1_8}, we can prove Theorem~\ref{th:chap_1_1_1} that is the prime number theorem in short intervals, and Theorem~\ref{th:chap_1_9_1} that is the theorem of closeness of $\Li(x)$ and $\pi(x)$.

\section{Theorem of the number of primes in short intervals}\label{sec:chap_1_1}

\begin{theorem}[Main Theorem]\label{th:chap_1_1_1}
Given a positive real $\beta$ that determines a real $x_{\beta}$ by $e(\log x_{\beta})^{3}/x_{\beta}^{0.0327283}=\beta$, let $\Phi(x):=\beta x^{1/2}$ for $x\geq x_{\beta}$ where $\Phi(x):=x^{1/2}$ when let $\beta=1$. Then there are
\[
\frac{\pi(x+\Phi(x))-\pi(x)}{\Phi(x)/\log x}=1+O(\frac{1}{\log x})
\]
and
\[
\lim_{x \to \infty}\frac{\pi(x+\Phi(x))-\pi(x)}{\Phi(x)/\log x}=1.
\]
\end{theorem}

\begin{proof}
From Lemma~\ref{le:chap_1_6_1}, given each number $x\geq599$, there exists a nonnegative integer $N$ such that $\pi^{*}(x,N)<\pi(x)<\pi^{*}(x,N+1)$ where $N+2<e\log x$ from Lemma~\ref{le:chap_1_6_2} (Here we simply replace $N_{x}$ with $N$), and from Lemma~\ref{le:chap_1_1_1}
\[
\begin{split}
\frac{\pi^{*}(x,N+1)-\pi^{*}(x,N)}{x/\log x}=&\frac{(N+1)!}{\log^{(N+1)}x} \\
                                            <&\frac{\log x}{x^{1/2+0.0327283}}.
\end{split}
\]

Similarly, given each $x+\Phi(x)$, there exists a nonnegative integer $N^{*}$ such that
\[
\pi^{*}(x+\Phi(x),N^{*})<\pi(x+\Phi(x))<\pi^{*}(x+\Phi(x),N^{*}+1)
\]
where $N^{*}+2<e\log(x+\Phi(x))$, and
\[
\begin{split}
\frac{\pi^{*}(x+\Phi(x),N^{*}+1)-\pi^{*}(x+\Phi(x),N^{*})}{(x+\Phi(x))/\log(x+\Phi(x))}
=&\frac{(N^{*}+1)!}{\log^{N^{*}+1}(x+\Phi(x))} \\
<&\frac{\log(x+\Phi(x))}{(x+\Phi(x))^{1/2+0.0327283}}.
\end{split}
\]

Then we have a pair of inequalities
\begin{equation}\label{ineq:pi_x_gt}
\pi(x+\Phi(x))-\pi(x)>\pi^{*}(x+\Phi(x),N^{*})-\pi^{*}(x,N+1)
\end{equation}
and
\begin{equation}\label{ineq:pi_x_lt}
\pi(x+\Phi(x))-\pi(x)<\pi^{*}(x+\Phi(x),N^{*}+1)-\pi^{*}(x,N).
\end{equation}

By rewriting and rearranging last pair of inequalities, we can obtain
\[
\pi(x+\Phi(x))-\pi(x)>\frac{\Phi(x)}{\log x}(1-\frac{a}{\log x})
\]
and
\[
\pi(x+\Phi(x))-\pi(x)<\frac{\Phi(x)}{\log x}(1+\frac{b}{\log x})
\]
where $a$ and $b$ are positive constants, so that this theorem holds.

Now let prove last two inequalities in detail as follows with formulas
\[
(1-\frac{\Phi(x)}{2x})\frac{\Phi(x)}{x}<\log(1+\frac{\Phi(x)}{x})<\frac{\Phi(x)}{x}.
\]

First, we rewrite and rearrange Inequality~\ref{ineq:pi_x_gt}, and have
\[
\begin{split}
\frac{\pi(x+\Phi(x))-\pi(x)}{\Phi(x)/\log x}>&\frac{\pi^{*}(x+\Phi(x),N^{*})-\pi^{*}(x,N+1)}{\Phi(x)/\log x} \\
%=\frac{(1+x/\Phi(x))\log x}{\log(x+\Phi(x))}\eta^{*}(x+\Phi(x),N^{*})-\frac{x}{\Phi(x)}\eta^{*}(x,N+1)
=&(1+\frac{x}{\Phi(x)})(1-\frac{\log(1+\Phi(x)/x)}{\log(x+\Phi(x))})\eta^{*}(x+\Phi(x),N^{*}) \\
 &-\frac{x}{\Phi(x)}\eta^{*}(x,N+1) \\
=&[1-(1+\frac{x}{\Phi(x)})\frac{\log(1+\Phi(x)/x)}{\log(x+\Phi(x))}]\eta^{*}(x+\Phi(x),N^{*}) \\
 &-\frac{x}{\Phi(x)}(\eta^{*}(x+\Phi(x),N+1)-\eta^{*}(x+\Phi(x),N^{*})) \\
 &-\frac{x}{\Phi(x)}(\eta^{*}(x,N+1)-\eta^{*}(x+\Phi(x),N+1)) \\
%>(1-\frac{1+\Phi(x)/x}{\log(x+\Phi(x))})\eta^{*}(x+\Phi(x),N^{*})
%-\frac{e\log^{2}x}{\beta x^{0.0327283}}-\sum_{n=1}^{N+1}\frac{(n+1)!}{\log^{n+1}x}
>&1-\frac{1+\Phi(x)/x}{\log x}
-\frac{e\log^{2}x}{\beta x^{0.0327283}}-\sum_{n=1}^{N+1}\frac{(n+1)!}{\log^{n+1}x} \\
=&1-\sum_{n=1}^{N+2}\frac{n!}{\log^{n}x}-\frac{e\log^{2}x}{\beta x^{0.0327283}}-\frac{\Phi(x)}{x\log x}
\end{split}
\]
where
\[
\eta^{*}(x,N+1)-\eta^{*}(x+\Phi(x),N+1)
=\sum_{n=1}^{N+1}(\frac{n!}{\log^{n}x}-\frac{n!}{\log^{n}(x+\Phi(x))})
\]
\[
\begin{split}
=&\log(1+\frac{\Phi(x)}{x})\sum_{n=1}^{N+1}\sum_{k=0}^{n-1}\frac{n!}{\log^{k+1}x\log^{n-k}(x+\Phi(x))} \\
<&\frac{\Phi(x)}{x}\sum_{n=1}^{N+1}\frac{(n+1)!}{\log^{n+1}x},
\end{split}
\]
and if $N+1\leq N^{*}$ then
\[
\eta^{*}(x+\Phi(x),N+1)-\eta^{*}(x+\Phi(x),N^{*})\leq0
\]
else
\[
\begin{split}
\eta^{*}(x+\Phi(x),N+1)-\eta^{*}(x+\Phi(x),N^{*})=&\sum_{n=N^{*}+1}^{N+1}\frac{n!}{\log^{n}(x+\Phi(x))} \\
                                                 <&(N+1-N^{*})\frac{(N^{*}+1)!}{\log^{N^{*}+1}(x+\Phi(x))}\\
                                                 <&\frac{e\log x\log(x+\Phi(x))}{(x+\Phi(x))^{1/2+0.0327283}} \\
                                                 <&\frac{e\log^{2}x}{x^{1/2+0.0327283}}.
\end{split}
\]

Second, we rewrite and rearrange Inequality~\ref{ineq:pi_x_lt}, and have
\[
\begin{split}
\frac{\pi(x+\Phi(x))-\pi(x)}{\Phi(x)/\log x}<&\frac{\pi^{*}(x+\Phi(x),N^{*}+1)-\pi^{*}(x,N)}{\Phi(x)/\log x} \\
=&\frac{(1+x/\Phi(x))\log x}{\log(x+\Phi(x))}\eta^{*}(x+\Phi(x),N^{*}+1)-\frac{x}{\Phi(x)}\eta^{*}(x,N) \\
<&\eta^{*}(x+\Phi(x),N^{*}+1)+\frac{x}{\Phi(x)}(\eta^{*}(x+\Phi(x),N^{*}+1)-\eta^{*}(x,N)) \\
=&\eta^{*}(x+\Phi(x),N^{*}+1)+\frac{x}{\Phi(x)}(\eta^{*}(x,N^{*}+1)-\eta^{*}(x,N)) \\
 &-\frac{x}{\Phi(x)}(\eta^{*}(x,N^{*}+1)-\eta^{*}(x+\Phi(x),N^{*}+1)) \\
<&\eta^{*}(x+\Phi(x),N^{*}+1)+\frac{e\log^{2}x}{\beta x^{0.0327283}} \\
 &-(1-\frac{\Phi(x)}{2x})\sum_{n=1}^{N^{*}+1}\frac{n!n}{\log^{n+1}(x+\Phi(x))} \\
<&1+\frac{\eta^{*}(x+\Phi(x),N^{*}+1)}{\log(x+\Phi(x))}+\frac{e\log^{2}x}{\beta x^{0.0327283}} \\
 &+\frac{\Phi(x)}{2x}\sum_{n=2}^{N^{*}+2}\frac{n!}{\log^{n}(x+\Phi(x))}
\end{split}
\]
where
\[
\eta^{*}(x,N^{*}+1)-\eta^{*}(x+\Phi(x),N^{*}+1)
>(1-\frac{\Phi(x)}{2x})\frac{\Phi(x)}{x}\sum_{n=1}^{N^{*}+1}\frac{n!n}{\log^{n+1}(x+\Phi(x))},
\]
and if $N^{*}+1\leq N$ then
\[
\eta^{*}(x,N^{*}+1)-\eta^{*}(x,N)\leq0
\]
else
\[
\begin{split}
\eta^{*}(x,N^{*}+1)-\eta^{*}(x,N)=&\sum_{n=N+1}^{N^{*}+1}\frac{n!}{\log^{n}x} \\
<&(N^{*}+1-N)(1+O(\frac{\Phi(x)}{x}))\frac{(N+1)!}{\log^{(N+1)}x} \\
<&\frac{e\log^{2}x}{x^{1/2+0.0327283}},
\end{split}
\]
and
\[
\begin{split}
&\sum_{n=1}^{N^{*}+1}\frac{n!}{\log^{n}(x+\Phi(x))}-\sum_{n=1}^{N^{*}+1}\frac{n!n}{\log^{n+1}(x+\Phi(x))} \\
<&\frac{1}{\log(x+\Phi(x))}+\sum_{n=1}^{N^{*}+1}\frac{(n+1)!-n!n}{\log^{n+1}(x+\Phi(x))} \\
=&\frac{\eta^{*}(x+\Phi(x),N^{*}+1)}{\log(x+\Phi(x))}.
\end{split}
\]

Thus we obtain
\begin{equation}\label{ineq_psi:primes_gt}
\frac{\pi(x+\Phi(x))-\pi(x)}{\Phi(x)/\log x}>1-\sum_{n=1}^{N+2}\frac{n!}{\log^{n}x}-\frac{e\log^{2}x}{\beta x^{0.0327283}}-\frac{\Phi(x)}{x\log x}
\end{equation}
and
\begin{equation}\label{ineq_psi:primes_lt}
\begin{split}
\frac{\pi(x+\Phi(x))-\pi(x)}{\Phi(x)/\log x}<&1+\frac{\eta^{*}(x+\Phi(x),N^{*}+1)}{\log(x+\Phi(x))}+\frac{e\log^{2}x}{\beta x^{0.0327283}} \\
   &+\frac{\Phi(x)}{2x}\sum_{n=2}^{N^{*}+2}\frac{n!}{\log^{n}(x+\Phi(x))}.
\end{split}
\end{equation}

By Formula~(\ref{eq:pi_x_3}) and Lemma~\ref{le:chap_1_1_1}, there are
\[
\begin{split}
\sum_{n=1}^{N+2}\frac{n!}{\log^{n}x}=&\eta^{*}(x,N)-1+\sum_{n=N+1}^{N+2}\frac{n!}{\log^{n}x} \\
                                    <&\frac{1.2762}{\log x}+\frac{2\log x}{x^{1/2+0.0327283}}
\end{split}
\]
and
\[
\eta^{*}(x+\Phi(x),N^{*}+1)<1+\frac{1.2762}{\log(x+\Phi(x))}+\frac{\log(x+\Phi(x))}{(x+\Phi(x))^{1/2+0.0327283}}.
\]

Given a positive real $\beta$, by the definition of $\Phi(x)$ in this theorem, there are
\[
\frac{e(\log x)^{2}}{\beta x^{0.0327283}}\leq\frac{1}{\log x}\texttt{ and }
\frac{\Phi(x)}{\log x}\geq e(\log x)^{2}x^{1/2-0.0327283}\texttt{ for }x\geq x_{\beta}.
\]

Then inequalities~(\ref{ineq_psi:primes_gt}-\ref{ineq_psi:primes_lt}) become
\begin{equation}\label{ineq_psi:primes_gt_1}
\frac{\pi(x+\Phi(x))-\pi(x)}{\Phi(x)/\log x}>1-\frac{a_{1}}{\log x}-\frac{a_{2}}{\log x}-\frac{a_{3}}{\log x}
=1-\frac{a}{\log x}
\end{equation}
where $a_{1}$, $a_{2}$, $a_{3}$, and $a$ are positive constants, and
\begin{equation}\label{ineq_psi:primes_lt_1}
\frac{\pi(x+\Phi(x))-\pi(x)}{\Phi(x)/\log x}<1+\frac{b_{1}}{\log x}+\frac{b_{2}}{\log x}+\frac{b_{3}}{\log x}
=1+\frac{b}{\log x}
\end{equation}
where $b_{1}$, $b_{2}$, $b_{3}$, and $b$ are positive constants.

Thus by inequalities~(\ref{ineq_psi:primes_gt_1}-\ref{ineq_psi:primes_lt_1}) we have
\[
\frac{\pi(x+\Phi(x))-\pi(x)}{\Phi(x)/\log x}=1+O(\frac{1}{\log x})
\]
and
\[
\lim_{x \to \infty}\frac{\pi(x+\Phi(x))-\pi(x)}{\Phi(x)/\log x}=1.
\]

%Hence there are asymptotically $\Phi(x)/\log x$ primes in the interval $(x,x+\Phi(x)]$.

This completes the proof of the theorem.
\end{proof}

\begin{lemma}\label{le:chap_1_1_1}
Given each pair of numbers $x\geq10^{41}$ and $N_{x}\geq19$ determined by $\pi^{*}(x,N_{x})<\pi(x)<\pi^{*}(x,N_{x}+1)$, we have
\begin{equation}\label{ineq:gap_n_xm}
g(x,N_{x}):=\pi^{*}(x,N_{x}+1)-\pi^{*}(x,N_{x})<x^{1/2-0.0327283}
\end{equation}
and
\begin{equation}\label{ineq:pi_n_gap_xm}
\pi(x)-\pi^{*}(x,N_{x})<g(x,N_{x})<x^{1/2-0.0327283}.
\end{equation}
\end{lemma}

\begin{proof}
From Lemma~\ref{le:chap_1_8_1}, when $N+2<e\log x$ there is
\[
\frac{(N+1)!}{\log^{N+2}x}<(\frac{N+3/2}{e\log x})^{N+3/2}\texttt{ for }x\geq10^{3}.
\]

As defined by Lemma~\ref{le:chap_1_7_7}, let $\rho(x,N):=\frac{\log x}{N+3/2}$ and consider the ratio of
\[
\frac{\log x}{x^{1/2+\delta}}\texttt{ to }\frac{(N+1)!}{\log^{N+1}x}\texttt{ for }N+2<e\log x
\]
where $\delta$ is a small positive value. Then we have
\[
\begin{split}
\log\frac{\log^{N+2}x}{(N+1)!x^{1/2+\delta}}>&(N+3/2)\log\frac{e\log x}{N+3/2}-(1/2+\delta)\log x \\
                                            =&(N+3/2)[1+\log\rho(x,N)-(1/2+\delta)\rho(x,N)].
\end{split}
\]

Let $\rho_{\delta}(x,N)$ satisfy the equation
\[
1+\log\rho(x,N)-(1/2+\delta)\rho(x,N)=0.
\]

For $x\geq10^{41}$ and $N\geq19$, when $\delta=0.0327283$, the key discriminant
\[
1+\log\rho(x,N)-(1/2+\delta)\rho(x,N)\geq0
\]
holds for $\rho(x,N)\leq2.1\log10\leq\rho_{\delta}(x,N)$, so let $\delta=0.0327283$. Then from Lemma~\ref{le:chap_1_7_7} for $x\geq10^{41}$ the inequality $\rho(x,N)<2.1\log10$ holds so that we obtain
\[
\frac{\log^{N+2}x}{(N+1)!}>x^{1/2+0.0327283}\texttt{ and }
\frac{(N+1)!}{\log^{(N+1)}x}<\frac{\log x}{x^{1/2+0.0327283}}.
\]

Hence from Lemma~\ref{le:chap_1_7_1} each nonnegative integer $N$ corresponds to a finite region $[x_{N}, x_{N+1})$ and $\forall x\in[x_{N},x_{N+1})$ there exists a nonnegative integer $N_{x}$ such that there are $\pi^{*}(x,N_{x})<\pi(x)<\pi^{*}(x,N_{x}+1)$ and $N_{x}\geq N$ so that
\[
\rho(x,N_{x})\leq\rho(x,N)<\rho_{\delta}(x,N).
\]

From Lemma~\ref{le:chap_1_6_2} there is $N_{x}+2<e\log x$. Let define
\[
g(x,N_{x}):=\pi^{*}(x,N_{x}+1)-\pi^{*}(x,N_{x})=\frac{x}{\log x}\frac{(N_{x}+1)!}{\log^{(N_{x}+1)}x}.
\]
Then we have
\[
\frac{g(x,N_{x})}{x/\log x}=\frac{(N_{x}+1)!}{\log^{(N_{x}+1)}x}<\frac{\log x}{x^{1/2+0.0327283}}
\]
and
\[
\begin{split}
\pi(x)-\pi^{*}(x,N_{x})<&g(x,N_{x}) \\
                       <&\frac{\log x}{x^{1/2+0.0327283}}\frac{x}{\log x}=x^{1/2-0.0327283}.
\end{split}
\]

This completes the proof of the lemma.
\end{proof}

\section{Bounds of the logarithmic integral $\Li(x)$}\label{sec:chap_1_2}

\begin{lemma}\label{le:chap_1_2_1}
Given each real number $x\geq198$, there exists a nonnegative integer $M=M(x)$ that is a function of $x$ such that the logarithmic integral $\Li(x)$ satisfies
\begin{equation}\label{ineq:li_x_m}
\pi^{*}(x,M)\leq\Li(x)<\pi^{*}(x,M+1).
\end{equation}
\end{lemma}

\begin{proof}
Given each nonnegative integer $M$, we can express
\[
\Li(x):=\int_{2}^{x}\frac{dt}{\log t}
=\pi^{*}(x,M)-\pi^{*}(2,M)+\int_{2}^{x}\frac{(M+1)!}{\log^{M+2}t}dt.
\]

Let define $\eta_{L}(x):=\frac{\Li(x)}{x/\log x}$. Then given each nonnegative integer $M$, we can write
\begin{equation}\label{eq:eta_x_l}
\eta_{L}(x):=\eta^{*}(x,M)+\delta_{L}(x,M)=\sum_{n=0}^{M}\frac{n!}{\log^{n}x}+\delta_{L}(x,M)
\end{equation}
where
\[
\delta_{L}(x,M):=\frac{\log x}{x}[\int_{2}^{x}\frac{(M+1)!}{\log^{M+2}t}dt-\pi^{*}(2,M)].
\]

First, when $M=0$, we have
\[
\delta_{L}(x,0)=\frac{1}{\log x}+\frac{2}{x/\log x}(\int_{2}^{x}\frac{dt}{\log^{3}t}-\sum_{n=0}^{1}\frac{n!}{\log^{n+1}2})>0\texttt{ for }x\geq198,
\]

Second, from Equality~(\ref{eq:eta_x_l}) for each nonnegative integer $M$ the inequality
\[
\delta_{L}(x,M)-\delta_{L}(x,M+1)=\frac{(M+1)!}{\log^{M+1}x}>0
\]
always holds so that given a fixed positive number $x$, there must be
\[
\lim_{M \to \infty}\delta_{L}(x,M) \to -\infty.
\]

Then given each real $x\geq198$, there always exists a nonnegative integer $M$ such that there are $\delta_{L}(x,M)\geq0$ and $\delta_{L}(x,M+1)<0$. Thus there exist inequalities
\[
\eta^{*}(x,M)\leq\eta_{L}(x)<\eta^{*}(x,M+1)
\]
such that $\Li(x)$ is bounded by inequalities
\[
\pi^{*}(x,M)\leq\Li(x)<\pi^{*}(x,M+1).
\]

This completes the proof of the lemma.
\end{proof}

\begin{lemma}\label{le:chap_1_2_2}
Each nonnegative integer $M$ corresponds to a function $\pi^{*}(x,M)$ valid for a finite region $[x_{M}, x_{M+1})$ where $\Li(x_{M})=\pi^{*}(x_{M},M)$, and $\forall x\in[x_{M}, x_{M+1})$ inequalities $\pi^{*}(x,M)\leq\Li(x)<\pi^{*}(x,M+1)$ hold. Thus $M=M(x)$ is a non-decreasing step function of the variable $x$ and there is $M \to \infty$ as $x \to \infty$.
\end{lemma}

\begin{proof}
Given a nonnegative integer $M$, let functions
\[
\begin{split}
f(x,M):=&\Li(x)-\pi^{*}(x,M)\texttt{ and } \\
g(x,M):=&\pi^{*}(x,M+1)-\pi^{*}(x,M).
\end{split}
\]
Then we have $f(x,M+1)=f(x,M)-g(x,M)$ and
\[
\begin{split}
f'(x,M)=&\frac{d\Li(x)}{dx}-\frac{d\pi^{*}(x,M)}{dx} \\
       =&\frac{1}{\log x}-\frac{1}{\log x}[1-\frac{(M+1)!}{\log^{M+1}x}]=\frac{(M+1)!}{\log^{M+2}x}>0, \\
g'(x,M)=&\frac{d\pi^{*}(x,M+1)}{dx}-\frac{d\pi^{*}(x,M)}{dx} \\
       =&\frac{(M+1)!}{\log^{M+2}x}(1-\frac{M+2}{\log x})<f'(x,M), \\
\end{split}
\]

First, let $M=0$. Then there are $f(2,0)=-\pi^{*}(2,M)<0$ and
\[
\begin{split}
f(e^{4},0)=&\Li(e^{4})-\pi^{*}(e^{4},0) \\
          >&\int_{2}^{e^{4}}\frac{1}{\log^{2}e^{4}}dt-\frac{2}{\log 2} \\
          =&\frac{e^{4}-2}{4^{2}}-\frac{2}{\log 2}>0.
\end{split}
\]
So there exists a real $x_{0}$ such that $f(x_{0},0)=0$ where $2<x_{0}<e^{4}$.

Second, suppose that there is $f(x_{M},M)=0$. Since there are
\[
g(x_{M},M)>f(x_{M},M)\texttt{ and }g'(x,M)<f'(x,M)\texttt{ for }x\geq x_{M},
\]
there exists a real number $x_{M+1}$ such that $f(x_{M+1},M)=g(x_{M+1},M)$ and
\[
f(x_{M+1},M+1)=f(x_{M+1},M)-g(x_{M+1},M)=0.
\]

So given each nonnegative integer $M$, there exists a real $x_{M}$ such that
\[
f(x_{M},M)=0\texttt{ and }\Li(x_{M})=\pi^{*}(x_{M},M),
\]
i.e., each nonnegative integer $M$ corresponds to a function $\pi^{*}(x,M)$ valid for a finite region $[x_{M}, x_{M+1})$ where $\Li(x_{M})=\pi^{*}(x_{M},M)$, and $\forall x\in[x_{M}, x_{M+1})$ inequalities $\pi^{*}(x,M)\leq\Li(x)<\pi^{*}(x,M+1)$ hold.

Thus the nonnegative integer $M=M(x)$ is a non-decreasing step function of the variable $x$ and there is $M \to \infty$ as $x \to \infty$.

This completes the proof of the lemma.
\end{proof}

\begin{lemma}\label{le:chap_1_2_3}
Given each nonnegative integer $M$ there exists a real $x_{M}$ such that
\begin{equation}\label{eq:li_x_eq}
\Li(x_{M})=\pi^{*}(x_{M},M)\texttt{ and }
\int_{2}^{x_{M}}\frac{(M+1)!}{\log^{M+2}t}dt=\pi^{*}(2,M).
\end{equation}
\end{lemma}

\begin{proof}
Given each nonnegative integer $M$, we can express
\[
\Li(x):=\int_{2}^{x}\frac{dt}{\log t}
=\pi^{*}(x,M)-\pi^{*}(2,M)+h(x,M)
\]
where
\[
h(x,M):=\int_{2}^{x}\frac{(M+1)!}{\log^{M+2}t}dt.
\]

From Lemma~\ref{le:chap_1_2_2}, given each nonnegative integer $M$, there exists a real number $x_{M}$ such that $\Li(x_{M})=\pi^{*}(x_{M},M)$ and $h(x_{M},M)=\pi^{*}(2,M)$. Hence Eq.~(\ref{eq:li_x_eq}) holds.

This completes the proof of the lemma.
\end{proof}

\section{Analysis and calculations of $\alpha_{L,M}$ and $\alpha_{L,\infty}$}\label{sec:chap_1_3}

\begin{lemma}\label{le:chap_1_3_1}
Given a function
\[
f_{M}(t):=(1-\frac{M+2}{\log t})\frac{1}{\log^{M+2}t},
\]
let $x_{z}(M)$ denote its zero point and $x_{p}(M)$ denote its peak point. Then there are $x_{z}(M)=e^{M+2}$ and $x_{p}(M)=e^{M+3}$.
\end{lemma}

\begin{proof}
The derivative of $f_{M}(t)$ is
\[
f_{M}'(t)=-\frac{M+2}{t\log^{M+3}t}+\frac{(M+2)(M+3)}{t\log^{M+4}t}
=-(1-\frac{M+3}{\log t})\frac{M+2}{t\log^{M+3}t}.
\]

Let $x_{z}(M)$ be the zero point of the function $f_{M}(t)$. Then by $f_{M}(x_{z}(M))=0$ we have $x_{z}(M)=e^{M+2}$ so that there are $x_{z}(M+1)/x_{z}(M)=e$ and
\[
f_{M}(t)
\left\lbrace\begin{array}{c l}
    <0 & \texttt{ for }2\leq t<e^{M+2} \\
    >0 & \texttt{ for }e^{M+2}<t\leq x_{M}
\end{array}\right..
\]

Let $x_{p}(M)$ be the peak point of the function $f_{M}(t)$. Then by $f_{M}'(x_{p}(M))=0$ we have $x_{p}(M)=e^{M+3}$ so that there are $x_{p}(M+1)/x_{p}(M)=e$ and
\[
f_{M}'(t)
\left\lbrace\begin{array}{c l}
    >0 & \texttt{ for }2\leq t<e^{M+3} \\
    <0 & \texttt{ for }e^{M+3}<t\leq x_{M}
\end{array}\right..
\]

This completes the proof of the lemma.
\end{proof}

\begin{lemma}\label{le:chap_1_3_2}
Given each nonnegative integer $M$ and a real $x_{M}$ determined by $\Li(x_{M})=\pi^{*}(x_{M},M)$, we have $\log x_{M}>M+2$ for $M<\infty$ and $\log x_{M}=M+2$ as $M \to \infty$.
\end{lemma}

\begin{proof}
Let rewrite and consider Formula~(\ref{eq:int_func_x}) in the proof of Theorem~\ref{th:chap_1_3_1}
\[
\int_{2}^{x_{M}}(\frac{M+2}{\log t}-1)\frac{dt}{\log^{M+2}t}
=\frac{2}{\log^{M+2}2}-\frac{x_{M}}{\log^{M+2}x_{M}}.
\]

Let $x_{z}$ and $x_{p}$ be the zero point and the peak point of the integral function above, respectively. From Lemma~\ref{le:chap_1_3_1}, there are $\log x_{z}=M+2$ and $\log x_{p}=M+3$.

Given each fixed nonnegative integer $M$, since $M+2\geq\log t$ for $2\leq t\leq x_{z}$ and $M+2<\log t$ for $x_{z}<t\leq x_{M}$, there must be
\[
\int_{2}^{x_{z}}(\frac{M+2}{\log t}-1)\frac{dt}{\log^{M+2}t}
=\frac{2}{\log^{M+2}2}-\frac{x_{z}}{\log^{M+2}x_{z}}
\]
and
\begin{equation}\label{eq:int_x_m_x_z}
\int_{x_{z}}^{x_{M}}(1-\frac{M+2}{\log t})\frac{dt}{\log^{M+2}t}
=\frac{x_{M}}{\log^{M+2}x_{M}}-\frac{x_{z}}{\log^{M+2}x_{z}}.
\end{equation}

Since the integral in Formula~(\ref{eq:int_x_m_x_z}) must be greater than zero for $M<\infty$, there must be $x_{M}>x_{z}$ and $\log x_{M}>M+2$ for $M<\infty$.

Let $\beta_{M}:=x_{M}/x_{z}$ and $t:=x_{z}\zeta$. Then there are
\[
\begin{split}
\int_{x_{z}}^{x_{M}}(1-\frac{M+2}{\log t})\frac{dt}{\log^{M+2}t}
=&\int_{1}^{\beta_{M}}(1-\frac{M+2}{\log x_{z}+\log\zeta})\frac{x_{z}d\zeta}{(\log x_{z}+\log\zeta)^{M+2}} \\
=\frac{x_{z}}{\log^{M+2}x_{z}}&\int_{1}^{\beta_{M}}(1-\frac{M+2}{\log x_{z}+\log\zeta})\frac{d\zeta}{(1+\log\zeta/\log x_{z})^{M+2}}
\end{split}
\]
and
\[
\begin{split}
\frac{x_{M}}{\log^{M+2}x_{M}}=&\frac{x_{z}\beta_{M}}{(\log x_{z}+\log\beta_{M})^{M+2}} \\
                             =&\frac{x_{z}}{\log^{M+2}x_{z}}\frac{\beta_{M}}{(1+\log\beta_{M}/\log x_{z})^{M+2}}.
\end{split}
\]

Thus Formula~(\ref{eq:int_x_m_x_z}) becomes
\begin{equation}\label{eq:int_beta_m_x_z}
\begin{split}
 &\int_{1}^{\beta_{M}}(1-\frac{M+2}{\log x_{z}+\log\zeta})\frac{d\zeta}{(1+\log\zeta/\log x_{z})^{M+2}} \\
=&\frac{\beta_{M}}{(1+\log\beta_{M}/\log x_{z})^{M+2}}-1.
\end{split}
\end{equation}

As $M \to \infty$, since $\log x_{z}=M+2$, we can write
\[
\lim_{M \to \infty}(1+\frac{\log\zeta}{\log x_{z}})^{M+2}
=\lim_{M \to \infty}(1+\frac{\log\zeta}{M+2})^{\frac{M+2}{\log\zeta}\log\zeta}
=e^{\log\zeta}=\zeta
\]
and
\[
\lim_{M \to \infty}(1+\frac{\log\beta_{M}}{\log x_{z}})^{M+2}
=\lim_{M \to \infty}(1+\frac{\log\beta_{M}}{M+2})^{\frac{M+2}{\log\beta_{M}}\log\beta_{M}}
=e^{\log\beta_{M}}=\beta_{M}.
\]
%where we make use of formulas
%\[
%lim_{x \to \infty}(1+\frac{1}{x})^{x}=e\texttt{ and }x=e^{\log x}.
%\]

Thus as $M \to \infty$, Formula~(\ref{eq:int_beta_m_x_z}) becomes
\[
\int_{1}^{\beta_{M}}(1-\frac{M+2}{\log x_{z}+\log\zeta})\frac{d\zeta}{\zeta}=0\texttt{ as }M \to \infty,
\]
so that we obtain
\[
\frac{\log x_{M}}{\log x_{z}}-1-\log\frac{\log x_{M}}{\log x_{z}}=0\texttt{ as }M \to \infty
\]
where since $\log x_{z}=M+2$
\[
\begin{split}
\int_{1}^{\beta_{M}}(1-\frac{M+2}{\log x_{z}+\log\zeta})\frac{d\zeta}{\zeta}
=&[\log\zeta-(M+2)\log\log(x_{z}\zeta)]|_{1}^{\beta_{M}} \\
=&\log\beta_{M}-(M+2)\log\log x_{M}+(M+2)\log\log x_{z} \\
%=\log x_{M}-\log x_{z}-(M+2)\log\frac{\log x_{M}}{\log x_{z}}
=&(\frac{\log x_{M}}{\log x_{z}}-1-\log\frac{\log x_{M}}{\log x_{z}})\log x_{z}.
\end{split}
\]

Therefore there must be
\[
\lim_{M \to \infty}\frac{\log x_{M}}{\log x_{z}}=1
\]
such that
\[
\lim_{M \to \infty}(\frac{\log x_{M}}{\log x_{z}}-1-\log\frac{\log x_{M}}{\log x_{z}})=0,
\]
which implies $x_{M}=x_{z}$ and $\log x_{M}=M+2$ as $M \to \infty$.

This completes the proof of the lemma.
\end{proof}

\begin{lemma}\label{le:chap_1_3_3}
Given each nonnegative integer $M$ and a real $x_{M}$ determined by $\Li(x_{M})=\pi^{*}(x_{M},M)$, we have
\[
\log\alpha_{L,M}=\frac{\log x_{M}}{M+2}+o(\frac{\log x_{M}}{M+2})\texttt{ and }
\lim_{M \to \infty}\alpha_{L,M}=\alpha_{L,\infty}=e
\]
where $\alpha_{L,M}:=x_{M+1}/x_{M}$ and $\alpha_{L,\infty}$ is a constant equal to $e$.
\end{lemma}

\begin{proof}
Let rewrite and consider Formula~(\ref{eq:int_func_x}) in the proof of Theorem~\ref{th:chap_1_3_1}
\[
\int_{2}^{x_{M}}(\frac{M+2}{\log t}-1)\frac{dt}{\log^{M+2}t}
=\frac{2}{\log^{M+2}2}-\frac{x_{M}}{\log^{M+2}x_{M}}.
\]

Let $x_{z}$ and $x_{p}$ be the zero point and the peak point of the integral function above, respectively. From Lemma~\ref{le:chap_1_3_1}, there are $\log x_{z}=M+2$ and $\log x_{p}=M+3$.

Since $\log x_{z}$ is proportional to $M+2$, $\log x_{p}$ is proportional to $M+3$, and from Lemma~\ref{le:chap_1_3_2} the number $x_{M}$ tends to $x_{z}$ as $M \to \infty$, then $\log x_{M}$ must also be approximately proportional to $M+2$.

First, let $\log x_{M}\sim\rho(M+2)$ where $\rho$ is a positive constant. Then there are
\[
\log\alpha_{L,M}=\log x_{M+1}-\log x_{M}
\sim\rho(M+1+2)-\rho(M+2)=\rho.
\]

Thus as $M \to \infty$, there must be
\[
\log\alpha_{L,M}\sim\frac{\log x_{M}}{M+2}\texttt{ and }\log\alpha_{L,M}=\frac{\log x_{M}}{M+2}+o(\frac{\log x_{M}}{M+2}).
\]

Second, as $M \to \infty$, we can write
\[
\begin{split}
\lim_{M \to \infty}(1+\frac{\log\xi}{\log x_{M}})^{M+3}
=&\lim_{M \to \infty}(1+\frac{\log\xi}{\log x_{M}})^{(M+3)(\log x_{M}/\log\xi)(\log\xi/\log x_{M})} \\
=&\lim_{M \to \infty}(1+\frac{\log\xi}{\log x_{M}})^{(\log x_{M}/\log\xi)(\log\xi/\log\alpha_{L,M})(M+3)/(M+2)} \\
=&e^{\log\xi/\log\alpha_{L,\infty}}
=\xi^{1/\log\alpha_{L,\infty}}
\end{split}
\]
and get
\[
\lim_{M \to \infty}\int_{1}^{\alpha_{L,M}}\frac{d\xi}{(1+\log\xi/\log x_{M})^{M+3}}
=\int_{1}^{\alpha_{L,\infty}}\frac{d\xi}{\xi^{1/\log\alpha_{L,\infty}}}.
\]
%where we make use of formulas
%\[
%\lim_{x \to \infty}(1+\frac{1}{x})^{x}=e\texttt{ and }x=e^{\log x}.
%\]

So as $M \to \infty$, Formula~(\ref{eq:int_alpha_m}) in the proof of Theorem~\ref{th:chap_1_3_1} becomes
\[
\int_{1}^{\alpha_{L,\infty}}\frac{d\xi}{\xi^{1/\log\alpha_{L,\infty}}}=\log\alpha_{L,\infty}\texttt{ or }
\alpha_{L,\infty}=e\log\alpha_{L,\infty}.
\]

Only when $\alpha_{L,\infty}=e$ the above equalities hold. Then we obtain
\[
\lim_{M \to \infty}\alpha_{L,M}=\alpha_{L,\infty}=e
\]
where $\alpha_{L,\infty}$ is a positive constant and must be equal to $e$.

This completes the proof of the lemma.
\end{proof}

\begin{theorem}\label{th:chap_1_3_1}
Each nonnegative integer $M$ corresponds to a real number $x_{M}$ determined by $\Li(x_{M})=\pi^{*}(x_{M},M)$. Let define $\alpha_{L,M}:=x_{M+1}/x_{M}$. Then there are
\begin{equation}\label{eq:lim_alpha_m}
\log\alpha_{L,M}=\log\alpha_{L,M+1}+o(\log\alpha_{L,M+1})\texttt{ and }
\lim_{M \to \infty}\alpha_{L,M}=\alpha_{L,\infty}=e
\end{equation}
where $\alpha_{L,\infty}$ is a constant equal to $e$.
\end{theorem}

\begin{proof}
From Lemma~\ref{le:chap_1_2_3}, given a nonnegative integer $M$, there is a real $x_{M}$ such that
\[
\Li(x_{M})=\pi^{*}(x_{M},M)\texttt{ and }
\int_{2}^{x_{M}}\frac{(M+1)!}{\log^{M+2}t}dt=\pi^{*}(2,M).
\]

By formulas
\[
\int_{2}^{x_{M}}\frac{(M+2)!}{\log^{M+3}t}dt
-\int_{2}^{x_{M}}\frac{(M+1)!}{\log^{M+2}t}dt
=\int_{2}^{x_{M}}(\frac{M+2}{\log t}-1)\frac{(M+1)!}{\log^{M+2}t}dt
\]
and
\[
\int_{2}^{x_{M}}\frac{(M+1)!}{\log^{M+2}t}dt
=\frac{(M+1)!t}{\log^{M+2}t}|_{2}^{x_{M}}
+\int_{2}^{x_{M}}\frac{(M+2)!}{\log^{M+3}t}dt,
\]
we have
\begin{equation}\label{eq:int_func_x}
\int_{2}^{x_{M}}(\frac{M+2}{\log t}-1)\frac{dt}{\log^{M+2}t}
=\frac{2}{\log^{M+2}2}-\frac{x_{M}}{\log^{M+2}x_{M}}.
\end{equation}

Since there are
\[
\pi^{*}(2,M+1)-\pi^{*}(2,M)=\frac{2(M+1)!}{\log^{M+2}2}
\]
and
\[
\int_{2}^{x_{M+1}}\frac{(M+2)!}{\log^{M+3}t}dt
-\int_{2}^{x_{M}}\frac{(M+1)!}{\log^{M+2}t}dt
=\int_{x_{M}}^{x_{M+1}}\frac{(M+2)!}{\log^{M+3}t}dt
-\frac{(M+1)!t}{\log^{M+2}t}|_{2}^{x_{M}},
\]
we have
\[
\int_{x_{M}}^{x_{M+1}}\frac{M+2}{\log^{M+3}t}dt=\frac{x_{M}}{\log^{M+2}x_{M}}.
\]

Let define $\alpha_{L,M}:=x_{M+1}/x_{M}$ and $t:=x_{M}\xi$, and write
\[
\begin{split}
\int_{x_{M}}^{x_{M+1}}\frac{M+2}{\log^{M+3}t}dt
=&\int_{x_{M}}^{\alpha_{L,M}x_{M}}\frac{M+2}{\log^{M+3}t}dt \\
=&\int_{1}^{\alpha_{L,M}}\frac{(M+2)x_{M}}{(\log x_{M}+\log\xi)^{M+3}}d\xi.
\end{split}
\]
Then we have
\begin{equation}\label{eq:int_alpha_m}
\int_{1}^{\alpha_{L,M}}\frac{d\xi}{(1+\log\xi/\log x_{M})^{M+3}}=\frac{\log x_{M}}{M+2}.
\end{equation}

Thus from Lemma~\ref{le:chap_1_3_3} we have
\[
\log\alpha_{L,M}=\frac{\log x_{M}}{M+2}+o(\frac{\log x_{M}}{M+2})\texttt{ and }
\lim_{M \to \infty}\alpha_{L,M}=\alpha_{L,\infty}=e
\]
where $\alpha_{L,\infty}$ is a constant equal to $e$.

With the definition $\alpha_{L,M}:=x_{M+1}/x_{M}$, let rewrite
\[
\log\alpha_{L,M}=\frac{\log x_{M+1}-\log\alpha_{L,M}}{M+2}+o(\frac{\log x_{M}}{M+2}).
\]
Rearranging the above equality we obtain
\[
\begin{split}
\log\alpha_{L,M}=&\frac{\log x_{M+1}}{M+3}+\frac{M+2}{M+3}o(\frac{\log x_{M}}{M+2}) \\
                =&\log\alpha_{L,M+1}+o(\frac{\log x_{M+1}}{M+3})+o(\frac{\log x_{M}}{M+3}) \\
                =&\log\alpha_{L,M+1}+o(\log\alpha_{L,M+1}).
\end{split}
\]

This completes the proof of the theorem.
\end{proof}

\section{Lemmas of the piecewise function $\pi^{*}(x,M)$}\label{sec:chap_1_4}

\begin{lemma}\label{le:chap_1_4_1}
Each nonnegative integer $M$ corresponds to a finite region $[x_{M}, x_{M+1})$ with $\alpha_{L,M}:=x_{M+1}/x_{M}$. There exists a nonnegative integer $M_{\alpha}$ such that for $M>M_{\alpha}$ there are a constant $\alpha_{L,\infty}$ equal to $e$ and
\begin{equation}\label{eq:chap_1_4_1}
\log\alpha_{L,M}=(1+o(1))\log\alpha_{L,\infty}\texttt{ and }\lim_{M \to \infty}\alpha_{L,M}=\alpha_{L,\infty}=e,
\end{equation}
and $\forall x\in[x_{M}, x_{M+1})$ with $M>M_{\alpha}$, inequalities $\pi^{*}(x,M)\leq\Li(x)<\pi^{*}(x,M+1)$ hold and there are
\begin{equation}\label{eq:rho_m_x}
M=M(x)=\frac{\log x}{\log\alpha_{L,\infty}}+O(1)=\log x+O(1).
\end{equation}
\end{lemma}

\begin{proof}
From Lemma~\ref{le:chap_1_2_2}, there is $M \to \infty$ as $x \to \infty$, and each nonnegative integer $M$ corresponds to a finite region $[x_{M}, x_{M+1})$ where $\Li(x_{M})=\pi^{*}(x_{M},M)$, and also $\forall x\in[x_{M}, x_{M+1})$ we have
\[
\pi^{*}(x,M)\leq\Li(x)<\pi^{*}(x,M+1)
\]
and
\[
\eta^{*}(x,M)\leq\eta_{L}(x)<\eta^{*}(x,M+1).
\]

Given a nonnegative integer $M$, let $\alpha_{L,M}:=x_{M+1}/x_{M}$. From Theorem~\ref{th:chap_1_3_1}, we have a constant $\alpha_{L,\infty}$ equal to $e$ and
\[
\log\alpha_{L,M}=\log\alpha_{L,M+1}+o(\log\alpha_{L,M+1})\texttt{ and }\lim_{M \to \infty}\alpha_{L,M}=\alpha_{L,\infty}=e.
\]
So there exists a nonnegative integer $M_{\alpha}$ such that for $M\geq M_{\alpha}$, there must be
\[
\log\alpha_{L,M}=(1+o(1))\log\alpha_{L,\infty}\texttt{ and }\lim_{M \to \infty}\alpha_{L,M}=\alpha_{L,\infty}=e
\]

Then for $M>M_{\alpha}$ there are
\[
\begin{split}
\log x_{M}=&\log x_{M_{\alpha}}+\sum_{k=M_{\alpha}}^{M-1}\log\alpha_{L,k} \\
%=\log x_{M_{\alpha}}+\sum_{k=M_{\alpha}}^{M-1}(1+o(1))\log\alpha_{L,\infty}
          =&\log x_{M_{\alpha}}+(M-M_{\alpha})(1+o(1))\log\alpha_{L,\infty},
\end{split}
\]
so that there must be
\[
\frac{\log(x_{M}/x_{M_{\alpha}})}{M-M_{\alpha}}=(1+o(1))\log\alpha_{L,\infty}.
\]

Then we obtain
\[
\begin{split}
M-M_{\alpha}=&\frac{\log x_{M}-\log x_{M_{\alpha}}}{(1+o(1))\log\alpha_{L,\infty}} \\
            =&\frac{\log x_{M}-\log x_{M_{\alpha}}}{\log\alpha_{L,\infty}}(1+o(1))
\end{split}
\]
and
\[
M=\frac{\log x_{M}}{\log\alpha_{L,\infty}}(1+o(1))-\frac{\log x_{M_{\alpha}}}{\log\alpha_{L,\infty}}(1+o(1))+M_{\alpha}.
\]

Since each nonnegative integer $M$ corresponds to a finite region $[x_{M},x_{M+1})$ of the variable $x$, $\forall x\in[x_{M}, x_{M+1})$ there are
\[
M\leq\frac{\log x}{\log\alpha_{L,\infty}}-\frac{\log x_{M_{\alpha}}}{\log\alpha_{L,\infty}}+M_{\alpha}
+\frac{\log x-\log x_{M_{\alpha}}}{\log\alpha_{L,\infty}}o(1)<M+1.
\]

By the analysis above, if $M \to \infty$ then there must be $1/\log x \to 0$ so that there are $1/\log x=o(1)$ and $o(\log x)=O(1)$. Thus for $x\geq x_{M_{\alpha}+1}$ there must be
\begin{equation}\label{eq:m_alpha_k}
\begin{split}
M=M(x)=&\frac{\log x}{\log\alpha_{L,\infty}}-\frac{\log x_{M_{\alpha}}}{\log\alpha_{L,\infty}}+M_{\alpha}+O(1) \\
      =&\frac{\log x}{\log\alpha_{L,\infty}}+O(1).
\end{split}
\end{equation}

This completes the proof of the lemma.
\end{proof}

\begin{lemma}\label{le:chap_1_4_2}
Let define
\[
\rho_{L}(x,M):=\frac{\log x}{M+2}\texttt{ where }\pi^{*}(x,M)\leq\Li(x)<\pi^{*}(x,M+1).
\]
Then there are
\begin{equation}\label{eq:rho_x_m_eq}
\begin{split}
\rho_{L}(x,M)=&(1+O(\frac{1}{M+2}))\log\alpha_{L,\infty} \\
             =&1+O(\frac{1}{M+2})\texttt{ for }x\geq198
\end{split}
\end{equation}
and
\begin{equation}\label{eq:rho_x_m_inf}
\lim_{M \to \infty}\rho_{L}(x,M)=\log\alpha_{L,\infty}=1.
\end{equation}
\end{lemma}

\begin{proof}
From Lemma~\ref{le:chap_1_2_2} and Lemma~\ref{le:chap_1_4_1} $\rho_{L}(x,M)$ must be a piecewise increasing function of the variable $\log x$ and must tend to a limiting value as $\log x \to \infty$.

First, Formula~(\ref{eq:rho_m_x}) can be rewritten as the following formula
\[
\frac{\log x}{M+2}=(1+O(\frac{1}{M+2}))\log\alpha_{L,\infty}.
\]
Then Formula~(\ref{eq:rho_x_m_eq}) holds for $M>M_{\alpha}$, i.e.,
\[
\rho_{L}(x,M):=\frac{\log x}{M+2}=(1+O(\frac{1}{M+2}))\log\alpha_{L,\infty}.
\]

Second, as $M \to \infty$, Formula~(\ref{eq:rho_x_m_eq}) becomes
\[
\lim_{M \to \infty}\rho_{L}(x,M)=(1+O(\lim_{M \to \infty}\frac{1}{M+2}))\log\alpha_{L,\infty}=\log\alpha_{L,\infty}=1.
\]

This completes the proof of the lemma.
\end{proof}

\section{Theorem of upper bounds of truncation gaps}\label{sec:chap_1_5}

\begin{theorem}\label{th:chap_1_5_1}
Given each pair of numbers $x\geq730$ and $M$ determined by $\pi^{*}(x,M)\leq\Li(x)<\pi^{*}(x,M+1)$, there are
\begin{equation}\label{eq:chap_1_5_1}
g(x,M):=\pi^{*}(x,M+1)-\pi^{*}(x,M)<x^{1/64}.
\end{equation}
\end{theorem}

\begin{proof}
From Lemma~\ref{le:chap_1_5_1} and Lemma~\ref{le:chap_1_3_2} we have
\[
\frac{\log^{M+2}x}{(M+1)!}>(\frac{e\log x}{M+2})^{M+2}\texttt{ for }M\geq5\texttt{ or }x\geq2250.
\]

As defined by Lemma~\ref{le:chap_1_4_2}, let $\rho_{L}(x,M):=\frac{\log x}{M+2}$. Then we consider the ratio of
\[
\frac{\log x}{x^{63/64}}\texttt{ to }\frac{(M+1)!}{\log^{M+1}x}
\]
and have
\[
\begin{split}
\log\frac{\log^{M+2}x}{(M+1)!x^{63/64}}
>&(M+2)\log(\frac{e\log x}{M+2})-\frac{63}{64}\log x \\
=&(M+2)[1+\log\rho_{L}(x,M)-(63/64)\rho_{L}(x,M)].
\end{split}
\]

Let $\rho_{0}(x,M)$ satisfy the equation $1+\log\rho_{L}(x,M)-(63/64)\rho_{L}(x,M)=0$. By calculation, we have $1.20698278<\rho_{0}(x,M)<1.20698279$.

Since the key discriminant
\[
1+\log\rho_{L}(x,M)-(63/64)\rho_{L}(x,M)\geq0
\]
holds for $\rho_{L}(x,M)\leq\rho_{0}(x,M)$, thus from Lemma~\ref{le:chap_1_4_2} for $x\geq198$ the inequality $\rho_{L}(x,M)<\rho_{0}(x,M)$ holds so that for $M\geq5$ or $x\geq2250$ we obtain
\[
\frac{\log^{M+2}x}{(M+1)!}>x^{63/64}\texttt{ and }
\frac{(M+1)!}{\log^{(M+1)}x}<\frac{\log x}{x^{63/64}}.
\]

By calculation, the inequality $\frac{\log^{M+2}x}{(M+1)!}>x^{63/64}$ holds for $2250>x\geq730$.

From Lemma~\ref{le:chap_1_2_1}, given each real number $x\geq198$, there exists a nonnegative integer $M=M(x)$ such that $\pi^{*}(x,M)\leq\Li(x)<\pi^{*}(x,M+1)$. Let define
\[
g(x,M):=\pi^{*}(x,M+1)-\pi^{*}(x,M)=\frac{x}{\log x}\frac{(M+1)!}{\log^{(M+1)}x}.
\]

Then for each positive number $x\geq730$, we have
\[
\frac{g(x,M)}{x/\log x}=\frac{(M+1)!}{\log^{(M+1)}x}<\frac{\log x}{x^{63/64}}\texttt{ and }
g(x,M)<\frac{\log x}{x^{63/64}}\frac{x}{\log x}=x^{1/64}.
\]

This completes the proof of the theorem.
\end{proof}

\begin{lemma}\label{le:chap_1_5_1}
\[
(M+1)!<(\frac{M+2}{e})^{M+2}\texttt{ for }M\geq5.
\]
\end{lemma}

\begin{proof}
Based upon the Stirling's formula
\[
n!=\sqrt{2\pi n}(\frac{n}{e})^{n}e^{\frac{\theta}{12n}},0<\theta<1,
\]
for $M\geq45$, there are
\[
\sqrt{2\pi}e^{\frac{1}{2}+\frac{1}{12(M+1)}}<(\frac{M+2}{e})^{1/2}
\]
and
\[
(M+1)!<\sqrt{2\pi}e^{\frac{1}{2}+\frac{1}{12(M+1)}}(\frac{M+1}{e})^{M+3/2}<(\frac{M+2}{e})^{M+2}.
\]

By calculation, the inequality $(M+1)!<(\frac{M+2}{e})^{M+2}$ holds for $45>M\geq5$.

This completes the proof of the lemma.
\end{proof}

\section{Lemmas of bounds of the prime counting function $\pi(x)$}\label{sec:chap_1_6}

\begin{lemma}\label{le:chap_1_6_1}
Given each $x\geq599$, there exists a nonnegative integer function $N_{x}=N(x)$ such that the prime counting function $\pi(x)$ satisfies
\begin{equation}\label{eq:chap_1_6_1}
\pi^{*}(x,N_{x})<\pi(x)<\pi^{*}(x,N_{x}+1).
\end{equation}
\end{lemma}

\begin{proof}
Let $\eta(x):=\frac{\pi(x)}{x/\log x}$. Then for each nonnegative integer $N$, we can write
\begin{equation}\label{eq:pi_n}
\eta(x)=\eta^{*}(x,N)+\delta(x,N)=\sum_{n=0}^{N}\frac{n!}{\log^{n}x}+\delta(x,N).
\end{equation}

First, when $N=0$, by inequalities~(\ref{eq:pi_x_2}-\ref{eq:pi_x_3}) there must be $\delta(x,0)>0$.

Second, from Equality~(\ref{eq:pi_n}) for each nonnegative integer $N$ the inequality
\[
\delta(x,N)-\delta(x,N+1)=\frac{(N+1)!}{\log^{N+1}x}>0
\]
always holds so that given a fixed positive number $x$, there must be
\[
\lim_{N \to \infty}\delta(x,N) \to -\infty.
\]

Then for each number $x\geq599$, there always exists a nonnegative integer $N_{x}=N(x)$ such that there are $\delta(x,N_{x})\geq0$ and $\delta(x,N_{x}+1)<0$. Thus there exist inequalities
\[
\eta^{*}(x,N_{x})\leq\eta(x)<\eta^{*}(x,N_{x}+1)
\]
such that $\pi(x)$ is bounded by inequalities
\[
\pi^{*}(x,N_{x})\leq\pi(x)<\pi^{*}(x,N_{x}+1).
\]

On the other hand, when $x$ is a positive integer, the value of $\pi^{*}(x,N_{x})$ must be an irrational number and can not be equal to a positive integer $\pi(x)$. So $\pi(x)$ must be bounded by inequalities
\[
\pi^{*}(x,N_{x})<\pi(x)<\pi^{*}(x,N_{x}+1)
\]
and there must be
\[
\eta^{*}(x,N_{x})<\eta(x)<\eta^{*}(x,N_{x}+1).
\]

This completes the proof of the lemma.
\end{proof}

We simply replace $N_{x}$ with $N$ in the following lemma.
\begin{lemma}\label{le:chap_1_6_2}
Given two numbers $x$ and $N$ determined by
\[
\pi^{*}(x,N)<\pi(x)<\pi^{*}(x,N+1),
\]
if $\pi(x)\leq\Li(x)$ then $N+2<\log x$ else $\log x\leq N+2\leq c_{n}\log x$ where $c_{n}$ is a constant greater than $1$ and less than $e$.
\end{lemma}

\begin{proof}
From Lemma~\ref{le:chap_1_3_2}, given two numbers $x_{M}$ and $M$ determined by inequalities $\pi^{*}(x_{M},M)\leq\Li(x_{M})<\pi^{*}(x_{M},M+1)$, there are $\log x_{M}>M+2$ for $M<\infty$ and $\log x_{M}=M+2$ as $M \to \infty$.

From Lemma~\ref{le:chap_1_2_2}, given a nonnegative integer $M$ and $\forall x\in[x_{M}, x_{M+1})$, all pairs of a positive number $x$ and the same nonnegative integer $M$ are determined by $\pi^{*}(x,M)\leq\Li(x)<\pi^{*}(x,M+1)$. Then $\forall x\in[x_{M}, x_{M+1})$ the inequality $M+2<\log x$ holds.

From Lemma~\ref{le:chap_1_2_1} and Lemma~\ref{le:chap_1_6_1}, given each number $x\geq10^{3}$, there always exist a unique pair of nonnegative integers $M=M(x)$ and $N=N(x)$ such that
\[
\pi^{*}(x,M)\leq\Li(x)<\pi^{*}(x,M+1)\texttt{ and }\pi^{*}(x,N)<\pi(x)<\pi^{*}(x,N+1).
\]

First, for all $x$ satisfying $\pi(x)\leq\Li(x)$, there must be
\[
\pi^{*}(x,N)<\pi(x)\leq\Li(x)<\pi^{*}(x,M+1),
\]
i.e., there must be $N\leq M$ and $N+2\leq M+2<\log x$.

Second, for all $x$ satisfying $\pi(x)>\Li(x)$, there must be
\[
\pi^{*}(x,M)\leq\Li(x)<\pi(x)<\pi^{*}(x,N+1),
\]
i.e., there must be $M<N+1$ and $N\geq M=\log x+O(1)$ from Lemma~\ref{le:chap_1_4_2}.

On the other hand, based upon the Stirling's formula
\[
n!=\sqrt{2\pi n}(\frac{n}{e})^{n}e^{\frac{\theta}{12n}},0<\theta<1,
\]
when $N+2\geq e\log x$, there shall be
\[
\begin{split}
\frac{N!}{\log^{N}x}>&\sqrt{2\pi N}(\frac{N}{e\log x})^{N}\geq\sqrt{2\pi N}(\frac{N}{N+2})^{N} \\
                    =&\sqrt{2\pi N}(1+\frac{2}{N})^{-N}\geq\sqrt{2\pi N}e^{-2}
\end{split}
\]
where we make use of the formula $\lim_{x \to \infty}(1+1/x)^{x}=e$, and
\[
\begin{split}
\eta^{*}(x,N)=&\eta^{*}(x,N-1)+\frac{N!}{\log^{N}x} \\
             >&\eta^{*}(x,N-1)+\sqrt{2\pi N}e^{-2}>2,
\end{split}
\]
since the upper bound for the first sign change of $\pi(x)-li(x)$ was improved to $1.6\cdot10^{1165}$~\cite{AC,AM} such that $\log x>1165\log10$ and $N>e^{4}/(2\pi)$.

But this result conflicts the fact that there is $\eta^{*}(x,N)<\eta(x)<2$ for $x\geq59$ by Formula~(\ref{eq:pi_x_2}). Thus there must be $N+2<e\log x$ when a finite number $x$ and a pair of nonnegative integers $M$ and $N\geq M$ satisfy
\[
\pi^{*}(x,M)\leq\Li(x)<\pi^{*}(x,M+1)\texttt{ and }\pi^{*}(x,N)<\pi(x)<\pi^{*}(x,N+1)
\]
where $M=M(x)=\log x+O(1)$ from Lemma~\ref{le:chap_1_4_2}. Therefore we can obtain
\[
M+2\leq N+2\leq c_{n}\log x\texttt{ for }x\geq59\texttt{ where }1<c_{n}<e.
\]

So given a pair of numbers $x$ and $N$ satisfying $\pi^{*}(x,N)<\pi(x)<\pi^{*}(x,N+1)$, if $\pi(x)\leq\Li(x)$ then $N+2<\log x$ else $\log x\leq N+2\leq c_{n}\log x$ where $c_{n}$ is a constant greater than $1$ and less than $e$.

This completes the proof of the lemma.
\end{proof}

\begin{lemma}\label{le:chap_1_6_3}
Let $\lambda(x)$ depend on the variable $x$ and $\lambda(x) \to \infty$ as $x \to \infty$. Then
\begin{equation}\label{eq:pi_x_big_o_1}
\frac{\pi(x)}{x/\log x}=\frac{\Li(x)}{x/\log x}+O(\frac{1}{\log^{\lambda(x)}x}).
\end{equation}
\end{lemma}

\begin{proof}
Let $\lambda(x)=a\sqrt{\log x}/\log\log x-1$. Then Formula~(\ref{eq:pi_x_4}) becomes
\[
\pi(x)=\Li(x)+O(x/\log^{\lambda(x)+1}x)
\]
where $\lambda(x) \to \infty$ as $x \to \infty$. The closeness of $\Li(x)$ and $\pi(x)$ determined by Formula~(\ref{eq:pi_x_4}) depends on the power $\lambda(x)$ of $\log x$ and is bounded by $O(x/\log^{\lambda(x)+1}x)$.

Similarly, for large numbers $x$, formulas~(\ref{eq:pi_x_5}-\ref{eq:pi_x_7}) become
\[
\pi(x)=\Li(x)+O(x/\log^{\lambda(x)+1}x)
\]
where by Formula~(\ref{eq:pi_x_5})
\[
\lambda(x)=a\sqrt{\log x}/\sqrt{\log\log x}-1 \to \infty\texttt{ as }x \to \infty
\]
or by Formula~(\ref{eq:pi_x_6})
\[
\lambda(x)=a(\log x)^{3/5}/(\log\log x)^{6/5}-1 \to \infty\texttt{ as }x \to \infty
\]
or by Formula~(\ref{eq:pi_x_7})
\[
\lambda(x)=A(x)-1 \to \infty\texttt{ as }x \to \infty
\]
where $A(x)$ implies that each integer $A>1$ corresponds to a sufficiently large $x$.

The closeness of $\Li(x)$ and $\pi(x)$ determined by formulas~(\ref{eq:pi_x_5}-\ref{eq:pi_x_7}) depends on the power $\lambda(x)$ of $\log x$ and is also bounded by $O(x/\log^{\lambda(x)+1}x)$.

This completes the proof of the lemma.
\end{proof}

We simply replace $N_{x}$ with $N$ in the following lemma.
\begin{lemma}\label{le:chap_1_6_4}
When a pair of numbers $x$ and $N$ are determined by $\pi^{*}(x,N)<\pi(x)<\pi^{*}(x,N+1)$, there is $N \to \infty$ as $x \to \infty$.
\end{lemma}

\begin{proof}
From Lemma~\ref{le:chap_1_6_3}, the closeness of $\Li(x)$ and $\pi(x)$ depends on the power $\lambda(x)$ of $\log x$ and $\lambda(x) \to \infty$ as $x \to \infty$.

From Lemma~\ref{le:chap_1_2_1} and Lemma~\ref{le:chap_1_6_1}, given each number $x\geq10^{3}$, there always exist a unique pair of nonnegative integers $M=M(x)$ and $N=N(x)$ such that
\[
\pi^{*}(x,M)\leq\Li(x)<\pi^{*}(x,M+1)\texttt{ and }\pi^{*}(x,N)<\pi(x)<\pi^{*}(x,N+1).
\]

If $\pi(x)\leq\Li(x)$ then from Lemma~\ref{le:chap_1_6_2} there are $N+2\leq M+2<\log x$ and
\begin{equation}\label{eq:pi_x_big_o_2}
\begin{split}
\frac{\Li(x)}{x/\log x}-\frac{\pi(x)}{x/\log x}<&\eta^{*}(x,M+1)-\eta^{*}(x,N)=\sum_{n=N+1}^{M+1}\frac{n!}{\log^{n}x} \\
\leq&(M+1-N)\frac{(N+1)!}{\log^{N+1}x}.
\end{split}
\end{equation}

If $\pi(x)>\Li(x)$ then from Lemma~\ref{le:chap_1_6_2} there are $M+2<N+2<e\log x$ and
\begin{equation}\label{eq:pi_x_big_o_3}
\begin{split}
\frac{\pi(x)}{x/\log x}-\frac{\Li(x)}{x/\log x}<&\eta^{*}(x,N+1)-\eta^{*}(x,M)=\sum_{n=M+1}^{N+1}\frac{n!}{\log^{n}x} \\
\leq&(N+1-M)\max\{\frac{(M+1)!}{\log^{M+1}x},\frac{(N+1)!}{\log^{N+1}x}\}.
\end{split}
\end{equation}

If the difference $|M-N|$ is finite, then from Lemma~\ref{le:chap_1_2_2} there are $M \to \infty$ and $N \to \infty$ as $x \to \infty$.

Otherwise, by comparing Formula~(\ref{eq:pi_x_big_o_1}) with Formula~(\ref{eq:pi_x_big_o_2}), there must be $N\geq\lambda(x)$, and by comparing Formula~(\ref{eq:pi_x_big_o_1}) with Formula~(\ref{eq:pi_x_big_o_3}), there must be $N\geq\lambda(x)$ or $N>M\geq\lambda(x)$. Hence there is $N \to \infty$ as $\lambda(x) \to \infty$ or $x \to \infty$.

This completes the proof of the lemma.
\end{proof}

\section{Lemmas of the piecewise function $\pi^{*}(x,N)$}\label{sec:chap_1_7}

\begin{lemma}\label{le:chap_1_7_1}
Each nonnegative integer $N$ corresponds to a function $\pi^{*}(x,N)$ valid for a finite region $[x_{N},x_{N+1})$, and $\forall x\in[x_{N},x_{N+1})$ there exists a nonnegative integer $N_{x}$ such that the pair of numbers $x$ and $N_{x}$ satisfy
\[
\pi^{*}(x,N_{x})<\pi(x)<\pi^{*}(x,N_{x}+1)
\]
and there is always $N_{x}\geq N$.
\end{lemma}

\begin{proof}
First, for $N=0,1$, by Formula~(\ref{eq:pi_x_3}), we can set $x_{N}=10^{2N+3}$ and when each pair of numbers $x\geq x_{N}$ and $N_{x}$ satisfy $\pi^{*}(x,N_{x})<\pi(x)<\pi^{*}(x,N_{x}+1)$ there is always $N_{x}\geq N$.

Second, suppose that given a nonnegative integer $N$, there exists a positive integer $x_{N}$ such that when each pair of numbers $x\geq x_{N}$ and $N_{x}$ are determined by $\pi^{*}(x,N_{x})<\pi(x)<\pi^{*}(x,N_{x}+1)$ there is always $N_{x}\geq N$. Then $N_{x}$ must tend to $\infty$ as $x \to \infty$ from Lemma~\ref{le:chap_1_6_4}. So there must exist a positive integer $x_{N+1}$ such that when each pair of numbers $x\geq x_{N+1}$ and $N_{x}$ satisfy
\[
\pi^{*}(x,N_{x})<\pi(x)<\pi^{*}(x,N_{x}+1),
\]
there is always $N_{x}\geq N+1$.

Hence each nonnegative integer $N$ corresponds to a function $\pi^{*}(x,N)$ valid for a finite region $[x_{N},x_{N+1})$, and $\forall x\in[x_{N},x_{N+1})$ there exists a nonnegative integer $N_{x}$ such that the pair of numbers $x$ and $N_{x}$ satisfy
\[
\pi^{*}(x,N_{x})<\pi(x)<\pi^{*}(x,N_{x}+1)
\]
and there is always $N_{x}\geq N$.

This completes the proof of the lemma.
\end{proof}

\begin{lemma}\label{le:chap_1_7_2}
Each nonnegative integer $N$ in the function $\pi^{*}(x,N)$ corresponds to a finite region $[x_{N}, x_{N+1})$ and a real number $\alpha_{N}:=x_{N+1}/x_{N}$. Then there exists a nonnegative integer $N_{\alpha}$ such that given each nonnegative integer $N\geq N_{\alpha}$ we have
\begin{equation}\label{eq:alpha_n}
\log\alpha_{N}=\log\alpha_{\infty}+O(\frac{1}{\log x_{N}})\texttt{ and }\lim_{N \to \infty}\alpha_{N}=\alpha_{\infty}
\end{equation}
where $\alpha_{\infty}$ is a positive constant.
\end{lemma}

\begin{proof}
From Lemma~\ref{le:chap_1_7_1}, each nonnegative integer $N$ corresponds to a function $\pi^{*}(x,N)$ valid for a finite region $[x_{N},x_{N+1})$, and $\forall x\in[x_{N},x_{N+1})$ there always exists a nonnegative integer $N_{x}\geq N$ such that
\[
\pi^{*}(x,N_{x})<\pi(x)<\pi^{*}(x,N_{x}+1)
\]
and
\[
\eta^{*}(x,N_{x})<\eta(x)<\eta^{*}(x,N_{x}+1).
\]

Given each fixed nonnegative integer $N$, let $y:=\log x$ and $x$ be a continuous real variable. Then $\eta^{*}(x,N)$ becomes a monotonically decreasing smooth continuous function $\eta^{*}(y,N)$ of the variable $y=\log x$. The second derivative of $\eta^{*}(y,N)$ is greater than zero for $y<\infty$ and tends to zero as $y \to \infty$.

With a difference $\Delta N=1$, since $N \to \infty$ as $y \to \infty$, the second difference of $\eta^{*}(y,N)$ on the variable $N$ is equivalent to
\[
\frac{\Delta^{2}\eta^{*}(y,N)}{(\Delta N)^{2}}\sim\frac{d^{2}\eta^{*}(y,N)}{(dy)^{2}}(\frac{\Delta y}{\Delta N})^{2}\geq0\texttt{ as }N \to \infty.
\]
So the second difference of $\eta^{*}(y,N)$ on the variable $N$ is also greater than zero for $y<\infty$ and tends to zero as $y \to \infty$.

Thus the first difference $\Delta\eta^{*}(y,N)/\Delta N$ of $\eta^{*}(y,N)$ on the variable $N$ must be approximately equal to a constant and tend to zero as $N \to \infty$ in a finite local region $[\log x_{N}, \log x_{N+2}]$ for the variable $y=\log x$ with a difference $\Delta N=1$ where the ratio of $\Delta N$ to $N$ is $1/N$ and tends to zero as $N \to \infty$. Hence there must be
\begin{equation}\label{eq:diff_x_n}
\frac{\eta^{*}(x_{N},N)-\eta^{*}(x_{N+1},N+1)}{\eta^{*}(x_{N+1},N+1)-\eta^{*}(x_{N+2},N+2)}>1\texttt{ for }N<\infty
\end{equation}
and
\begin{equation}\label{eq:diff_lim}
\lim_{N \to \infty}\frac{\eta^{*}(x_{N},N)-\eta^{*}(x_{N+1},N+1)}{\eta^{*}(x_{N+1},N+1)-\eta^{*}(x_{N+2},N+2)}=1.
\end{equation}

Then, for $N\geq2$, there are
\[
\begin{split}
\eta^{*}(x_{N},N)-\eta^{*}(x_{N+1},N+1)
=&\sum_{n=1}^{N}\frac{n!}{\log^{n}x_{N}}-\sum_{n=1}^{N+1}\frac{n!}{\log^{n}x_{N+1}} \\
=&A_{1}(N)\frac{\log x_{N+1}-\log x_{N}}{\log x_{N}\log x_{N+1}}-\frac{(N+1)!}{\log^{N+1}x_{N+1}}
\end{split}
\]
and
\[
\eta^{*}(x_{N+1},N+1)-\eta^{*}(x_{N+2},N+2)
=\sum_{n=1}^{N+1}\frac{n!}{\log^{n}x_{N+1}}-\sum_{n=1}^{N+2}\frac{n!}{\log^{n}x_{N+2}}
\]
\[
=A_{1}(N+1)\frac{\log x_{N+2}-\log x_{N+1}}{\log x_{N+1}\log x_{N+2}}
-\frac{(N+2)!}{\log^{N+2}x_{N+2}}
\]
where
\[
A_{1}(N)=\sum_{n=1}^{N}\sum_{k=0}^{n-1}\frac{n!}{\log^{k}x_{N}\log^{n-1-k}x_{N+1}}
\]
and for $N<\infty$
\[
\begin{split}
A_{1}(N)-A_{1}(N+1)=&\sum_{n=2}^{N}\frac{n!}{\log^{n-1}x_{N+1}}\sum_{k=0}^{n-1}(\frac{\log^{k}x_{N+1}}{\log^{k}x_{N}}-\frac{\log^{k}x_{N+1}}{\log^{k}x_{N+2}})\\
                    &-\frac{(N+1)!}{\log^{N}x_{N+1}}\sum_{k=0}^{N}\frac{\log^{k}x_{N+1}}{\log^{k}x_{N+2}}>0\texttt{ or }\frac{A_{1}(N+1)}{A_{1}(N)}<1.
\end{split}
\]

Let define $\alpha_{N}:=x_{N+1}/x_{N}$ for $N>1$. Then formulas~(\ref{eq:diff_x_n}-\ref{eq:diff_lim}) become
\[
\frac{\log\alpha_{N}}{\log\alpha_{N+1}}>\frac{A_{1}(N+1)}{A_{1}(N)}\frac{\log x_{N}}{\log x_{N+2}}\texttt{ for }N<\infty
\]
and
\[
\lim_{N \to \infty}\frac{\log\alpha_{N}}{\log\alpha_{N+1}}=1.
\]

Since the second difference of $\eta^{*}(y,N)$ on the variable $N$ is greater than zero for $y<\infty$ and tends to zero as $y \to \infty$, and since there are
\[
\frac{A_{1}(N+1)}{A_{1}(N)}\frac{\log x_{N}}{\log x_{N+2}}<1\texttt{ for }N<\infty
\]
and
\[
\lim_{N \to \infty}\frac{A_{1}(N+1)}{A_{1}(N)}\frac{\log x_{N}}{\log x_{N+2}}=1,
\]
there must be $\log\alpha_{N}=\log\alpha_{N+1}+O(1/\log x_{N})$ for $N>1$.

Thus there exists a nonnegative integer $N_{\alpha}$ such that given each integer $N\geq N_{\alpha}$, there must be
\[
\log\alpha_{N}=\log\alpha_{\infty}+O(\frac{1}{\log x_{N}})
\]
and
\[
\lim_{N \to \infty}\alpha_{N}=\alpha_{\infty}
\]
where $\alpha_{\infty}$ is a positive constant that can be determined by Lemma~\ref{le:chap_1_7_5}.

This completes the proof of the lemma.
\end{proof}

\begin{lemma}\label{le:chap_1_7_3}
The nonnegative integer $N$ in the piecewise function $\pi^{*}(x,N)$ is a non-decreasing step function of the variable $\log x$ and there exists a nonnegative integer $N_{\alpha}$ such that for $N\geq N_{\alpha}$ and $\forall x\in[x_{N},x_{N+1})$ there are
\begin{equation}\label{eq:n_alpha}
N=\frac{\log x}{\log\alpha_{\infty}}+O(1)\texttt{ and }x=(\alpha_{\infty}^{1/2})^{2N+3+O(1)}.
\end{equation}
\end{lemma}

\begin{proof}
From Lemma~\ref{le:chap_1_7_2}, each nonnegative integer $N$ in the function $\pi^{*}(x,N)$ corresponds to a finite region $[x_{N}, x_{N+1})$ and a real number $\alpha_{N}:=x_{N+1}/x_{N}$. Then there exists a nonnegative integer $N_{\alpha}$ such that given each nonnegative integer $N\geq N_{\alpha}$, we have
\[
\log\alpha_{N}=\log\alpha_{\infty}+O(\frac{1}{\log x_{N}})\texttt{ and }\lim_{N \to \infty}\alpha_{N}=\alpha_{\infty}
\]
where $\alpha_{\infty}$ is a positive constant.

Then given each nonnegative integer $N\geq N_{\alpha}$, there are
\[
\begin{split}
\log x_{N}=&\log x_{N_{\alpha}}+\sum_{k=N_{\alpha}}^{N-1}\log\alpha_{k} \\
          =&\log x_{N_{\alpha}}+\sum_{k=N_{\alpha}}^{N-1}(\log\alpha_{\infty}+O(\frac{1}{\log x_{k}})) \\
          =&\log x_{N_{\alpha}}+(N-N_{\alpha})(\log\alpha_{\infty}+o(1))
\end{split}
\]
where since there are $\log x_{N}=o(x_{N})$ as $N \to \infty$ and
\[
\begin{split}
\sum_{k=N_{\alpha}}^{N-1}\frac{1}{\log x_{k}}\sim&\sum_{k=N_{\alpha}}^{N-1}\int_{\log x_{k}}^{\log x_{k+1}}\frac{dt}{t}=\int_{\log x_{N_{\alpha}}}^{\log x_{N}}\frac{dt}{t} \\
=&\log\frac{\log x_{N}}{\log x_{N_{\alpha}}}=o(\log\frac{x_{N}}{x_{N_{\alpha}}}),
\end{split}
\]
and the ratio of the real $\log x_{N}$ to the number $N$ must be finite, we can also write
\[
\begin{split}
\frac{\log(x_{N}/x_{N_{\alpha}})}{N-N_{\alpha}}-\log\alpha_{\infty}
=&O(\frac{1}{N-N_{\alpha}}\sum_{k=N_{\alpha}}^{N-1}\frac{1}{\log x_{k}}) \\
=&o(\frac{\log(x_{N}/x_{N_{\alpha}})}{N-N_{\alpha}})=o(1),
\end{split}
\]
so that there must be
\[
\frac{\log(x_{N}/x_{N_{\alpha}})}{N-N_{\alpha}}=\log\alpha_{\infty}+o(1).
\]

Then we obtain
\[
N-N_{\alpha}=\frac{\log x_{N}-\log x_{N_{\alpha}}}{\log\alpha_{\infty}+o(1)}=\frac{\log x_{N}-\log x_{N_{\alpha}}}{\log\alpha_{\infty}}(1+o(1))
\]
and
\[
N=\frac{\log x_{N}}{\log\alpha_{\infty}}(1+o(1))-\frac{\log x_{N_{\alpha}}}{\log\alpha_{\infty}}(1+o(1))+N_{\alpha}.
\]

Since each nonnegative integer $N$ corresponds to a finite region $[x_{N},x_{N+1})$ of the variable $x$, $\forall x\in[x_{N},x_{N+1})$ there are
\[
N\leq\frac{\log x}{\log\alpha_{\infty}}-\frac{\log x_{N_{\alpha}}}{\log\alpha_{\infty}}+N_{\alpha}
+\frac{\log x-\log x_{N_{\alpha}}}{\log\alpha_{\infty}}o(1)<N+1.
\]

By the analysis above, if $N \to \infty$ then there must be $1/\log x \to 0$ such that there are $1/\log x=o(1)$ and $o(\log x)=O(1)$. Thus $\forall x\in[x_{N},x_{N+1})$ there are
\begin{equation}\label{eq:n_alpha_k}
N=\frac{\log x}{\log\alpha_{\infty}}-\frac{\log x_{N_{\alpha}}}{\log\alpha_{\infty}}+N_{\alpha}+O(1)
=\frac{\log x}{\log\alpha_{\infty}}+O(1).
\end{equation}
It implies that the nonnegative integer $N$ is a non-decreasing step function of the variable $\log x$.

On the other hand, by the definition $\rho(x,N):=\frac{2\log x}{2N+3}$ in Lemma~\ref{le:chap_1_7_7}, we can write Formula~(\ref{eq:n_alpha_k}) as $\forall x\in[x_{N},x_{N+1})$
\[
\log x=(N+3/2)\log\alpha_{\infty}+O(1)\texttt{ or }x=(\alpha_{\infty}^{1/2})^{2N+3+O(1)}.
\]

This completes the proof of the lemma.
\end{proof}

\begin{lemma}\label{le:chap_1_7_4}
There exists a positive integer $x_{0}$ such that each nonnegative integer $N$ corresponds to a finite region $[x_{N},x_{N+1})$ with $\alpha_{N}:=x_{N+1}/x_{N}=\alpha_{\infty}$ for $N=0,1,2,\ldots$ and $\forall x\in[x_{N},x_{N+1})$ there exists a nonnegative integer $N_{x}$ such that
\[
\pi^{*}(x,N_{x})<\pi(x)<\pi^{*}(x,N_{x}+1)\texttt{ and }N_{x}\geq N.
\]
\end{lemma}

\begin{proof}
From Lemma~\ref{le:chap_1_7_1}, each nonnegative integer $N$ corresponds to a function $\pi^{*}(x,N)$ valid for a finite region $[x_{N},x_{N+1})$ such that  $\forall x\in[x_{N},x_{N+1})$ there are
\[
\pi^{*}(x,N_{x})<\pi(x)<\pi^{*}(x,N_{x}+1)\texttt{ and }N_{x}\geq N.
\]
On the other hand, from Lemma~\ref{le:chap_1_7_1}, there are $N_{x}\geq N$ for $x\geq x_{N}$ and $N_{x}\geq N+1$ for $x\geq x_{N+1}$.

From Lemma~\ref{le:chap_1_7_2} for $N\geq N_{\alpha}$, there are
\[
\log\alpha_{N}=\log\alpha_{\infty}+O(\frac{1}{\log x_{N}})\texttt{ and }\lim_{N \to \infty}\alpha_{N}=\alpha_{\infty}.
\]

Then by the proof of Lemma~\ref{le:chap_1_7_2} there must be
\begin{equation}\label{eq:alpha_n_lt}
\alpha_{N}\leq\alpha_{N+1}\leq\alpha_{\infty}\texttt{ for }N<\infty\texttt{ and }\lim_{N \to \infty}\alpha_{N}=\alpha_{\infty}
\end{equation}
or
\begin{equation}\label{eq:alpha_n_gt}
\alpha_{N}\geq\alpha_{N+1}\geq\alpha_{\infty}\texttt{ for }N<\infty\texttt{ and }\lim_{N \to \infty}\alpha_{N}=\alpha_{\infty}.
\end{equation}

First, in case of Formula~\ref{eq:alpha_n_lt}, if $\alpha_{N}<\alpha_{\infty}$ then $\forall x\in[x_{N},x_{N+1})$ there are $\alpha_{\infty}x>\alpha_{N}x_{N}=x_{N+1}$ for $N<\infty$.

Let $x_{N+1}^{*}=\alpha_{\infty}x_{N}$. Then there are $x_{N+1}=\alpha_{N}x_{N}<\alpha_{\infty}x_{N}=x_{N+1}^{*}$ since $\alpha_{N}<\alpha_{\infty}$. We can enlarge the old region $[x_{N},x_{N+1})$ to a new region $[x_{N},x_{N+1}^{*})$, so that we have a new $\alpha_{N}^{*}=x_{N+1}^{*}/x_{N}=\alpha_{\infty}$.

Thus for each nonnegative integer $N=0,1,2,\ldots$, if $\alpha_{N}<\alpha_{\infty}$ then given each $x$ in the new enlarged region $[x_{N},x_{N+1}^{*})$, there exist a pair of numbers $N_{x}$ and $N_{x}^{*}$ such that the pair of numbers $x$ and $N_{x}$ satisfy
\[
\pi^{*}(x,N_{x})<\pi(x)<\pi^{*}(x,N_{x}+1)\texttt{ and }N_{x}\geq N,
\]
and the pair of numbers $\alpha_{\infty}x$ and $N_{x}^{*}$ satisfy
\[
\pi^{*}(\alpha_{\infty}x,N_{x}^{*})<\pi(\alpha_{\infty}x)<\pi^{*}(\alpha_{\infty}x,N_{x}^{*}+1)\texttt{ and }N_{x}^{*}\geq N+1
\]
where whether $N_{x}\leq N_{x}^{*}$ or $N_{x}>N_{x}^{*}$ it is important that because of $x\in[x_{N},x_{N+1}^{*})$ there are $\alpha_{\infty}x\geq x_{N+1}^{*}$ and $N_{x}^{*}\geq N+1$.

Second, in case of Formula~\ref{eq:alpha_n_gt}, if $\alpha_{N}>\alpha_{\infty}$ then $\forall x\in[x_{N+1}/\alpha_{\infty},x_{N+1})$ there is $\alpha_{\infty}x\geq x_{N+1}$ for $N<\infty$.

Let $x_{N}^{*}=x_{N+1}/\alpha_{\infty}$ such that we have a new $\alpha_{N}^{*}=x_{N+1}/x_{N}^{*}=\alpha_{\infty}$. Then there is $x_{N}<x_{N}^{*}$. We can reduce the old region $[x_{N},x_{N+1})$ to a new region $[x_{N}^{*},x_{N+1})$.

Thus for each nonnegative integer $N=k,k-1,k-2,\ldots,2,1,0$ where $k<\infty$, if $\alpha_{N}>\alpha_{\infty}$ then given each $x$ in the new reduced region $[x_{N}^{*},x_{N+1})$, there exist a pair of numbers $N_{x}$ and $N_{x}^{*}$ such that the pair of numbers $x$ and $N_{x}$ satisfy
\[
\pi^{*}(x,N_{x})<\pi(x)<\pi^{*}(x,N_{x}+1)\texttt{ and }N_{x}\geq N,
\]
and the pair of numbers $\alpha_{\infty}x$ and $N_{x}^{*}$ satisfy
\[
\pi^{*}(\alpha_{\infty}x,N_{x}^{*})<\pi(\alpha_{\infty}x)<\pi^{*}(\alpha_{\infty}x,N_{x}^{*}+1)\texttt{ and }N_{x}^{*}\geq N+1
\]
where whether $N_{x}\leq N_{x}^{*}$ or $N_{x}>N_{x}^{*}$ it is important that because of $x\in[x_{N}^{*},x_{N+1})$ there are $\alpha_{\infty}x\geq x_{N+1}$ and $N_{x}^{*}\geq N+1$.

In general, for each nonnegative integer $N=0,1,2,\ldots$, if $\alpha_{N}<\alpha_{\infty}$ then given each $x$ in the new enlarged region $[x_{N},x_{N+1}^{*})$ where $x_{N+1}^{*}=\alpha_{\infty}x_{N}$, there exist a pair of numbers $N_{x}$ and $N_{x}^{*}$ such that the pair of numbers $x$ and $N_{x}$ satisfy
\[
\pi^{*}(x,N_{x})<\pi(x)<\pi^{*}(x,N_{x}+1)\texttt{ and }N_{x}\geq N,
\]
and the pair of numbers $\alpha_{\infty}x$ and $N_{x}^{*}$ satisfy
\[
\pi^{*}(\alpha_{\infty}x,N_{x}^{*})<\pi(\alpha_{\infty}x)<\pi^{*}(\alpha_{\infty}x,N_{x}^{*}+1)\texttt{ and }N_{x}^{*}\geq N+1,
\]
else if $\alpha_{N}>\alpha_{\infty}$ then for each natural number $K=N,N-1,N-2,\ldots,2,1,0$, if $\alpha_{K}>\alpha_{\infty}$ then given each $x$ in the new reduced region $[x_{K}^{*},x_{K+1})$ where $x_{K}^{*}=x_{K+1}/\alpha_{\infty}$, there exist a pair of numbers $N_{x}$ and $N_{x}^{*}$ such that the pair of numbers $x$ and $N_{x}$ satisfy
\[
\pi^{*}(x,N_{x})<\pi(x)<\pi^{*}(x,N_{x}+1)\texttt{ and }N_{x}\geq K,
\]
and the pair of numbers $\alpha_{\infty}x$ and $N_{x}^{*}$ satisfy
\[
\pi^{*}(\alpha_{\infty}x,N_{x}^{*})<\pi(\alpha_{\infty}x)<\pi^{*}(\alpha_{\infty}x,N_{x}^{*}+1)\texttt{ and }N_{x}^{*}\geq K+1.
\]

Therefore there exists a positive integer $x_{0}$ such that each nonnegative integer $N$ corresponds to a finite region $[x_{N},x_{N+1})$ with $\alpha_{N}:=x_{N+1}/x_{N}=\alpha_{\infty}$ for $N=0,1,2,\ldots$ and $\forall x\in[x_{N},x_{N+1})$ there exists a nonnegative integer $N_{x}$ such that
\[
\pi^{*}(x,N_{x})<\pi(x)<\pi^{*}(x,N_{x}+1)\texttt{ and }N_{x}\geq N.
\]

This completes the proof of the lemma.
\end{proof}

\begin{lemma}\label{le:chap_1_7_5}
We can define a piecewise function $\pi^{*}(x,N)$ with $x_{N}:=10^{2N+3}$ for $N=0,1,2,\ldots$. Then not only each nonnegative integer $N$ corresponds to a finite region $[x_{N},x_{N+1})$ with $\alpha_{N}:=x_{N+1}/x_{N}=10^{2}$ so that there are $\alpha_{\infty}=10^{2}$ and $N_{\alpha}=0$, but also the piecewise function $\pi^{*}(x,N)$ is a lower bound of the prime counting function $\pi(x)$ since $\forall x\in[x_{N},x_{N+1})$ there exists a nonnegative integer $N_{x}$ such that
\[
\pi^{*}(x,N_{x})<\pi(x)<\pi^{*}(x,N_{x}+1)\texttt{ and }N_{x}\geq N.
\]
\end{lemma}

\begin{proof}
Now, let determine the constants $\alpha_{\infty}$ and $N_{\alpha}$ in Lemma~\ref{le:chap_1_7_2} and Lemma~\ref{le:chap_1_7_3}.

From Lemma~\ref{le:chap_1_7_4}, there exists a positive integer $x_{0}$ such that each nonnegative integer $N$ corresponds to a finite region $[x_{N},x_{N+1})$ with $\alpha_{N}:=x_{N+1}/x_{N}=\alpha_{\infty}$ for $N=0,1,2,\ldots$ and $\forall x\in[x_{N},x_{N+1})$ there exists a nonnegative integer $N_{x}$ such that
\[
\pi^{*}(x,N_{x})<\pi(x)<\pi^{*}(x,N_{x}+1)\texttt{ and }N_{x}\geq N.
\]

Since Lemma~\ref{le:chap_1_8_1} and Theorem~\ref{th:chap_1_8_1} hold for $x\geq10^{3}$, we can determine the constant $\alpha_{\infty}$ by some data for $x\geq10^{3}$.

First, let $x_{N}:=10^{2N+3}$ for $N=0,1$. Then by Formula~(\ref{eq:pi_x_3}) and Lemma~\ref{le:chap_1_6_1}, given a number $x\geq599$, there exists a nonnegative integer $N_{x}$ such that there are
\[
\pi^{*}(x,1)\leq\pi^{*}(x,N_{x})<\pi(x)<\pi^{*}(x,N_{x}+1)\texttt{ and }
N_{x}\geq1\texttt{ for }x\geq599,
\]
so that there are $\alpha_{0}:=x_{1}/x_{0}=10^{2}$ and $N_{x}>N=0$ for $x_{0}\leq x<x_{1}$.

Second, let $x_{N}:=10^{2N+3}$ for $N=2,3,\ldots,12$. Then given a number $x=x_{N}$, we obtain a nonnegative integer $N_{x}$ by calculation~\cite{AA,AB,AD}, which satisfies
\[
\pi^{*}(x,N_{x})<\pi(x)<\pi^{*}(x,N_{x}+1)\texttt{ and }
N_{x}\geq N\texttt{ for }N=2,3,\ldots,12,
\]
so that there are $\alpha_{N}:=x_{N+1}/x_{N}=10^{2}$ for $N=1,2,\ldots,11$.

By Formula~(\ref{eq:n_alpha}), let a real number
\[
r=2N+3+O(1)=2N+3+\xi\texttt{ where }0\leq\xi<2.
\]
Then $\forall x\in[x_{N},x_{N+1})$ there is $x=\alpha_{\infty}^{r/2}$. When $r$ is a positive odd integer and $x$ is a positive integer, then $\alpha_{\infty}^{1/2}$ must be a positive integer. When $r$ is a positive odd integer and $x=10^{r}$, then $\alpha_{\infty}^{1/2}$ must be equal to 10.

So for each number $x=10^{r}=x_{N}$ where $x_{N}=10^{2N+3}$ for $N=0,1,2,\ldots,12$, the number $r=2N+3$ is a positive odd integer. Hence there must be $\alpha_{\infty}=10^{2}$.

Thus we can define a piecewise function $\pi^{*}(x,N)$ with $x_{N}:=10^{2N+3}$ for $N=0,1,2,\ldots$. Then not only each nonnegative integer $N$ corresponds to a finite region $[x_{N},x_{N+1})$ with $\alpha_{N}:=x_{N+1}/x_{N}=10^{2}$ so that there are $\alpha_{\infty}=10^{2}$ and $N_{\alpha}=0$, but also the piecewise function $\pi^{*}(x,N)$ is a lower bound of the prime counting function $\pi(x)$ since $\forall x\in[x_{N},x_{N+1})$ there exists a nonnegative integer $N_{x}$ such that
\[
\pi^{*}(x,N_{x})<\pi(x)<\pi^{*}(x,N_{x}+1)\texttt{ and }N_{x}\geq N.
\]

This completes the proof of the lemma.
\end{proof}

\begin{lemma}\label{le:chap_1_7_6}
We can define a piecewise function $\pi^{*}(x,N)$ with $x_{N}:=10^{2N+4}$ for $N=0,1,2,\ldots$. Then not only each nonnegative integer $N$ corresponds to a finite region $[x_{N},x_{N+1})$ with $\alpha_{N}:=x_{N+1}/x_{N}=10^{2}$ so that there are $\alpha_{\infty}=10^{2}$ and $N_{\alpha}=0$, but also the piecewise function $\pi^{*}(x,N)$ is a lower bound of the prime counting function $\pi(x)$ since $\forall x\in[x_{N},x_{N+1})$ there exists a nonnegative integer $N_{x}$ such that
\[
\pi^{*}(x,N_{x})<\pi(x)<\pi^{*}(x,N_{x}+1)\texttt{ and }N_{x}\geq N.
\]
\end{lemma}

\begin{proof}
From Lemma~\ref{le:chap_1_5_1}, there is
\[
\frac{(N+1)!}{\log^{N+2}x}<(\frac{N+2}{e\log x})^{N+2}\texttt{ for }N\geq5.
\]

Let define
\[
\rho(x,N):=\frac{\log x}{N+2}
\]
and
\[
x_{N}:=10^{2N+4}\texttt{ for }N=0,1,2,\ldots,11.
\]
Then given a number $x=x_{N}$, there exists a nonnegative integer $N_{x}$ such that
\[
\pi^{*}(x,N_{x})<\pi(x)<\pi^{*}(x,N_{x}+1)
\]
and
\[
N_{x}\geq N\texttt{ for }N=0,1,2,\ldots,11
\]
by calculation~\cite{AA,AB,AD}, so that there are
\[
\alpha_{N}:=x_{N+1}/x_{N}=10^{2}\texttt{ for }N=1,2,\ldots,10.
\]

By Formula~(\ref{eq:n_alpha}), let a real number
\[
r=2N+4+O(1)=2N+4+\xi\texttt{ where }0\leq\xi<2.
\]
Then $\forall x\in[x_{N},x_{N+1})$ there is $x=\alpha_{\infty}^{r/2}$. When $r$ is a positive even integer and $x$ is a positive integer, then $\alpha_{\infty}^{1/2}$ must be a positive integer. When $r$ is a positive even integer and $x=10^{r}$, then $\alpha_{\infty}^{1/2}$ must be equal to 10.

So for each number $x=10^{r}=x_{N}$ where $x_{N}=10^{2N+4}$ for $N=0,1,2,\ldots,10$, the number $r=2N+4$ is a positive even integer. Hence there must also be $\alpha_{\infty}=10^{2}$.

Thus we can also define a piecewise function $\pi^{*}(x,N)$ with $x_{N}:=10^{2N+4}$ for $N=0,1,2,\ldots$. Then not only each nonnegative integer $N$ corresponds to a finite region $[x_{N},x_{N+1})$ with $\alpha_{N}:=x_{N+1}/x_{N}=10^{2}$ so that there are $\alpha_{\infty}=10^{2}$ and $N_{\alpha}=0$, but also the piecewise function $\pi^{*}(x,N)$ is a lower bound of the prime counting function $\pi(x)$ since $\forall x\in[x_{N},x_{N+1})$ there exists a nonnegative integer $N_{x}$ such that
\[
\pi^{*}(x,N_{x})<\pi(x)<\pi^{*}(x,N_{x}+1)\texttt{ and }N_{x}\geq N.
\]

This completes the proof of the lemma.
\end{proof}

\begin{lemma}\label{le:chap_1_7_7}
For each nonnegative integer $N$ and $\forall x\in[x_{N},x_{N+1})$ let define
\[
\rho(x,N):=\frac{2\log x}{2N+3}.
\]
Then there are
\begin{equation}\label{eq:rho_x_n_eq}
\begin{split}
\rho(x,N)=&\log\alpha_{\infty}*(1+O(\frac{1}{2N+3})) \\
         =&2\log10*(1+O(\frac{1}{2N+3}))
\end{split}
\end{equation}
and
\begin{equation}\label{eq:rho_x_n_inf}
\lim_{N \to \infty}\rho(x,N)=\log\alpha_{\infty}=2\log10
\end{equation}
where the positive constant $\alpha_{\infty}=10^{2}$ so that $\forall x\in[x_{N},x_{N+1})$
\begin{equation}\label{eq:n_alpha_1}
N=\frac{\log x}{2\log10}+O(1)\texttt{ and }x=10^{2N+3+O(1)},
\end{equation}
and
\[
\rho(x,N)<\{
\begin{array}{cc}
  2.16\log10 & \texttt{ for }x\geq10^{3}, \\
  2.10\log10 & \texttt{ for }x\geq10^{41}.
\end{array}
\]

When each pair of numbers $x\in[x_{N},x_{N+1})$ and $N_{x}$ are determined by
\[
\pi^{*}(x,N_{x})<\pi(x)<\pi^{*}(x,N_{x}+1)
\]
there are
\[
\rho(x,N_{x})\leq\rho(x,N)\texttt{ and }\rho(x,N_{x})\leq2\log10\texttt{ as }N \to \infty.
\]
\end{lemma}

\begin{proof}
By the definition of $\eta^{*}(x,N)$ we have
\[
\eta^{*}(x,N+1)-\eta^{*}(x,N)=\frac{(N+1)!}{\log^{(N+1)}x}.
\]

With $N+2\leq N_{x}+2<e\log x$, from Lemma~\ref{le:chap_1_8_1}, there is
\[
\frac{(N+1)!}{\log^{N+2}x}<(\frac{N+3/2}{e\log x})^{N+3/2}\texttt{ for }x\geq10^{3}.
\]

Since $N+3/2$ is a key expression in the inequality above, let define
\[
\rho(x,N):=\frac{2\log x}{2N+3}
\]
to estimate $\eta^{*}(x,N+1)-\eta^{*}(x,N)$ as that done in the proof of Theorem~\ref{th:chap_1_8_1}.

From Lemma~\ref{le:chap_1_7_1}, each nonnegative integer $N$ corresponds to a function $\pi^{*}(x,N)$ valid for a finite region $[x_{N},x_{N+1})$. Then $\forall x\in[x_{N},x_{N+1})$ there are
\[
\rho(x,N):=\frac{2\log x}{2N+3}=\log\alpha_{\infty}*(1+O(\frac{1}{2N+3})).
\]

As $N \to \infty$, Formula~(\ref{eq:rho_x_n_eq}) becomes
\[
\lim_{N \to \infty}\rho(x,N)=\log\alpha_{\infty}*(1+O(\lim_{N \to \infty}\frac{1}{2N+3}))=\log\alpha_{\infty}.
\]

From Lemma~\ref{le:chap_1_7_5} with the positive constant $\alpha_{\infty}=10^{2}$, Formula~(\ref{eq:n_alpha}) becomes Formula~(\ref{eq:n_alpha_1}) and the last equalities in formulas~(\ref{eq:rho_x_n_eq}-\ref{eq:rho_x_n_inf}) also hold.

Thus there are
\[
\rho(x,N)<\{
\begin{array}{cc}
  2.16\log10 & \texttt{ for }x\geq10^{3} \\
  2.10\log10 & \texttt{ for }N\geq19
\end{array}
\]
and $x\in[x_{N},x_{N+1})$.

From Lemma~\ref{le:chap_1_7_1}, given each nonnegative integer $N$, when each pair of numbers $x\in[x_{N},x_{N+1})$ and $N_{x}$ are determined by
\[
\pi^{*}(x,N_{x})<\pi(x)<\pi^{*}(x,N_{x}+1)
\]
there is always $N_{x}\geq N$. Thus $\forall x\in[x_{N},x_{N+1})$ there are
\[
\rho(x,N_{x})\leq\rho(x,N)\texttt{ and }\rho(x,N_{x})\leq2\log10\texttt{ as }N \to \infty.
\]

This completes the proof of the lemma.
\end{proof}

\begin{lemma}\label{le:chap_1_7_8}
For each nonnegative integer $N$ and $\forall x\in[x_{N},x_{N+1})$ let define
\[
\rho(x,N):=\frac{\log x}{N+2}.
\]
Then there are
\begin{equation}\label{eq:rho_x_n_eq_2}
\begin{split}
\rho(x,N)=&\log\alpha_{\infty}*(1+O(\frac{1}{N+2})) \\
         =&2\log10*(1+O(\frac{1}{N+2}))
\end{split}
\end{equation}
and
\begin{equation}\label{eq:rho_x_n_inf_2}
\lim_{N \to \infty}\rho(x,N)=\log\alpha_{\infty}=2\log10
\end{equation}
where the positive constant $\alpha_{\infty}=10^{2}$ so that $\forall x\in[x_{N},x_{N+1})$
\begin{equation}\label{eq:n_alpha_1_2}
N=\frac{\log x}{2\log10}+O(1)\texttt{ and }x=10^{2N+4+O(1)},
\end{equation}
and
\[
\rho(x,N)<\{
\begin{array}{cc}
  2.16\log10 & \texttt{ for }x\geq10^{4} \\
  2.10\log10 & \texttt{ for }x\geq10^{42}
\end{array}.
\]

When each pair of numbers $x\in[x_{N},x_{N+1})$ and $N_{x}$ are determined by
\[
\pi^{*}(x,N_{x})<\pi(x)<\pi^{*}(x,N_{x}+1)
\]
there are
\[
\rho(x,N_{x})\leq\rho(x,N)\texttt{ and }\rho(x,N_{x})\leq2\log10\texttt{ as }N \to \infty.
\]
\end{lemma}

\begin{proof}
By the definition of $\eta^{*}(x,N)$ we have
\[
\eta^{*}(x,N+1)-\eta^{*}(x,N)=\frac{(N+1)!}{\log^{(N+1)}x}.
\]

With $N+2\leq N_{x}+2<e\log x$, from Lemma~\ref{le:chap_1_5_1}, there is
\[
\frac{(N+1)!}{\log^{N+2}x}<(\frac{N+2}{e\log x})^{N+2}\texttt{ for }N\geq5.
\]

Since $N+2$ is a key expression in the inequality above, let define
\[
\rho(x,N):=\frac{\log x}{N+2}
\]
to estimate $\eta^{*}(x,N+1)-\eta^{*}(x,N)$ as that done in the proof of Theorem~\ref{th:chap_1_8_2}.

From Lemma~\ref{le:chap_1_7_1}, each nonnegative integer $N$ corresponds to a function $\pi^{*}(x,N)$ valid for a finite region $[x_{N},x_{N+1})$. Then $\forall x\in[x_{N},x_{N+1})$ there are
\[
\rho(x,N):=\frac{\log x}{N+2}=\log\alpha_{\infty}*(1+O(\frac{1}{N+2})).
\]

As $N \to \infty$, Formula~(\ref{eq:rho_x_n_eq_2}) becomes
\[
\lim_{N \to \infty}\rho(x,N)=\log\alpha_{\infty}*(1+O(\lim_{N \to \infty}\frac{1}{N+2}))=\log\alpha_{\infty}.
\]

From Lemma~\ref{le:chap_1_7_6} with the positive constant $\alpha_{\infty}=10^{2}$, Formula~(\ref{eq:n_alpha_1_2}) holds and the last equalities in formulas~(\ref{eq:rho_x_n_eq_2}-\ref{eq:rho_x_n_inf_2}) also hold.

Thus there are
\[
\rho(x,N)<\{
\begin{array}{cc}
  2.16\log10 & \texttt{ for }x\geq10^{4} \\
  2.10\log10 & \texttt{ for }N\geq19
\end{array}
\]
and $x\in[x_{N},x_{N+1})$.

From Lemma~\ref{le:chap_1_7_1}, given each nonnegative integer $N$, when each pair of numbers $x\in[x_{N},x_{N+1})$ and $N_{x}$ are determined by
\[
\pi^{*}(x,N_{x})<\pi(x)<\pi^{*}(x,N_{x}+1)
\]
there is always $N_{x}\geq N$. Thus $\forall x\in[x_{N},x_{N+1})$ there are
\[
\rho(x,N_{x})\leq\rho(x,N)\texttt{ and }\rho(x,N_{x})\leq2\log10\texttt{ as }N \to \infty.
\]

This completes the proof of the lemma.
\end{proof}

\section{Theorems of upper bounds of prime number gaps}\label{sec:chap_1_8}

\begin{theorem}\label{th:chap_1_8_1}
For each pair of numbers $x\geq10^{3}$ and $N_{x}$ determined by
\[
\pi^{*}(x,N_{x})<\pi(x)<\pi^{*}(x,N_{x}+1),
\]
there are
\begin{equation}\label{ineq:gap_n}
g(x,N_{x}):=\pi^{*}(x,N_{x}+1)-\pi^{*}(x,N_{x})<\sqrt{x}.
\end{equation}
\end{theorem}

\begin{proof}
From Lemma~\ref{le:chap_1_8_1}, when $N+2<e\log x$ there is
\[
\frac{(N+1)!}{\log^{N+2}x}<(\frac{N+3/2}{e\log x})^{N+3/2}\texttt{ for }x\geq10^{3}.
\]

As defined by Lemma~\ref{le:chap_1_7_7}, let
\[
\rho(x,N):=\frac{\log x}{N+3/2}.
\]
Then we consider the ratio of
\[
\frac{\log x}{x^{1/2}}\texttt{ to }\frac{(N+1)!}{\log^{N+1}x}\texttt{ for }N+2<e\log x
\]
and have
\[
\begin{split}
\log\frac{\log^{N+2}x}{(N+1)!x^{1/2}}
>&(N+3/2)\log\frac{e\log x}{N+3/2}-\frac{1}{2}\log x \\
=&(N+3/2)(1+\log\rho(x,N)-\rho(x,N)/2).
\end{split}
\]

Let $\rho_{0}(x,N)$ satisfy the equation
\[
1+\log\rho(x,N)-\rho(x,N)/2=0.
\]
By calculation $5.35669398<\rho_{0}(x,N)<5.35669399$.

Since the key discriminant
\[
1+\log\rho(x,N)-\rho(x,N)/2\geq0
\]
holds for $\rho(x,N)\leq\rho_{0}(x,N)$, thus, from Lemma~\ref{le:chap_1_7_7}, $\forall x\in[x_{N},x_{N+1})$ the inequality
\[
\rho(x,N)<2.16\log10<\rho_{0}(x,N)
\]
holds so that we obtain
\[
\frac{\log^{N+2}x}{(N+1)!}>x^{1/2}\texttt{ and }
\frac{(N+1)!}{\log^{(N+1)}x}<\frac{\log x}{x^{1/2}}.
\]

Hence from Lemma~\ref{le:chap_1_7_1}, given each pair of numbers $x\in[x_{N},x_{N+1})$ and $N_{x}$ determined by
\[
\pi^{*}(x,N_{x})<\pi(x)<\pi^{*}(x,N_{x}+1),
\]
we have $N_{x}\geq N$ so that
\[
\rho(x,N_{x})\leq\rho(x,N)<\rho_{0}(x,N).
\]

From Lemma~\ref{le:chap_1_6_2} there is $N_{x}+2<e\log x$. Let define
\[
\begin{split}
g(x,N_{x}):=&\pi^{*}(x,N_{x}+1)-\pi^{*}(x,N_{x}) \\
           =&\frac{x}{\log x}\frac{(N_{x}+1)!}{\log^{(N_{x}+1)}x}.
\end{split}
\]
Then we have
\[
\frac{g(x,N_{x})}{x/\log x}=\frac{(N_{x}+1)!}{\log^{(N_{x}+1)}x}<\frac{\log x}{x^{1/2}}
\]
and
\[
g(x,N_{x})<\frac{\log x}{x^{1/2}}\frac{x}{\log x}=\sqrt{x}.
\]

This completes the proof of the theorem.
\end{proof}

\begin{theorem}\label{th:chap_1_8_2}
For each pair of numbers $x\geq10^{4}$ and $N_{x}\geq5$ determined by
\[
\pi^{*}(x,N_{x})<\pi(x)<\pi^{*}(x,N_{x}+1),
\]
there are
\begin{equation}\label{ineq:gap_n}
g(x,N_{x}):=\pi^{*}(x,N_{x}+1)-\pi^{*}(x,N_{x})<\sqrt{x}.
\end{equation}
\end{theorem}

\begin{proof}
From Lemma~\ref{le:chap_1_5_1}, there is
\[
\frac{(N+1)!}{\log^{N+2}x}<(\frac{N+2}{e\log x})^{N+2}\texttt{ for }N\geq5.
\]

As defined by Lemma~\ref{le:chap_1_7_8}, let
\[
\rho(x,N):=\frac{\log x}{N+2}.
\]
Then we consider the ratio of
\[
\frac{\log x}{x^{1/2}}\texttt{ to }\frac{(N+1)!}{\log^{N+1}x}\texttt{ for }N\geq5
\]
and have
\[
\begin{split}
\log\frac{\log^{N+2}x}{(N+1)!x^{1/2}}
>&(N+2)\log\frac{e\log x}{N+2}-\frac{1}{2}\log x \\
=&(N+2)(1+\log\rho(x,N)-\rho(x,N)/2).
\end{split}
\]

Let $\rho_{0}(x,N)$ satisfy the equation
\[
1+\log\rho(x,N)-\rho(x,N)/2=0.
\]
By calculation $5.35669398<\rho_{0}(x,N)<5.35669399$.

Since the key discriminant
\[
1+\log\rho(x,N)-\rho(x,N)/2\geq0
\]
holds for $\rho(x,N)\leq\rho_{0}(x,N)$, thus, from Lemma~\ref{le:chap_1_7_8}, $\forall x\in[x_{N},x_{N+1})$ the inequality
\[
\rho(x,N)<2.16\log10<\rho_{0}(x,N)
\]
holds so that we obtain
\[
\frac{\log^{N+2}x}{(N+1)!}>x^{1/2}\texttt{ and }
\frac{(N+1)!}{\log^{(N+1)}x}<\frac{\log x}{x^{1/2}}.
\]

Hence from Lemma~\ref{le:chap_1_7_1}, given each pair of numbers $x\in[x_{N},x_{N+1})$ and $N_{x}$ determined by
\[
\pi^{*}(x,N_{x})<\pi(x)<\pi^{*}(x,N_{x}+1),
\]
we have $N_{x}\geq N$ so that
\[
\rho(x,N_{x})\leq\rho(x,N)<\rho_{0}(x,N).
\]

From Lemma~\ref{le:chap_1_6_2} there is $N_{x}+2<e\log x$. Let define
\[
\begin{split}
g(x,N_{x}):=&\pi^{*}(x,N_{x}+1)-\pi^{*}(x,N_{x}) \\
           =&\frac{x}{\log x}\frac{(N_{x}+1)!}{\log^{(N_{x}+1)}x}.
\end{split}
\]
Then for $N_{x}\geq5$ we have
\[
\frac{g(x,N_{x})}{x/\log x}=\frac{(N_{x}+1)!}{\log^{(N_{x}+1)}x}<\frac{\log x}{x^{1/2}}
\]
and
\[
g(x,N_{x})<\frac{\log x}{x^{1/2}}\frac{x}{\log x}=\sqrt{x}.
\]

This completes the proof of the theorem.
\end{proof}

\begin{lemma}\label{le:chap_1_8_1}
When $N+2<e\log x$, there is
\[
\frac{(N+1)!}{\log^{N+2}x}<(\frac{N+3/2}{e\log x})^{N+3/2}\texttt{ for }x\geq10^{3}.
\]
\end{lemma}

\begin{proof}
Based upon the Stirling's formula
\[
n!=\sqrt{2\pi n}(\frac{n}{e})^{n}e^{\frac{\theta}{12n}},0<\theta<1,
\]
for $N\geq46$, there are
\[
\sqrt{2\pi}e^{\frac{1}{2}+\frac{1}{12(N+1)}}<(\frac{N+3/2}{e})^{1/2}
\]
and
\[
\begin{split}
(N+1)!<&\sqrt{2\pi}e^{\frac{1}{2}+\frac{1}{12(N+1)}}(\frac{N+1}{e})^{N+3/2} \\
      <&(\frac{N+3/2}{e})^{N+2}.
\end{split}
\]

By calculation, the inequality
\[
(N+1)!<(\frac{N+3/2}{e})^{N+2}
\]
holds for $46>N\geq16$.

When $N+2<e\log x$ there are
\[
\begin{split}
\frac{(N+1)!}{\log^{N+2}x}<&(\frac{N+3/2}{e\log x})^{N+2} \\
                          <&(\frac{N+3/2}{e\log x})^{N+3/2}\texttt{ for }N\geq16.
\end{split}
\]

By calculation, the inequality
\[
\frac{(N+1)!}{\log^{N+2}x}<(\frac{N+3/2}{e\log x})^{N+3/2}
\]
holds for $x\geq10^{3}$.

Hence when $N+2<e\log x$ we have
\[
\frac{(N+1)!}{\log^{N+2}x}<(\frac{N+3/2}{e\log x})^{N+3/2}\texttt{ for }x\geq10^{3}.
\]

This completes the proof of the lemma.
\end{proof}

\section{Theorem of closeness of $\Li(x)$ and $\pi(x)$}\label{sec:chap_1_9}

In 1901, H. Koch~\cite{AP} showed that if and only if the Riemann Hypothesis is true, then $\pi(x)=\Li(x)+O(\sqrt{x}\log x)$. So the Riemann Hypothesis is equivalent to the assertion that $|\Li(x)-\pi(x)|\leq c\sqrt{x}\log x$ for some value of $c$.

Based on results in sections~\ref{sec:chap_1_2}~-~\ref{sec:chap_1_8}, we shall prove the following theorem.

\begin{theorem}[Theorem of closeness of $\Li(x)$ and $\pi(x)$]\label{th:chap_1_9_1}
For large numbers $x$, the prime counting function $\pi(x)$ satisfies
\begin{equation}\label{eq:pi_Li_ref}
\pi(x)=\Li(x)+O(\sqrt{x}\log x).
\end{equation}
\end{theorem}

\begin{proof}
First, from Lemma~\ref{le:chap_1_2_1}, given each real number $x\geq198$, there exists a nonnegative integer $M=M(x)$ that is a function of the variable $x$ such that
\[
\pi^{*}(x,M)\leq\Li(x)<\pi^{*}(x,M+1)
\]
where $M=M(x)=\log x+O(1)$ by Formula~(\ref{eq:rho_m_x}), and from Theorem~\ref{th:chap_1_5_1}
\[
\pi^{*}(x,M+1)-\pi^{*}(x,M)<x^{1/64}\texttt{ for }x\geq730.
\]

Second, from Lemma~\ref{le:chap_1_6_1}, given each real number $x\geq599$, there exists a nonnegative integer $N_{x}=N(x)$ that is a function of the variable $x$ such that
\[
\pi^{*}(x,N_{x})<\pi(x)<\pi^{*}(x,N_{x}+1)
\]
where $N_{x}+2<e\log x$ from Lemma~\ref{le:chap_1_6_2}, and from Theorem~\ref{th:chap_1_8_1}
\[
\pi^{*}(x,N_{x}+1)-\pi^{*}(x,N_{x})<\sqrt{x}.
\]

So given each real number $x\geq10^{3}$, there exist a unique pair of nonnegative integers $M=M(x)$ and $N_{x}=N(x)$ such that
\[
\pi^{*}(x,M)\leq\Li(x)<\pi^{*}(x,M+1)
\]
and
\[
\pi^{*}(x,N_{x})<\pi(x)<\pi^{*}(x,N_{x}+1).
\]

If $\pi(x)<\Li(x)$ then we have $N_{x}+2<\log x$ from Lemma~\ref{le:chap_1_6_2} and
\[
\begin{split}
|\Li(x)-\pi(x)|=&\Li(x)-\pi(x) \\
               <&\pi^{*}(x,M+1)-\pi^{*}(x,N_{x}) \\
               <&c_{m}\sqrt{x}\log x
\end{split}
\]
where $M+1-N_{x}<M=M(x)=\log x+O(1)$ and
\[
\begin{split}
\pi^{*}(x,M+1)-\pi^{*}(x,N_{x})=&\sum_{n=N_{x}+1}^{M+1}(\pi^{*}(x,n)-\pi^{*}(x,n-1)) \\
<&(M+1-N_{x})(\pi^{*}(x,N_{x}+1)-\pi^{*}(x,N_{x})) \\
<&M\sqrt{x}\leq c_{m}\sqrt{x}\log x
\end{split}
\]
where $c_{m}$ is a constant greater than 1.

If $\pi(x)>\Li(x)$ then we have $N_{x}+2<e\log x$ and
\[
\begin{split}
|\Li(x)-\pi(x)|=&\pi(x)-\Li(x) \\
<&\pi^{*}(x,N_{x}+1)-\pi^{*}(x,M) \\
<&c_{n}\sqrt{x}\log x
\end{split}
\]
where $N_{x}+1-M<N_{x}=N(x)\leq c_{n}\log x$ and
\[
\begin{split}
&\pi^{*}(x,N_{x}+1)-\pi^{*}(x,M)=\sum_{n=M+1}^{N_{x}+1}(\pi^{*}(x,n)-\pi^{*}(x,n-1)) \\
<&(N_{x}+1-M)\max\{\pi^{*}(x,M+1)-\pi^{*}(x,M),\pi^{*}(x,N_{x}+1)-\pi^{*}(x,N_{x})\} \\
<&N_{x}\max\{x^{1/64},x^{1/2}\}\leq c_{n}x^{1/2}\log x
\end{split}
\]
where $c_{n}$ is a constant greater than $c_{m}$ and less than $e$.

Hence for each positive number $x\geq10^{3}$ we have
\[
|\Li(x)-\pi(x)|\leq c\sqrt{x}\log x
\]
and
\[
\pi(x)=\Li(x)+O(\sqrt{x}\log x)
\]
where $c$ is a constant greater than $c_{m}$ and less than $c_{n}$.

This completes the proof of the theorem.
\end{proof}

%% file: chap2.tex
%-----------------------------------------------------------------------
% Beginning of chap2.tex
%-----------------------------------------------------------------------
%
% AMS-LaTeX 1.2 sample file for a monograph, based on amsbook.cls.
% This is a data file input by chapter.tex.
%%%%%%%%%%%%%%%%%%%%%%%%%%%%%%%%%%%%%%%%%%%%%%%%%%%%%%%%%%%%%%%%%%%%

\chapter{Prime values of Bunyakovsky's polynomials}\label{ch:chap_2}

\section*{Introduction}

The Bunyakovsky conjecture stated in 1857 by the Russian mathematician Viktor Bunyakovsky, asserts when a polynomial $f(x)$ in one variable with positive degree and integer coefficients should have infinitely many prime values for positive integer inputs. Three necessary conditions are

(1) the leading coefficient of $f(x)$ is positive,

(2) the polynomial is irreducible over the integers,

(3) as $m$ runs over the positive integers, the numbers $f(m)$ should not share a common factor greater than 1.

Bunyakovsky's conjecture is that these three conditions are sufficient: if $f(x)$ satisfies these three conditions, which may be called a Bunyakovsky's polynomial, then $f(m)$ is prime for infinitely many positive integers $m$~\cite{AD,FB}.

If $f(x)$ has degree $n\geq2$ and with integer coefficients, for $X\geq1$ let
\[
\pi_{f}(x)=\pi_{f(x)}(X)=\#\{m\geq1||f(m)|\leq X\texttt{ and }|f(m)|\texttt{ is a prime}\}.
\]

In 1922, Nagell showed that $\lim_{X \to \infty}\pi_{f(x)}(X)/X=0$, i.e., $\lim_{x \to \infty}\pi_{f}(x)/x=0$, so there are few prime values. In 1931, Heilbronn proved the more precise statement~\cite{AD}:

There exists a positive constant $C$ (dependent on $f(x)$), such that
\[
\pi_{f(x)}(X)\leq C\frac{X^{1/n}}{\log X},\texttt{ for every }X\geq1,
\]
i.e., there exists a positive constant $C^{*}$ (dependent on $f(x)$), such that
\[
\pi_{f}(x)\leq C^{*}\frac{x}{\log x},\texttt{ for every }x\geq1.
\]

From Theorem~\ref{th:chap_1_1_1}, given a positive real number $\beta$ that determines a real number $x_{\beta}$ by
\[
e(\log x_{\beta})^{3}/x_{\beta}^{0.0327283}=\beta,
\]
let $g(x):=\beta x^{1/2}$ for $x\geq x_{\beta}$ where $g(x):=x^{1/2}$ when let $\beta=1$. Then there are
\[
\frac{\pi(x+g(x))-\pi(x)}{g(x)/\log x}=1+O(\frac{1}{\log x})
\]
and
\[
\lim_{x \to \infty}\frac{\pi(x+g(x))-\pi(x)}{g(x)/\log x}=1.
\]

Hence as $x \to \infty$ the number of odd primes in the interval $(x, x+g(x)]$ is asymptotically equal to
\[
\pi(x+g(x))-\pi(x)=g(x)/\log x+o(g(x)/\log x)
\]
and for each real number $\beta$ there exists a finite positive number $x_{\beta}$ such that for every number $x$ greater than $x_{\beta}$, there are asymptotically $g(x)/\log x$ primes in the interval $(x,x+g(x)]$.

Based on Theorem~\ref{th:chap_1_1_1}, it is proved in this chapter that given a polynomial
\[
f(x)=a_{n}x^{n}+a_{n-1}x^{n-1}+\cdots+a_{1}x+a_{0}=\bar{f}(x)+a_{0}
\]
in one variable with positive degree and integer coefficients, if its integer coefficients satisfy the following three necessary conditions:

(1) $a_{0}>0$, $a_{n}>0$, and $a_{i}\geq0$ for $i=1,2,\ldots,n-1$,

(2) the polynomial is irreducible over the integers,

(3) as $m$ runs over the positive integers, the numbers $f(m)$ should not share a common factor greater than 1,\\
or it can be converted to a new polynomial $f(x)$ such that the integer coefficients of the new polynomial $f(x)$ satisfy these three necessary conditions, then $f(m)$ is prime for infinitely many positive integers $m$.

The key ideas of the proof of the main theorem are as follows: For each $r=1,2,\ldots$,

(1) define a modulus $M_{r}$ on the range $(\bar{f}(x_{r}),\bar{f}(x_{r+1}))$,

(2) make use of Theorem~\ref{th:chap_1_1_1} to estimate the number of primes in a table formed by some values of integral polynomials $\bar{f}(x)+h$ where $h=1,2,\ldots,M_{r}$ in the range $(\bar{f}(x_{r}),\bar{f}(x_{r+1}))$,

(3) make use of the orthogonality relation for Dirichlet characters to estimate the number of primes $p=f(x)$ in the range $(\bar{f}(x_{r}),\bar{f}(x_{r+1}))$, which is estimated greater than $c_{f}x_{r}M_{r}/\log x_{r}$ where $c_{f}$ is a positive constant dependent on the polynomial $f(x)$, \\
so that there exists a large positive number $r_{n}$ such that for each large positive number $r=r_{n},r_{n}+1,r_{n}+2,\ldots$ there are more than $c_{f}x_{r}M_{r}/\log x_{r}$ primes $p=f(mx_{r})$ in the range $(\bar{f}(x_{r}),\bar{f}(x_{r+1}))$, i.e., $f(m)$ is prime for infinitely many positive integers $m$. These results are stated in Theorem~\ref{th:chap_2_4_1}. Then the main theorem~\ref{th:chap_2_4_2} is proved.

%%Hence, based on the main theorem, it is proved that there are infinitely many primes of the form $a^{2}+1$.

\section{Definitions}\label{sec:chap_2_1_10}

\begin{definition}\label{def:chap_2_1_1}
When a polynomial
\[
f(x)=a_{n}x^{n}+a_{n-1}x^{n-1}+\cdots+a_{1}x+a_{0}=\bar{f}(x)+a_{0}
\]
in one variable with positive degree $n\geq2$ and integer coefficients that satisfy the following necessary conditions:
\[
a_{0}>0,\texttt{ }a_{n}>0,\texttt{ and }a_{i}\geq0\texttt{ for }i=1,2,\ldots,n-1,
\]
there exists a positive integer $x_{1}$ such that $a_{i}\leq a_{n}x_{1}^{n-i}$ for $i=1,2,\ldots,n-1$,
\[
\sqrt{na_{n}}x_{1}^{n}>a_{0},\texttt{ }gcd(a_{0},x_{1})=1,\texttt{ and }e\log^{3}(na_{n}x_{1}^{n})\leq(na_{n}x_{1}^{n})^{0.0327283}.
\]

Let $x_{r+1}=x_{r}^{2}$ and each positive integer $M_{r}$ satisfy $na_{n}x_{r}^{n}>M_{r}\geq\sqrt{na_{n}}x_{r}^{n}$ and $gcd(a_{0},M_{r})=1$ for $r=1,2,\ldots$. Then there are
\[
\bar{f}(x_{r})=a_{n}x_{r}^{n}+a_{n-1}x_{r}^{n-1}+\cdots+a_{1}x_{r}\leq na_{n}x_{r}^{n}
\]
and
\[
e\log^{3}(na_{n}x_{r+1}^{n})\leq(na_{n}x_{r+1}^{n})^{0.0327283}\texttt{ for }r=1,2,\ldots.
\]
\end{definition}

Since $x_{r+1}=x_{r}^{2}$, there are $x_{r+1}/x_{r}=x_{r}$ and $x_{r+1}/x_{r} \to \infty$ as $x_{r} \to \infty$. Thus $x_{r+1}/x_{r}$ is with a same order of infinity as $x_{r}$.

Hence for any arbitrary small positive value $\epsilon$, there exists a number $r_{\epsilon}$ such that there are $x_{r}\geq 1/\epsilon$ and $x_{r+1}\geq x_{r}/\epsilon$ for $r\geq r_{\epsilon}$, which implies $x_{r} \to \infty$ and $x_{r+1}/x_{r} \to \infty$ as $\epsilon \to 0$.

\begin{definition}\label{def:chap_2_1_2}
With $0<h<M_{r}$ and $gcd(h,M_{r})=1$, let denote
\[
\pi(x;x_{r})=\sum_{\bar{f}(\xi)<p<\bar{f}(\xi)+M_{r},x_{r}\leq\xi=mx_{r}<x}1
\]
and
\[
\pi(x;x_{r},h)=\sum_{p=\bar{f}(\xi)+h,0<h<M_{r}, x_{r}\leq\xi=mx_{r}<x}1.
\]
Then the function $\pi(x;x_{r})$ counts the number of primes $p=\bar{f}(\xi)+h$ where
\[
0<h<M_{r},\texttt{ }x_{r}\leq\xi=mx_{r}<x\texttt{ for }m=1,2,\ldots,
\]
and $h$ varies in the range $[1,M_{r})$, and the function $\pi(x;x_{r},h)$ counts the number of primes
\[
p=\bar{f}(\xi)+h\texttt{ where }0<h<M_{r},\texttt{ }x_{r}\leq\xi=mx_{r}<x\texttt{ for }m=1,2,\ldots,
\]
and $h$ is fixed.
\end{definition}

\begin{definition}\label{def:chap_2_1_3}
Let define prime sets $\{p\}$ on polynomials $f(\xi)=\bar{f}(\xi)+h$ by
\[
Q_{i}(r,x)=\{p|p=\bar{f}(q)+h_{i},0<h_{i}<M_{r}, x_{r}\leq q=mx_{r}<x,gcd(h_{i},M_{r})=1\}
\]
for $i=1,2,\ldots,\varphi(M_{r})$, and with $0<h<M_{r}$ and $gcd(h,M_{r})=1$
\[
Q_{f}(r,x)=\bigcup_{i=1}^{\varphi(M_{r})}Q_{i}(r,x)=\{p|p=\bar{f}(q)+h,0<h<M_{r}, x_{r}\leq q=mx_{r}<x\}
\]
where $m=1,2,\ldots$ and $\varphi(M_{r})$ is the Euler's totient function.
\end{definition}

Let $\pi_{f}(r,x)$ count the number of primes $p\in Q_{f}(r,x)$ and $\pi_{f}(r,x,h_{i})$ count the number of primes $p\in Q_{i}(r,x)$ for $i=1,2,\ldots,\varphi(M_{r})$. Then there are
\[
\pi(x;x_{r})=\pi_{f}(r,x)\texttt{ and }\pi(x;x_{r},h_{i})=\pi_{f}(r,x,h_{i})\texttt{ for }i=1,2,\ldots,\varphi(M_{r}).
\]

\begin{definition}\label{def:chap_2_1_4}
With complex variable $s$ and primes $p$, for $i=1,2,\ldots,\varphi(M_{r})$, let define complex functions on prime sets $Q_{i}(r,x)$ by
\[
\zeta_{i}(r,x,s)=\prod_{p\in Q_{i}(r,x)}(1-p^{-s})^{-1}.
\]

Let define a complex function $\zeta_{f}(r,x,s)$ on prime sets $Q_{i}(r,x)$ and $Q_{f}(r,x)$ by
\[
\zeta_{f}(r,x,s)=\prod_{i=1}^{\varphi(M_{r})}\zeta_{i}(r,x,s)=\prod_{p\in Q_{f}(r,x)}(1-p^{-s})^{-1}.
\]
\end{definition}

\begin{definition}\label{def:chap_2_1_5}
Let define functions on prime sets $Q_{i}(r,x)$ and $Q_{f}(r,x)$ by
\[
\Lambda_{i}(r,\xi)=\{
       \begin{array}{cc}
         \log p & \texttt{ if }\xi=p^{m},\texttt{ for any }m>0\texttt{ and }p\in Q_{i}(r,x), \\
         0 & otherwise.
       \end{array}
\]
for $i=1,2,\ldots,\varphi(M_{r})$, and
\[
\Lambda_{f}(r,\xi)=\{
       \begin{array}{cc}
         \log p & \texttt{ if }\xi=p^{m},\texttt{ for any }m>0\texttt{ and }p\in Q_{f}(r,x), \\
         0 & otherwise.
       \end{array}
\]
\end{definition}

Differentiating the relations
\[
\log\zeta_{i}(r,x,s)=\sum_{m,p\in Q_{i}(r,x)}\frac{1}{mp^{ms}}\texttt{ for }i=1,2,\ldots,\varphi(M_{r}),
\]
and
\[
\log\zeta_{f}(r,x,s)=\sum_{m,p\in Q_{f}(r,x)}\frac{1}{mp^{ms}},
\]
we obtain
\[
\frac{\zeta_{i}'(r,x,s)}{\zeta_{i}(r,x,s)}
=-\sum_{m,p\in Q_{i}(r,x)}\frac{\log p}{p^{ms}}
=-\sum_{\xi=1}^{\infty}\frac{\Lambda_{i}(r,\xi)}{\xi^{s}}\texttt{ for }i=1,2,\ldots,\varphi(M_{r}),
\]
and
\[
\frac{\zeta_{f}'(r,x,s)}{\zeta_{f}(r,x,s)}
=-\sum_{m,p\in Q_{f}(r,x)}\frac{\log p}{p^{ms}}
=-\sum_{\xi=1}^{\infty}\frac{\Lambda_{f}(r,\xi)}{\xi^{s}}.
\]

\begin{definition}\label{def:chap_2_1_6}
For each Dirichlet character $\chi$ and $i=1,2,\ldots,\varphi(M_{r})$, let define functions on prime sets $Q_{i}(r,x)$ by
\[
L_{i}(r,x,s,\chi)=\prod_{p\in Q_{i}(r,x)}(1-\chi(p)p^{-s})^{-1}.
\]

For each Dirichlet character $\chi$, let define a function $L_{f}(r,x,s,\chi)$ on prime sets $Q_{i}(r,x)$ and $Q_{f}(r,x)$ by
\[
L_{f}(r,x,s,\chi)=\prod_{i=1}^{\varphi(M_{r})}L_{i}(r,x,s,\chi)=\prod_{p\in Q_{f}(r,x)}(1-\chi(p)p^{-s})^{-1}.
\]
\end{definition}

\begin{definition}\label{def:chap_2_1_7}
For primes $p\in Q_{f}(r,x_{r+1})$, similar to the proof of the prime number lemma for arithmetic progressions in~\cite{FD}, our object is to get an estimate for
\begin{equation}
\pi_{f}(r,x_{r+1})=\sum_{p\in Q_{f}(r,x_{r+1})}1=\sum_{\xi=f(x_{r})}^{f(x_{r+1})}P_{f}(r,\xi)
\end{equation}
where $P_{f}$ is the characteristic function of the numbers:
\[
P_{f}(r,\xi)=\{
       \begin{array}{cc}
         1 & \texttt{ if }\xi\in Q_{f}(r,x_{r+1}), \\
         0 & otherwise.
       \end{array}
\]
\end{definition}
While $P_{f}$ itself does not arise in a natural way, the function $P_{f}^{*}$ such that
\[
P_{f}^{*}(r,\xi)=\{
       \begin{array}{cc}
         \frac{1}{m} & \texttt{ if $\xi=p^{m}$ for some $m$ and $p\in Q_{f}(r,x_{r+1})$}, \\
         0 & otherwise,
       \end{array}
\]
occurs in the Dirichlet series for $\log\zeta_{f}(r,x_{r+1},s)$:
\begin{equation}\label{eq:def_zeta}
\log\zeta_{f}(r,x_{r+1},s)=\sum_{m,p\in Q_{f}(r,x_{r+1})}\frac{1}{mp^{ms}}
=\sum_{\xi=1}^{\infty}\frac{P_{f}^{*}(r,\xi)}{\xi^{s}}.
\end{equation}
For each fixed number $m$, the number of $m$th powers of primes $p\in Q_{f}(r,x_{r+1})$ is less than or equal to the number of primes $p\in Q_{f}(r,x_{r+1})$ satisfying $f(x_{r})<p<\sqrt[m]{f(x_{r+1})}$, so that with $f^{1/m}(x)<f(x^{1/m})$ for $m\geq2$
\[
\begin{split}
\sum_{\xi=f(x_{r})}^{f(x_{r+1})}P_{f}^{*}(r,\xi)
=&\sum_{\xi=f(x_{r})}^{f(x_{r+1})}P_{f}(r,\xi)
+\frac{1}{2}\sum_{\xi=f(x_{r})}^{\sqrt{f(x_{r+1})}}P_{f}(r,\xi)
+\frac{1}{3}\sum_{\xi=f(x_{r})}^{\sqrt[3]{f(x_{r+1})}}P_{f}(r,\xi)+\cdots \\
\leq&\sum_{\xi=f(x_{r})}^{f(x_{r+1})}P_{f}(r,\xi)
+\frac{1}{2}\sum_{\xi=f(x_{r})}^{f(\sqrt{x_{r+1}})}P_{f}(r,\xi)
+\frac{1}{3}\sum_{\xi=f(x_{r})}^{f(\sqrt[3]{x_{r+1}})}P_{f}(r,\xi)+\cdots \\
=&\pi_{f}(r,x_{r+1})+\frac{\pi_{f}(r,\sqrt{x_{r+1}})}{2}+\frac{\pi_{f}(r,\sqrt[3]{x_{r+1}})}{3}+\cdots.
\end{split}
\]

For $m=2,3,\ldots$, by Definition~\ref{def:chap_2_1_3} there are
\[
\pi_{f}(r,x_{r+1}^{1/m})\leq\pi_{f}(r,x_{r}^{2/2})=\pi_{f}(r,x_{r})\sim0
\]
such that we can write
\[
\pi_{f}(r,x_{r+1}^{1/m})=O(x_{r+1}^{-m/2})\texttt{ for }m=2,3,\ldots.
\]

Thus there are
\[
\frac{\pi_{f}(r,\sqrt{x_{r+1}})}{2}+\frac{\pi_{f}(r,\sqrt[3]{x_{r+1}})}{3}+\cdots=O(\frac{1}{x_{r+1}}).
\]

Since $x_{r+1}=x_{r}^{2}$, there are $x_{r+1}/x_{r}=x_{r}$ and $x_{r+1}/x_{r} \to \infty$ as $x_{r} \to \infty$. Thus $x_{r+1}/x_{r}$ is with a same order of infinity as $x_{r}$ and for any arbitrary small positive value $\epsilon$, there should exist a number $r_{\epsilon}$ such that for $r\geq r_{\epsilon}$, there are
\[
\frac{\pi_{f}(r,\sqrt{x_{r+1}})}{2}+\frac{\pi_{f}(r,\sqrt[3]{x_{r+1}})}{3}+\cdots
<\epsilon\pi_{f}(r,x_{r+1})
\]
and
\[
\sum_{\xi=f(x_{r})}^{f(x_{r+1})}P_{f}^{*}(r,\xi)\sim\pi_{f}(r,x_{r+1})\texttt{ as }\epsilon \to 0.
\]

\section{Lemmas of Bunyakovsky's polynomials}\label{sec:chap_2_2_11}

\begin{lemma}\label{le:chap_2_2_1}
When a polynomial
\begin{equation}\label{eq:bunyakovsky_col_1}
f(x)=a_{n}x^{n}+a_{n-1}x^{n-1}+\cdots+a_{1}x+a_{0}
\end{equation}
in one variable with positive degree and integer coefficients that satisfy $a_{0}\neq0$ and $a_{n}>0$, then $f(x)$ can always be converted to a new polynomial
\begin{equation}\label{eq:bunyakovsky_col_1_y}
g(y)=b_{n}y^{n}+b_{n-1}y^{n-1}+\cdots+b_{1}y+b_{0}
\end{equation}
such that $b_{0}>0$, $b_{n}>0$, and $b_{i}\geq0$ for $i=1,2,\ldots,n-1$.
\end{lemma}

\begin{proof}
If some of $a_{i}$ for $i=1,2,\ldots,n-1$ are negative, then by choosing a suitable positive number $k$ and setting $x=y+k$, $f(x)$ becomes
\[
\begin{split}
f(y+k)=&a_{n}(y+k)^{n}+a_{n-1}(y+k)^{n-1}+\cdots+a_{1}(y+k)+a_{0} \\
      =&a_{n}[y^{n}+C(n,1)ky^{n-1}+\cdots+C(n,n-1)k^{n-1}y+k^{n}] \\
       &+a_{n-1}[y^{n-1}+C(n-1,1)ky^{n-2}+\cdots+C(n-1,n-2)k^{n-2}y+k^{n-1}] \\
       &+\cdots+a_{1}(y+k)+a_{0} \\
      =&a_{n}y^{n}+[a_{n}C(n,1)k+a_{n-1}]y^{n-1}+\cdots \\
       &+[a_{n}C(n,1)k^{n-1}+a_{n-1}C(n-1,1)k^{n-2}+\cdots+a_{1}]y \\
       &+(a_{n}k^{n}+a_{n-1}k^{n-1}+\cdots+a_{1}k+a_{0}) \\
      =&g(y)=b_{n}y^{n}+b_{n-1}y^{n-1}+\cdots+b_{1}y+b_{0}
\end{split}
\]
where
\[
b_{i}=\sum_{j=i}^{n}a_{j}C(j,j-i)k^{j-i}\texttt{ for }i=0,1,2,\ldots,n\texttt{ and }C(n,i)=\frac{n!}{i!(n-i)!},
\]
such that the inequalities $b_{0}>0$, $b_{n}>0$, and $b_{i}\geq0$ for $i=1,2,\ldots,n-1$ always holds.

Hence without loss of generality, we can assume that the inequalities $a_{0}>0$, $a_{n}>0$, and $a_{i}\geq0$ for $i=1,2,\ldots,n-1$ always holds.

This completes the proof of the lemma.
\end{proof}

\begin{lemma}[Primes in rows]\label{le:chap_2_2_2}
When a polynomial
\begin{equation}\label{eq:bunyakovsky_row_p}
f(x)=a_{n}x^{n}+a_{n-1}x^{n-1}+\cdots+a_{1}x+a_{0}=\bar{f}(x)+a_{0}
\end{equation}
in one variable with positive degree and integer coefficients that satisfy the following necessary conditions:
\[
a_{0}>0,\texttt{ }a_{n}>0,\texttt{ and }a_{i}\geq0\texttt{ for }i=1,2,\ldots,n-1,
\]
by Definition~(\ref{def:chap_2_1_1}), if numbers
\[
\bar{f}(mx_{r})+1,\bar{f}(mx_{r})+2,\ldots,\bar{f}(mx_{r})+M_{r}
\]
for $m=1,2,\ldots,m_{r}$ are arranged in $m_{r}$ rows each containing $M_{r}$ numbers:
\begin{equation}\label{tab:bunyakovsky_row_p}
\begin{array}{cccc}
  \bar{f}(x_{r})+1, & \bar{f}(x_{r})+2, & \cdots, & \bar{f}(x_{r})+M_{r} \\
  \bar{f}(2x_{r})+1, & \bar{f}(2x_{r})+2, & \cdots, & \bar{f}(2x_{r})+M_{r} \\
  \cdots, & \cdots, & \cdots, & \cdots \\
  \bar{f}(m_{r}x_{r})+1, & \bar{f}(m_{r}x_{r})+2, & \cdots, & \bar{f}(m_{r}x_{r})+M_{r}
\end{array}
\end{equation}
where $m_{r}=x_{r+1}/x_{r}=x_{r}$, $na_{n}x_{r}^{n}>M_{r}\geq\sqrt{na_{n}}x_{r}^{n}$ and $gcd(a_{0},M_{r})=1$, then Table~\ref{tab:bunyakovsky_row_p} contains asymptotically more than $c_{f}x_{r}M_{r}/\log x_{r}$ primes $p=\bar{f}(mx_{r})+h_{c}$ and $\pi_{f}(r,x_{r+1})$ is asymptotically greater than $c_{f}x_{r}M_{r}/\log x_{r}$ where $c_{f}$ is a positive constant dependent on the polynomial $f(x)$ and $h_{c}$ is an integer satisfying $0<h_{c}<M_{r}$ and $gcd(h_{c},M_{r})=1$.
\end{lemma}

\begin{proof}
From Theorem~\ref{th:chap_1_1_1}, given a positive real number $\beta$ that determines a real number $y_{\beta}$ by $e(\log y_{\beta})^{3}/y_{\beta}^{0.0327283}=\beta$, let $g(y):=\beta y^{1/2}$ for $y\geq y_{\beta}$ where $g(y):=y^{1/2}$ when let $\beta=1$. Then for every $y\geq y_{\beta}$, there are $e(\log y)^{3}/y^{0.0327283}\leq\beta$ and asymptotically $g(y)/\log y$ primes in the interval $(y,y+g(y)]$, i.e., there must be
\[
\pi(y+g(y))-\pi(y)=\frac{g(y)}{\log y}(1+O(\frac{1}{\log y}))\sim\frac{g(y)}{\log y}.
\]

Thus for $y=\bar{f}(x_{r+1})$ there must be
\[
g(\bar{f}(x_{r+1}))\geq(\bar{f}(x_{r+1}))^{1/2-0.0327283}e(\log\bar{f}(x_{r+1}))^{3}.
\]

It can be satisfied. Let $g(\bar{f}(x_{r+1}))=M_{r}$. Then by Definition~(\ref{def:chap_2_1_1}) we have
\[
\begin{split}
g(\bar{f}(x_{r+1}))\geq&\sqrt{na_{n}}x_{r}^{n}=\sqrt{na_{n}x_{r}^{2n}}=\sqrt{na_{n}x_{r+1}^{n}} \\
\geq&(na_{n}x_{r+1}^{n})^{1/2-0.0327283}e(\log(na_{n}x_{r+1}^{n}))^{3} \\
\geq&(\bar{f}(x_{r+1}))^{1/2-0.0327283}e(\log\bar{f}(x_{r+1}))^{3}.
\end{split}
\]

For $m=1,2,\ldots,m_{r}$, let $y=\bar{f}(mx_{r})$ and $g(y)=M_{r}$. Then there are
\[
\begin{split}
g(y)\geq\sqrt{na_{n}}x_{r}^{n}\geq&(\bar{f}(x_{r+1}))^{1/2-0.0327283}e(\log\bar{f}(x_{r+1}))^{3} \\
\geq&(\bar{f}(mx_{r}))^{1/2-0.0327283}e(\log\bar{f}(mx_{r}))^{3} \\
=&y^{1/2-0.0327283}e(\log y)^{3}
\end{split}
\]
such that for large values of $x_{r}$ we have
\[
\begin{split}
\pi(y+g(y))-\pi(y)=&\frac{g(y)}{\log y}(1+O(\frac{1}{\log y}))\sim\frac{g(y)}{\log y}
                  =\frac{M_{r}}{\log(\bar{f}(mx_{r}))} \\
                  \geq&\frac{M_{r}}{\log(na_{n}(mx_{r})^{n})}
                  \geq\frac{\sqrt{na_{n}}x_{r}^{n}}{\log(na_{n}(mx_{r})^{n})}.
\end{split}
\]

Especially, for large values of $x_{r}$ we have
\[
\begin{split}
\pi(\bar{f}(x_{r})+g(\bar{f}(x_{r})))-\pi(\bar{f}(x_{r}))
\sim&\frac{g(\bar{f}(x_{r}))}{\log(\bar{f}(x_{r}))}=\frac{M_{r}}{\log(\bar{f}(x_{r}))} \\
\geq&\frac{M_{r}}{\log(na_{n}x_{r}^{n})}
\geq\frac{\sqrt{na_{n}}x_{r}^{n}}{\log(na_{n}x_{r}^{n})}
\end{split}
\]
and
\[
\begin{split}
\pi(\bar{f}(x_{r+1})+g(\bar{f}(x_{r+1})))-\pi(\bar{f}(x_{r+1}))
\sim&\frac{g(\bar{f}(x_{r+1}))}{\log(\bar{f}(x_{r+1}))}=\frac{M_{r+1}}{\log(\bar{f}(x_{r+1}))} \\
\geq&\frac{M_{r}}{\log(na_{n}x_{r}^{2n})}
\geq\frac{\sqrt{na_{n}}x_{r}^{n}}{\log(na_{n}x_{r}^{2n})}.
\end{split}
\]

Finally, for each $m=1,2,\ldots,m_{r}$, there are
\[
\begin{split}
\bar{f}(mx_{r}+x_{r})-\bar{f}(mx_{r})=&a_{n}[(mx_{r}+x_{r})^{n}-(mx_{r})^{n}] \\
                                      &+a_{n-1}[(mx_{r}+x_{r})^{n-1}-(mx_{r})^{n-1}] \\
                                      &+\cdots+a_{1}[(mx_{r}+x_{r})-(mx_{r})] \\
\geq&a_{n}x_{r}^{n}\sum_{k=0}^{n-1}(m+1)^{n-1-k}m^{k}>a_{n}x_{r}^{n}\sum_{k=0}^{n-1}m^{n-1} \\
=&na_{n}x_{r}^{n}m^{n-1}\geq na_{n}x_{r}^{n}>M_{r}.
\end{split}
\]
Thus Table~\ref{tab:bunyakovsky_row_p} contains no repeated primes. On the other hand, for each prime $p$ in Table~\ref{tab:bunyakovsky_row_p}, there is $gcd(p,M_{r})=1$, which implies that the prime $p$ belongs to a column $h_{p}$ satisfying $0<h_{p}<M_{r}$ and $gcd(h_{p},M_{r})=1$. Thus no prime is missed out when $\varphi(M_{r})$ columns $h$ satisfying $gcd(h,M_{r})=1$ are getting out of Table~\ref{tab:bunyakovsky_row_p}.

Hence for $m=1,2,\ldots,m_{r}$ each $m$th row of Table~\ref{tab:bunyakovsky_row_p} contains asymptotically $M_{r}/\log(na_{n}(mx_{r})^{n})$ primes $p=\bar{f}(mx_{r})+h_{c}$ where $h_{c}$ is an integer satisfying $0<h_{c}<M_{r}$ and $gcd(h_{c},M_{r})=1$ in the interval $(\bar{f}(mx_{r}),\bar{f}(mx_{r})+M_{r})$ such that Table~\ref{tab:bunyakovsky_row_p} contains asymptotically more than $c_{f}x_{r}M_{r}/\log x_{r}$ primes and $\pi_{f}(r,x_{r+1})$ is asymptotically greater than $c_{f}x_{r}M_{r}/\log x_{r}$ where $c_{f}$ is a positive constant dependent on the polynomial $f(x)$, since
\[
\begin{split}
\pi_{f}(r,x_{r+1})\sim&\sum_{m=1}^{m_{r}-1}\frac{M_{r}}{\log(na_{n}(mx_{r})^{n})} \\
=&\frac{M_{r}}{\log(na_{n}x_{r}^{2n})}\sum_{m=1}^{m_{r}-1}\frac{\log(na_{n}x_{r}^{2n})}{\log(na_{n}(mx_{r})^{n})}\\
=&\frac{M_{r}}{\log(na_{n}x_{r}^{2n})}\sum_{m=1}^{m_{r}-1}(1+\frac{\log(x_{r}/m)^{n}}{\log(na_{n}(mx_{r})^{n})})\\
>&\frac{m_{r}M_{r}}{\log(na_{n}x_{r}^{2n})}\geq c_{f}\frac{x_{r}M_{r}}{\log x_{r}}
\end{split}
\]
where for $\log(na_{n}x_{r}^{n})>3\log(na_{n}3^{2n})$
\[
\begin{split}
\sum_{m=1}^{m_{r}-1}\frac{\log(x_{r}/m)^{n}}{\log(na_{n}(mx_{r})^{n})}
>&\frac{\log x_{r}^{n}}{\log(na_{n}x_{r}^{n})}+\frac{\log(x_{r}/2)^{n}}{\log(na_{n}(2x_{r})^{n})}+\frac{\log(x_{r}/3)^{n}}{\log(na_{n}(3x_{r})^{n})}\\
=&3-\frac{\log(na_{n})}{\log(na_{n}x_{r}^{n})}-\frac{\log(na_{n}2^{2n})}{\log(na_{n}(2x_{r})^{n})}-\frac{\log(na_{n}3^{2n})}{\log(na_{n}(3x_{r})^{n})}\\
>&3-3\frac{\log(na_{n}3^{2n})}{\log(na_{n}x_{r}^{n})}>2.
\end{split}
\]

This completes the proof of the lemma.
\end{proof}

\section{Lemmas of polynomial progressions}\label{sec:chap_2_3_12}

\subsection{Lemmas to be referenced}

Some lemmas proved in~\cite{FD} are stated in this subsection.

\begin{lemma}[Theorem 6-7 in~\cite{FD}]\label{le:chap_2_3_1_th_6_7}
If $(a,m)=(b,m)=1$, then
\[
\sum_{\chi\in X(m)}\frac{\chi(a)}{\chi(b)}=\{
\begin{array}{cc}
  \varphi(m) & \texttt{ if }a\equiv b(mod\texttt{ }m), \\
  0 & \texttt{ otherwise}.
\end{array}
\]
\end{lemma}

If $f$ is regular on the vertical line $\sigma_{0}+it$, and if
\[
\lim_{a \to \infty,b \to \infty}\int_{\sigma_{0}-ai}^{\sigma_{0}+bi}f(s)ds=
\lim_{a \to \infty,b \to \infty}\int_{-a}^{b}f(\sigma_{0}+it)idt
\]
exists, then we abbreviate this limit to $\int_{(\sigma_{0})}f(s)ds$.

\begin{lemma}[Theorem 7-9 in~\cite{FD}]\label{le:chap_2_3_2_th_7_9}
\[
\frac{1}{2\pi i}\int_{(2)}\frac{y^{s}}{s^{2}}ds=\{
\begin{array}{cc}
  0 & \texttt{ for }0<y<1, \\
  \log y & \texttt{ for }y\geq1.
\end{array}
\]
\end{lemma}

\subsection{Lemmas to be proved}

\begin{lemma}\label{le:chap_2_3_3}
\[
\sum_{p\in Q_{f}(r,x_{r+1}),m\geq2,p^{m}\leq f(x_{r+1})}\frac{1}{m}\log\frac{f(x_{r+1})}{p^{m}}=O(\frac{\log x_{r+1}}{x_{r+1}}).
\]
\end{lemma}

\begin{proof}
For $m=2,3,\ldots$, by Definition~\ref{def:chap_2_1_3} there are
\[
\pi_{f}(r,x_{r+1}^{1/m})\leq\pi_{f}(r,x_{r}^{2/2})=\pi_{f}(r,x_{r})\sim0
\]
such that we can write
\[
\pi_{f}(r,x_{r+1}^{1/m})=O(x_{r+1}^{-m/2})\texttt{ for }m=2,3,\ldots.
\]

Thus there are
\[
\frac{\pi_{f}(r,\sqrt{x_{r+1}})}{2}+\frac{\pi_{f}(r,\sqrt[3]{x_{r+1}})}{3}+\cdots=O(\frac{1}{x_{r+1}})
\]
and
\[
\begin{split}
&\sum_{p\in Q_{f}(r,x_{r+1}),m\geq2,p^{m}\leq f(x_{r+1})}\frac{1}{m}\log\frac{f(x_{r+1})}{p^{m}} \\
\leq&\sum_{p\in Q_{f}(r,x_{r+1}),m\geq2,p^{m}\leq f(x_{r+1})}\frac{1}{m}\log f(x_{r+1}) \\
\leq&(\frac{\pi_{f}(r,x_{r+1}^{1/2})}{2}+\frac{\pi_{f}(r,x_{r+1}^{1/3})}{3}+\cdots)\log f(x_{r+1}) \\
=&O(\frac{\log f(x_{r+1})}{x_{r+1}})=O(\frac{\log x_{r+1}}{x_{r+1}}).
\end{split}
\]

This completes the proof of the lemma.
\end{proof}

\begin{lemma}[Theorem 7-15 in~\cite{FD}]\label{le:chap_2_3_4_th_7_15}
\begin{equation}\label{eq:th_7_15}
\begin{split}
\sum_{p\in Q_{f}(r,x_{r+1})}\log\frac{f(x_{r+1})}{p}
=&\frac{1}{2\pi i}\int_{(2)}\frac{f^{s}(x_{r+1})}{s^{2}}\log\zeta_{f}(r,x_{r+1},s)ds \\
&+O(\frac{\log x_{r+1}}{x_{r+1}}).
\end{split}
\end{equation}
\end{lemma}

\begin{proof}
Using Lemma~\ref{le:chap_2_3_2_th_7_9}, similar to the proof of Theorem 7-15 in~\cite{FD}, we have
\[
\begin{split}
\frac{1}{2\pi i}\int_{(2)}\frac{f^{s}(x_{r+1})}{s^{2}}\log\zeta_{f}(r,x_{r+1},s)ds
=&\frac{1}{2\pi i}\int_{(2)}\frac{1}{s^{2}}\sum_{\xi=1}^{\infty}\frac{\Lambda_{f}(r,\xi)}{\log\xi}(\frac{f(x_{r+1})}{\xi})^{s}ds\\
=&\frac{1}{2\pi i}\sum_{\xi=1}^{\infty}\frac{\Lambda_{f}(r,\xi)}{\log\xi}\int_{(2)}\frac{1}{s^{2}}(\frac{f(x_{r+1})}{\xi})^{s}ds\\
=&\sum_{\xi\leq f(x_{r+1})}\frac{\Lambda_{f}(r,\xi)}{\log\xi}\log\frac{f(x_{r+1})}{\xi} \\
=&\sum_{m,p\in Q_{f}(r,x_{r+1}),p^{m}\leq f(x_{r+1})}\frac{1}{m}\log\frac{f(x_{r+1})}{p^{m}} \\
=&\sum_{p\in Q_{f}(r,x_{r+1})}\log\frac{f(x_{r+1})}{p} \\
 &+\sum_{p\in Q_{f}(r,x_{r+1}),m\geq2,p^{m}\leq f(x_{r+1})}\frac{1}{m}\log\frac{f(x_{r+1})}{p^{m}}.
\end{split}
\]

Thus from Lemma~\ref{le:chap_2_3_3} we obtain
\[
\sum_{p\in Q_{f}(r,x_{r+1})}\log\frac{f(x_{r+1})}{p}
=\frac{1}{2\pi i}\int_{(2)}\frac{f^{s}(x_{r+1})}{s^{2}}\log\zeta_{f}(r,x_{r+1},s)ds+O(\frac{\log x_{r+1}}{x_{r+1}}).
\]

This completes the proof of the lemma.
\end{proof}

\begin{lemma}[Theorem 7-24 in~\cite{FD}]\label{le:chap_2_3_5_th_7_24}
For $(h,M_{r})=1$ and $p=\bar{f}(q)+h\in Q_{f}(r,x_{r+1})$,
\begin{equation}\label{eq:th_7_24}
\begin{split}
\sum_{p\in Q_{f}(r,x_{r+1}),p\equiv h(mod\texttt{ }M_{r})}\log\frac{f(x_{r+1})}{p}
=&\frac{1}{\varphi(M_{r})}\sum_{p\in Q_{f}(r,x_{r+1})}\log\frac{f(x_{r+1})}{p} \\
 &+O(\frac{\log x_{r+1}}{x_{r+1}}).
\end{split}
\end{equation}
\end{lemma}

\begin{proof}
Using Lemma~\ref{le:chap_2_3_2_th_7_9} and the series expansion for
\[
\begin{split}
\log L_{f}(r,x_{r+1},s,\chi)
=&\sum_{i=1}^{\varphi(M_{r})}\log L_{i}(r,x_{r+1},s,\chi) \\
=&\sum_{m,p\in Q_{f}(r,x_{r+1})}\frac{\chi(p^{m})}{mp^{ms}},
\end{split}
\]
we obtain
\[
\begin{split}
&\frac{1}{2\pi i}\int_{(2)}\frac{f^{s}(x_{r+1})}{s^{2}}\log L_{f}(r,x_{r+1},s,\chi)ds \\
=&\frac{1}{2\pi i}\int_{(2)}\frac{f^{s}(x_{r+1})}{s^{2}}\sum_{m,p\in Q_{f}(r,x_{r+1})}\frac{\chi(p^{m})}{mp^{ms}}ds \\
=&\frac{1}{2\pi i}\sum_{m,p\in Q_{f}(r,x_{r+1})}\frac{\chi(p^{m})}{m}\int_{(2)}\frac{(f(x_{r+1})/p^{m})^{s}}{s^{2}}ds \\
=&\sum_{m,p\in Q_{f}(r,x_{r+1}),p^{m}\leq f(x_{r+1})}\frac{\chi(p^{m})}{m}\log\frac{f(x_{r+1})}{p^{m}} \\
=&\sum_{p\in Q_{f}(r,x_{r+1})}\chi(p)\log\frac{f(x_{r+1})}{p} \\
 &+\sum_{p\in Q_{f}(r,x_{r+1}),m\geq2,p^{m}\leq f(x_{r+1})}\frac{\chi(p^{m})}{m}\log\frac{f(x_{r+1})}{p^{m}} \\
=&\sum_{p\in Q_{f}(r,x_{r+1})}\chi(p)\log\frac{f(x_{r+1})}{p}
+O(\frac{\log x_{r+1}}{x_{r+1}})
\end{split}
\]
where since
\[
|\sum_{p\in Q_{f}(r,x_{r+1}),m\geq2,p^{m}\leq f(x_{r+1})}\frac{\chi(p^{m})}{m}\log\frac{f(x_{r+1})}{p^{m}}|
\]
\[
\leq\sum_{p\in Q_{f}(r,x_{r+1}),m\geq2,p^{m}\leq f(x_{r+1})}\frac{1}{m}\log\frac{f(x_{r+1})}{p^{m}},
\]
from Lemma~\ref{le:chap_2_3_3} there is
\[
\sum_{p\in Q_{f}(r,x_{r+1}),m\geq2,p^{m}\leq f(x_{r+1})}\frac{\chi(p^{m})}{m}\log\frac{f(x_{r+1})}{p^{m}}
=O(\frac{\log x_{r+1}}{x_{r+1}}).
\]
Multiplying by $1/\chi(h)$ and summing over all characters modulo $M_{r}$, we deduce with the help of Lemma~\ref{le:chap_2_3_1_th_6_7} that
\[
\begin{split}
 &\sum_{\chi}\frac{1}{\chi(h)}\sum_{p\in Q_{f}(r,x_{r+1})}\chi(p)\log\frac{f(x_{r+1})}{p} \\
=&\varphi(M_{r})\sum_{p\in Q_{f}(r,x_{r+1}),p\equiv h(mod\texttt{ }M_{r})}\log\frac{f(x_{r+1})}{p} \\
=&\frac{1}{2\pi i}\sum_{\chi}\frac{1}{\chi(h)}
\int_{(2)}\frac{f^{s}(x_{r+1})}{s^{2}}\log L_{f}(r,x_{r+1},s,\chi)ds
+\varphi(M_{r})O(\frac{\log x_{r+1}}{x_{r+1}}).
\end{split}
\]
Thus we have
\begin{equation}\label{eq:th_7_24_15}
\begin{split}
&\sum_{p\in Q_{f}(r,x_{r+1}),p\equiv h(mod\texttt{ }M_{r})}\log\frac{f(x_{r+1})}{p} \\
=&\frac{1}{2\pi i\varphi(M_{r})}\sum_{\chi}\frac{1}{\chi(h)}\int_{(2)}\frac{f^{s}(x_{r+1})}{s^{2}}\log L_{f}(r,x_{r+1},s,\chi)ds+O(\frac{\log x_{r+1}}{x_{r+1}}).
\end{split}
\end{equation}

To estimate the integrals appearing in Eq.~(\ref{eq:th_7_24_15}), we must distinguish two cases. First consider the case $\chi=\chi_{0}$. Since there is $p\dagger M_{r}$ for each $p\in Q_{f}(r,x_{r+1})$, there are
\[
\begin{split}
L_{f}(r,x_{r+1},s,\chi_{0})=&\prod_{p\dagger M_{r},p\in Q_{f}(r,x_{r+1})}(1-p^{-s})^{-1} \\
=&\prod_{p\in Q_{f}(r,x_{r+1})}(1-p^{-s})^{-1}=\zeta_{f}(r,x_{r+1},s),
\end{split}
\]
so that from Lemma~\ref{le:chap_2_3_4_th_7_15} there are
\[
\begin{split}
 &\frac{1}{2\pi i}\frac{1}{\chi_{0}(h)}\int_{(2)}\frac{f^{s}(x_{r+1})}{s^{2}}\log L_{f}(r,x_{r+1},s,\chi_{0})ds\\
=&\frac{1}{2\pi i}\int_{(2)}\frac{f^{s}(x_{r+1})}{s^{2}}\log\zeta_{f}(r,x_{r+1},s)ds \\
=&\sum_{p\in Q_{f}(r,x_{r+1})}\log\frac{f(x_{r+1})}{p}+O(\frac{\log x_{r+1}}{x_{r+1}}).
\end{split}
\]
On the other hand, if $\chi\neq\chi_{0}$, then $L_{f}(r,x_{r+1},s,\chi)$ has no pole at $s=1$, but the other properties used earlier still obtain. Hence similar to the proof of Theorem 7-24 in~\cite{FD}, the function
\[
\psi(s,\chi)=\frac{f^{s}(x_{r+1})}{s^{2}}\log L_{f}(r,x_{r+1},s,\chi)
\]
is regular in the region of the integral, and the integral along the bounders of the region either tends to zero or is
\[
O(\frac{\log x_{r+1}}{x_{r+1}}).
\]
It follows that
\[
\sum_{p\in Q_{f}(r,x_{r+1}),p\equiv h(mod\texttt{ }M_{r})}\log\frac{f(x_{r+1})}{p}
\]
\[
=\frac{1}{\varphi(M_{r})}\sum_{p\in Q_{f}(r,x_{r+1})}\log\frac{f(x_{r+1})}{p}+O(\frac{\log x_{r+1}}{x_{r+1}}).
\]

Then for each $(h,M_{r})=1$ where $h=h_{i}$, $(h_{i},M_{r})=1$, and $i=1,2,\ldots,\varphi(M_{r})$, we obtain
\[
\sum_{p\in Q_{f}(r,x_{r+1}),p\equiv h(mod\texttt{ }M_{r})}\log\frac{f(x_{r+1})}{p}
\]
\[
=\frac{1}{\varphi(M_{r})}\sum_{p\in Q_{f}(r,x_{r+1})}\log\frac{f(x_{r+1})}{p}
+O(\frac{\log x_{r+1}}{x_{r+1}}).
\]

This completes the proof of the lemma.
\end{proof}

\begin{lemma}\label{le:chap_2_3_6_th_7_26}
Given a fixed $x_{r}$ and $0<h<M_{r}$, if $(h,M_{r})=1$ and $p=f(q)=\bar{f}(q)+h\in Q_{f}(r,x_{r+1})$, then
\begin{equation}\label{eq:th_7_26}
\pi_{f}(r,x_{r+1},h)=\frac{\pi_{f}(r,x_{r+1})}{\varphi(M_{r})}+O(1).
\end{equation}
\end{lemma}

\begin{proof}
From Lemma~\ref{le:chap_2_3_4_th_7_15} and Lemma~\ref{le:chap_2_3_5_th_7_24}, there is
\[
\sum_{p\in Q_{f}(r,x_{r+1}),p\equiv h(mod\texttt{ }M_{r})}\log\frac{f(x_{r+1})}{p}
\]
\[
=\frac{1}{\varphi(M_{r})}\sum_{p\in Q_{f}(r,x_{r+1})}\log\frac{f(x_{r+1})}{p}
+O(\frac{\log x_{r+1}}{x_{r+1}}).
\]

In the equality, put
\[
\delta=\delta(x)=x^{-1/2}.
\]
Then since $f(x(1+\delta))\sim f(x)(1+n\delta)$ and $\log(1+n\delta)\sim n\delta$ as $x \to \infty$, we have
\[
\begin{split}
&\sum_{p\in Q_{f}(r,x(1+\delta))}\log\frac{f(x)(1+n\delta)}{p}
-\sum_{p\in Q_{f}(r,x)}\log\frac{f(x)}{p} \\
=&\sum_{p\in Q_{f}(r,x)}\log(1+n\delta)
+\sum_{p=f(mx_{r}),m_{x}<m\leq m_{x}(1+\delta)}\log\frac{f(x)(1+n\delta)}{p} \\
=&\log(1+n\delta)\pi_{f}(r,x)+O(\log(1+n\delta)*\delta m_{x})
\end{split}
\]
where $m_{x}=x/x_{r}$, and
\[
\begin{split}
&\sum_{p\in Q_{f}(r,x(1+\delta)),p\equiv h(mod\texttt{ }M_{r})}\log\frac{f(x)(1+n\delta)}{p}
-\sum_{p\in Q_{f}(r,x),p\equiv h(mod\texttt{ }M_{r})}\log\frac{f(x)}{p} \\
=&\sum_{p\in Q_{f}(r,x),p\equiv h(mod\texttt{ }M_{r})}\log(1+n\delta) \\
 &+\sum_{p=f(mx_{r}),m_{x}<m\leq m_{x}(1+\delta),p\equiv h(mod\texttt{ }M_{r})}\log\frac{f(x)(1+n\delta)}{p} \\
=&\log(1+n\delta)\pi_{f}(r,x,h)+O(\log(1+n\delta)*\delta m_{x}).
\end{split}
\]
So that there are
\[
\begin{split}
\pi_{f}(r,x,h)-\frac{\pi_{f}(r,x)}{\varphi(M_{r})}=&O(\delta m_{x})+O(\frac{\log x}{\delta x}) \\
                                                  =&O(\frac{m_{x}}{\sqrt{x}})+O(\frac{\log x}{\sqrt{x}})
                                                   =O(\frac{m_{x}+\log x}{\sqrt{x}}).
\end{split}
\]

Thus for $x=x_{r+1}=x_{r}^{2}$ we obtain
\[
\pi_{f}(r,x_{r+1},h)=\frac{\pi_{f}(r,x_{r+1})}{\varphi(M_{r})}+O(1).
\]

This completes the proof of the lemma.
\end{proof}

\section{Theorems of Bunyakovsky's polynomials}\label{sec:chap_2_4_13}

\begin{theorem}[Prime number theorem of Bunyakovsky's polynomials]\label{th:chap_2_4_1}
If a polynomial
\begin{equation}\label{eq:bunyakovsky_col_p}
f(x)=a_{n}x^{n}+a_{n-1}x^{n-1}+\cdots+a_{1}x+a_{0}=\bar{f}(x)+a_{0}
\end{equation}
in one variable with positive degree and integer coefficients satisfies the following three necessary conditions:

(1) $a_{0}>0$, $a_{n}>0$, and $a_{i}\geq0$ for $i=1,2,\ldots,n-1$,

(2) the polynomial is irreducible over the integers,

(3) as $m$ runs over the positive integers the numbers $f(m)$ should not share a common factor greater than 1, \\ then $f(m)$ is prime for infinitely many positive integers $m$.
\end{theorem}

\begin{proof}
From Lemma~\ref{le:chap_2_2_2}, there are asymptotically more than $c_{f}x_{r}M_{r}/\log x_{r}$ primes $p=\bar{f}(mx_{r})+h$ for $m=1,2,\ldots,m_{r}$ in the range $(\bar{f}(x_{r}),\bar{f}(x_{r+1}))$ where $c_{f}$ is a positive constant dependent on the polynomial $f(x)$. So there must be
\begin{equation}\label{eq:dirichlet_c_4}
\pi(x_{r+1};x_{r})=\sum_{\bar{f}(\xi)<p<\bar{f}(\xi)+M_{r},x_{r}\leq\xi=mx_{r}<x_{r+1}}1
>c_{f}\frac{x_{r}M_{r}}{\log x_{r}}
\end{equation}
and
\begin{equation}\label{eq:dirichlet_c_5}
\frac{\pi(x_{r+1};x_{r})}{\varphi(M_{r})}
>c_{f}\frac{M_{r}}{\varphi(M_{r})}\frac{x_{r}}{\log x_{r}}
>c_{f}\frac{x_{r}}{\log x_{r}}.
\end{equation}

For any positive number $n\geq2$, let constrain primes $p=\bar{f}(\xi)+h$ where $0<h<M_{r}$ and $gcd(h,M_{r})=1$. By Formulas~(\ref{eq:dirichlet_c_4}-\ref{eq:dirichlet_c_5}) we have
\[
\frac{x_{r+1}}{x_{r}} \to \infty\texttt{ and }
\frac{\pi(x_{r+1};x_{r})}{\varphi(M_{r})} \to \infty\texttt{ as }x_{r} \to \infty.
\]

Since $x_{r+1}/x_{r}=x_{r}$, there is $x_{r+1}/x_{r} \to \infty$ as $x_{r} \to \infty$. Thus $x_{r+1}/x_{r}$ is with a same order of infinity as $x_{r}$.

Hence for any arbitrary small positive value $\epsilon$, there exists a number $r_{\epsilon}$ such that for $r\geq r_{\epsilon}$, there is $x_{r}\geq 1/\epsilon$ and $x_{r+1}\geq x_{r}/\epsilon$, which implies that as $\epsilon \to 0$ there are $x_{r} \to \infty$, $x_{r+1}/x_{r} \to \infty$, and $\pi(x_{r+1};x_{r}) \to \infty$ such that the number of elements of the prime set $Q_{f}(r,x_{r+1})$ is infinite.

Let $\pi_{f}(r,x_{r+1},h)$ count the number of primes $p\in Q_{i}(r,x_{r+1})$ where $i$ is fixed and $gcd(h,M_{r})=1$. By Definitions~(\ref{def:chap_2_1_2}-\ref{def:chap_2_1_3}), for $m=1,2,\ldots$ there are
\[
\pi_{f}(r,x_{r+1})=\pi(x_{r+1};x_{r})=\sum_{\bar{f}(\xi)<p<\bar{f}(\xi)+M_{r},x_{r}\leq\xi=mx_{r}<x_{r+1}}1
\]
and
\[
\pi_{f}(r,x_{r+1},h)=\pi(x_{r+1};x_{r},h)=\sum_{p=\bar{f}(\xi)+h,0<h<M_{r},x_{r}\leq\xi=mx_{r}<x_{r+1}}1.
\]

Then from Lemma~\ref{le:chap_2_3_6_th_7_26} proved in Sec.~(\ref{sec:chap_2_3_12}) we can obtain
\begin{equation}\label{eq:dirichlet_c_6}
\pi_{f}(r,x_{r+1},h)=\frac{\pi_{f}(r,x_{r+1})}{\varphi(M_{r})}+O(1)
\end{equation}
where $0<h<M_{r}$ and $gcd(h,M_{r})=1$ for any arbitrary small positive value $\epsilon$ and a corresponding number $r_{\epsilon}$ when considering a fixed modulus $x_{r}$ for the variable $r\geq r_{\epsilon}$.

By definitions of $\pi(x_{r+1};x_{r})$ and $\pi(x_{r+1};x_{r},h)$ in Sec.~(\ref{sec:chap_2_1_10}), Formula~(\ref{eq:dirichlet_c_6}) becomes
\[
\pi(x_{r+1};x_{r},h)=\frac{\pi(x_{r+1};x_{r})}{\varphi(M_{r})}+O(1).
\]

Thus by Formula~(\ref{eq:dirichlet_c_5}), we obtain that for $0<h<M_{r}$ if $gcd(h,M_{r})=1$ then
\begin{equation}\label{eq:dirichlet_c_7}
\pi(x_{r+1};x_{r},h)\sim\frac{\pi(x_{r+1};x_{r})}{\varphi(M_{r})}
>c_{f}\frac{x_{r}}{\log x_{r}}.
\end{equation}

Hence for any positive number $n$, there exists a large positive number $r_{n}$ such that for each large positive number $r=r_{n},r_{n}+1,r_{n}+2,\ldots$ there are more than $c_{f}x_{r}M_{r}/\log x_{r}$ primes $p=f(mx_{r})$ in the range $(\bar{f}(x_{r}),\bar{f}(x_{r+1}))$. Then $f(m)$ is prime for infinitely many positive integers $m$.

This completes the proof of the theorem.
\end{proof}

\begin{theorem}[Bunyakovsky's Conjecture]\label{th:chap_2_4_2}
Given a polynomial
\begin{equation}\label{eq:bunyakovsky}
f(x)=a_{n}x^{n}+a_{n-1}x^{n-1}+\cdots+a_{1}x+a_{0}
\end{equation}
in one variable with positive degree and integer coefficients, if its integer coefficients satisfy the following three necessary conditions:

(1) $a_{0}>0$, $a_{n}>0$, and $a_{i}\geq0$ for $i=1,2,\ldots,n-1$,

(2) the polynomial is irreducible over the integers,

(3) as $m$ runs over the positive integers the numbers $f(m)$ should not share a common factor greater than 1, \\
or it can be converted to a new polynomial $f(x)$ such that the integer coefficients of the new polynomial $f(x)$ satisfy these three necessary conditions, then $f(m)$ is prime for infinitely many positive integers $m$.
\end{theorem}

\begin{proof}
When the integer coefficients of the polynomial $f(x)$ satisfy these three necessary conditions, or from Lemma~\ref{le:chap_2_2_1} $f(x)$ can be converted to a new polynomial $f(x)$ such that the integer coefficients of the new polynomial $f(x)$ satisfy these three necessary conditions stated in the theorem, then from Theorem~\ref{th:chap_2_4_1} $f(m)$ is prime for infinitely many positive integers $m$.

This completes the proof of the theorem.
\end{proof}

\section{Corollaries of Bunyakovsky's polynomials}\label{sec:chap_2_5_14}

\begin{corollary}[Landau's Problem]\label{co:chap_2_5_1}
There are infinitely many primes of the form $a^{2}+1$.
\end{corollary}

\begin{proof}
Consider $f(x)=a_{2}x^{2}+a_{0}=(2k)^{2}x^{2}+1$ where $k\geq1$, $a_{2}=(2k)^{2}$ and $a_{0}=1$. Then by Definition~(\ref{def:chap_2_1_1}) there are $M_{r}=3kx_{r}^{n}>\sqrt{na_{n}}x_{r}^{n}=\sqrt{2a_{2}}x_{r}^{n}$.

Since $a_{0}$ satisfies $0<a_{0}<M_{r}$ and only $a_{0}$ is odd, from Theorem~\ref{th:chap_2_4_2}, $f(m)=(2k)^{2}m^{2}+1$ is prime for infinitely many positive integers $m$.

Since $f(m)$ is of the form $a^{2}+1$, there are infinitely many primes of the form $a^{2}+1$.

This completes the proof of the corollary.
\end{proof}

A conjecture posed by Hardy $\&$ Littlewood in 1923~\cite{FB} is proved as following:
\begin{corollary}[Hardy $\&$ Littlewood's Conjecture]\label{co:chap_2_5_2}
Let $n$ be an odd integer, $n>1$; let $h$ be a positive integer which is not an $e$th power of an integer; for any factor $e>1$ of $n$. Then there exist infinitely many nature numbers $m$ such that $m^{n}+h$ is a prime.
\end{corollary}

\begin{proof}
Let $f(x)=a_{n}x^{n}+a_{0}=x^{n}+h$ where $a_{n}=1$ and $a_{0}=h$. Since $a_{0}>0$ and $a_{n}>0$, the first necessary condition of Theorem~\ref{th:chap_2_4_2} is satisfied. Since $h$ is not an $e$th power of an integer for any factor $e>1$ of $n$ so that the polynomial is irreducible over the integers, the second necessary condition of Theorem~\ref{th:chap_2_4_2} is satisfied. Since $gcd(a_{n},a_{0})=1$, and $a_{n}$, $a_{0}$ are not both even, the third necessary condition of Theorem~\ref{th:chap_2_4_2} is satisfied.

Then from Theorem~\ref{th:chap_2_4_2} $f(m)$ is prime for infinitely many positive integers $m$. Hence there exist infinitely many nature numbers $m$ such that $m^{n}+h$ is a prime.

This completes the proof of the corollary.
\end{proof}

A conjecture~\cite{FB} of quadratic polynomials is proved as following:
\begin{corollary}[Conjecture F]\label{co:chap_2_5_3}
Let $a>0,b>0,c>0$ be integers such that $gcd(a,b,c)=1$, $b^{2}-4ac$ is not a square and $a+b$, $c$ are not both even. Then there are infinitely many primes of the form $am^{2}+bm+c$.
\end{corollary}

\begin{proof}
Let a polynomial $f(x)=ax^{2}+bx+c$. Since $a>0,b>0,c>0$, the first necessary condition of Theorem~\ref{th:chap_2_4_2} is satisfied. Since $b^{2}-4ac$ is not a square so that the polynomial is irreducible over the integers, the second necessary condition of Theorem~\ref{th:chap_2_4_2} is satisfied. Since $gcd(a,b,c)=1$, and $a+b$, $c$ are not both even, the third necessary condition of Theorem~\ref{th:chap_2_4_2} is satisfied.

Then from Theorem~\ref{th:chap_2_4_2} $f(m)$ is prime for infinitely many positive integers $m$. Since $f(m)$ is of the form $am^{2}+bm+c$, there are infinitely many primes of the form $am^{2}+bm+c$.

This completes the proof of the corollary.
\end{proof}

A conjecture of a kind of special polynomials is proved as following:
\begin{corollary}\label{co:chap_2_5_4}
For each integer $n>1$, a polynomial
\[
f(m)=1+2m+\cdots+(n-1)m^{n-2}+nm^{n-1}
\]
is prime for infinitely many positive integers $m$.
\end{corollary}

\begin{proof}
Since the polynomial is irreducible over the integers and all of its integer coefficients are positive, the first and second necessary conditions of Theorem~\ref{th:chap_2_4_2} are satisfied. Since $gcd(1,2,\ldots,n)=1$, and the number $1$ and the number $2+3+\cdots+n$ are not both even, the third necessary condition of Theorem~\ref{th:chap_2_4_2} is satisfied.

Then from Theorem~\ref{th:chap_2_4_2} $f(m)$ is prime for infinitely many positive integers $m$.

This completes the proof of the corollary.
\end{proof}

A conjecture of another kind of special polynomials is proved as following:
\begin{corollary}\label{co:chap_2_5_5}
For any group of consecutive primes $p_{m},p_{m+1},\ldots,p_{n}$ where $p_{m}$ is the $m$th prime and $m<n$, a polynomial
\[
f(b)=\sum_{k=m}^{n}p_{k}b^{n-k}=p_{m}b^{n-m}+p_{m+1}b^{n-m-1}+\cdots+p_{n-1}b+p_{n}
\]
is prime and greater than $p_{n}$ for infinitely many positive integers $b$.
\end{corollary}

\begin{proof}
Since the polynomial is irreducible over the integers and all of its integer coefficients are positive, the first and second necessary conditions of Theorem~\ref{th:chap_2_4_2} are satisfied. Since $gcd(p_{m},p_{m+1},\ldots,p_{n})=1$, and the number $p_{m}+p_{m+1}+\cdots+p_{n-1}$ and the number $p_{n}$ are not both even, the third necessary condition of Theorem~\ref{th:chap_2_4_2} is satisfied.

Then from Theorem~\ref{th:chap_2_4_2} $f(b)$ is prime and greater than $p_{n}$ for infinitely many positive integers $b$.

This completes the proof of the corollary.
\end{proof}

%% file: chap3.tex
%-----------------------------------------------------------------------
% Beginning of chap3.tex
%-----------------------------------------------------------------------
%
% AMS-LaTeX 1.2 sample file for a monograph, based on amsbook.cls.
% This is a data file input by chapter.tex.
%%%%%%%%%%%%%%%%%%%%%%%%%%%%%%%%%%%%%%%%%%%%%%%%%%%%%%%%%%%%%%%%%%%%

\chapter{Prime number theorem of the Bunyakovsky's polynomials}\label{ch:chap_3}

\section*{Introduction}

The Bunyakovsky conjecture stated in 1857 by the Russian mathematician Viktor Bunyakovsky, asserts when a polynomial $f(x)$ in one variable with positive degree and integer coefficients should have infinitely many prime values for positive integer inputs. Three necessary conditions are

(1) the leading coefficient of $f(x)$ is positive,

(2) the polynomial is irreducible over the integers,

(3) as $m$ runs over the positive integers the numbers $f(m)$ should not share a common factor greater than 1.

Bunyakovsky's conjecture is that these three conditions are sufficient: if $f(x)$ satisfies the three conditions, which may be called a Bunyakovsky's polynomial, then $f(m)$ is prime for infinitely many positive integers $m$~\cite{AD,FB}.

In 1922, Nagell showed that $\lim_{X \to \infty}\pi_{f(x)}(X)/X=0$, i.e.  $\lim_{x \to \infty}\pi_{f}(x)/x=0$, so there are few prime values. In 1931, Heilbronn proved the more precise statement~\cite{AD}:

There exists a positive constant $C$ (dependent on $f(x)$), such that
\[
\pi_{f(x)}(X)\leq C\frac{X^{1/n}}{\log X},\texttt{ for every }X\geq1,
\]
i.e. there exists a positive constant $C^{*}$ (dependent on $f(x)$), such that
\[
\pi_{f}(x)\leq C^{*}\frac{x}{\log x},\texttt{ for every }x\geq1.
\]

From Theorem~\ref{th:chap_2_4_2}, we know that if $f(x)$ is a Bunyakovsky's polynomial, then $f(m)$ is prime for infinitely many positive integers $m$.

Let $\pi_{f}(x)$ denote the number of primes $f(q)$ where $f$ is a Bunyakovsky's polynomial with a positive degree $n$ and $0<q\leq x$. It is proved in this chapter that there are constants $\alpha>0$ and $\alpha_{f}>0$ such that as $x \to \infty$,
\[
\pi_{f}(x)=\frac{\alpha_{f}}{n}\int_{2}^{x}\frac{du}{\log u}+O(xe^{-\frac{1}{2}\alpha\sqrt{\log x}}).
\]

Hence, based on the prime number theorem for Bunyakovsky's polynomials, the prime number theorem for primes of the form $a^{2}+1$ is proved.

\section{Definitions}\label{sec:chap_3_1_15}

\begin{definition}\label{def:chap_3_1_1}
When a polynomial $f(x)$ in one variable with positive degree and integer coefficients which satisfy the following three necessary conditions:

(1) the leading coefficient of $f(x)$ is positive,

(2) the polynomial is irreducible over the integers,

(3) as $m$ runs over the positive integers the numbers $f(m)$ should not share a common factor greater than 1, \\ then let define a number set $Q_{f}(x)$ on the polynomial $f(x)$ by
\[
Q_{f}(x)=\{q|f(q)\texttt{ is prime}, 0<q\leq x\}
\]
and let $\pi_{f}(x)$ count the number of primes $p=f(q),q\in Q_{f}(x)$.
\end{definition}

\begin{definition}\label{def:chap_3_1_2}
With complex variable $s$ and primes $p$, let define a complex function $\zeta_{f}(s)$ on the number set $Q_{f}(\infty)$ by
\[
\zeta_{f}(s)=\prod_{p=f(q),q\in Q_{f}(\infty)}(1-p^{-s})^{-1}.
\]
\end{definition}

Because of the Riemann zeta function
\[
\zeta(s)=\zeta_{f}(r,s)\prod_{p\neq f(q),q\in Q_{f}(r,\infty)}(1-p^{-s})^{-1}
\]
and the relation
\[
\log\zeta(s)=\log\zeta_{f}(r,s)+\sum_{m,p\neq f(q),q\in Q_{f}(r,\infty)}\frac{1}{mp^{ms}},
\]
we have
\[
0<|\zeta_{f}(r,s)/\zeta(s)|<1
\]
and
\[
0<|\log\zeta_{f}(r,s)/\log\zeta(s)|<1.
\]

\begin{definition}\label{def:chap_3_1_3}
Let define a function $\Lambda_{f}(\xi)$ on the number set $Q_{f}(\infty)$ by
\[
\Lambda_{f}(\xi)=\{
       \begin{array}{cc}
         \log p & \texttt{ if }\xi=p^{m},\texttt{ for any }m>0\texttt{ and }p=f(q),q\in Q_{f}(\infty), \\
         0 & otherwise.
       \end{array}
\]
\end{definition}

Differentiating the relation
\[
\log\zeta_{f}(s)=\sum_{m,p=f(q),q\in Q_{f}(\infty)}\frac{1}{mp^{ms}},
\]
we obtain
\[
\frac{\zeta_{f}'(s)}{\zeta_{f}(s)}
=-\sum_{m,p=f(q),q\in Q_{f}(\infty)}\frac{\log p}{p^{ms}}
=-\sum_{\xi=1}^{\infty}\frac{\Lambda_{f}(\xi)}{\xi^{s}}.
\]

\begin{definition}\label{def:chap_3_1_4}
For numbers $q\in Q_{f}(x)$, similar to the proof of the prime number theorem for arithmetic progressions in~\cite{FD}, our object is to get an estimate for
\begin{equation}
\pi_{f}(x)=\sum_{p=f(q),q\in Q_{f}(x)}1=\sum_{\xi=1}^{f(x)}P_{f}(\xi)
\end{equation}
where $P_{f}$ is the characteristic function of the numbers:
\[
P_{f}(\xi)=\{
       \begin{array}{cc}
         1 & \texttt{ if }\xi=f(q),q\in Q_{f}(\infty), \\
         0 & otherwise.
       \end{array}
\]
\end{definition}
While $P_{f}$ itself does not arise in a natural way, the function $P_{f}^{*}$ such that
\[
P_{f}^{*}(\xi)=\{
       \begin{array}{cc}
         \frac{1}{m} & \texttt{ if $\xi=p^{m}$ for some $m$ and $p=f(q),q\in Q_{f}(\infty)$}, \\
         0 & otherwise,
       \end{array}
\]
occurs in the Dirichlet series for $\log\zeta_{f}(s)$:
\begin{equation}\label{eq_pnt:bunyakovsky_zeta}
\log\zeta_{f}(s)=\sum_{m,p=f(q),q\in Q_{f}(\infty)}\frac{1}{mp^{ms}}
=\sum_{\xi=1}^{\infty}\frac{P_{f}^{*}(\xi)}{\xi^{s}}.
\end{equation}
For each fixed $m$, the number of $m$th powers of numbers $p=f(q),q\in Q_{f}(x)$, which do not exceed $f(x)$ is equal to the number of numbers $q\in Q_{f}(x)$ which do not exceed $\sqrt[m]{x}$, so that
\[
\begin{split}
\sum_{\xi=1}^{f(x)}P_{f}^{*}(\xi)=&\sum_{\xi=1}^{f(x)}P_{f}(\xi)
+\frac{1}{2}\sum_{\xi=1}^{\sqrt{f(x)}}P_{f}(\xi)
+\frac{1}{3}\sum_{\xi=1}^{\sqrt[3]{f(x)}}P_{f}(\xi)+\cdots \\
\leq&\sum_{\xi=1}^{f(x)}P_{f}(\xi)
+\frac{1}{2}\sum_{\xi=1}^{f(\sqrt{x})}P_{f}(\xi)
+\frac{1}{3}\sum_{\xi=1}^{f(\sqrt[3]{x})}P_{f}(\xi)+\cdots \\
=&\pi_{f}(x)+\frac{\pi_{f}(\sqrt{x})}{2}+\frac{\pi_{f}(\sqrt[3]{x})}{3}+\cdots,
\end{split}
\]
and since, for $m\geq2$,
\[
\pi_{f}(x^{1/m})<cx^{1/m}\leq c\sqrt{x}=o(\frac{x}{\log x}),
\]
it is to be expected that
\[
\sum_{\xi=1}^{f(x)}P_{f}^{*}(\xi)\sim\pi_{f}(x).
\]

\section{Theorems of Bunyakovsky's polynomials}\label{sec:chap_3_2_16}

\begin{theorem}\label{th:chap_3_2_1}
Given a Bunyakovsky's polynomial $f(x)$, by Definition (\ref{def:chap_3_1_2}), let
\[
\log\zeta_{f}(s)=(x/f(x))^{s}\theta_{f}(x,s)\log\zeta(s).
\]
Then, for $0<s=\sigma\leq1$, there are
\[
0<|\theta_{f}(x,s)|\leq C^{*}n\texttt{ and }\frac{|\theta_{f}(x,s)|}{n\log x}\leq\frac{\pi_{f}(x)}{x}\leq \frac{C^{*}}{\log x}
\]
where $n$ is the positive degree of the polynomial $f(x)$ and $C^{*}$ is a positive constant (dependent on $f(x)$).
\end{theorem}

\begin{proof}
First, by Definition (\ref{def:chap_3_1_2}), there is $0<|\log\zeta_{f}(r,s)/\log\zeta(s)|<1$ such that there is
\[
|\theta_{f}(x,s)|=|\frac{\log\zeta_{f}(s)}{\log\zeta(s)}(\frac{f(x)}{x})^{s}|>0.
\]

Next, let define a function $\zeta_{f}(s,x)$ on the number set $Q_{f}(x)$ by
\[
\log\zeta_{f}(s,x)=\sum_{p=f(q),q\in Q_{f}(x)}(1-p^{-s})^{-1}=\frac{\pi_{f}(x)}{1-\rho_{f}^{-s}(s)}
\]
and define a function $\zeta(s,f(x))$ on the function $\zeta(s)$ by
\[
\log\zeta(s,f(x))=\sum_{p\leq f(x)}(1-p^{-s})^{-1}=\frac{\pi(f(x))}{1-\rho_{p}^{-s}(s)}.
\]
Then there are
\[
\frac{1}{1-\rho_{f}^{-s}(s)}=\frac{\log\zeta_{f}(s,x)}{\pi_{f}(x)}=\frac{1}{\pi_{f}(x)}\sum_{p=f(q),q\in Q_{f}(x)}(1-p^{-s})^{-1}
\]
and
\[
\frac{1}{1-\rho_{p}^{-s}(s)}=\frac{\log\zeta(s,f(x))}{\pi(f(x))}=\frac{1}{\pi(f(x))}\sum_{p\leq f(x)}(1-p^{-s})^{-1}.
\]
It shows that $(1-\rho_{f}^{-s}(s))^{-1}$ and $(1-\rho_{p}^{-s}(s))^{-1}$ are the average values of $(1-p^{-s})^{-1}$ for $\log\zeta_{f}(s,x)$ and $\log\zeta(s,f(x))$, respectively.

Let the elements of the number set $Q_{f}(x)$ be in a sequence
\[
0<q_{1}<q_{2}<\cdots<q_{i-1}<q_{i}<\cdots<q_{i_{x}}\leq x
\]
where $i_{x}=\pi_{f}(x)$, and for $i=2,3,\cdots,i_{x}$ let $m_{i}$ denote the number of primes $p$ satisfying $f(q_{i-1})<p\leq f(q_{i})$. Then for $0<s=\sigma\leq1$ the inequality
\[
(1-f(q_{i})^{-\sigma})^{-1}\leq\frac{1}{m_{i}}\sum_{f(q_{i-1})<p\leq f(q_{i})}(1-p^{-\sigma})^{-1}
\]
always holds for $i=2,3,\cdots,i_{x}$. Thus the inequalities
\[
|\rho_{p}^{s}(s)|\leq|\rho_{f}^{s}(s)|\texttt{ and }|1-\rho_{p}^{-s}(s)|\leq|1-\rho_{f}^{-s}(s)|
\]
always hold. On the other hand, as $x \to \infty$ there are
\[
|\rho_{p}^{s}(s)| \to \infty\texttt{, }
|\rho_{f}^{s}(s)| \to \infty\texttt{ and }\lim_{x \to \infty}|\frac{1-\rho_{p}^{-s}(s)}{1-\rho_{f}^{-s}(s)}|=1.
\]
Thus there should be
\[
|\frac{\log\zeta_{f}(s,x)}{\log\zeta(s,f(x))}|=|\frac{1-\rho_{p}^{-s}(s)}{1-\rho_{f}^{-s}(s)}|\frac{\pi_{f}(x)}{\pi(f(x))}\leq\frac{\pi_{f}(x)}{\pi(f(x))}.
\]

Since as $x \to \infty$ and $0<s=\sigma\leq1$
\[
\begin{split}
|\theta_{f}(x,s)|=&|\frac{\log\zeta_{f}(s)/x^{s}}{\log\zeta(s)/f^{s}(x)}|=\lim_{x \to \infty}|\frac{\log\zeta_{f}(s,x)/x^{s}}{\log\zeta(s,f(x))/f^{s}(x)}| \\
              \leq&\frac{\pi_{f}(x)/x}{\pi(f(x))/f(x)}\leq\frac{(C^{*}x/\log x)/x}{(f(x)/\log f(x))/f(x)}=C^{*}n
\end{split}
\]
where $n$ is the positive degree of the polynomial $f(x)$, the inequalities
\[
0<|\theta_{f}(x,s)|\leq C^{*}n\texttt{ and }\frac{|\theta_{f}(x,s)|}{n\log x}\leq\frac{\pi_{f}(x)}{x}\leq \frac{C^{*}}{\log x}
\]
hold for $0<s=\sigma\leq1$.

Hence for $0<s=\sigma\leq1$ $\theta_{f}(x,s)$ is the ratio of the average density $\log\zeta_{f}(s,x)/x^{s}$ of the distribution of primes $p=f(q),q\in Q_{f}(x)$ in the region $[1,f(x)]$ to the average density $\log\zeta(s,f(x))/f^{s}(x)$ of the distribution of primes $p$ in the range $[1,f(x)]$.

This completes the proof of the theorem.
\end{proof}

\begin{theorem}\label{th:chap_3_2_2_7_15}
There are constants $\alpha>0$ and $\alpha_{f}>0$ such that as $x \to \infty$,
\begin{equation}\label{eq_pnt:bunyakovsky_7_15}
\sum_{p=f(q),q\in Q_{f}(x)}\log\frac{f(x)}{p}
=\alpha_{f}\int_{c}^{1}\frac{x^{s}}{s^{2}}ds+O(xe^{-\alpha\sqrt{\log x}}),
\end{equation}
for some $c$ with $0<c<1$.
\end{theorem}

\begin{proof}
Using Theorem~7-9 in~\cite{FD}:
\[
\frac{1}{2\pi i}\int_{(2)}\frac{y^{s}}{s^{2}}ds=\{
\begin{array}{cc}
  0 & \texttt{ for }0<y<1, \\
  \log y & \texttt{ for }y\geq1,
\end{array}
\]
similar to the proof of Theorem~7-15 in~\cite{FD}, we have
\[
\begin{split}
\frac{1}{2\pi i}\int_{(2)}\frac{f^{s}(x)}{s^{2}}\log\zeta_{f}(s)ds
=&\frac{1}{2\pi i}\int_{(2)}\frac{1}{s^{2}}\sum_{\xi=1}^{\infty}\frac{\Lambda_{f}(\xi)}{\log\xi}(\frac{f(x)}{\xi})^{s}ds \\
=&\frac{1}{2\pi i}\sum_{\xi=1}^{\infty}\frac{\Lambda_{f}(\xi)}{\log\xi}\int_{(2)}\frac{1}{s^{2}}(\frac{f(x)}{\xi})^{s}ds \\
=&\sum_{\xi\leq f(x)}\frac{\Lambda_{f}(\xi)}{\log\xi}\log\frac{f(x)}{\xi} \\
=&\sum_{m,p=f(q),q\in Q_{f}(x),p^{m}\leq f(x)}\frac{1}{m}\log\frac{f(x)}{p^{m}} \\
=&\sum_{p=f(q),q\in Q_{f}(x)}\log\frac{f(x)}{p} \\
 &+\sum_{p=f(q),q\in Q_{f}(x),m\geq2,p^{m}\leq f(x)}\frac{1}{m}\log\frac{f(x)}{p^{m}}.
\end{split}
\]

As noted earlier, the number of terms in the last sum is
\[
\pi_{f}(x^{1/2})+\pi_{f}(x^{1/3})+\cdots<x^{1/2}+x^{1/3}+\cdots+x^{1/u}<ux^{1/2},
\]
where $u$ is the smallest number such that $x^{1/u}<2$. Thus
\[
\pi_{f}(x^{1/2})+\pi_{f}(x^{1/3})+\cdots=O(x^{1/2}\log x),
\]
so that
\[
\sum_{p=f(q),q\in Q_{f}(x)}\log\frac{f(x)}{p}
=\frac{1}{2\pi i}\int_{(2)}\frac{f^{s}(x)}{s^{2}}\log\zeta_{f}(s)ds+O(xe^{-\sqrt{\log x}}),
\]
since
\[
\begin{split}
\sum_{p=f(q),q\in Q_{f}(x),m\geq2,p^{m}\leq f(x)}\frac{1}{m}\log\frac{f(x)}{p^{m}}
\leq&\sum_{p=f(q),q\in Q_{f}(x),m\geq2,p^{m}\leq f(x)}\log f(x) \\
=&O(\sqrt{x}\log x\log f(x))=O(\sqrt{x}\log^{2}x) \\
=&O(xe^{-\sqrt{\log x}}).
\end{split}
\]

By the analysis of the integral $\int_{(2)}(x^{s}/s^{2})\log\zeta(s)ds$ in~\cite{FD}, there is
\[
\int_{(2)}\frac{x^{s}}{s^{2}}\log\zeta(s)ds
=2\pi i\int_{c}^{1}\frac{x^{s}}{s^{2}}ds+O(xe^{-\alpha\sqrt{\log x}}).
\]

By the analysis of the integral $\int_{(2)}(x^{s}/s^{2})\theta_{f}(x,s)\log\zeta(s)ds$, which is similar to the analysis of the integral $\int_{(2)}(x^{s}/s^{2})\log\zeta(s)ds$ in~\cite{FD}, we have
\[
\int_{(2)}\frac{x^{s}}{s^{2}}\theta_{f}(x,s)\log\zeta(s)ds
=2\pi i\int_{c}^{1}\frac{x^{s}}{s^{2}}\theta_{f}(x,s)ds+O(xe^{-\alpha\sqrt{\log x}}),
\]
so that there are
\[
\begin{split}
\sum_{p=f(q),q\in Q_{f}(x)}\log\frac{f(x)}{p}
=&\frac{1}{2\pi i}\int_{(2)}\frac{f^{s}(x)}{s^{2}}\log\zeta_{f}(s)ds+O(xe^{-\sqrt{\log x}}) \\
=&\frac{1}{2\pi i}\int_{(2)}\frac{x^{s}}{s^{2}}\theta_{f}(x,s)\log\zeta(s)ds+O(xe^{-\sqrt{\log x}}) \\
=&\int_{c}^{1}\frac{x^{s}}{s^{2}}\theta_{f}(x,s)ds+O(xe^{-\alpha\sqrt{\log x}})
\end{split}
\]
where by Theorem~\ref{th:chap_3_2_1} for $0<s=\sigma\leq1$, the inequality $0<|\theta_{f}(x,s)|\leq C^{*}n$ holds such that there is $|(x^{s}/s^{2})\theta_{f}(x,s)|\leq C^{*}n|x^{s}/s^{2}|$ where $n$ is the positive degree of the polynomial $f(x)$ and $C^{*}$ is a positive constant (dependent on $f(x)$).

On the other hand, by the proof of Theorem~7-15 in~\cite{FD}, we have
\[
\sum_{p\leq x}\log\frac{x}{p}
=\int_{c}^{1}\frac{x^{s}}{s^{2}}ds+O(xe^{-\alpha\sqrt{\log x}}).
\]

By denoting
\[
\frac{\alpha_{f}}{n}=\lim_{x \to \infty}(\sum_{p=f(q),0<q\leq x}\log\frac{x}{q})/(\sum_{p\leq x}\log\frac{x}{p})
\]
where $\alpha_{f}>0$, there are
\[
\begin{split}
\lim_{x \to \infty}\sum_{p=f(q),q\in Q_{f}(x)}\log\frac{f(x)}{p}=&\lim_{x \to \infty}\sum_{q\in Q_{f}(x)}n\log\frac{x}{q} \\
=&n\lim_{x \to \infty}\sum_{p=f(q),0<q\leq x}\log\frac{x}{q} \\
=&\alpha_{f}\lim_{x \to \infty}\sum_{p\leq x}\log\frac{x}{p}.
\end{split}
\]

Thus there should be
\[
\int_{c}^{1}\frac{x^{s}}{s^{2}}\theta_{f}(x,s)ds\sim\alpha_{f}\int_{c}^{1}\frac{x^{s}}{s^{2}}ds\texttt{ as }x \to \infty,
\]
so that we have
\[
\sum_{p=f(q),q\in Q_{f}(x)}\log\frac{f(x)}{p}
=\alpha_{f}\int_{c}^{1}\frac{x^{s}}{s^{2}}ds+O(xe^{-\alpha\sqrt{\log x}}).
\]

This completes the proof of the theorem.
\end{proof}

\begin{theorem}[Prime number theorem for Bunyakovsky's polynomials]\label{th:chap_3_2_3_7_16}
Let $\pi_{f}(x)$ denote the number of primes $f(q)$ where $f$ is a Bunyakovsky's polynomial with a positive degree $n$ and $0<q\leq x$. There are constants $\alpha>0$ and $\alpha_{f}>0$ such that as $x \to \infty$,
\[
\pi_{f}(x)=\frac{\alpha_{f}}{n}\int_{2}^{x}\frac{du}{\log u}+O(xe^{-\frac{1}{2}\alpha\sqrt{\log x}}).
\]
\end{theorem}

\begin{proof}
In Eq. (\ref{eq_pnt:bunyakovsky_7_15}), put
\[
\delta=\delta(x)=e^{-\frac{1}{2}\alpha\sqrt{\log x}}.
\]
Then since $f(x(1+\delta))\sim f(x)(1+n\delta)$ and $\log(1+n\delta)\sim n\delta$ as $x \to \infty$, we have
\[
\begin{split}
&\sum_{p=f(q),q\in Q_{f}(x(1+\delta))}\log\frac{f(x)(1+n\delta)}{p}
-\sum_{p=f(q),q\in Q_{f}(x)}\log\frac{f(x)}{p} \\
=&\sum_{p=f(q),q\in Q_{f}(x)}\log(1+n\delta)
+\sum_{p=f(q),x<q\leq x(1+\delta)}\log\frac{f(x)(1+n\delta)}{p} \\
=&\log(1+n\delta)\pi_{f}(x)+O(\log(1+n\delta)*\delta x) \\
=&\alpha_{f}\int_{c}^{1}\frac{x^{s}}{s^{2}}((1+\delta)^{s}-1)ds+O(xe^{-\alpha\sqrt{\log x}}),
\end{split}
\]
so that
\[
\pi_{f}(x)
=\frac{\alpha_{f}}{\log(1+n\delta)}\int_{c}^{1}\frac{(1+\delta)^{s}-1}{s^{2}}x^{s}ds
+O(\delta x)+O(\frac{xe^{-\alpha\sqrt{\log x}}}{\delta}).
\]

Now
\[
(1+\delta)^{s}-1=s\delta+\frac{s(s-1)}{2!}(1+\vartheta\delta)^{s-2}\delta^{2},
\]
where $0<\vartheta<1$, so that for $0<s<1$,
\[
|(1+\delta)^{s}-1-s\delta|\leq\frac{s|s-1|}{2}\delta^{2}<\delta^{2}.
\]

Thus, making the change of variable $x^{s}=u$, we obtain
\[
\begin{split}
\int_{c}^{1}\frac{(1+\delta)^{s}-1}{s^{2}}x^{s}ds
=&\delta\int_{c}^{1}\frac{x^{s}}{s}ds+O(\delta^{2}\int_{c}^{1}\frac{x^{s}}{s^{2}}ds) \\
=&\delta\int_{x^{c}}^{x}\frac{du}{\log u}+O(\delta^{2}\int_{c}^{1}x^{s}ds) \\
=&\delta\int_{2}^{x}\frac{du}{\log u}+O(\delta^{2}x).
\end{split}
\]

Finally,
\[
\begin{split}
\pi_{f}(x)
=&\frac{\alpha_{f}\delta}{\log(1+n\delta)}\int_{2}^{x}\frac{du}{\log u}
+O(\delta x)+O(\frac{xe^{-\alpha\sqrt{\log x}}}{\delta}) \\
=&(1+O(\delta))\frac{\alpha_{f}}{n}\int_{2}^{x}\frac{du}{\log u}+O(\delta x)+O(xe^{-\frac{1}{2}\alpha\sqrt{\log x}}) \\
=&\frac{\alpha_{f}}{n}\int_{2}^{x}\frac{du}{\log u}+O(\delta x)
=\frac{\alpha_{f}}{n}\int_{2}^{x}\frac{du}{\log u}+O(xe^{-\frac{1}{2}\alpha\sqrt{\log x}}).
\end{split}
\]

This completes the proof of the theorem.
\end{proof}

\begin{theorem}[Bunyakovsky's Conjecture]\label{th:chap_3_2_4}
If a polynomial $f(x)$ in one variable with positive degree $n$ and integer coefficients which satisfy the following three necessary conditions:

(1) the leading coefficient of $f(x)$ is positive,

(2) the polynomial is irreducible over the integers,

(3) as $m$ runs over the positive integers the numbers $f(m)$ should not share a common factor greater than 1, \\ then there are constants $\alpha>0$ and $\alpha_{f}>0$ such that as $x \to \infty$,
\[
\pi_{f}(x)=\frac{\alpha_{f}}{n}\int_{2}^{x}\frac{du}{\log u}+O(xe^{-\frac{1}{2}\alpha\sqrt{\log x}})
\]
where $\pi_{f}(x)$ denote the number of primes $f(m)$ for all positive integers $m\leq x$, and $f(m)$ is prime for infinitely many positive integers $m$.
\end{theorem}

\begin{proof}
By Theorem~\ref{th:chap_3_2_3_7_16}, there are constants $\alpha>0$ and $\alpha_{f}>0$ such that the number of primes $f(m)$ for all $m\leq x$ is equal to
\[
\pi_{f}(x)=\frac{\alpha_{f}}{n}\int_{2}^{x}\frac{du}{\log u}+O(xe^{-\frac{1}{2}\alpha\sqrt{\log x}}).
\]

Since $\pi_{f}(x) \to \infty$ as $x \to \infty$.  Hence $f(m)$ is prime for infinitely many positive integers $m$.

This completes the proof of the theorem.
\end{proof}

\section{Corollaries of Bunyakovsky's polynomials}\label{sec:chap_3_3_17}

\begin{corollary}[Landau's Problem]\label{co:chap_3_3_1}
For a polynomial $f(x)=x^{2}+1$, let $\pi_{f}(x)$ denote the number of primes $f(m)=m^{2}+1$ for all positive integers $m\leq x$. Then there are constants $\alpha>0$ and $\alpha_{f}>0$ such that as $x \to \infty$,
\[
\pi_{f}(x)=\frac{\alpha_{f}}{2}\int_{2}^{x}\frac{du}{\log u}+O(xe^{-\frac{1}{2}\alpha\sqrt{\log x}})
\]
and there are infinitely many primes of the form $a^{2}+1$.
\end{corollary}

\begin{proof}
Since $f(x)=x^{2}+1$ is a Bunyakovsky's polynomial, then by Theorem~\ref{th:chap_3_2_4}, there are constants $\alpha>0$ and $\alpha_{f}>0$ such that as $x \to \infty$,
\[
\pi_{f}(x)=\frac{\alpha_{f}}{2}\int_{2}^{x}\frac{du}{\log u}+O(xe^{-\frac{1}{2}\alpha\sqrt{\log x}})
\]
and $f(m)=m^{2}+1$ is prime for infinitely many positive integers $m$.

Thus there are infinitely many primes of the form $a^{2}+1$.

This completes the proof of the corollary.
\end{proof}

A conjecture posed by Hardy $\&$ Littlewood in 1923~\cite{FB} is proved and expressed as the following theorem:
\begin{corollary}[Hardy $\&$ Littlewood's Conjecture]\label{co:chap_3_3_2}
Let $n$ be an odd integer, $n>1$; let $h$ be a positive integer which is not an $e$th power of an integer; for any factor $e>1$ of $n$. For a polynomial $f(x)=x^{n}+h$, let $\pi_{f}(x)$ denote the number of primes $f(m)=m^{n}+h$ for all positive integers $m\leq x$. Then there are constants $\alpha>0$ and $\alpha_{f}>0$ such that as $x \to \infty$,
\[
\pi_{f}(x)=\frac{\alpha_{f}}{n}\int_{2}^{x}\frac{du}{\log u}+O(xe^{-\frac{1}{2}\alpha\sqrt{\log x}})
\]
and there exist infinitely many nature numbers $m$ such that $m^{n}+h$ is a prime.
\end{corollary}

\begin{proof}
Let $f(x)=a_{n}x^{n}+a_{0}=x^{n}+h$ where $a_{n}=1$ and $a_{0}=h$. Since $x^{n}+h$ is a Bunyakovsky's polynomial, then by Theorem~\ref{th:chap_3_2_4}, there are constants $\alpha>0$ and $\alpha_{f}>0$ such that as $x \to \infty$,
\[
\pi_{f}(x)=\frac{\alpha_{f}}{n}\int_{2}^{x}\frac{du}{\log u}+O(xe^{-\frac{1}{2}\alpha\sqrt{\log x}})
\]
and $f(m)=m^{n}+h$ is prime for infinitely many positive integers $m$.

Hence there exist infinitely many nature numbers $m$ such that $m^{n}+h$ is a prime.

This completes the proof of the corollary.
\end{proof}

A conjecture~\cite{FB} of quadratic polynomials is proved and expressed as the following theorem:
\begin{corollary}[Conjecture F]\label{co:chap_3_3_3}
Let $a>0,b>0,c>0$ be integers such that $gcd(a,b,c)=1$, $b^{2}-4ac$ is not a square and $a+b$, $c$ are not both even. For a polynomial $f(x)=ax^{2}+bx+c$, let $\pi_{f}(x)$ denote the number of primes $f(m)=am^{2}+bm+c$ for all positive integers $m\leq x$. Then there are constants $\alpha>0$ and $\alpha_{f}>0$ such that as $x \to \infty$,
\[
\pi_{f}(x)=\frac{\alpha_{f}}{2}\int_{2}^{x}\frac{du}{\log u}+O(xe^{-\frac{1}{2}\alpha\sqrt{\log x}})
\]
and there are infinitely many primes of the form $am^{2}+bm+c$.
\end{corollary}

\begin{proof}
Since $a>0,b>0,c>0$, the first necessary condition of Theorem~\ref{th:chap_3_2_4} is satisfied. Since $b^{2}-4ac$ is not a square so that the polynomial is irreducible over the integers, the second necessary condition of Theorem~\ref{th:chap_3_2_4} is satisfied. Since $gcd(a,b,c)=1$, and $a+b$, $c$ are not both even, then the third necessary condition of Theorem~\ref{th:chap_3_2_4} is satisfied.

Then by Theorem~\ref{th:chap_3_2_4}, there are constants $\alpha>0$ and $\alpha_{f}>0$ such that as $x \to \infty$,
\[
\pi_{f}(x)=\frac{\alpha_{f}}{2}\int_{2}^{x}\frac{du}{\log u}+O(xe^{-\frac{1}{2}\alpha\sqrt{\log x}})
\]
and $f(m)$ is prime for infinitely many positive integers $m$. Since $f(m)$ is of the form $am^{2}+bm+c$, thus there are infinitely many primes of the form $am^{2}+bm+c$.

This completes the proof of the corollary.
\end{proof}

A conjecture of a kind of special polynomials is proved and expressed as the following theorem:
\begin{corollary}\label{co:chap_3_3_4}
For each integer $n>1$, a polynomial
\[
f(m)=1+2m+\cdots+(n-1)m^{n-1}+nm^{n-1}+nm^{n-1}.
\]
Let $\pi_{f}(x)$ denote the number of primes $f(m)$ for all positive integers $m\leq x$. Then there are constants $\alpha>0$ and $\alpha_{f}>0$ such that as $x \to \infty$,
\[
\pi_{f}(x)=\frac{\alpha_{f}}{n-1}\int_{2}^{x}\frac{du}{\log u}+O(xe^{-\frac{1}{2}\alpha\sqrt{\log x}})
\]
and $f(m)$ is prime for infinitely many positive integers $m$.
\end{corollary}

\begin{proof}
Since the leading coefficient of the polynomial is positive, the first necessary condition of Theorem~\ref{th:chap_3_2_4} is satisfied. Since the polynomial is irreducible over the integers, the second necessary condition of Theorem~\ref{th:chap_3_2_4} is satisfied. Since $gcd(1,2,\cdots,n)=1$, and $1$, $2+3+\cdots+n$ are not both even, then the third necessary condition of Theorem~\ref{th:chap_3_2_4} is satisfied.

Then by Theorem~\ref{th:chap_3_2_4}, there are constants $\alpha>0$ and $\alpha_{f}>0$ such that as $x \to \infty$,
\[
\pi_{f}(x)=\frac{\alpha_{f}}{n-1}\int_{2}^{x}\frac{du}{\log u}+O(xe^{-\frac{1}{2}\alpha\sqrt{\log x}})
\]
and $f(m)$ is prime for infinitely many positive integers $m$.

This completes the proof of the corollary.
\end{proof}

A conjecture of another kind of special polynomials is proved and expressed as the following theorem:
\begin{corollary}\label{co:chap_3_3_5}
For any group of consecutive primes $p_{m},p_{m+1},\cdots,p_{n}$ where $p_{m}$ is the $m$th prime and $m<n$, a polynomial
\[
f(b)=\sum_{k=m}^{n}p_{k}b^{n-k}=p_{m}b^{n-m}+p_{m+1}b^{n-m-1}+\cdots+p_{n-1}b+p_{n}.
\]
Let $\pi_{f}(x)$ denote the number of primes $f(b)$ for all positive integers $b\leq x$. Then there are constants $\alpha>0$ and $\alpha_{f}>0$ such that as $x \to \infty$,
\[
\pi_{f}(x)=\frac{\alpha_{f}}{n-m}\int_{2}^{x}\frac{du}{\log u}+O(xe^{-\frac{1}{2}\alpha\sqrt{\log x}})
\]
and $f(b)$ is prime and greater than $p_{n}$ for infinitely many positive integers $b$.
\end{corollary}

\begin{proof}
Since the leading coefficient of the polynomial is positive, the first necessary condition of Theorem~\ref{th:chap_3_2_4} is satisfied. Since the polynomial is irreducible over the integers, the second necessary condition of Theorem~\ref{th:chap_3_2_4} is satisfied. Since $gcd(p_{m},p_{m+1},\cdots,p_{n})=1$, and $p_{m}+p_{m+1}+\cdots+p_{n-1}$, $p_{n}$ are not both even, then the third necessary condition of Theorem~\ref{th:chap_3_2_4} is satisfied.

Then by Theorem~\ref{th:chap_3_2_4}, there are constants $\alpha>0$ and $\alpha_{f}>0$ such that as $x \to \infty$,
\[
\pi_{f}(x)=\frac{\alpha_{f}}{n-m}\int_{2}^{x}\frac{du}{\log u}+O(xe^{-\frac{1}{2}\alpha\sqrt{\log x}})
\]
and $f(b)$ is prime and greater than $p_{n}$ for infinitely many positive integers $b$.

This completes the proof of the corollary.
\end{proof}

%% file: chap4.tex
%-----------------------------------------------------------------------
% Beginning of chap4.tex
%-----------------------------------------------------------------------
%
% AMS-LaTeX 1.2 sample file for a monograph, based on amsbook.cls.
% This is a data file input by chapter.tex.
%%%%%%%%%%%%%%%%%%%%%%%%%%%%%%%%%%%%%%%%%%%%%%%%%%%%%%%%%%%%%%%%%%%%

\chapter{The number of integers with $k$ prime factors in short intervals}\label{ch:chap_4}

\section*{Introduction}

Landau discussed products of $k$ prime factors and proved the following theorem in 1900~\cite{SKA,SKB}:
\begin{theorem}[Theorem 437 in~\cite{SKA}]\label{th:chap_4_0_1}
Let $k\geq2$ and consider a positive integer $P_{k}$ which is the product of just $k$ distinct prime factors, i.e.
\[
P_{k}=p_{1}p_{2}\cdots p_{k}.
\]
So there are $\omega(P_{k})=k$ and $\Omega(P_{k})=k$ where $\omega(P_{k})$ and $\Omega(P_{k})$ are the number of different prime factors of $P_{k}$ and the total number of prime factors of $P_{k}$, respectively. Let write $\pi_{k}(x)$ for the number of these (squarefree) $P_{k}\leq x$. Then there is asymptotically
\[
\pi_{k}(x)\sim\frac{x(\log\log x)^{k-1}}{(k-1)!\log x}\texttt{ for }k\geq2.
\]
\end{theorem}

From Theorem~\ref{th:chap_1_1_1}, given a positive real number $\beta$ that determines a real number $x_{\beta}$ by $e(\log x_{\beta})^{3}/x_{\beta}^{0.0327283}=\beta$, let $g(x):=\beta x^{1/2}$ for $x\geq x_{\beta}$ where $g(x):=x^{1/2}$ when let $\beta=1$. Then there are
\[
\frac{\pi(x+g(x))-\pi(x)}{g(x)/\log x}=1+O(\frac{1}{\log x})
\]
and
\[
\lim_{x \to \infty}\frac{\pi(x+g(x))-\pi(x)}{g(x)/\log x}=1.
\]

Based on Theorem~\ref{th:chap_4_0_1} and Theorem~\ref{th:chap_1_1_1}, a theorem for the numbers of integers with $k$ prime factors in short intervals is proved in next section.

\section{Theorem of the number of integers with $k$ prime factors}\label{sec:chap_4_1_18}

\begin{theorem}\label{th:chap_4_1_1}
Given a positive real number $\beta$ that determines a real number $x_{\beta}$ by $e(\log x_{\beta})^{3}/x_{\beta}^{0.0327283}=\beta$, let $g(x):=\beta x^{1/2}$ for $x\geq x_{\beta}$ where $g(x):=x^{1/2}$ when let $\beta=1$. Let $\pi_{k}(x)$ denote the number of positive integers not exceeding $x$, which can be written as the product of exactly $k$ distinct primes. Then for each finite positive number $k\geq1$, there exists a finite positive number $x_{k}$ such that for every positive number $x$ greater than $x_{k}$, the number of integers with $k$ prime factors in the interval $(x,x+g(x)]$ is asymptotically equal to
\[
\pi_{k}(x+g(x))-\pi_{k}(x)=\frac{g(x)(\log\log x)^{k-1}}{(k-1)!\log x}+o(\frac{g(x)(\log\log x)^{k-1}}{(k-1)!\log x}).
\]
\end{theorem}

\begin{proof}
Given a positive real number $\beta$ that determines a real number $x_{\beta}$ by $e(\log x_{\beta})^{3}/x_{\beta}^{0.0327283}=\beta$, let $g(x):=\beta x^{1/2}$ for $x\geq x_{\beta}$ where $g(x):=x^{1/2}$ when let $\beta=1$. Then $x$ satisfies $e(\log x)^{3}/x^{0.0327283}\leq\beta$. From Theorem~\ref{th:chap_1_1_1} we have
\[
\pi(x+g(x))-\pi(x)=\frac{g(x)}{\log x}+o(\frac{g(x)}{\log x}),
\]
so that Theorem~\ref{th:chap_4_1_1} holds for $k=1$ with the definition $0!=1$.

As defined in Theorem~\ref{th:chap_4_0_1} for any finite positive number $k\geq2$, let denote
\[
P_{k}=\prod_{i=1}^{k}p_{i}.
\]

We write $\tau_{k}(x)$ for the number of all products of just $k$ prime factors, which do not exceed $x$ and include those (not squarefree) products. By introducing three auxiliary functions
\[
L_{k}(x)=\sum\frac{1}{p_{1}p_{2}\cdots p_{k}},
\Pi_{k}(x)=\sum1=\sum_{m\leq x}c_{m}
\]
and
\[
\vartheta_{k}(x)=\sum\log(p_{1}p_{2}\cdots p_{k})=\sum_{m\leq x}c_{m}\log m,
\]
in the proof of Theorem~\ref{th:chap_4_0_1}~\cite{SKA,SKB}, Landau obtained
\[
L_{k}(x)\sim(\log\log x)^{k}\texttt{ for }k\geq1,
\]
\[
\vartheta_{k}(x)=\Pi_{k}(x)\log x+\int_{2}^{x}\frac{\Pi_{k}(t)}{t}dt\texttt{ and }\Pi_{k}(t)=O(t),
\]
so that for $k\geq2$ there are
\[
\vartheta_{k}(x)=kx(\log\log x)^{k-1}+o(kx(\log\log x)^{k-1}),
\]
\[
\begin{split}
\Pi_{k}(x)=&\frac{\vartheta_{k}(x)}{\log x}+O(\frac{x}{\log x}) \\
          =&\frac{kx(\log\log x)^{k-1}}{\log x}+o(\frac{kx(\log\log x)^{k-1}}{\log x})
\end{split}
\]
and asymptotically
\[
\pi_{k}(x)=\frac{x(\log\log x)^{k-1}}{(k-1)!\log x}+o(\frac{x(\log\log x)^{k-1}}{(k-1)!\log x})
\]
since there are
\[
k!\pi_{k}(x)\leq\Pi_{k}(x)\leq k!\tau_{k}(x)\texttt{ for }k\geq1,
\]
and
\[
\pi_{k}(x)\sim\tau_{k}(x)\texttt{ as }x \to \infty
\]
where since considering those integers $n\leq x$ in the form of $P_{k}$ with at least two of the $p$ equal, the number of those integers $n$ is $\tau_{k}(x)-\pi_{k}(x)$, and every such $n$ can be expressed in the form of $P_{k}$ with $p_{k-1}=p_{k}$, and so
\[
\tau_{k}(x)-\pi_{k}(x)\leq\sum_{p_{1}p_{2}\cdots p_{k-1}^{2}\leq x}1\leq\sum_{p_{1}p_{2}\cdots p_{k-1}\leq x}1=\Pi_{k-1}(x)\texttt{ for }k\geq2.
\]

By these results proved by Landau~\cite{SKA,SKB} and mentioned above, we obtain
\[
\begin{split}
&\vartheta_{k}(x+g(x))-\vartheta_{k}(x) \\
=&\Pi_{k}(x+g(x))\log(x+g(x))-\Pi_{k}(x)\log x+\int_{x}^{x+g(x)}\frac{\Pi_{k}(t)}{t}dt \\
=&[\Pi_{k}(x+g(x))-\Pi_{k}(x)]\log x+O(g(x)).
\end{split}
\]
Similar to the proof of
\[
F_{k}(x)=\vartheta_{k}(x)-kxL_{k-1}(x)
        =o(x(\log\log x)^{k-1})\texttt{ for }k\geq1
\]
in~\cite{SKA,SKB} by induction with
\[
kF_{k+1}(x)=(k+1)\sum_{p\leq x}F_{k}(\frac{x}{p}),
\]
we use this to prove by induction that
\begin{equation}\label{eq_ksi:landau_o}
F_{k}(x+g(x))-F_{k}(x)=o(g(x)(\log\log x)^{k-1})\texttt{ for }k\geq1.
\end{equation}

First from Theorem~\ref{th:chap_1_1_1} we have
\[
\pi(x+g(x))-\pi(x)=\frac{g(x)}{\log x}+o(\frac{g(x)}{\log x}),
\]
so that as $x \to \infty$ there are
\[
\begin{split}
F_{1}(x+g(x))-F_{1}(x)=&[\vartheta_{1}(x+g(x))-(x+g(x))]-[\vartheta_{1}(x)-x] \\
                      =&\vartheta_{1}(x+g(x))-\vartheta_{1}(x)-g(x) \\
                      =&\vartheta(x+g(x))-\vartheta(x)-g(x) \\
                      =&\sum_{p\leq x+g(x)}\log p-\sum_{p\leq x}\log p-g(x) \\
                      =&\sum_{x<p\leq x+g(x)}\log p-g(x) \\
                   \leq&[\pi(x+g(x))-\pi(x)]\log(x+g(x))-g(x)=o(g(x)).
\end{split}
\]
Thus Eq.~(\ref{eq_ksi:landau_o}) holds for $k=1$.

Next let us suppose that Eq.~(\ref{eq_ksi:landau_o}) holds for $k=K\geq1$ so that
\[
\begin{split}
K[F_{K+1}(x+g(x))-F_{K+1}(x)]=&(K+1)[\sum_{p\leq x+g(x)}F_{K}(\frac{x+g(x)}{p})-\sum_{p\leq x}F_{K}(\frac{x}{p})] \\
=&(K+1)\{\sum_{p\leq x}[F_{K}(\frac{x+g(x)}{p})-F_{K}(\frac{x}{p})] \\
 &+\sum_{x<p\leq x+g(x)}F_{K}(\frac{x+g(x)}{p})\} \\
=&(K+1)[\sum_{p\leq x}o(\frac{g(x)}{p}(\log\log x)^{K-1})+O(g(x))] \\
=&(K+1)[o(g(x)(\log\log x)^{K-1})\sum_{p\leq x}\frac{1}{p}+O(g(x))] \\
=&(K+1)[o(g(x)(\log\log x)^{K})+O(g(x))] \\
=&(K+1)o(g(x)(\log\log x)^{K})
\end{split}
\]
where $O(g(x))=o(g(x)(\log\log x)^{K})$ as $x \to \infty$ for $K\geq1$. Then, for any $\epsilon>0$, there is an $x_{0}=x_{0}(K,\epsilon)$ such that
\[
|F_{K}(x)|<\epsilon x(\log\log x)^{K-1}
\]
for all $x\geq x_{0}$, and by the definition of $F_{K}(x)$, we see that $|F_{K}(x)|\leq D$ for $1\leq x<x_{0}$, where $D$ depends only on $K$ and $\epsilon$, so that since
\[
\frac{x+g(x)}{x_{0}}<x<p\leq x+g(x),
\]
we have
\[
\sum_{x<p\leq x+g(x)}|F_{K}(\frac{x+g(x)}{p})|<D[\pi(x+g(x))-\pi(x)]<Dg(x).
\]

Then since $K+1\leq2K$, there is
\[
F_{K+1}(x+g(x))-F_{K+1}(x)=o(g(x)(\log\log x)^{K}).
\]
This implies Eq.~(\ref{eq_ksi:landau_o}) for $k=K+1$ and it follows for all $k\geq1$ by induction.

Hence for $k\geq2$ as $x \to \infty$, by using the definition of $F_{k}(x)$ there are
\[
\begin{split}
\vartheta_{k}(x+g(x))-\vartheta_{k}(x)=&kg(x)L_{k-1}(x)+F_{k}(x+g(x))-F_{k}(x) \\
                                      =&kg(x)(\log\log x)^{k-1}+o(kg(x)(\log\log x)^{k-1}), \\
\Pi_{k}(x+g(x))-\Pi_{k}(x)=&\frac{\vartheta_{k}(x+g(x))-\vartheta_{k}(x)}{\log x}+O(\frac{g(x)}{\log x}) \\
                          =&\frac{kg(x)(\log\log x)^{k-1}}{\log x}+o(\frac{kg(x)(\log\log x)^{k-1}}{\log x}),
\end{split}
\]
and asymptotically
\[
\pi_{k}(x+g(x))-\pi_{k}(x)=\frac{(\log\log x)^{k-1}}{(k-1)!}\frac{g(x)}{\log x}+o(\frac{(\log\log x)^{k-1}}{(k-1)!}\frac{g(x)}{\log x})
\]
since for $k\geq1$ there are
\[
\begin{split}
k![\pi_{k}(x+g(x))-\pi_{k}(x)]\leq&\Pi_{k}(x+g(x))-\Pi_{k}(x) \\
                              \leq&k![\tau_{k}(x+g(x))-\tau_{k}(x)].
\end{split}
\]

Thus for $k\geq2$ the number of integers with $k$ prime factors in the interval $(x,x+g(x)]$ is asymptotically equal to
\[
\pi_{k}(x+g(x))-\pi_{k}(x)=\frac{(\log\log x)^{k-1}}{(k-1)!}\frac{g(x)}{\log x}+o(\frac{(\log\log x)^{k-1}}{(k-1)!}\frac{g(x)}{\log x}),
\]
since there are
\[
\begin{split}
k![\pi_{k}(x+g(x))-\pi_{k}(x)]\leq&\Pi_{k}(x+g(x))-\Pi_{k}(x) \\
                              \leq&k![\tau_{k}(x+g(x))-\tau_{k}(x)]
\end{split}
\]
for $k\geq1$, and
\[
\pi_{k}(x+g(x))-\pi_{k}(x)\sim\tau_{k}(x+g(x))-\tau_{k}(x)\texttt{ as }x \to \infty
\]
where since considering those integers $n$ in the form of $P_{k}$ with at least two of the $p$ equal in the interval $(x,x+g(x)]$, the number of those integers $n$ is
\[
[\tau_{k}(x+g(x))-\tau_{k}(x)]-[\pi_{k}(x+g(x))-\pi_{k}(x)],
\]
and every such $n$ in the interval $(x,x+g(x)]$ can be expressed in the form of $P_{k}$ with $p_{k-1}=p_{k}$, and so
\[
\begin{split}
&[\tau_{k}(x+g(x))-\tau_{k}(x)]-[\pi_{k}(x+g(x))-\pi_{k}(x)] \\
\leq&\sum_{x<p_{1}p_{2}\cdots p_{k-1}^{2}\leq x+g(x)}1\leq\sum_{x<p_{1}p_{2}\cdots p_{k-1}\leq x+g(x)}1 \\
=&\Pi_{k-1}(x+g(x))-\Pi_{k-1}(x)\texttt{ for }k\geq2.
\end{split}
\]

Thus, for any finite positive number $k\geq1$, the number of integers with $k$ prime factors in the interval $(x,x+g(x)]$ is asymptotically equal to
\[
\pi_{k}(x+g(x))-\pi_{k}(x)=\frac{g(x)(\log\log x)^{k-1}}{(k-1)!\log x}+o(\frac{g(x)(\log\log x)^{k-1}}{(k-1)!\log x}).
\]

Given a positive real number $\beta$ that determines a real number $x_{\beta}$ by
\[
e(\log x_{\beta})^{3}/x_{\beta}^{0.0327283}=\beta,
\]
let $g(x):=\beta x^{1/2}$ for $x\geq x_{\beta}$ where $g(x):=x^{1/2}$ when let $\beta=1$. Then $x$ satisfies $(\log x)^{5/2}/x^{0.0327283}\leq\beta$, and the number of integers with $k$ prime factors in the interval $(x,x+\beta x^{1/2}]$ is asymptotically equal to
\[
\pi_{k}(x+\beta x^{1/2})-\pi_{k}(x)=\frac{\beta x^{1/2}(\log\log x)^{k-1}}{(k-1)!\log x}+o(\frac{\beta x^{1/2}(\log\log x)^{k-1}}{(k-1)!\log x}).
\]

Hence, for any finite positive number $k\geq1$, there exists a finite positive number $x_{k}$ such that for every number $x$ greater than $x_{k}$, there is always at least one product $P_{k}$ in the interval $(x,x+\beta x^{1/2}]$.

This completes the proof of the theorem.
\end{proof}

\section{Extension of Theorem~\ref{th:chap_4_1_1} with an exclusive prime $p=2$}\label{sec:chap_4_2_19}

\begin{theorem}\label{th:chap_4_2_1}
Let $k\geq2$ and consider a positive integer $P_{k,2}$ which is the product of just $k$ distinct prime factors greater than 2, i.e.,
\[
P_{k,2}=p_{1}p_{2}\cdots p_{k}\texttt{ where }2<p_{1}<p_{2}<\cdots<p_{k}.
\]
So there are $\omega(P_{k,2})=k$ and $\Omega(P_{k,2})=k$ where $\omega(P_{k,2})$ and $\Omega(P_{k,2})$ are the number of different prime factors of $P_{k,2}$ and the total number of prime factors of $P_{k,2}$, respectively. Let write $\pi_{k,2}(x)$ for the number of these (squarefree) $P_{k,2}\leq x$. Then there is asymptotically
\[
\pi_{k,2}(x)\sim\frac{x(\log\log x)^{k-1}}{(k-1)!\log x}\texttt{ for }k\geq2.
\]
\end{theorem}

\begin{proof}
The theorem can be proved in two ways.

First, as the proof of Theorem~\ref{th:chap_4_0_1} in~\cite{SKA,SKB}, we can similarly prove this theorem.

Second, consider the product $P_{k}$ of just $k$ distinct prime factors with $p_{1}=2$, i.e.,
\[
P_{k}=p_{1}p_{2}\cdots p_{k}=2P_{k-1,2}.
\]

Then, from Theorem~\ref{th:chap_4_0_1} and definitions in Theorem~\ref{th:chap_4_2_1}, there are
\[
\begin{split}
\pi_{k,2}(x)=&\pi_{k}(x)-\pi_{k-1,2}(x/2) \\
\sim&\frac{x(\log\log x)^{k-1}}{(k-1)!\log x}-\frac{(x/2)(\log\log(x/2))^{k-2}}{(k-2)!\log(x/2)} \\
\sim&\frac{x(\log\log x)^{k-1}}{(k-1)!\log x}(1-\frac{k-1}{2\log\log x}) \\
\sim&\frac{x(\log\log x)^{k-1}}{(k-1)!\log x}.
\end{split}
\]

Then there is asymptotically
\[
\pi_{k,2}(x)\sim\frac{x(\log\log x)^{k-1}}{(k-1)!\log x}\texttt{ for }k\geq2.
\]

This completes the proof of the theorem.
\end{proof}

\begin{theorem}\label{th:chap_4_2_2}
Given a positive real number $\beta$ that determines a real number $x_{\beta}$ by $e(\log x_{\beta})^{3}/x_{\beta}^{0.0327283}=\beta$, let $g(x):=\beta x^{1/2}$ for $x\geq x_{\beta}$ where $g(x):=x^{1/2}$ when let $\beta=1$. Let $\pi_{k,2}(x)$ denote the number of positive integers not exceeding $x$ that can be written as the product $P_{k,2}$ of exactly $k$ distinct primes greater than 2. Then for each finite positive number $k\geq2$, there exists a finite positive number $x_{k}$ such that for every positive number $x$ greater than $x_{k}$, the number of integers with $k$ prime factors greater than 2 in the interval $(x,x+g(x)]$ is asymptotically equal to
\[
\pi_{k,2}(x+g(x))-\pi_{k,2}(x)=\frac{g(x)(\log\log x)^{k-1}}{(k-1)!\log x}+o(\frac{g(x)(\log\log x)^{k-1}}{(k-1)!\log x}).
\]
\end{theorem}

\begin{proof}
The theorem can be proved in two ways.

First, as the proof of Theorem~\ref{th:chap_4_1_1} above, we can similarly prove this theorem.

Second, consider the product $P_{k}$ of just $k$ distinct prime factors with $p_{1}=2$, i.e.,
\[
P_{k}=p_{1}p_{2}\cdots p_{k}=2P_{k-1,2}.
\]

Then, from Theorem~\ref{th:chap_4_1_1} and definitions in Theorem~\ref{th:chap_4_2_1}, there are asymptotically
\[
\begin{split}
&\pi_{k,2}(x+g(x))-\pi_{k,2}(x) \\
=&[\pi_{k}(x+g(x))-\pi_{k-1,2}(x/2+g(x)/2)]-[\pi_{k}(x)-\pi_{k-1,2}(x/2)] \\
=&[\pi_{k}(x+g(x))-\pi_{k}(x)]-[\pi_{k-1,2}(x/2+g(x)/2)-\pi_{k-1,2}(x/2)] \\
\sim&\frac{g(x)(\log\log x)^{k-1}}{(k-1)!\log x}-\frac{g(x)(\log\log(x/2))^{k-2}}{2(k-2)!\log(x/2)} \\
\sim&\frac{g(x)(\log\log x)^{k-1}}{(k-1)!\log x}(1-\frac{k-1}{2\log\log x}) \\
\sim&\frac{g(x)(\log\log x)^{k-1}}{(k-1)!\log x}.
\end{split}
\]

Then there is asymptotically
\[
\pi_{k,2}(x+g(x))-\pi_{k,2}(x)=\frac{g(x)(\log\log x)^{k-1}}{(k-1)!\log x}+o(\frac{g(x)(\log\log x)^{k-1}}{(k-1)!\log x}).
\]

This completes the proof of the theorem.
\end{proof}

\section{Extension of Theorem~\ref{th:chap_4_1_1} with an exclusive prime $p$}\label{sec:chap_4_3_20}

\begin{theorem}\label{th:chap_4_3_1}
Let $k\geq2$ and consider a positive integer $P_{k,p}$ which is the product of just $k$ distinct prime factors not equal to an exclusive prime $p$, i.e.,
\[
P_{k,p}=p_{1}p_{2}\cdots p_{k}\texttt{ where }p_{i}\neq p\texttt{ for }i=1,2,\ldots,k.
\]
So there are $\omega(P_{k,p})=k$ and $\Omega(P_{k,p})=k$ where $\omega(P_{k,p})$ and $\Omega(P_{k,p})$ are the number of different prime factors of $P_{k,p}$ and the total number of prime factors of $P_{k,p}$, respectively. Let write $\pi_{k,p}(x)$ for the number of these (squarefree) $P_{k,p}\leq x$. Then there is asymptotically
\[
\pi_{k,p}(x)\sim\frac{x(\log\log x)^{k-1}}{(k-1)!\log x}\texttt{ for }k\geq2.
\]
\end{theorem}

\begin{proof}
The theorem can be proved in two ways.

First, as the proof of Theorem~\ref{th:chap_4_0_1} in~\cite{SKA,SKB}, we can similarly prove this theorem.

Second, consider the product $P_{k}$ of just $k$ distinct prime factors with $p_{1}=p$, i.e.,
\[
P_{k}=p_{1}p_{2}\cdots p_{k}=pP_{k-1,p}.
\]

Then, from Theorem~\ref{th:chap_4_0_1} and definitions in Theorem~\ref{th:chap_4_3_1}, there are
\[
\begin{split}
\pi_{k,p}(x)=&\pi_{k}(x)-\pi_{k-1,p}(x/p) \\
\sim&\frac{x(\log\log x)^{k-1}}{(k-1)!\log x}-\frac{(x/p)(\log\log(x/p))^{k-2}}{(k-2)!\log(x/p)} \\
\sim&\frac{x(\log\log x)^{k-1}}{(k-1)!\log x}(1-\frac{k-1}{p\log\log x}) \\
\sim&\frac{x(\log\log x)^{k-1}}{(k-1)!\log x}.
\end{split}
\]

Then there is asymptotically
\[
\pi_{k,p}(x)\sim\frac{x(\log\log x)^{k-1}}{(k-1)!\log x}\texttt{ for }k\geq2.
\]

This completes the proof of the theorem.
\end{proof}

\begin{theorem}\label{th:chap_4_3_2}
Given a positive real number $\beta$ that determines a real number $x_{\beta}$ by $e(\log x_{\beta})^{3}/x_{\beta}^{0.0327283}=\beta$, let $g(x):=\beta x^{1/2}$ for $x\geq x_{\beta}$ where $g(x):=x^{1/2}$ when let $\beta=1$. Let $\pi_{k,p}(x)$ denote the number of positive integers not exceeding $x$ that can be written as the product $P_{k,p}$ of exactly $k$ distinct primes not equal to an exclusive prime $p$. Then for each finite positive number $k\geq2$, there exists a finite positive number $x_{k}$ such that for every positive number $x$ greater than $x_{k}$, the number of integers with $k$ prime factors not equal to an exclusive prime $p$ in the interval $(x,x+g(x)]$ is asymptotically equal to
\[
\pi_{k,p}(x+g(x))-\pi_{k,p}(x)=\frac{g(x)(\log\log x)^{k-1}}{(k-1)!\log x}+o(\frac{g(x)(\log\log x)^{k-1}}{(k-1)!\log x}).
\]
\end{theorem}

\begin{proof}
The proof can be in two ways.

First, as the proof of Theorem~\ref{th:chap_4_1_1} above, we can similarly prove this theorem.

Second, consider the product $P_{k}$ of just $k$ distinct prime factors with $p_{1}=p$, i.e.,
\[
P_{k}=p_{1}p_{2}\cdots p_{k}=pP_{k-1,p}.
\]

Then, from Theorem~\ref{th:chap_4_1_1} and definitions in Theorem~\ref{th:chap_4_3_1}, there are asymptotically
\[
\begin{split}
&\pi_{k,p}(x+g(x))-\pi_{k,p}(x) \\
=&[\pi_{k}(x+g(x))-\pi_{k-1,p}(x/p+g(x)/p)]-[\pi_{k}(x)-\pi_{k-1,p}(x/p)] \\
=&[\pi_{k}(x+g(x))-\pi_{k}(x)]-[\pi_{k-1,p}(x/p+g(x)/p)-\pi_{k-1,p}(x/p)] \\
\sim&\frac{g(x)(\log\log x)^{k-1}}{(k-1)!\log x}-\frac{g(x)(\log\log(x/p))^{k-2}}{p(k-2)!\log(x/p)} \\
\sim&\frac{g(x)(\log\log x)^{k-1}}{(k-1)!\log x}(1-\frac{k-1}{p\log\log x}) \\
\sim&\frac{g(x)(\log\log x)^{k-1}}{(k-1)!\log x}.
\end{split}
\]

Then there is asymptotically
\[
\pi_{k,p}(x+g(x))-\pi_{k,p}(x)=\frac{g(x)(\log\log x)^{k-1}}{(k-1)!\log x}+o(\frac{g(x)(\log\log x)^{k-1}}{(k-1)!\log x}).
\]

This completes the proof of the theorem.
\end{proof}

%% file: chap5.tex
%-----------------------------------------------------------------------
% Beginning of chap5.tex
%-----------------------------------------------------------------------
%
% AMS-LaTeX 1.2 sample file for a monograph, based on amsbook.cls.
% This is a data file input by chapter.tex.
%%%%%%%%%%%%%%%%%%%%%%%%%%%%%%%%%%%%%%%%%%%%%%%%%%%%%%%%%%%%%%%%%%%%

\chapter{Prime products of Schinzel's polynomials}\label{ch:chap_5}

\section*{Introduction}

In his joint paper with W. Sierpinski, Schinzel proposed the following conjectures~\cite{AD,FB}:

(1) Let $k\geq2$, let $f_{1}(x),f_{2}(x),\ldots,f_{k}(x)$ be irreducible polynomials over the integers, each with integer coefficients and a positive leading coefficient. Assume that the following condition holds:
(*) there does not exist any integer $d>1$ dividing all the products $f_{1}(j)f_{2}(j)\cdots f_{k}(j)$, for every integer $j$.

Then there exist infinitely many natural numbers $m$ such that all numbers $f_{1}(m),f_{2}(m),\ldots,f_{k}(m)$ are primes.

(2) Under the assumptions for $f_{1}(x),f_{2}(x),\ldots,f_{k}(x)$, there exists a natural number $m$ such that the numbers $f_{1}(m),f_{2}(m)$,$\cdots,f_{k}(m)$ are primes.

Landau discussed products of $k$ prime factors and proved Theorem~\ref{th:chap_4_0_1} in 1900~\cite{SKA,SKB}. Similarly, as the proof of Theorem~\ref{th:chap_4_0_1} in~\cite{SKA}, we have proved Theorem~\ref{th:chap_4_2_1} in Chap.~\ref{ch:chap_4}, which is the extension of Theorem~\ref{th:chap_4_0_1} with an exclusive prime $p=2$. With $P_{k}$ and $\pi_{k}(x)$ replacing $P_{k,2}$ and $\pi_{k,2}(x)$, respectively,   it states as follows: Let $k\geq2$ and consider a positive integer $P_{k}$ which is the product of just $k$ distinct prime factors greater than 2, i.e.,
\[
P_{k}=p_{1}p_{2}\cdots p_{k}\texttt{ where }2<p_{1}<p_{2}<\cdots<p_{k}.
\]
So there are $\omega(P_{k})=k$ and $\Omega(P_{k})=k$ where $\omega(P_{k})$ and $\Omega(P_{k})$ are the number of different prime factors of $P_{k}$ and the total number of prime factors of $P_{k}$, respectively. Let write $\pi_{k}(x)$ for the number of these (squarefree) $P_{k}\leq x$. Then there is asymptotically
\[
\pi_{k}(x)\sim\frac{x(\log\log x)^{k-1}}{(k-1)!\log x}\texttt{ for }k\geq2.
\]

From Theorem~\ref{th:chap_4_1_1}, given a positive real number $\beta$ that determines a real number $x_{\beta}$ by $e(\log x_{\beta})^{3}/x_{\beta}^{0.0327283}=\beta$, let $g(x):=\beta x^{1/2}$ for $x\geq x_{\beta}$ where $g(x):=x^{1/2}$ when let $\beta=1$. Let $\pi_{k}(x)$ denote the number of positive integers not exceeding $x$, which can be written as the product of exactly $k$ distinct primes greater than 2. Then for each finite positive number $k\geq1$, there exists a finite positive number $x_{k}$ such that for every positive number $x$ greater than $x_{k}$, the number of integers with $k$ prime factors in the interval $(x,x+g(x)]$ is asymptotically equal to
\[
\pi_{k}(x+g(x))-\pi_{k}(x)=\frac{g(x)(\log\log x)^{k-1}}{(k-1)!\log x}+o(\frac{g(x)(\log\log x)^{k-1}}{(k-1)!\log x}).
\]

Based on Theorem~\ref{th:chap_4_2_1} and Theorem~\ref{th:chap_4_1_1}, a theorem is proved in this chapter. It states as follows: For any finite positive number $k\geq2$, let $f_{1}(x),f_{2}(x),\ldots,f_{k}(x)$ be irreducible polynomials over the integers, each with integer coefficients and a positive leading coefficient. Assume that the following condition holds: there does not exist any integer $d>1$ dividing all the products $f_{1}(j)f_{2}(j)\cdots f_{k}(j)$, for every integer $j$, and the polynomial
\[
f(x)=\prod_{i=1}^{k}f_{i}(x)=a_{n}x^{n}+a_{n-1}x^{n-1}+\cdots+a_{1}x+a_{0}=\bar{f}(x)+a_{0}
\]
which may be called a Schinzel's polynomial, if its integer coefficients satisfy the above necessary conditions and the following conditions or from Lemma~\ref{le:chap_5_2_1} it can be converted to a new polynomial $f(x)$ such that the integer coefficients of the new polynomial $f(x)$ satisfy the above necessary conditions and the following conditions: $a_{0}>0$ and $a_{0}$ is odd, $a_{n}>0$, and $a_{i}\geq0$ for $i=1,2,\ldots,n-1$, then there exist infinitely many natural numbers $m$ such that all numbers $f_{1}(m),f_{2}(m),\ldots,f_{k}(m)$ are primes.

The key ideas of the proof of the main theorem are as follows:
For each $r=1,2,\ldots$,

(1) define a modulus $M_{r}$ on the range $(\bar{f}(x_{r}),\bar{f}(x_{r+1}))$,

(2) make use of Theorem~\ref{th:chap_4_1_1} to estimate the number of products of $k$ distinct primes greater than 2 in a table formed by some values of integral polynomials $\bar{f}(x)+h$ where $h=1,2,\ldots,M_{r}$ in the range $(\bar{f}(x_{r}),\bar{f}(x_{r+1}))$, then Lemma~\ref{le:chap_5_2_2} states that each row of the table contains at least one product $P_{k}$;

(3) make use of the orthogonality relation for Dirichlet characters to estimate the number of products of $k$ distinct primes $\prod_{i}^{k}p_{i}=f(x)$ in the range $(\bar{f}(x_{r}),\bar{f}(x_{r+1}))$, which is estimated greater than $c_{f}x_{r}M_{r}/\log x_{r}$ where $c_{f}$ is a positive constant dependent on the polynomial $f(x)$, \\
so that there exists a large positive number $r_{n}$ such that for each large positive number $r=r_{n},r_{n}+1,r_{n}+2,\ldots$ there are more than $c_{f}x_{r}M_{r}/\log x_{r}$ products $P_{k}=f(mx_{r})$ of $k$ distinct primes greater than 2 in the range $(\bar{f}(x_{r}),\bar{f}(x_{r+1}))$. These results are stated in Theorem~\ref{th:chap_5_4_1}. Then the main theorem~\ref{th:chap_5_4_2} is proved.

Hence, based on the main theorem, the Schinzel's hypothesis H and the Dickson's conjecture can be proved, so that some conjectures in the number theory can be proved. For instance, it can be proved that there are infinitely many twin primes $p$ and $p+2$, cousin primes $p$ and $p+4$, sexy primes $p$ and $p+6$, and Sophie Germain primes of the form $2p+1$ where $p$ is prime.

\section{Definitions}\label{sec:chap_5_1_21}

Given any finite positive number $k\geq2$, let $P_{k}$ denote a product of $k$ distinct primes greater than 2.

\begin{definition}\label{def:chap_5_1_1}
When a polynomial
\[
f(x)=a_{n}x^{n}+a_{n-1}x^{n-1}+\cdots+a_{1}x+a_{0}=\bar{f}(x)+a_{0}
\]
in one variable with positive degree $n\geq2$ and integer coefficients that satisfy the following necessary conditions:
\[
a_{0}>0\texttt{ and }a_{0}\texttt{ is odd, }a_{n}>0,\texttt{ and }a_{i}\geq0\texttt{ for }i=1,2,\ldots,n-1,
\]
there exists a positive integer $x_{1}$ such that $a_{i}\leq a_{n}x_{1}^{n-i}$ for $i=1,2,\ldots,n-1$,
\[
\sqrt{na_{n}}x_{1}^{n}>a_{0},\texttt{ }gcd(a_{0},x_{1})=1,\texttt{ and }e\log^{3}(na_{n}x_{1}^{n})\leq(na_{n}x_{1}^{n})^{0.0327283}.
\]

For $r=1,2,\ldots$ let $x_{r+1}=x_{r}^{2}$ and each integer $M_{r}=2^{m}\geq\sqrt{na_{n}}x_{r}^{n}>2^{m-1}$. Then for $r=1,2,\ldots$ there are $gcd(a_{0},M_{r})=1$,
\[
\bar{f}(x_{r})=a_{n}x_{r}^{n}+a_{n-1}x_{r}^{n-1}+\cdots+a_{1}x_{r}\leq na_{n}x_{r}^{n},
\]
and
\[
e\log^{3}(na_{n}x_{r+1}^{n})\leq(na_{n}x_{r+1}^{n})^{0.0327283}.
\]
\end{definition}

Since $x_{r+1}=x_{r}^{2}$, there are $x_{r+1}/x_{r}=x_{r}$ and $x_{r+1}/x_{r} \to \infty$ as $x_{r} \to \infty$. Thus $x_{r+1}/x_{r}$ is with a same order of infinity as $x_{r}$.

Hence for any arbitrary small positive value $\epsilon$, there exists a number $r_{\epsilon}$ such that for $r\geq r_{\epsilon}$, there are $x_{r}\geq 1/\epsilon$ and $x_{r+1}\geq x_{r}/\epsilon$, which implies $x_{r} \to \infty$ and $x_{r+1}/x_{r} \to \infty$ as $\epsilon \to 0$.

\begin{definition}\label{def:chap_5_1_2}
With $0<h<M_{r}$ and $gcd(h,M_{r})=1$, let denote
\[
\pi(x;x_{r})=\sum_{\bar{f}(\xi)<P_{k}<\bar{f}(\xi)+M_{r},x_{r}\leq\xi=mx_{r}<x}1
\]
and
\[
\pi(x;x_{r},h)=\sum_{P_{k}=\bar{f}(\xi)+h,0<h<M_{r},x_{r}\leq\xi=mx_{r}<x}1.
\]

The function $\pi(x;x_{r})$ counts the number of products $P_{k}=\bar{f}(\xi)+h$ where $0<h<M_{r},x_{r}\leq\xi=mx_{r}<x$ for $m=1,2,\ldots$ and $h$ varies in the range $[1,M_{r})$, and the function $\pi(x;x_{r},h)$ counts the number of products $P_{k}=\bar{f}(\xi)+h$ where $0<h<M_{r},x_{r}\leq\xi=mx_{r}<x$ for $m=1,2,\ldots$ and $h$ is fixed.
\end{definition}

\begin{definition}\label{def:chap_5_1_3}
Let define number sets $\{P_{k}\}$ on polynomials $\bar{f}(\xi)+h$ by
\[
Q_{i}(r,x)=\{P_{k}|P_{k}=\bar{f}(q)+h_{i},0<h_{i}<M_{r},x_{r}\leq q=mx_{r}<x,gcd(h_{i},M_{r})=1\}
\]
for $i=1,2,\ldots,\varphi(M_{r})$, and with $0<h<M_{r}$ and $gcd(h,M_{r})=1$
\[
Q_{f}(r,x)=\bigcup_{i=1}^{\varphi(M_{r})}Q_{i}(r,x)=\{P_{k}|P_{k}=\bar{f}(q)+h,0<h<M_{r},x_{r}\leq q=mx_{r}<x\}
\]
where $m=1,2,\ldots$ and $\varphi(M_{r})$ is the Euler's totient function.
\end{definition}

Let $\pi_{f}(r,x)$ count the number of products $P_{k}\in Q_{f}(r,x)$ and $\pi_{f}(r,x,h_{i})$ count the number of products $P_{k}\in Q_{i}(r,x)$ for $i=1,2,\ldots,\varphi(M_{r})$. Then there are $\pi(x;x_{r})=\pi_{f}(r,x)$ and $\pi(x;x_{r},h_{i})=\pi_{f}(r,x,h_{i})$ for $i=1,2,\ldots,\varphi(M_{r})$.

\begin{definition}\label{def:chap_5_1_4}
With complex variable $s$ and products $P_{k}$, for $i=1,2,\ldots,\varphi(M_{r})$, let define complex functions on number sets $Q_{i}(r,x)$ by
\[
\zeta_{i}(r,x,s)=\prod_{P_{k}\in Q_{i}(r,x)}(1-P_{k}^{-s})^{-1}.
\]

Let define a complex function $\zeta_{f}(r,x,s)$ on number sets $Q_{i}(r,x)$ and $Q_{f}(r,x)$ by
\[
\zeta_{f}(r,x,s)=\prod_{i=1}^{\varphi(M_{r})}\zeta_{i}(r,x,s)=\prod_{P_{k}\in Q_{f}(r,x)}(1-P_{k}^{-s})^{-1}.
\]
\end{definition}

\begin{definition}\label{def:chap_5_1_5}
Let define functions on number sets $Q_{i}(r,x)$ and $Q_{f}(r,x)$ by
\[
\Lambda_{i}(r,\xi)=\{
       \begin{array}{cc}
         \log P_{k} & \texttt{ if }\xi=P_{k}^{m},\texttt{ for any }m>0\texttt{ and }P_{k}\in Q_{i}(r,x), \\
         0 & otherwise.
       \end{array}
\]
for $i=1,2,\ldots,\varphi(M_{r})$, and
\[
\Lambda_{f}(r,\xi)=\{
       \begin{array}{cc}
         \log P_{k} & \texttt{ if }\xi=P_{k}^{m},\texttt{ for any }m>0\texttt{ and }P_{k}\in Q_{f}(r,x), \\
         0 & otherwise.
       \end{array}
\]
\end{definition}

Differentiating the relations
\[
\log\zeta_{i}(r,x,s)=\sum_{m,P_{k}\in Q_{i}(r,x)}\frac{1}{mP_{k}^{ms}}\texttt{ for }i=1,2,\ldots,\varphi(M_{r})
\]
and
\[
\log\zeta_{f}(r,x,s)=\sum_{m,P_{k}\in Q_{f}(r,x)}\frac{1}{mP_{k}^{ms}},
\]
we obtain
\[
\frac{\zeta_{i}'(r,x,s)}{\zeta_{i}(r,x,s)}
=-\sum_{m,P_{k}\in Q_{i}(r,x)}\frac{\log P_{k}}{P_{k}^{ms}}
=-\sum_{\xi=1}^{\infty}\frac{\Lambda_{i}(r,\xi)}{\xi^{s}}\texttt{ for }i=1,2,\ldots,\varphi(M_{r})
\]
and
\[
\frac{\zeta_{f}'(r,x,s)}{\zeta_{f}(r,x,s)}
=-\sum_{m,P_{k}\in Q_{f}(r,x)}\frac{\log P_{k}}{P_{k}^{ms}}
=-\sum_{\xi=1}^{\infty}\frac{\Lambda_{f}(r,\xi)}{\xi^{s}}.
\]

\begin{definition}\label{def:chap_5_1_6}
For each Dirichlet character $\chi$ and $i=1,2,\ldots,\varphi(M_{r})$, let define functions on number sets $Q_{i}(r,x)$ by
\[
L_{i}(r,x,s,\chi)=\prod_{P_{k}\in Q_{i}(r,x)}(1-\chi(P_{k})P_{k}^{-s})^{-1}.
\]

For each Dirichlet character $\chi$, let define a function $L_{f}(r,x,s,\chi)$ on number sets $Q_{i}(r,x)$ and $Q_{f}(r,x)$ by
\[
L_{f}(r,x,s,\chi)=\prod_{i=1}^{\varphi(M_{r})}L_{i}(r,x,s,\chi)=\prod_{P_{k}\in Q_{f}(r,x)}(1-\chi(P_{k})P_{k}^{-s})^{-1}.
\]
\end{definition}

\begin{definition}\label{def:chap_5_1_7}
For numbers $P_{k}\in Q_{f}(r,x_{r+1})$, similar to the proof of the prime number lemma for arithmetic progressions in~\cite{FD}, our object is to get an estimate for
\begin{equation}
\pi_{f}(r,x_{r+1})=\sum_{P_{k}\in Q_{f}(r,x_{r+1})}1=\sum_{\xi=f(x_{r})}^{f(x_{r+1})}P_{f}(r,\xi)
\end{equation}
where $P_{f}$ is the characteristic function of the numbers:
\[
P_{f}(r,\xi)=\{
       \begin{array}{cc}
         1 & \texttt{ if }\xi\in Q_{f}(r,x_{r+1}), \\
         0 & otherwise.
       \end{array}
\]
\end{definition}
While $P_{f}$ itself does not arise in a natural way, the function $P_{f}^{*}$ such that
\[
P_{f}^{*}(r,\xi)=\{
       \begin{array}{cc}
         \frac{1}{m} & \texttt{ if $\xi=P_{k}^{m}$ for some $m$ and $P_{k}\in Q_{f}(r,x_{r+1})$}, \\
         0 & otherwise,
       \end{array}
\]
occurs in the Dirichlet series for $\log\zeta_{f}(r,x_{r+1},s)$:
\begin{equation}\label{eq_h:def_zeta}
\log\zeta_{f}(r,x_{r+1},s)=\sum_{m,P_{k}\in Q_{f}(r,x_{r+1})}\frac{1}{mP_{k}^{ms}}
=\sum_{\xi=1}^{\infty}\frac{P_{f}^{*}(r,\xi)}{\xi^{s}}.
\end{equation}
For each fixed $m$, the number of $m$th powers of numbers $P_{k}\in Q_{f}(r,x_{r+1})$ is less than or equal to the number of numbers $P_{k}\in Q_{f}(r,x_{r+1})$ satisfying $f(x_{r})<P_{k}<\sqrt[m]{f(x_{r+1})}$, so that with $f^{1/m}(x)<f(x^{1/m})$ for $m\geq2$
\[
\begin{split}
\sum_{\xi=1}^{f(x_{r+1})}P_{f}^{*}(r,\xi)
=&\sum_{\xi=f(x_{r})}^{f(x_{r+1})}P_{f}(r,\xi)
+\frac{1}{2}\sum_{\xi=f(x_{r})}^{\sqrt{f(x_{r+1})}}P_{f}(r,\xi)
+\frac{1}{3}\sum_{\xi=f(x_{r})}^{\sqrt[3]{f(x_{r+1})}}P_{f}(r,\xi)+\cdots \\
\leq&\sum_{\xi=f(x_{r})}^{f(x_{r+1})}P_{f}(r,\xi)
+\frac{1}{2}\sum_{\xi=f(x_{r})}^{f(\sqrt{x_{r+1}})}P_{f}(r,\xi)
+\frac{1}{3}\sum_{\xi=f(x_{r})}^{f(\sqrt[3]{x_{r+1}})}P_{f}(r,\xi)+\cdots \\
=&\pi_{f}(r,x_{r+1})+\frac{\pi_{f}(r,\sqrt{x_{r+1}})}{2}+\frac{\pi_{f}(r,\sqrt[3]{x_{r+1}})}{3}+\cdots.
\end{split}
\]

By Definition~\ref{def:chap_5_1_3} there are
\[
\pi_{f}(r,x_{r+1}^{1/m})\leq\pi_{f}(r,x_{r}^{2/2})=\pi_{f}(r,x_{r})\sim0\texttt{ for }m=2,3,\ldots
\]
such that we can write
\[
\pi_{f}(r,x_{r+1}^{1/m})=O(x_{r+1}^{-m/2})\texttt{ for }m=2,3,\ldots.
\]

Thus there are
\[
\frac{\pi_{f}(r,\sqrt{x_{r+1}})}{2}+\frac{\pi_{f}(r,\sqrt[3]{x_{r+1}})}{3}+\cdots=O(\frac{1}{x_{r+1}}).
\]

Since $x_{r+1}=x_{r}^{2}$, there are $x_{r+1}/x_{r}=x_{r}$ and $x_{r+1}/x_{r} \to \infty$ as $x_{r} \to \infty$. Thus $x_{r+1}/x_{r}$ is with a same order of infinity as $x_{r}$ and for any arbitrary small positive value $\epsilon$, there should exist a number $r_{\epsilon}$ such that for $r\geq r_{\epsilon}$, there are
\[
\frac{\pi_{f}(r,\sqrt{x_{r+1}})}{2}+\frac{\pi_{f}(r,\sqrt[3]{x_{r+1}})}{3}+\cdots
<\epsilon\pi_{f}(r,x_{r+1})
\]
and
\[
\sum_{\xi=f(x_{r})}^{f(x_{r+1})}P_{f}^{*}(r,\xi)\sim\pi_{f}(r,x_{r+1})\texttt{ as }\epsilon \to 0.
\]

\section{Lemmas of Schinzel's polynomials}\label{sec:chap_5_2_22}

\begin{lemma}\label{le:chap_5_2_1}
When a polynomial
\begin{equation}\label{eq_h:schinzel_col_1}
f(x)=a_{n}x^{n}+a_{n-1}x^{n-1}+\cdots+a_{1}x+a_{0}
\end{equation}
in one variable with positive degree and integer coefficients that satisfy $a_{0}\neq0$ and $a_{n}>0$, then $f(x)$ can always be converted to a new polynomial
\begin{equation}\label{eq_h:schinzel_col_1_y}
g(y)=b_{n}y^{n}+b_{n-1}y^{n-1}+\cdots+b_{1}y+b_{0}
\end{equation}
such that $b_{0}>0$ and $b_{0}$ is odd, $b_{n}>0$, and $b_{i}\geq0$ for $i=1,2,\ldots,n-1$.
\end{lemma}

\begin{proof}
If some of $a_{i}$ for $i=1,2,\ldots,n-1$ are negative, then by choosing a suitable positive number $k$ and setting $x=y+k$, $f(x)$ becomes
\[
\begin{split}
f(y+k)=&a_{n}(y+k)^{n}+a_{n-1}(y+k)^{n-1}+\cdots+a_{1}(y+k)+a_{0} \\
      =&a_{n}[y^{n}+C(n,1)ky^{n-1}+\cdots+C(n,n-1)k^{n-1}y+k^{n}] \\
       &+a_{n-1}[y^{n-1}+C(n-1,1)ky^{n-2}+\cdots+C(n-1,n-2)k^{n-2}y+k^{n-1}] \\
       &+\cdots+a_{1}(y+k)+a_{0} \\
      =&a_{n}y^{n}+[a_{n}C(n,1)k+a_{n-1}]y^{n-1}+\cdots \\
       &+[a_{n}C(n,1)k^{n-1}+a_{n-1}C(n-1,1)k^{n-2}+\cdots+a_{1}]y \\
       &+(a_{n}k^{n}+a_{n-1}k^{n-1}+\cdots+a_{1}k+a_{0}) \\
      =&g(y)=b_{n}y^{n}+b_{n-1}y^{n-1}+\cdots+b_{1}y+b_{0}
\end{split}
\]
where
\[
b_{i}=\sum_{j=i}^{n}a_{j}C(j,j-i)k^{j-i}\texttt{ for }i=0,1,2,\ldots,n\texttt{ and }C(n,i)=\frac{n!}{i!(n-i)!},
\]
such that the inequalities $b_{0}>0$ and $b_{0}$ is odd, $b_{n}>0$, and $b_{i}\geq0$ for $i=1,2,\ldots,n-1$ always holds.

Hence without loss of generality, we can assume that the inequalities $a_{0}>0$ and $a_{0}$ is odd, $a_{n}>0$, and $a_{i}\geq0$ for $i=1,2,\ldots,n-1$ always holds.

This completes the proof of the lemma.
\end{proof}

\begin{lemma}[Products of $k$ distinct primes in rows]\label{le:chap_5_2_2}
When a polynomial
\begin{equation}\label{eq_h:schinzel_row_p}
f(x)=\prod_{i=1}^{k}f_{i}(x)=a_{n}x^{n}+a_{n-1}x^{n-1}+\cdots+a_{1}x+a_{0}=\bar{f}(x)+a_{0}
\end{equation}
in one variable with positive degree and integer coefficients that satisfy the following necessary conditions: $a_{0}>0$ and $a_{0}$ is odd, $a_{n}>0$, and $a_{i}\geq0$ for $i=1,2,\ldots,n-1$, by Definition (\ref{def:chap_5_1_1}), if numbers
\[
\bar{f}(mx_{r})+1,\bar{f}(mx_{r})+2,\ldots,\bar{f}(mx_{r})+M_{r}
\]
for $m=1,2,\ldots,m_{r}$ are arranged in $m_{r}$ rows each containing $M_{r}$ numbers:
\begin{equation}\label{tab:schinzel_row_p}
\begin{array}{cccc}
  \bar{f}(x_{r})+1, & \bar{f}(x_{r})+2, & \cdots, & \bar{f}(x_{r})+M_{r} \\
  \bar{f}(2x_{r})+1, & \bar{f}(2x_{r})+2, & \cdots, & \bar{f}(2x_{r})+M_{r} \\
  \cdots, & \cdots, & \cdots, & \cdots \\
  \bar{f}(m_{r}x_{r})+1, & \bar{f}(m_{r}x_{r})+2, & \cdots, & \bar{f}(m_{r}x_{r})+M_{r}
\end{array}
\end{equation}
where $m_{r}=x_{r+1}/x_{r}=x_{r}$ and $gcd(a_{0},M_{r})=1$, then Table~\ref{tab:schinzel_row_p} contains asymptotically more than $c_{f}x_{r}M_{r}/\log x_{r}$ products
\[
P_{k}=\prod_{i=1}^{k}p_{i}=\bar{f}(mx_{r})+h_{c}
\]
and $\pi_{f}(r,x_{r+1})$ is asymptotically greater than $c_{f}x_{r}M_{r}/\log x_{r}$ where $c_{f}$ is a positive constant dependent on the polynomial $f(x)$, $p_{i}$ is prime, and $h_{c}$ is an integer satisfying $0<h_{c}<M_{r}$ and $gcd(h_{c},M_{r})=1$.
\end{lemma}

\begin{proof}
From Theorem~\ref{th:chap_4_1_1}, given a positive real number $\beta$ that determines a real number $y_{\beta}$ by $e(\log y_{\beta})^{3}/y_{\beta}^{0.0327283}=\beta$, let $g(y):=\beta y^{1/2}$ for $y\geq y_{\beta}$ where $g(y):=y^{1/2}$ when let $\beta=1$. Then for every $y\geq y_{\beta}$, there are $e(\log y)^{3}/y^{0.0327283}\leq\beta$ and asymptotically $g(y)/\log y$ product $P_{k}$ of $k$ distinct primes greater than 2 in the interval $(y, y+g(y)]$, i.e., there must be
\[
\pi_{k}(y+g(y))-\pi_{k}(y)=\frac{g(y)}{\log y}(1+O(\frac{1}{\log y}))\sim\frac{g(y)}{\log y}.
\]

Thus for $y=\bar{f}(x_{r+1})$ there must be
\[
g(\bar{f}(x_{r+1}))\geq(\bar{f}(x_{r+1}))^{1/2-0.0327283}e(\log\bar{f}(x_{r+1}))^{3}.
\]

It can be satisfied. Let $g(\bar{f}(x_{r+1}))=M_{r}$. Then by Definition (\ref{def:chap_5_1_1}) we have
\[
\begin{split}
g(\bar{f}(x_{r+1}))\geq&\sqrt{na_{n}}x_{r}^{n}=\sqrt{na_{n}x_{r}^{2n}}=\sqrt{na_{n}x_{r+1}^{n}} \\
\geq&(na_{n}x_{r+1}^{n})^{1/2-0.0327283}e(\log(na_{n}x_{r+1}^{n}))^{3} \\
\geq&(\bar{f}(x_{r+1}))^{1/2-0.0327283}e(\log\bar{f}(x_{r+1}))^{3}.
\end{split}
\]

For $m=1,2,\ldots,m_{r}$, let $y=\bar{f}(mx_{r})$ and $g(y)=M_{r}$. Then there are
\[
\begin{split}
g(y)\geq\sqrt{na_{n}}x_{r}^{n}\geq&(\bar{f}(x_{r+1}))^{1/2-0.0327283}e(\log\bar{f}(x_{r+1}))^{3} \\
\geq&(\bar{f}(mx_{r}))^{1/2-0.0327283}e(\log\bar{f}(mx_{r}))^{3} \\
=&y^{1/2-0.0327283}e(\log y)^{3}
\end{split}
\]
such that for large values of $x_{r}$ we have
\[
\begin{split}
\pi_{k}(y+g(y))-\pi_{k}(y)=&\frac{g(y)}{\log y}(1+O(\frac{1}{\log y})) \\
                       \sim&\frac{g(y)}{\log y}=\frac{M_{r}}{\log(\bar{f}(mx_{r}))} \\
                       \geq&\frac{M_{r}}{\log(na_{n}(mx_{r})^{n})}
                       \geq\frac{\sqrt{na_{n}}x_{r}^{n}}{\log(na_{n}(mx_{r})^{n})}.
\end{split}
\]

Especially, for large values of $x_{r}$ we have
\[
\begin{split}
\pi_{k}(\bar{f}(x_{r})+g(\bar{f}(x_{r})))-\pi_{k}(\bar{f}(x_{r}))
\sim&\frac{g(\bar{f}(x_{r}))}{\log(\bar{f}(x_{r}))}=\frac{M_{r}}{\log(\bar{f}(x_{r}))} \\
\geq&\frac{M_{r}}{\log(na_{n}x_{r}^{n})}\geq\frac{\sqrt{na_{n}}x_{r}^{n}}{\log(na_{n}x_{r}^{n})}
\end{split}
\]
and
\[
\begin{split}
\pi_{k}(\bar{f}(x_{r+1})+g(\bar{f}(x_{r+1})))-\pi_{k}(\bar{f}(x_{r+1}))
\sim&\frac{g(\bar{f}(x_{r+1}))}{\log(\bar{f}(x_{r+1}))}=\frac{M_{r+1}}{\log(\bar{f}(x_{r+1}))} \\
\geq&\frac{M_{r}}{\log(na_{n}x_{r}^{2n})}\geq\frac{\sqrt{na_{n}}x_{r}^{n}}{\log(na_{n}x_{r}^{2n})}.
\end{split}
\]

Finally, for each $m=1,2,\ldots,m_{r}$, there are
\[
\begin{split}
\bar{f}(mx_{r}+x_{r})-\bar{f}(mx_{r})=&a_{n}[(mx_{r}+x_{r})^{n}-(mx_{r})^{n}] \\
                                      &+a_{n-1}[(mx_{r}+x_{r})^{n-1}-(mx_{r})^{n-1}] \\
                                      &+\cdots+a_{1}[(mx_{r}+x_{r})-(mx_{r})] \\
                                  \geq&a_{n}x_{r}^{n}\sum_{k=0}^{n-1}(m+1)^{n-1-k}m^{k}.
\end{split}
\]

Then for $m=1$ there are
\[
\begin{split}
\bar{f}(x_{r}+x_{r})-\bar{f}(x_{r})\geq&a_{n}x_{r}^{n}\sum_{k=0}^{n-1}2^{n-1-k} \\
                                      =&a_{n}x_{r}^{n}(2^{n-1}+2^{n-2}+\cdots+2+1) \\
                                      =&a_{n}x_{r}^{n}(2^{n}-1)>2\sqrt{na_{n}}x_{r}^{n}>M_{r},
\end{split}
\]
and for $m=2,3,\ldots,m_{r}$ there are
\[
\begin{split}
\bar{f}(mx_{r}+x_{r})-\bar{f}(mx_{r})\geq&a_{n}x_{r}^{n}\sum_{k=0}^{n-1}(m+1)^{n-1-k}m^{k} \\
                                        >&a_{n}x_{r}^{n}\sum_{k=0}^{n-1}m^{n-1}=na_{n}x_{r}^{n}m^{n-1} \\
                                        >&2\sqrt{na_{n}}x_{r}^{n}>M_{r}.
\end{split}
\]
Thus Table~\ref{tab:schinzel_row_p} contains no repeated products $P_{k}$. On the other hand, for each product $P_{k}$ in Table~\ref{tab:schinzel_row_p}, there is $gcd(P_{k},M_{r})=1$, which implies that the product $P_{k}$ belongs to a column $h_{P}$ satisfying $0<h_{P}<M_{r}$ and $gcd(h_{P},M_{r})=1$. Thus no product $P_{k}$ is missed out when $\varphi(M_{r})$ columns $h$ satisfying $gcd(h,M_{r})=1$ are getting out of Table~\ref{tab:schinzel_row_p}.

Hence for $m=1,2,\ldots,m_{r}$ each $m$th row of Table~\ref{tab:schinzel_row_p} contains asymptotically $M_{r}/\log(na_{n}(mx_{r})^{n})$ products
\[
P_{k}=\prod_{i=1}^{k}p_{i}=\bar{f}(mx_{r})+h_{c}
\]
where $p_{i}$ is prime and $h_{c}$ is an integer satisfying $0<h_{c}<M_{r}$ and $gcd(h_{c},M_{r})=1$ in the interval $(\bar{f}(mx_{r}),\bar{f}(mx_{r})+M_{r})$ such that Table~\ref{tab:schinzel_row_p} contains asymptotically more than $c_{f}x_{r}M_{r}/\log x_{r}$ products $P_{k}$ and $\pi_{f}(r,x_{r+1})$ is asymptotically greater than $c_{f}x_{r}M_{r}/\log x_{r}$ where $c_{f}$ is a positive constant dependent on the polynomial $f(x)$, since
\[
\begin{split}
\pi_{f}(r,x_{r+1})\sim&\sum_{m=1}^{m_{r}-1}\frac{M_{r}}{\log(na_{n}(mx_{r})^{n})} \\
=&\frac{M_{r}}{\log(na_{n}x_{r}^{2n})}\sum_{m=1}^{m_{r}-1}\frac{\log(na_{n}x_{r}^{2n})}{\log(na_{n}(mx_{r})^{n})}\\
=&\frac{M_{r}}{\log(na_{n}x_{r}^{2n})}\sum_{m=1}^{m_{r}-1}(1+\frac{\log(x_{r}/m)^{n}}{\log(na_{n}(mx_{r})^{n})})\\
>&\frac{m_{r}M_{r}}{\log(na_{n}x_{r}^{2n})}\geq c_{f}\frac{x_{r}M_{r}}{\log x_{r}}
\end{split}
\]
where for $\log(na_{n}x_{r}^{n})>3\log(na_{n}3^{2n})$
\[
\begin{split}
\sum_{m=1}^{m_{r}-1}\frac{\log(x_{r}/m)^{n}}{\log(na_{n}(mx_{r})^{n})}
>&\frac{\log x_{r}^{n}}{\log(na_{n}x_{r}^{n})}+\frac{\log(x_{r}/2)^{n}}{\log(na_{n}(2x_{r})^{n})}+\frac{\log(x_{r}/3)^{n}}{\log(na_{n}(3x_{r})^{n})}\\
=&3-\frac{\log(na_{n})}{\log(na_{n}x_{r}^{n})}-\frac{\log(na_{n}2^{2n})}{\log(na_{n}(2x_{r})^{n})}-\frac{\log(na_{n}3^{2n})}{\log(na_{n}(3x_{r})^{n})}\\
>&3-3\frac{\log(na_{n}3^{2n})}{\log(na_{n}x_{r}^{n})}>2.
\end{split}
\]

The proof of the lemma is completed.
\end{proof}

\section{Lemmas of $k$ prime factors for polynomial progressions}\label{sec:chap_5_3_23}

\subsection{Lemmas to be referenced}

Some lemmas proved in~\cite{FD} are stated in this subsection.

\begin{lemma}[Theorem 6-7 in~\cite{FD}]\label{le:chap_5_3_1_th_6_7}
If $(a,m)=(b,m)=1$, then
\[
\sum_{\chi\in X(m)}\frac{\chi(a)}{\chi(b)}=\{
\begin{array}{cc}
  \varphi(m) & \texttt{ if }a\equiv b(mod\texttt{ }m), \\
  0 & \texttt{ otherwise}.
\end{array}
\]
\end{lemma}

If $f$ is regular on the vertical line $\sigma_{0}+it$, and if
\[
\lim_{a \to \infty,b \to \infty}\int_{\sigma_{0}-ai}^{\sigma_{0}+bi}f(s)ds=
\lim_{a \to \infty,b \to \infty}\int_{-a}^{b}f(\sigma_{0}+it)idt
\]
exists, then we abbreviate this limit to $\int_{(\sigma_{0})}f(s)ds$.

\begin{lemma}[Theorem 7-9 in~\cite{FD}]\label{le:chap_5_3_2_th_7_9}
\[
\frac{1}{2\pi i}\int_{(2)}\frac{y^{s}}{s^{2}}ds=\{
\begin{array}{cc}
  0 & \texttt{ for }0<y<1, \\
  \log y & \texttt{ for }y\geq1.
\end{array}
\]
\end{lemma}

\subsection{Lemmas to be proved}

\begin{lemma}\label{le:chap_5_3_3}
\[
\sum_{P_{k}\in Q_{f}(r,x_{r+1}),m\geq2,P_{k}^{m}\leq f(x_{r+1})}\frac{1}{m}\log\frac{f(x_{r+1})}{P_{k}^{m}}=O(\frac{\log x_{r+1}}{x_{r+1}}).
\]
\end{lemma}

\begin{proof}
By Definition~\ref{def:chap_5_1_3} there are
\[
\pi_{f}(r,x_{r+1}^{1/m})\leq\pi_{f}(r,x_{r}^{2/2})=\pi_{f}(r,x_{r})\sim0\texttt{ for }m=2,3,\ldots
\]
such that we can write
\[
\pi_{f}(r,x_{r+1}^{1/m})=O(x_{r+1}^{-m/2})\texttt{ for }m=2,3,\ldots.
\]

Thus there are
\[
\frac{\pi_{f}(r,\sqrt{x_{r+1}})}{2}+\frac{\pi_{f}(r,\sqrt[3]{x_{r+1}})}{3}+\cdots=O(\frac{1}{x_{r+1}})
\]
and
\[
\begin{split}
&\sum_{P_{k}\in Q_{f}(r,x_{r+1}),m\geq2,P_{k}^{m}\leq f(x_{r+1})}\frac{1}{m}\log\frac{f(x_{r+1})}{P_{k}^{m}} \\
\leq&\sum_{P_{k}\in Q_{f}(r,x_{r+1}),m\geq2,P_{k}^{m}\leq f(x_{r+1})}\frac{1}{m}\log f(x_{r+1}) \\
\leq&(\frac{\pi_{f}(r,x_{r+1}^{1/2})}{2}+\frac{\pi_{f}(r,x_{r+1}^{1/3})}{3}+\cdots)\log f(x_{r+1}) \\
=&O(\frac{\log f(x_{r+1})}{x_{r+1}})=O(\frac{\log x_{r+1}}{x_{r+1}}).
\end{split}
\]

This completes the proof of the lemma.
\end{proof}

\begin{lemma}[Theorem 7-15 in~\cite{FD}]\label{le:chap_5_3_4_th_7_15}
\begin{equation}\label{eq_h:th_7_15}
\begin{split}
\sum_{P_{k}\in Q_{f}(r,x_{r+1})}\log\frac{f(x_{r+1})}{P_{k}}
=&\frac{1}{2\pi i}\int_{(2)}\frac{f^{s}(x_{r+1})}{s^{2}}\log\zeta_{f}(r,x_{r+1},s)ds \\
 &+O(\frac{\log x_{r+1}}{x_{r+1}}).
\end{split}
\end{equation}
\end{lemma}

\begin{proof}
Using Lemma~\ref{le:chap_5_3_2_th_7_9}, similar to the proof of Theorem 7-15 in~\cite{FD}, we have
\[
\begin{split}
\frac{1}{2\pi i}\int_{(2)}\frac{f^{s}(x_{r+1})}{s^{2}}\log\zeta_{f}(r,x_{r+1},s)ds
=&\frac{1}{2\pi i}\int_{(2)}\frac{1}{s^{2}}\sum_{\xi=1}^{\infty}\frac{\Lambda_{f}(r,\xi)}{\log\xi}(\frac{f(x_{r+1})}{\xi})^{s}ds\\
=&\frac{1}{2\pi i}\sum_{\xi=1}^{\infty}\frac{\Lambda_{f}(r,\xi)}{\log\xi}\int_{(2)}\frac{1}{s^{2}}(\frac{f(x_{r+1})}{\xi})^{s}ds\\
=&\sum_{\xi\leq f(x_{r+1})}\frac{\Lambda_{f}(r,\xi)}{\log\xi}\log\frac{f(x_{r+1})}{\xi} \\
=&\sum_{m,P_{k}\in Q_{f}(r,x_{r+1}),P_{k}^{m}\leq f(x_{r+1})}\frac{1}{m}\log\frac{f(x_{r+1})}{P_{k}^{m}} \\
=&\sum_{P_{k}\in Q_{f}(r,x_{r+1})}\log\frac{f(x_{r+1})}{P_{k}} \\
+&\sum_{P_{k}\in Q_{f}(r,x_{r+1}),m\geq2,P_{k}^{m}\leq f(x_{r+1})}\frac{1}{m}\log\frac{f(x_{r+1})}{P_{k}^{m}}.
\end{split}
\]

Thus from Lemma~\ref{le:chap_5_3_3} we obtain
\[
\begin{split}
\sum_{P_{k}\in Q_{f}(r,x_{r+1})}\log\frac{f(x_{r+1})}{P_{k}}=&\frac{1}{2\pi i}\int_{(2)}\frac{f^{s}(x_{r+1})}{s^{2}}\log\zeta_{f}(r,x_{r+1},s)ds \\
&+O(\frac{\log x_{r+1}}{x_{r+1}}).
\end{split}
\]

This completes the proof of the lemma.
\end{proof}

\begin{lemma}[Theorem 7-24 in~\cite{FD}]\label{le:chap_5_3_5_th_7_24}
For $gcd(h,M_{r})=1$ and $P_{k}\in Q_{f}(r,x_{r+1})$,
\begin{equation}\label{eq_h:th_7_24}
\begin{split}
\sum_{P_{k}\in Q_{f}(r,x_{r+1}),P_{k}\equiv h(mod\texttt{ }M_{r})}\log\frac{f(x_{r+1})}{P_{k}}
=&\frac{1}{\varphi(M_{r})}\sum_{P_{k}\in Q_{f}(r,x_{r+1})}\log\frac{f(x_{r+1})}{P_{k}} \\
 &+O(\frac{\log x_{r+1}}{x_{r+1}}).
\end{split}
\end{equation}
\end{lemma}

\begin{proof}
Using Lemma~\ref{le:chap_5_3_2_th_7_9} and the series expansion for
\[
\log L_{f}(r,x_{r+1},s,\chi)
=\sum_{m,P_{k}\in Q_{f}(r,x_{r+1})}\frac{\chi(P_{k}^{m})}{mP_{k}^{ms}},
\]
we obtain
\[
\begin{split}
&\frac{1}{2\pi i}\int_{(2)}\frac{f^{s}(x_{r+1})}{s^{2}}\log L_{f}(r,x_{r+1},s,\chi)ds \\
=&\frac{1}{2\pi i}\int_{(2)}\frac{f^{s}(x_{r+1})}{s^{2}}\sum_{m,P_{k}\in Q_{f}(r,x_{r+1})}\frac{\chi(P_{k}^{m})}{mP_{k}^{ms}}ds \\
=&\frac{1}{2\pi i}\sum_{m,P_{k}\in Q_{f}(r,x_{r+1})}\frac{\chi(P_{k}^{m})}{m}\int_{(2)}\frac{(f(x)/P_{k}^{m})^{s}}{s^{2}}ds \\
=&\sum_{m,P_{k}\in Q_{f}(r,x_{r+1}),P_{k}^{m}\leq f(x_{r+1})}\frac{\chi(P_{k}^{m})\log(f(x)/P_{k}^{m})}{m} \\
=&\sum_{P_{k}\in Q_{f}(r,x_{r+1})}\chi(P_{k})\log\frac{f(x_{r+1})}{P_{k}} \\
 &+\sum_{P_{k}\in Q_{f}(r,x_{r+1}),m\geq2,P_{k}^{m}\leq f(x_{r+1})}\frac{\chi(P_{k}^{m})\log(f(x)/P_{k}^{m})}{m}\\
=&\sum_{P_{k}\in Q_{f}(r,x_{r+1})}\chi(P_{k})\log\frac{f(x_{r+1})}{P_{k}}
+O(\frac{\log x_{r+1}}{x_{r+1}})
\end{split}
\]
where since
\[
|\sum_{P_{k}\in Q_{f}(r,x_{r+1}),m\geq2,P_{k}^{m}\leq f(x_{r+1})}\frac{\chi(P_{k}^{m})}{m}\log\frac{f(x_{r+1})}{P_{k}^{m}}|
\]
\[
\leq\sum_{P_{k}\in Q_{f}(r,x_{r+1}),m\geq2,P_{k}^{m}\leq f(x_{r+1})}\frac{1}{m}\log\frac{f(x_{r+1})}{P_{k}^{m}}
\]
from Lemma~\ref{le:chap_5_3_3} there is
\[
\sum_{P_{k}\in Q_{f}(r,x_{r+1}),m\geq2,P_{k}^{m}\leq f(x_{r+1})}\frac{\chi(P_{k}^{m})}{m}\log\frac{f(x_{r+1})}{P_{k}^{m}}
=O(\frac{\log x_{r+1}}{x_{r+1}}).
\]
Multiplying by $1/\chi(h)$ and summing over all characters modulo $M_{r}$, we deduce with the help of Lemma~\ref{le:chap_5_3_1_th_6_7} that
\[
\begin{split}
&\sum_{\chi}\frac{1}{\chi(h)}\sum_{P_{k}\in Q_{f}(r,x_{r+1})}\chi(P_{k})\log\frac{f(x_{r+1})}{P_{k}} \\
=&\varphi(M_{r})\sum_{P_{k}\in Q_{f}(r,x_{r+1}),P_{k}\equiv h(mod\texttt{ }M_{r})}\log\frac{f(x_{r+1})}{P_{k}} \\
=&\frac{1}{2\pi i}\sum_{\chi}\frac{1}{\chi(h)}
\int_{(2)}\frac{f^{s}(x_{r+1})}{s^{2}}\log L_{f}(r,x_{r+1},s,\chi)ds
+\varphi(M_{r})O(\frac{\log x_{r+1}}{x_{r+1}}).
\end{split}
\]
Thus we have
\begin{equation}\label{eq_h:th_7_24_15}
\sum_{P_{k}\in Q_{f}(r,x_{r+1}),P_{k}\equiv h(mod\texttt{ }M_{r})}\log\frac{f(x_{r+1})}{P_{k}}
\end{equation}
\[
=\frac{1}{2\pi i\varphi(M_{r})}\sum_{\chi}\frac{1}{\chi(h)}\int_{(2)}\frac{f^{s}(x_{r+1})}{s^{2}}\log L_{f}(r,x_{r+1},s,\chi)ds+O(\frac{\log x_{r+1}}{x_{r+1}}).
\]

To estimate the integrals appearing in Eq. (\ref{eq_h:th_7_24_15}), we must distinguish two cases. First consider the case $\chi=\chi_{0}$. Since there is $gcd(P_{k},M_{r})=1$ for each $P_{k}\in Q_{f}(r,x_{r+1})$, there are
\[
\begin{split}
L_{f}(r,x_{r+1},s,\chi_{0})=&\prod_{gcd(P_{k},M_{r})=1,P_{k}\in Q_{f}(r,x_{r+1})}(1-P_{k}^{-s})^{-1} \\
                           =&\prod_{P_{k}\in Q_{f}(r,x_{r+1})}(1-P_{k}^{-s})^{-1}=\zeta_{f}(r,x_{r+1},s),
\end{split}
\]
so that from Lemmas~\ref{le:chap_5_3_4_th_7_15} there is
\[
\begin{split}
&\frac{1}{2\pi i}\frac{1}{\chi_{0}(h)}\int_{(2)}\frac{f^{s}(x_{r+1})}{s^{2}}\log L_{f}(r,x_{r+1},s,\chi_{0})ds \\
=&\frac{1}{2\pi i}\int_{(2)}\frac{f^{s}(x_{r+1})}{s^{2}}\log\zeta_{f}(r,x_{r+1},s)ds \\
=&\sum_{P_{k}\in Q_{f}(r,x_{r+1})}\log\frac{f(x_{r+1})}{P_{k}}+O(\frac{\log x_{r+1}}{x_{r+1}}).
\end{split}
\]
On the other hand, if $\chi\neq\chi_{0}$, then $L_{f}(r,x_{r+1},s,\chi)$ has no pole at $s=1$, but the other properties used earlier still obtain. Hence similar to the proof of Theorem 7-24 in~\cite{FD}, the function
\[
\psi(s,\chi)=\frac{f^{s}(x_{r+1})}{s^{2}}\log L_{f}(r,x_{r+1},s,\chi)
\]
is regular in the region of the integral, and the integral along the bounders of the region either tends to zero or is
\[
O(\frac{\log x_{r+1}}{x_{r+1}}).
\]
It follows that
\[
\begin{split}
\sum_{P_{k}\in Q_{f}(r,x_{r+1}),P_{k}\equiv h(mod\texttt{ }M_{r})}\log\frac{f(x_{r+1})}{P_{k}}
=&\frac{1}{\varphi(M_{r})}\sum_{P_{k}\in Q_{f}(r,x_{r+1})}\log\frac{f(x_{r+1})}{P_{k}} \\
 &+O(\frac{\log x_{r+1}}{x_{r+1}}).
\end{split}
\]

Then for each $(h,M_{r})=1$ where $h=h_{i}$, $(h_{i},M_{r})=1$, and $i=1,2,\ldots,\varphi(M_{r})$, we obtain
\[
\begin{split}
\sum_{P_{k}\in Q_{f}(r,x_{r+1}),P_{k}\equiv h(mod\texttt{ }M_{r})}\log\frac{f(x_{r+1})}{P_{k}}
=&\frac{1}{\varphi(M_{r})}\sum_{P_{k}\in Q_{f}(r,x_{r+1})}\log\frac{f(x_{r+1})}{P_{k}} \\
 &+O(\frac{\log x_{r+1}}{x_{r+1}}).
\end{split}
\]

This completes the proof of the lemma.
\end{proof}

\begin{lemma}\label{le:chap_5_3_6_th_7_26}
Given a fixed $x_{r}$ and $0<h<M_{r}$, if $(h,M_{r})=1$, then
\begin{equation}\label{eq_h:th_7_26}
\pi_{f}(r,x_{r+1},h)=\frac{\pi_{f}(r,x_{r+1})}{\varphi(M_{r})}+O(1).
\end{equation}
\end{lemma}

\begin{proof}
From Lemma~\ref{le:chap_5_3_4_th_7_15} and Lemma~\ref{le:chap_5_3_5_th_7_24}, there is
\[
\begin{split}
\sum_{P_{k}\in Q_{f}(r,x_{r+1}),P_{k}\equiv h(mod\texttt{ }M_{r})}\log\frac{f(x_{r+1})}{P_{k}}
=&\frac{1}{\varphi(M_{r})}\sum_{P_{k}\in Q_{f}(r,x_{r+1})}\log\frac{f(x_{r+1})}{P_{k}} \\
 &+O(\frac{\log x_{r+1}}{x_{r+1}}).
\end{split}
\]

In the equality, put
\[
\delta=\delta(x)=x^{-1/2}.
\]
Then assume that $f(x)(1+n\delta)$ corresponds to $x(1+\delta)$ and since $\log(1+n\delta)\sim n\delta$ as $x \to \infty$, we have
\[
\begin{split}
&\sum_{P_{k}\in Q_{f}(r,x(1+\delta))}\log\frac{f(x)(1+n\delta)}{P_{k}}
-\sum_{P_{k}\in Q_{f}(r,x)}\log\frac{f(x)}{P_{k}} \\
=&\sum_{P_{k}\in Q_{f}(r,x)}\log(1+n\delta)
+\sum_{P_{k}=f(mx_{r}),m_{x}<m\leq m_{x}(1+\delta)}\log\frac{f(x)(1+n\delta)}{P_{k}} \\
=&\log(1+n\delta)\pi_{f}(r,x)+O(\log(1+n\delta)*\delta m_{x})
\end{split}
\]
where $m_{x}=x/x_{r}$, and
\[
\begin{split}
&\sum_{P_{k}\in Q_{f}(r,x(1+\delta)),P_{k}\equiv h(mod\texttt{ }M_{r})}\log\frac{f(x)(1+n\delta)}{P_{k}}
-\sum_{P_{k}\in Q_{f}(r,x),P_{k}\equiv h(mod\texttt{ }M_{r})}\log\frac{f(x)}{P_{k}} \\
=&\sum_{P_{k}\in Q_{f}(r,x),P_{k}\equiv h(mod\texttt{ }M_{r})}\log(1+n\delta) \\
 &+\sum_{P_{k}=f(mx_{r}),m_{x}<m\leq m_{x}(1+\delta),P_{k}\equiv h(mod\texttt{ }M_{r})}\log\frac{f(x)(1+n\delta)}{P_{k}} \\
=&\log(1+n\delta)\pi_{f}(r,x,h)+O(\log(1+n\delta)*\delta m_{x}).
\end{split}
\]
So that there are
\[
\begin{split}
\pi_{f}(r,x,h)-\frac{\pi_{f}(r,x)}{\varphi(M_{r})}=&O(\delta m_{x})+O(\frac{\log x}{\delta x}) \\
=&O(\frac{m_{x}}{\sqrt{x}})+O(\frac{\log x}{\sqrt{x}})=O(\frac{m_{x}+\log x}{\sqrt{x}}).
\end{split}
\]

Thus for $x=x_{r+1}=x_{r}^{2}$ we obtain
\[
\pi_{f}(r,x_{r+1},h)=\frac{\pi_{f}(r,x_{r+1})}{\varphi(M_{r})}+O(1).
\]

This completes the proof of the lemma.
\end{proof}

\section{Theorems of Schinzel's polynomials}\label{sec:chap_5_4_24}

\begin{theorem}[Prime number theorem of Schinzel's polynomials]\label{th:chap_5_4_1}
Given any finite positive number $k\geq2$, when a polynomial
\begin{equation}\label{eq_h:schinzel_col_p}
f(x)=\prod_{i=1}^{k}f_{i}(x)=a_{n}x^{n}+a_{n-1}x^{n-1}+\cdots+a_{1}x+a_{0}=\bar{f}(x)+a_{0}
\end{equation}
in one variable with positive degree and integer coefficients that satisfy the following three necessary conditions:

(1) $a_{0}>0$ and $a_{0}$ is odd, $a_{n}>0$, and $a_{i}\geq0$ for $i=1,2,\ldots,n-1$,

(2) each polynomial $f_{i}(x)$ for $i=1,2,\ldots,k$ is irreducible over the integers,

(3) as $m$ runs over the positive integers the numbers $f(m)$ should not share a common factor greater than 1, \\ then $f(m)$ is a product of $k$ distinct primes greater than 2 for infinitely many positive integers $m$ so that there exist infinitely many natural numbers $m$ such that all numbers $f_{1}(m),f_{2}(m),\ldots,f_{k}(m)$ are primes.
\end{theorem}

\begin{proof}
From Lemma~\ref{le:chap_5_2_2}, there are asymptotically more than $c_{f}x_{r}M_{r}/\log x_{r}$ products $P_{k}=\bar{f}(mx_{r})+h$ for $m=1,2,\ldots,m_{r}$ in the range $(\bar{f}(x_{r}),\bar{f}(x_{r+1}))$ where $c_{f}$ is a positive constant dependent on the polynomial $f(x)$. So there must be
\begin{equation}\label{eq_h:dirichlet_c_4}
\pi(x_{r+1};x_{r})=\sum_{\bar{f}(\xi)<P_{k}<\bar{f}(\xi)+M_{r},x_{r}\leq\xi=mx_{r}<x_{r+1}}1
>c_{f}\frac{x_{r}M_{r}}{\log x_{r}}
\end{equation}
and
\begin{equation}\label{eq_h:dirichlet_c_5}
\frac{\pi(x_{r+1};x_{r})}{\varphi(M_{r})}
>c_{f}\frac{M_{r}}{\varphi(M_{r})}\frac{x_{r}}{\log x_{r}}
>c_{f}\frac{x_{r}}{\log x_{r}}.
\end{equation}

For any permissible pair of positive numbers $k\geq2$ and $n\geq k$, let constrain products $P_{k}=\bar{f}(x)+h$ where $0<h<M_{r}$ and $gcd(h,M_{r})=1$. By Formulas (\ref{eq_h:dirichlet_c_4}-\ref{eq_h:dirichlet_c_5}) we have
\[
\frac{x_{r+1}}{x_{r}} \to \infty\texttt{ and }
\frac{\pi(x_{r+1};x_{r})}{\varphi(M_{r})} \to \infty\texttt{ as }x_{r} \to \infty.
\]

Since $x_{r+1}=x_{r}^{2}$, there are $x_{r+1}/x_{r}=x_{r}$ and $x_{r+1}/x_{r} \to \infty$ as $x_{r} \to \infty$. Thus $x_{r+1}/x_{r}$ is with a same order of infinity as $x_{r}$.

Hence for any arbitrary small positive value $\epsilon$, there exists a number $r_{\epsilon}$ such that for $r\geq r_{\epsilon}$, there are $x_{r}\geq 1/\epsilon$ and $x_{r+1}\geq x_{r}/\epsilon$, which implies that as $\epsilon \to 0$ there are $x_{r} \to \infty$, $x_{r+1}/x_{r} \to \infty$, and $\pi(x;x_{r}) \to \infty$ such that the number of elements of the number set $Q_{f}(r,x_{r+1})$ is infinite.

Let $\pi_{f}(r,x_{r+1},h)$ count the number of products $P_{k}\in Q_{f}(r,x_{r+1})$ where $k$ is fixed and $gcd(h,M_{r})=1$. By Definitions (\ref{def:chap_5_1_2}-\ref{def:chap_5_1_3}), for $m=1,2,\ldots$ there are
\[
\pi_{f}(r,x_{r+1})=\pi(x;x_{r})=\sum_{\bar{f}(\xi)<P_{k}<\bar{f}(\xi)+M_{r},x_{r}\leq\xi=mx_{r}<x_{r+1}}1
\]
and
\[
\pi_{f}(r,x_{r+1},h)=\pi(x;x_{r},h)=\sum_{P_{k}=\bar{f}(\xi)+h,0<h<x_{r},x_{r}\leq\xi=mx_{r}<x_{r+1}}1.
\]

Then, for any arbitrary small positive value $\epsilon$ and a corresponding number $r_{\epsilon}$, when considering a fixed modulus $x_{r}$ where $r\geq r_{\epsilon}$, from Lemma~\ref{le:chap_5_3_6_th_7_26} proved in Sec. (\ref{sec:chap_5_3_23}) we can obtain that for $0<h<M_{r}$ if $gcd(h,M_{r})=1$ then we have
\begin{equation}\label{eq_h:dirichlet_c_6}
\pi_{f}(r,x_{r+1},h)=\frac{\pi_{f}(r,x_{r+1})}{\varphi(M_{r})}+O(1).
\end{equation}

By definitions of $\pi(x_{r+1};x_{r})$ and $\pi(x_{r+1};x_{r},h)$, Formula (\ref{eq_h:dirichlet_c_6}) becomes
\[
\pi(x_{r+1};x_{r},h)=\frac{\pi(x_{r+1};x_{r})}{\varphi(M_{r})}+O(1).
\]

Thus by Formula (\ref{eq_h:dirichlet_c_6}), we obtain that for $0<h<M_{r}$ if $gcd(h,M_{r})=1$ then
\begin{equation}\label{eq_h:dirichlet_c_7}
\pi(x_{r+1};x_{r},h)\sim\frac{\pi(x_{r+1};x_{r})}{\varphi(M_{r})}
>c_{f}\frac{x_{r}}{\log x_{r}}.
\end{equation}

Hence for any permissible pair of positive numbers $k\geq2$ and $n\geq k$, there exists a large positive number $r_{n}$ such that for each large positive number $r=r_{n},r_{n}+1,r_{n}+2,\ldots$ there are more than $c_{f}x_{r}M_{r}/\log x_{r}$ products $P_{k}=f(mx_{r})$ of $k$ distinct primes greater than 2 in the range $(\bar{f}(x_{r}),\bar{f}(x_{r+1}))$. Then $f(m)$ is a product of $k$ distinct primes greater than 2 for infinitely many positive integers $m$.

Since the number of polynomials $f_{i}(x)$ for $i=1,2,\ldots,k$ is equal to the number of prime factors of the product and each polynomial $f_{i}(x)$ is always greater than 1, by
\[
f(m)-P_{k}=\prod_{i=1}^{k}(f_{i}(m)-p_{i})=0,
\]
there must be $f_{i}(m)=p_{i}$ where $p_{i}$ is prime for $i=1,2,\ldots,k$, so that there exist infinitely many natural numbers $m$ such that all numbers $f_{1}(m),f_{2}(m),\ldots,f_{k}(m)$ are primes.

This completes the proof of the theorem.
\end{proof}

\begin{theorem}[Schinzel's hypothesis H]\label{th:chap_5_4_2}
Given any finite positive number $k\geq2$, let $f_{1}(x),f_{2}(x),\ldots,f_{k}(x)$ be irreducible polynomials over the integers, each with integer coefficients and a positive leading coefficient. Assume that the following condition holds: there does not exist any integer $d>1$ dividing all the products $f_{1}(j)f_{2}(j)\cdots f_{k}(j)$, for every integer $j$, and the polynomial
\begin{equation}\label{eq_h:schinzel}
f(x)=\prod_{i=1}^{k}f_{i}(x)=a_{n}x^{n}+a_{n-1}x^{n-1}+\cdots+a_{1}x+a_{0}
\end{equation}
can be converted to a new polynomial $f(x)$ such that the integer coefficients of the new polynomial $f(x)$ satisfy the above conditions and the following conditions: $a_{0}>0$ and $a_{0}$ is odd, $a_{n}>0$, and $a_{i}\geq0$ for $i=1,2,\ldots,n-1$, then there exist infinitely many natural numbers $m$ such that all numbers $f_{1}(m),f_{2}(m),\ldots,f_{k}(m)$ are primes.
\end{theorem}

\begin{proof}
If the integer coefficients of the polynomial $f(x)$ do not satisfy the necessary conditions stated in the theorem, then from Lemma~\ref{le:chap_5_2_1} $f(x)$ can be converted to a new polynomial $f(x)$ such that the integer coefficients of the new polynomial $f(x)$ satisfy the necessary conditions stated in the theorem.

Since $f_{1}(x),f_{2}(x),\ldots,f_{k}(x)$ are irreducible polynomials over the integers and there does not exist any integer $d>1$ dividing all the products $f_{1}(j)f_{2}(j)\cdots f_{k}(j)$, for every integer $j$, the second and third conditions of Theorem~\ref{th:chap_5_4_1} are satisfied, which are

(2) each polynomial $f_{i}(x)$ for $i=1,2,\ldots,k$ is irreducible over the integers,

(3) as $m$ runs over the positive integers the numbers $f(m)$ should not share a common factor greater than 1.

Since the polynomial $f(x)$ or the new converted polynomial $f(x)$ satisfies the above conditions and the other conditions of Theorem~\ref{th:chap_5_4_1}: $a_{0}>0$ and $a_{0}$ is odd, $a_{n}>0$, and $a_{i}\geq0$ for $i=1,2,\ldots,n-1$, from Theorem~\ref{th:chap_5_4_1} $f(m)$ is a product of $k$ distinct primes for infinitely many positive integers $m$. Since the number of polynomials $f_{i}(x)$ for $i=1,2,\ldots,k$ is equal to the number of prime factors of the product and each polynomial $f_{i}(x)$ is always greater than 1, there must be $f_{i}(m)=p_{i}$ where $p_{i}$ is prime for $i=1,2,\ldots,k$, so that there exist infinitely many natural numbers $m$ such that all numbers $f_{1}(m),f_{2}(m),\ldots,f_{k}(m)$ are primes.

This completes the proof of the theorem.
\end{proof}

\begin{theorem}[Dickson's Conjecture]\label{th:chap_5_4_3}
Given any finite positive number $k>1$, let $f_{i}(x)=b_{i}x+c_{i}$ with $b_{i}$, $c_{i}$ integers, $b_{i}\geq1$ (for $i=1,2,\ldots,k$). Assume that the following condition holds: there does not exist any integer $d>1$ dividing all the products $f_{1}(j)f_{2}(j)\cdots f_{k}(j)$, for every integer $j$, and the polynomial
\begin{equation}\label{eq_h:dickson}
f(x)=\prod_{i=1}^{k}f_{i}(x)=a_{n}x^{n}+a_{n-1}x^{n-1}+\cdots+a_{1}x+a_{0}
\end{equation}
can be converted to a new polynomial $f(x)$ such that the integer coefficients of the new polynomial $f(x)$ satisfy all conditions of Theorem~\ref{th:chap_5_4_1}, then there exist infinitely many natural numbers $m$ such that all numbers $f_{1}(m),f_{2}(m),\ldots,f_{k}(m)$ are primes.
\end{theorem}

\begin{proof}
If the integer coefficients of the polynomial $f(x)$ do not satisfy the necessary conditions stated in the theorem, then from Lemma~\ref{le:chap_5_2_1} $f(x)$ can be converted to a new polynomial $f(x)$ such that the integer coefficients of the new polynomial $f(x)$ satisfy the necessary conditions stated in the theorem.

Since $f_{1}(x),f_{2}(x),\ldots,f_{k}(x)$ are irreducible polynomials over the integers and there does not exist any integer $d>1$ dividing all the products $f_{1}(j)f_{2}(j)\cdots f_{k}(j)$, for every integer $j$, the second and third conditions of Theorem~\ref{th:chap_5_4_1} are satisfied, which are

(2) each polynomial $f_{i}(x)$ for $i=1,2,\ldots,k$ is irreducible over the integers,

(3) as $m$ runs over the positive integers the numbers $f(m)$ should not share a common factor greater than 1.

Since the polynomial $f(x)$ or the new converted polynomial $f(x)$ satisfies the above conditions and the other conditions of Theorem~\ref{th:chap_5_4_1}: $a_{0}>0$ and $a_{0}$ is odd, $a_{n}>0$, and $a_{i}\geq0$ for $i=1,2,\ldots,n-1$, from Theorem~\ref{th:chap_5_4_1} $f(m)$ is a product of $k$ distinct primes for infinitely many positive integers $m$. Since the number of polynomials $f_{i}(x)$ for $i=1,2,\ldots,k$ is equal to the number of prime factors of the product and each polynomial $f_{i}(x)$ is always greater than 1, there must be $f_{i}(m)=p_{i}$ where $p_{i}$ is prime for $i=1,2,\ldots,k$, so that there exist infinitely many natural numbers $m$ such that all numbers $f_{1}(m),f_{2}(m),\ldots,f_{k}(m)$ are primes.

This completes the proof of the theorem.
\end{proof}

\section{Corollaries of Schinzel's polynomials}\label{sec:chap_5_5_25}

A conjecture[2] of linear polynomials is proved and expressed as the following corollary:
\begin{corollary}[Conjecture $D_{1}$]\label{co:chap_5_5_1}
Given any finite positive number $k>1$, let $c_{1}<c_{2}<\cdots<c_{k}$ be nonzero integers and assume that
\[
f_{1}(x)=x+c_{1},f_{2}(x)=x+c_{2},\ldots,f_{k}(x)=x+c_{k}
\]
satisfy condition: there does not exist any integer $d>1$ dividing all the products $f_{1}(j)f_{2}(j)\cdots f_{k}(j)$, for every integer $j$, and the polynomial
\begin{equation}\label{eq_h:conjecture_d1}
f(x)=\prod_{i=1}^{k}f_{i}(x)=a_{n}x^{n}+a_{n-1}x^{n-1}+\cdots+a_{1}x+a_{0}
\end{equation}
can be converted to a new polynomial $f(x)$ such that the integer coefficients of the new polynomial $f(x)$ satisfy all conditions of Theorem~\ref{th:chap_5_4_1}, then there exist infinitely many natural numbers $m\geq1$ such that
\[
m+c_{1},m+c_{2},\ldots,m+c_{k}
\]
are consecutive primes.
\end{corollary}

\begin{proof}
If the integer coefficients of the polynomial $f(x)$ do not satisfy the necessary conditions stated in the lemma, then from Lemma~\ref{le:chap_5_2_1} $f(x)$ can be converted to a new polynomial $f(x)$ such that the integer coefficients of the new polynomial $f(x)$ satisfy the necessary conditions stated in the lemma.

From Theorem~\ref{th:chap_5_4_3}, $f(m)$ is a product of $k$ distinct primes for infinitely many positive integers $m$. Since the number of polynomials $f_{i}(x)$ for $i=1,2,\ldots,k$ is equal to the number of prime factors of the product and each polynomial $f_{i}(x)$ is always greater than 1, there must be $f_{i}(m)=p_{i}$ where $p_{i}$ is prime for $i=1,2,\ldots,k$, so that there exist infinitely many natural numbers $m$ such that all numbers $f_{1}(m),f_{2}(m),\ldots,f_{k}(m)$ are primes and
\[
m+c_{1},m+c_{2},\ldots,m+c_{k}
\]
are consecutive primes.

This completes the proof of the corollary.
\end{proof}

Hence, based on the main theorem and these corollaries, we can obtain the many interesting results which Schinzel $\&$ Sierpi$\acute{n}$ski proved, under the assumption that the Dickson's conjecture is valid[1-2]. Some of these results are as follows:

(1) For every even integer $2k\geq2$ there exist infinitely many pairs of consecutive primes with difference equal to $2k$. In particular, there exist infinitely many pairs of twin primes $p$ and $p+2$, cousin primes $p$ and $p+4$, and sexy primes $p$ and $p+6$;

(2) For every integer $m\geq1$ there exist $2m$ consecutive primes which are $m$ couples of twin primes;

(3) There are infinitely many prime triplets of the form $(p,p+2,p+6)$ and $(p,p+4,p+6)$;

(4) There are infinitely many prime quadruplets of the form $(p,p+2,p+6,p+8)$;

(5) There are infinitely many prime $k$-tuples in every admissible pattern for a prime constellation $(p,p+2c_{1},p+2c_{2},\ldots,p+2c_{k-1})$ for $k\geq2$ and $0<c_{1}<c_{2}<\cdots<c_{k-1}$;

(6) There are infinitely many Pythagorean triples of the form $(a,b,c)$ with $b$ even and $a$, $c$
primes, for example $(3,4,5)$, $(11,12,13)$, and many more.

\begin{corollary}[Sophie Germain primes]\label{co:chap_5_5_2}
There are infinitely many Sophie Germain primes of the form $2p+1$ where $p$ is prime.
\end{corollary}

\begin{proof}
Let $f_{1}(x)=2x+1$, $f_{2}(x)=2f_{1}(x)+1=4x+3$ and
\[
f(x)=\prod_{i=1}^{2}f_{i}(x)=(2x+1)(4x+3)
\]
\[
=8x^{2}+10x+3=a_{2}x^{2}+a_{1}x+a_{0}.
\]

Since $f_{1}(x)$ and $f_{2}(x)$ are irreducible over the integers, and since $f(1)=21$ and $f(2)=55$, $f_{1}(x)$ and $f_{2}(x)$ satisfy the first three conditions of Theorem~\ref{th:chap_5_4_1}. Since we can always choose a suitable positive number $k$ and let $x=ky$ to convert the polynomial $f(x)$ to a new polynomial
\[
g(y)=f(ky)=a_{2}k^{2}y^{2}+a_{1}ky+a_{0}=b_{2}y^{2}+b_{1}y+b_{0}
\]
such that the integer coefficients of the new polynomial $g(y)$ satisfy all conditions of Theorem~\ref{th:chap_5_4_1}, there exist infinitely many natural numbers $m\geq1$ such that $p=f_{1}(m)=2m+1$ and $2p+1=f_{2}(m)=4m+3$ are Sophie Germain primes.

This completes the proof of the corollary.
\end{proof}

Hence, based on the main theorem, there are infinitely many Sophie Germain primes of the form $2p+1$ where $p$ is prime.

%% file: chap6.tex
%-----------------------------------------------------------------------
% Beginning of chap7.tex
%-----------------------------------------------------------------------
%
% AMS-LaTeX 1.2 sample file for a monograph, based on amsbook.cls.
% This is a data file input by chapter.tex.
%%%%%%%%%%%%%%%%%%%%%%%%%%%%%%%%%%%%%%%%%%%%%%%%%%%%%%%%%%%%%%%%%%%%

\chapter{Prime product number theorem of the Schinzel's polynomials}\label{ch:chap_6}

\section*{Introduction}

In his joint paper with W. Sierpinski, Schinzel proposed the following conjectures~\cite{AD,FB}:

(1) Let $k\geq1$, let $f_{1}(x),f_{2}(x),\ldots,f_{k}(x)$ be irreducible polynomials over the integers, each with integer coefficients and a positive leading coefficient. Assume that the following condition holds:
(*) there does not exist any integer $n>1$ dividing all the products $f_{1}(j)f_{2}(j)\cdots f_{k}(j)$, for every integer $j$.

Then there exist infinitely many natural numbers $m$ such that all numbers $f_{1}(m),f_{2}(m),\ldots,f_{k}(m)$ are primes.

(2) Under the assumptions for $f_{1}(x),f_{2}(x),\ldots,f_{k}(x)$, there exists a natural number $m$ such that the numbers $f_{1}(m),f_{2}(m)$,$\cdots,f_{k}(m)$ are primes.

From Theorem~\ref{th:chap_5_4_1}, we know that there exist infinitely many natural numbers $m$ such that all numbers $f_{1}(m),f_{2}(m),\ldots,f_{k}(m)$ are primes.

Let $\pi_{k}(x)$ denote the number of natural numbers $m\leq x$ where all numbers $f_{1}(m),f_{2}(m),\ldots,f_{k}(m)$ are primes. The Bateman-Horn conjecture states that
\[
\pi_{k}(x)\sim\frac{C}{D}\int_{2}^{x}\frac{du}{\log^{k}u},\texttt{ }C=\prod_{p}\frac{1-N(p)/p}{(1-1/p)^{k}}
\]
where $D$ is the product of the degrees of the polynomials and where $C$ is the product over primes $p$ with $N(p)$ the number of solutions to $f(m)\equiv0(mod$ $p)$.

It is proved in this paper that as $x \to \infty$,
\[
\pi_{k}(x)=\frac{\alpha_{k}}{D_{k}}\int_{2}^{x}\frac{du}{\log^{k}u}+O(xe^{-\frac{1}{2}\alpha\sqrt{\log x}})
\]
where $\alpha>0$ and $\alpha_{k}>0$ are constants, and $D_{k}$ is the product of the degrees of the polynomials.

Hence, based on the prime number theorem for a group of polynomials of arbitrary degree, the prime number theorems for twin primes and Sophie Germain primes can be proved.

\section{Definitions}\label{sec:chap_6_1_26}

For any finite positive number $k\geq1$, let $f_{1}(x),f_{2}(x),\ldots,f_{k}(x)$ be irreducible polynomials over the integers, each with a positive leading coefficient. Assume that the following condition holds: there does not exist any integer $d>1$ dividing all the products $f_{1}(j)f_{2}(j)\cdots f_{k}(j)$, for every integer $j$.

\begin{definition}\label{def:chap_6_1_1}
Given a polynomial $f(x)=f_{1}(x)f_{2}(x)\cdots f_{k}(x)$, let define number sets $Q_{i}(x)$ on polynomials $f_{1}(x)$, $f_{2}(x)$, $\ldots$, and $f_{k}(x)$, respectively, by
\[
Q_{1}(x)=\{q|f_{1}(q)\texttt{ is prime},0<q\leq x\},
\]
\[
Q_{i}(x)=\{q|f_{i}(q)\texttt{ is prime},q\in Q_{i-1}(x)\}
\]
for $i=2,3,\ldots,k$, and define a number set $Q_{f}(x)$ on the polynomial $f(x)$ by
\[
Q_{f}(x)=\{q|f_{i}(q)\texttt{ is prime for }i=1,2,\ldots,k,0<q\leq x\}.
\]
Then there are
\[
Q_{1}(x)\supset Q_{2}(x)\supset\cdots\supset Q_{k-1}(x)\supset Q_{k}(x)=Q_{f}(x).
\]

Let $\pi_{i}(x)$ count the number of primes $f_{i}(q),q\in Q_{i}(x)$ for $i=1,2,\ldots,k$, respectively, and $\pi_{f}(x)$ count the number of prime products $f(q)=f_{1}(q)f_{2}(q)\cdots f_{k}(q),q\in Q_{f}(x)$. Then there is $\pi_{f}(x)=\pi_{k}(x)$.
\end{definition}

\begin{definition}\label{def:chap_6_1_2}
With complex variable $s$ and primes $p$, let define complex functions $\zeta_{i}(s)$ on number sets $Q_{i}(\infty)$ for $i=1,2,\ldots,k$, respectively, by
\[
\zeta_{i}(s)=\prod_{p=f_{i}(q),q\in Q_{i}(\infty)}(1-p^{-s})^{-1}.
\]
Then the complex function $\zeta_{f}(s)$ on the number set $Q_{f}(\infty)$ is define by
\[
\zeta_{f}(s)=\zeta_{k}(s)=\prod_{p=f_{k}(q),q\in Q_{k}(\infty)}(1-p^{-s})^{-1}.
\]
\end{definition}

For $i=1,2,\ldots,k$, because of the Riemann zeta function
\[
\zeta(s)=\zeta_{i}(s)\prod_{p\neq f_{i}(q),q\in Q_{i}(\infty)}(1-p^{-s})^{-1}
\]
and the relation
\[
\log\zeta(s)=\log\zeta_{i}(s)+\sum_{m,p\neq f_{i}(q),q\in Q_{i}(\infty)}\frac{1}{mp^{ms}},
\]
we have
\[
0<|\zeta_{i}(s)/\zeta(s)|<1
\]
and
\[
0<|\log\zeta_{i}(r,s)/\log\zeta(s)|<1.
\]
Also, we have
\[
0<|\zeta_{f}(s)/\zeta(s)|<1
\]
and
\[
0<|\log\zeta_{f}(r,s)/\log\zeta(s)|<1.
\]

\begin{definition}\label{def:chap_6_1_3}
Let define functions on the number set $Q_{i}(\infty)$ by
\[
\Lambda_{i}(\xi)=\{
       \begin{array}{cc}
         \log p & \texttt{ if }\xi=p^{m},\texttt{ for any }m>0\texttt{ and }p=f_{i}(q),q\in Q_{i}(\infty), \\
         0 & otherwise.
       \end{array}
\]
for $i=1,2,\ldots,k$, and a function on the number set $Q_{f}(\infty)$ by $\Lambda_{f}(\xi)=\Lambda_{k}(\xi)$, i.e.,
\[
\Lambda_{f}(\xi)=\{
       \begin{array}{cc}
         \log p & \texttt{ if }\xi=p^{m},\texttt{ for any }m>0\texttt{ and }p=f_{k}(q),q\in Q_{k}(\infty), \\
         0 & otherwise.
       \end{array}
\]
\end{definition}

For $i=1,2,\ldots,k$, differentiating the relation
\[
\log\zeta_{i}(s)=\sum_{m,p=f_{i}(q),q\in Q_{i}(\infty)}\frac{1}{mp^{ms}},
\]
we obtain
\[
\frac{\zeta_{i}'(s)}{\zeta_{i}(s)}
=-\sum_{m,p=f_{i}(q),q\in Q_{i}(\infty)}\frac{\log p}{p^{ms}}
=-\sum_{\xi=1}^{\infty}\frac{\Lambda_{i}(\xi)}{\xi^{s}}.
\]

Also there are
\[
\log\zeta_{f}(s)=\log\zeta_{k}(s)=\sum_{m,p=f_{k}(q),q\in Q_{k}(\infty)}\frac{1}{mp^{ms}},
\]
and
\[
\frac{\zeta_{f}'(s)}{\zeta_{f}(s)}
=\frac{\zeta_{k}'(s)}{\zeta_{k}(s)}
=-\sum_{m,p=f_{k}(q),q\in Q_{k}(\infty)}\frac{\log p}{p^{ms}}
=-\sum_{\xi=1}^{\infty}\frac{\Lambda_{k}(\xi)}{\xi^{s}}.
\]

\begin{definition}\label{def:chap_6_1_4}
For numbers $q\in Q_{i}(x)$ where $i=1,2,\ldots,k$, similar to the proof of the prime number theorem for arithmetic progressions in~\cite{FD}, our object is to get an estimate for
\begin{equation}
\pi_{i}(x)=\sum_{p=f_{i}(q),q\in Q_{i}(x)}1=\sum_{\xi=1}^{f_{i}(x)}P_{i}(\xi)
\end{equation}
where $P_{i}$ is the characteristic function of the numbers:
\[
P_{i}(\xi)=\{
       \begin{array}{cc}
         1 & \texttt{ if }\xi=f_{i}(q),q\in Q_{i}(\infty), \\
         0 & otherwise.
       \end{array}
\]
\end{definition}

While $P_{i}$ itself does not arise in a natural way, the function $P_{i}^{*}$ such that
\[
P_{i}^{*}(\xi)=\{
       \begin{array}{cc}
         \frac{1}{m} & \texttt{ if $\xi=p^{m}$ for some $m$ and $p=f_{i}(q),q\in Q_{i}(\infty)$}, \\
         0 & otherwise,
       \end{array}
\]
occurs in the Dirichlet series for $\log\zeta_{i}(s)$:
\begin{equation}\label{eq_pnt:schinzel_zeta}
\log\zeta_{i}(s)=\sum_{m,p=f_{i}(q),q\in Q_{i}(\infty)}\frac{1}{mp^{ms}}
=\sum_{\xi=1}^{\infty}\frac{P_{i}^{*}(\xi)}{\xi^{s}}.
\end{equation}
For each fixed $m$, the number of $m$th powers of numbers $p=f_{i}(q),q\in Q_{i}(x)$, which do not exceed $f_{i}(x)$ is equal to the number of numbers $q\in Q_{i}(x)$ which do not exceed $\sqrt[m]{x}$, so that
\[
\begin{split}
\sum_{\xi=1}^{f_{i}(x)}P_{i}^{*}(\xi)
=&\sum_{\xi=1}^{f_{i}(x)}P_{i}(\xi)
+\frac{1}{2}\sum_{\xi=1}^{\sqrt{f_{i}(x)}}P_{i}(\xi)
+\frac{1}{3}\sum_{\xi=1}^{\sqrt[3]{f_{i}(x)}}P_{i}(\xi)+\cdots \\
\leq&\sum_{\xi=1}^{f_{i}(x)}P_{i}(\xi)
+\frac{1}{2}\sum_{\xi=1}^{f_{i}(\sqrt{x})}P_{i}(\xi)
+\frac{1}{3}\sum_{\xi=1}^{f_{i}(\sqrt[3]{x})}P_{f}(\xi)+\cdots \\
=&\pi_{i}(x)+\frac{\pi_{i}(\sqrt{x})}{2}+\frac{\pi_{i}(\sqrt[3]{x})}{3}+\cdots,
\end{split}
\]
and since, for $m\geq2$,
\[
\pi_{i}(x^{1/m})<cx^{1/m}\leq c\sqrt{x}=o(\frac{x}{\log x}),
\]
it is to be expected that
\[
\sum_{\xi=1}^{f_{i}(x)}P_{i}^{*}(\xi)\sim\pi_{i}(x)\texttt{ for }i=1,2,\ldots,k.
\]

\section{Lemmas of Schinzel's polynomials}\label{sec:chap_6_2_27}

\begin{lemma}\label{le:chap_6_2_1}
By Definition (\ref{def:chap_6_1_2}), let $\log\zeta_{1}(s)=(x/f_{1}(x))^{s}\theta_{1}(x,s)\log\zeta(s)$. Then for $0<s=\sigma\leq1$
\[
0<|\theta_{1}(x,s)|\leq C_{1}^{*}n_{1}\texttt{ and }\frac{|\theta_{1}(x,s)|}{n_{1}\log x}\leq\frac{\pi_{1}(x)}{x}\leq \frac{C_{1}^{*}}{\log x}
\]
where $n_{1}$ is the positive degree of the polynomial $f_{1}(x)$ and $C_{1}^{*}$ is a positive constant (dependent on $f_{1}(x)$).
\end{lemma}

\begin{proof}
First, by Definition (\ref{def:chap_6_1_2}), there is $0<|\log\zeta_{1}(r,s)/\log\zeta(s)|<1$ such that there is
\[
|\theta_{1}(x,s)|=|\frac{\log\zeta_{1}(s)}{\log\zeta(s)}(\frac{f_{1}(x)}{x})^{s}|>0.
\]

Next, let define a function $\zeta_{1}(s,x)$ on the number set $Q_{1}(x)$ by
\[
\log\zeta_{1}(s,x)=\sum_{p=f_{1}(q),q\in Q_{1}(x)}(1-p^{-s})^{-1}=\frac{\pi_{1}(x)}{1-\rho_{1}^{-s}(s)}
\]
and define a function $\zeta(s,f_{1}(x))$ on the function $\zeta(s)$ by
\[
\log\zeta(s,f_{1}(x))=\sum_{p\leq f_{1}(x)}(1-p^{-s})^{-1}=\frac{\pi(f_{1}(x))}{1-\rho_{0}^{-s}(s)}.
\]
Then there are
\[
\frac{1}{1-\rho_{1}^{-s}(s)}=\frac{\log\zeta_{1}(s,x)}{\pi_{1}(x)}=\frac{1}{\pi_{1}(x)}\sum_{p=f_{1}(q),q\in Q_{1}(x)}(1-p^{-s})^{-1}
\]
and
\[
\frac{1}{1-\rho_{0}^{-s}(s)}=\frac{\log\zeta(s,f_{1}(x))}{\pi(f_{1}(x))}=\frac{1}{\pi(f_{1}(x))}\sum_{p\leq f_{1}(x)}(1-p^{-s})^{-1}.
\]
It shows that $(1-\rho_{1}^{-s}(s))^{-1}$ and $(1-\rho_{0}^{-s}(s))^{-1}$ are the average values of $(1-p^{-s})^{-1}$ for $\log\zeta_{1}(s,x)$ and $\log\zeta(s,f_{1}(x))$, respectively.

Let $q_{0}=0$ and the elements of the number set $Q_{1}(x)$ be in a sequence
\[
0<q_{1}<q_{2}<\cdots<q_{i-1}<q_{i}<\cdots<q_{i_{x}}\leq x
\]
where $i_{x}=\pi_{1}(x)$, and for $i=1,2,\ldots,i_{x}$ let $m_{i}$ denote the number of primes $p$ satisfying $f_{1}(q_{i-1})<p\leq f_{1}(q_{i})$. Then for $0<s=\sigma\leq1$ the inequality
\[
(1-f_{1}(q_{i})^{-\sigma})^{-1}\leq\frac{1}{m_{i}}\sum_{f_{1}(q_{i-1})<p\leq f_{1}(q_{i})}(1-p^{-\sigma})^{-1}
\]
always holds for $i=1,2,\ldots,i_{x}$. Thus the inequalities
\[
|\rho_{0}^{s}(s)|\leq|\rho_{1}^{s}(s)|\texttt{ and }|1-\rho_{0}^{-s}(s)|\leq|1-\rho_{1}^{-s}(s)|
\]
always hold. On the other hand, as $x \to \infty$ there are
\[
|\rho_{0}^{s}(s)| \to \infty,\texttt{ }
|\rho_{1}^{s}(s)| \to \infty,\texttt{ and }\lim_{x \to \infty}|\frac{1-\rho_{0}^{-s}(s)}{1-\rho_{1}^{-s}(s)}|=1.
\]
Thus there should be
\[
|\frac{\log\zeta_{1}(s,x)}{\log\zeta(s,f_{1}(x))}|=|\frac{1-\rho_{0}^{-s}(s)}{1-\rho_{1}^{-s}(s)}|\frac{\pi_{1}(x)}{\pi(f_{1}(x))}\leq\frac{\pi_{1}(x)}{\pi(f_{1}(x))}.
\]

In 1931, Heilbronn proved the statement~\cite{AD}:
There exists a positive constant $C_{1}^{*}$ (dependent on $f_{1}(x)$), such that
\[
\pi_{1}(x)\leq C_{1}^{*}\frac{x}{\log x},\texttt{ for every }x\geq1.
\]

Since as $x \to \infty$ and $0<s=\sigma\leq1$
\[
\begin{split}
|\theta_{1}(x,s)|=&|\frac{\log\zeta_{1}(s)}{\log\zeta(s)}(\frac{f_{1}(x)}{x})^{s}|=\lim_{x \to \infty}|\frac{\log\zeta_{1}(s,x)}{\log\zeta(s,f_{1}(x))}(\frac{f_{1}(x)}{x})^{s}| \\
\leq&\frac{\pi_{1}(x)}{\pi(f_{1}(x))}\cdot\frac{f_{1}(x)}{x}\leq\frac{(C_{1}^{*}x/\log x)/x}{(f_{1}(x)/\log f_{1}(x))/f_{1}(x)}=C_{1}^{*}n_{1}
\end{split}
\]
where $n_{1}$ is the positive degree of the polynomial $f_{1}(x)$, the inequalities
\[
0<|\theta_{1}(x,s)|\leq C_{1}^{*}n_{1}\texttt{ and }\frac{|\theta_{1}(x,s)|}{n_{1}\log x}\leq\frac{\pi_{1}(x)}{x}\leq \frac{C_{1}^{*}}{\log x}
\]
hold for $0<s=\sigma\leq1$.

Hence, for $0<s=\sigma\leq1$, $\theta_{1}(x,s)$ is the ratio of the average density $\log\zeta_{1}(s,x)/x^{s}$ of the distribution of primes $p=f_{1}(q),q\in Q_{1}(x)$ in the region $[1,f_{1}(x)]$ to the average density $\log\zeta(s,f_{1}(x))/f_{1}^{s}(x)$ of the distribution of primes $p$ in the range $[1,f_{1}(x)]$.

This completes the proof of the lemma.
\end{proof}

\begin{lemma}\label{le:chap_6_2_2}
By Definition (\ref{def:chap_6_1_2}), for each $i=1,2,\ldots,k-1$, let
\[
\log\zeta_{i+1}(s)=\theta_{i+1}(x,s)\log\zeta_{i}(s)=(x/f_{1}(x))^{s}\prod_{j=1}^{i+1}\theta_{j}(x,s)\log\zeta(s).
\]
Then for $0<s=\sigma\leq1$, there are
\[
0<|\theta_{i+1}(x,s)|=|\frac{\log\zeta_{i+1}(s)}{\log\zeta_{i}(s)}|
\sim\frac{\pi_{i+1}(x)}{\pi_{i}(x)}<1,
\]
and
\[
0<|\theta_{f}(x,s)|<C_{1}^{*}n_{1}\texttt{ where }\theta_{f}(x,s)=\prod_{j=1}^{k}\theta_{j}(x,s).
\]
\end{lemma}

\begin{proof}
First, by Definition (\ref{def:chap_6_1_2}), there is $0<|\log\zeta_{i+1}(r,s)/\log\zeta_{i}(s)|<1$ such that
\[
|\theta_{i+1}(x,s)|=|\frac{\log\zeta_{i+1}(s)}{\log\zeta_{i}(s)}|>0.
\]

Next, let define a function $\zeta_{i}(s,x)$ on the number set $Q_{i}(x)$ by
\[
\log\zeta_{i}(s,x)=\sum_{p=f_{i}(q),q\in Q_{i}(x)}(1-p^{-s})^{-1}=\frac{\pi_{i}(x)}{1-\rho_{i}^{-s}(s)}
\]
and define a function $\zeta_{i+1}(s,x)$ on the number set $Q_{i+1}(x)$ by
\[
\log\zeta_{i+1}(s,x)=\sum_{p=f_{i+1}(q),q\in Q_{i+1}(x)}(1-p^{-s})^{-1}=\frac{\pi_{i+1}(x)}{1-\rho_{i+1}^{-s}(s)}.
\]
Then there are
\[
\frac{1}{1-\rho_{i}^{-s}(s)}=\frac{\log\zeta_{i}(s,x)}{\pi_{i}(x)}=\frac{1}{\pi_{i}(x)}\sum_{p=f_{i}(q),q\in Q_{i}(x)}(1-p^{-s})^{-1}
\]
and
\[
\frac{1}{1-\rho_{i+1}^{-s}(s)}=\frac{\log\zeta_{i+1}(s,x)}{\pi_{i+1}(x)}=\frac{1}{\pi_{i+1}(x)}\sum_{p=f_{i+1}(q),q\in Q_{i+1}(x)}(1-p^{-s})^{-1}.
\]
It shows that $(1-\rho_{i}^{-s}(s))^{-1}$ and $(1-\rho_{i+1}^{-s}(s))^{-1}$ are the average values of $(1-p^{-s})^{-1}$ for $\log\zeta_{i}(s,x)$ and $\log\zeta_{i+1}(s,x)$, respectively.

Let $q_{0}=0$ and the elements of the number set $Q_{i+1}(x)$ be in a sequence
\[
0<q_{1}<q_{2}<\cdots<q_{j-1}<q_{j}<\cdots<q_{j_{x}}\leq x
\]
where $j_{x}=\pi_{i+1}(x)$, and for $j=1,2,\ldots,j_{x}$ let $m_{j}$ denote the number of primes $p=f_{i}(q)$ where $q\in Q_{i}(x)$ and $p$ satisfies $f_{i+1}(q_{j-1})<p\leq f_{i+1}(q_{j})$. Then for $0<s=\sigma\leq1$ the inequality
\[
(1-f_{i+1}(q_{j})^{-\sigma})^{-1}\leq\frac{1}{m_{j}}\sum_{f_{i+1}(q_{j-1})<p\leq f_{i+1}(q_{j})}(1-p^{-\sigma})^{-1}
\]
always holds for $j=1,2,\ldots,j_{x}$. Thus the inequalities
\[
|\rho_{i}^{s}(s)|\leq|\rho_{i+1}^{s}(s)|\texttt{ and }|1-\rho_{i}^{-s}(s)|\leq|1-\rho_{i+1}^{-s}(s)|
\]
always hold. On the other hand, as $x \to \infty$ there are
\[
|\rho_{i}^{s}(s)| \to \infty,\texttt{ }
|\rho_{i+1}^{s}(s)| \to \infty,\texttt{ and }\lim_{x \to \infty}|\frac{1-\rho_{i}^{-s}(s)}{1-\rho_{i+1}^{-s}(s)}|=1.
\]
Then as $x \to \infty$ there should be
\[
|\frac{\log\zeta_{i+1}(s,x)}{\log\zeta_{i}(s,x)}|=|\frac{1-\rho_{i}^{-s}(s)}{1-\rho_{i+1}^{-s}(s)}|\frac{\pi_{i+1}(x)}{\pi_{i}(x)}\sim\frac{\pi_{i+1}(x)}{\pi_{i}(x)}
\]
and
\[
|\theta_{i+1}(x,s)|=|\frac{\log\zeta_{i+1}(s)}{\log\zeta_{i}(s)}|
\sim\frac{\pi_{i+1}(x)}{\pi_{i}(x)}<1.
\]

Finally, from Lemma~\ref{le:chap_6_2_1}, there are
\[
\begin{split}
\theta_{f}(x,s)=&(\frac{f_{1}(x)}{x})^{s}\frac{\log\zeta_{k}(s)}{\log\zeta(s)} \\
=&(\frac{f_{1}(x)}{x})^{s}\frac{\log\zeta_{1}(s)}{\log\zeta(s)}\prod_{j=1}^{k-1}\frac{\log\zeta_{j+1}(s)}{\log\zeta_{j}(s)}
=\prod_{j=1}^{k}\theta_{j}(x,s)
\end{split}
\]
and
\[
0<|\theta_{f}(x,s)|=\prod_{j=1}^{k}|\theta_{j}(x,s)|<|\theta_{1}(x,s)|\leq C_{1}^{*}n_{1}.
\]

This completes the proof of the lemma.
\end{proof}

\begin{lemma}\label{le:chap_6_2_3_th_7_15}
There is a constant $\alpha>0$ such that as $x \to \infty$,
\begin{equation}\label{eq_pnt:schinzel_7_15}
\sum_{p=f_{1}(q),q\in Q_{1}(x)}\log\frac{f_{1}(x)}{p}
=\int_{c}^{1}\frac{x^{s}}{s^{2}}\theta_{1}(x,s)ds+O(xe^{-\alpha\sqrt{\log x}})
\end{equation}
for some $c$ with $0<c<1$.
\end{lemma}

\begin{proof}
Using Theorem~7-9 in~\cite{FD}:
\[
\frac{1}{2\pi i}\int_{(2)}\frac{y^{s}}{s^{2}}ds=\{
\begin{array}{cc}
  0 & \texttt{ for }0<y<1, \\
  \log y & \texttt{ for }y\geq1,
\end{array}
\]
similar to the proof of Theorem~7-15 in~\cite{FD}, we have
\[
\begin{split}
\frac{1}{2\pi i}\int_{(2)}\frac{f_{1}^{s}(x)}{s^{2}}\log\zeta_{1}(s)ds
=&\frac{1}{2\pi i}\int_{(2)}\frac{1}{s^{2}}\sum_{\xi=1}^{\infty}\frac{\Lambda_{1}(\xi)}{\log\xi}(\frac{f_{1}(x)}{\xi})^{s}ds \\
=&\frac{1}{2\pi i}\sum_{\xi=1}^{\infty}\frac{\Lambda_{1}(\xi)}{\log\xi}\int_{(2)}\frac{1}{s^{2}}(\frac{f_{1}(x)}{\xi})^{s}ds \\
=&\sum_{\xi\leq f_{1}(x)}\frac{\Lambda_{1}(\xi)}{\log\xi}\log\frac{f_{1}(x)}{\xi} \\
=&\sum_{m,p=f_{1}(q),q\in Q_{1}(x),p^{m}\leq f_{1}(x)}\frac{1}{m}\log\frac{f_{1}(x)}{p^{m}} \\
=&\sum_{p=f_{1}(q),q\in Q_{1}(x)}\log\frac{f_{1}(x)}{p} \\
 &+\sum_{p=f_{1}(q),q\in Q_{1}(x),m\geq2,p^{m}\leq f_{1}(x)}\frac{1}{m}\log\frac{f_{1}(x)}{p^{m}}.
\end{split}
\]

As noted earlier, the number of terms in the last sum is
\[
\pi_{1}(x^{1/2})+\pi_{1}(x^{1/3})+\cdots<x^{1/2}+x^{1/3}+\cdots+x^{1/u}<ux^{1/2},
\]
where $u$ is the smallest number such that $x^{1/u}<2$. Thus
\[
\pi_{1}(x^{1/2})+\pi_{1}(x^{1/3})+\cdots=O(x^{1/2}\log x),
\]
so that
\[
\begin{split}
\sum_{p=f_{1}(q),q\in Q_{1}(x),m\geq2,p^{m}\leq f_{1}(x)}\frac{1}{m}\log\frac{f_{1}(x)}{p^{m}}
\leq&\sum_{p=f_{1}(q),q\in Q_{1}(x),m\geq2,p^{m}\leq f_{1}(x)}\log f_{1}(x) \\
=&O(\sqrt{x}\log x\log f_{1}(x))=O(\sqrt{x}\log^{2}x) \\
=&O(xe^{-\sqrt{\log x}}).
\end{split}
\]
and
\[
\sum_{p=f_{1}(q),q\in Q_{1}(x)}\log\frac{f_{1}(x)}{p}
=\frac{1}{2\pi i}\int_{(2)}\frac{f_{1}^{s}(x)}{s^{2}}\log\zeta_{1}(s)ds+O(xe^{-\sqrt{\log x}}).
\]

Because from Lemma~\ref{le:chap_6_2_2} there is
\[
\log\zeta_{1}(s)=(x/f_{1}(x))^{s}\theta_{1}(x,s)\log\zeta(s)\texttt{ where }0<|\theta_{1}(x,s)|\leq C_{1}^{*}n_{1},
\]
the equation above becomes
\[
\sum_{p=f_{1}(q),q\in Q_{1}(x)}\log\frac{f_{1}(x)}{p}
=\frac{1}{2\pi i}\int_{(2)}\frac{x^{s}}{s^{2}}\theta_{1}(x,s)\log\zeta(s)ds+O(xe^{-\sqrt{\log x}})
\]
where $0<|\theta_{1}(x,s)|\leq C_{1}^{*}n_{1}$ and $n_{1}$ is the positive degree of the polynomial $f_{1}(s)$.

Similar to the analysis of the integral $\int_{(2)}(x^{s}/s^{2})\log\zeta(s)ds$ in~\cite{FD}, which results
\[
\int_{(2)}\frac{x^{s}}{s^{2}}\log\zeta(s)ds
=2\pi i\int_{c}^{1}\frac{x^{s}}{s^{2}}ds+O(xe^{-\alpha\sqrt{\log x}}),
\]
the analysis of the integral $\int_{(2)}(x^{s}/s^{2})\theta_{1}(x,s)\log\zeta(s)ds$ with $0<|\theta_{1}(x,s)|<C_{1}^{*}n_{1}$ results
\[
\int_{(2)}\frac{x^{s}}{s^{2}}\theta_{1}(x,s)\log\zeta(s)ds
=2\pi i\int_{c}^{1}\frac{x^{s}}{s^{2}}\theta_{1}(x,s)ds+O(xe^{-\alpha\sqrt{\log x}}),
\]
so that there are
\[
\begin{split}
\sum_{p=f_{1}(q),q\in Q_{1}(x)}\log\frac{f_{1}(x)}{p}
%=&\frac{1}{2\pi i}\int_{(2)}\frac{f_{1}^{s}(x)}{s^{2}}\log\zeta_{1}(s)ds+O(xe^{-\sqrt{\log x}}) \\
%=&\frac{1}{2\pi i}\int_{(2)}\frac{x^{s}}{s^{2}}\theta_{1}(x,s)\log\zeta(s)ds+O(xe^{-\sqrt{\log x}}) \\
=&\int_{c}^{1}\frac{x^{s}}{s^{2}}\theta_{1}(x,s)ds+O(xe^{-\alpha\sqrt{\log x}}).
\end{split}
\]

This completes the proof of the lemma.
\end{proof}

\begin{lemma}\label{le:chap_6_2_4_th_7_15}
There exist constants $\alpha>0$ and $\alpha_{1}>0$ such that as $x \to \infty$,
\begin{equation}\label{eq_pnt:schinzel_7_15_1}
\sum_{p=f_{1}(q),q\in Q_{1}(x)}\log\frac{f_{1}(x)}{p}
=\alpha_{1}\int_{c}^{1}\frac{x^{s}}{s^{2}}ds+O(xe^{-\alpha\sqrt{\log x}}).
\end{equation}
for some $c$ with $0<c<1$.
\end{lemma}

\begin{proof}
From Lemma~\ref{le:chap_6_2_3_th_7_15}, there is a constant $\alpha>0$ such that as $x \to \infty$,
\[
\sum_{p=f_{1}(q),q\in Q_{1}(x)}\log\frac{f_{1}(x)}{p}
=\int_{c}^{1}\frac{x^{s}}{s^{2}}\theta_{1}(x,s)ds+O(xe^{-\alpha\sqrt{\log x}})
\]
for some $c$ with $0<c<1$.

From the proof of Theorem~7-15 in~\cite{FD}, we have
\[
\sum_{p\leq x}\log\frac{x}{p}
=\int_{c}^{1}\frac{x^{s}}{s^{2}}ds+O(xe^{-\alpha\sqrt{\log x}}).
\]

By denoting
\[
\frac{\alpha_{1}}{n}=\lim_{x \to \infty}(\sum_{p=f_{1}(q),0<q\leq x}\log\frac{x}{q})/(\sum_{p\leq x}\log\frac{x}{p})
\]
where $\alpha_{1}>0$, there are
\[
\begin{split}
\lim_{x \to \infty}\sum_{p=f_{1}(q),q\in Q_{1}(x)}\log\frac{f_{1}(x)}{p}=&\lim_{x \to \infty}\sum_{q\in Q_{1}(x)}n\log\frac{x}{q} \\
=&n\lim_{x \to \infty}\sum_{p=f_{1}(q),0<q\leq x}\log\frac{x}{q} \\
=&\alpha_{1}\lim_{x \to \infty}\sum_{p\leq x}\log\frac{x}{p}.
\end{split}
\]

Thus there should be
\[
\int_{c}^{1}\frac{x^{s}}{s^{2}}\theta_{1}(x,s)ds\sim\alpha_{1}\int_{c}^{1}\frac{x^{s}}{s^{2}}ds\texttt{ as }x \to \infty,
\]
so that we have
\[
\sum_{p=f_{1}(q),q\in Q_{1}(x)}\log\frac{f_{1}(x)}{p}
=\alpha_{1}\int_{c}^{1}\frac{x^{s}}{s^{2}}ds+O(xe^{-\alpha\sqrt{\log x}}).
\]

This completes the proof of the lemma.
\end{proof}

\begin{lemma}\label{le:chap_6_2_5_th_7_16}
For a polynomial $f_{1}(x)$ in one variable with positive degree $n_{1}$ and integer coefficients satisfying three necessary conditions:

(1) the leading coefficient of $f_{1}(x)$ is positive,

(2) the polynomial is irreducible over the integers,

(3) as $m$ runs over the positive integers the numbers $f_{1}(m)$ should not share a common factor greater than 1, \\
let $\pi_{1}(x)$ denote the number of primes $f_{1}(q)$ where $0<q\leq x$. There is a constant $\alpha>0$ such that as $x \to \infty$,
\begin{equation}\label{eq_pnt:schinzel_7_16}
\pi_{1}(x)=\frac{1}{n_{1}}\int_{2}^{x}\theta_{1}(u)\frac{du}{\log u}+O(xe^{-\frac{1}{2}\alpha\sqrt{\log x}}).
\end{equation}
\end{lemma}

\begin{proof}
In Eq. (\ref{eq_pnt:schinzel_7_15}), put
\[
\delta=\delta(x)=e^{-\frac{1}{2}\alpha\sqrt{\log x}}.
\]
Then since $f_{1}(x(1+\delta))\sim f_{1}(x)(1+n_{1}\delta)$ and $\log(1+n_{1}\delta)\sim n_{1}\delta$ as $x \to \infty$, we have
\[
\begin{split}
&\sum_{p=f_{1}(q),q\in Q_{1}(x(1+\delta))}\log\frac{f_{1}(x)(1+n_{1}\delta)}{p}
-\sum_{p=f_{1}(q),q\in Q_{1}(x)}\log\frac{f_{1}(x)}{p} \\
=&\sum_{p=f_{1}(q),q\in Q_{1}(x)}\log(1+n_{1}\delta)
+\sum_{p=f_{1}(q),x<q\leq x(1+\delta)}\log\frac{f_{1}(x)(1+n_{1}\delta)}{p} \\
=&\log(1+n_{1}\delta)\pi_{1}(x)+O(\log(1+n_{1}\delta)*\delta x) \\
=&\int_{c}^{1}\frac{x^{s}}{s^{2}}((1+\delta)^{s}-1)\theta_{1}(x,s)ds+O(xe^{-\alpha\sqrt{\log x}}),
\end{split}
\]
so that
\[
\pi_{1}(x)
=\frac{1}{\log(1+n_{1}\delta)}\int_{c}^{1}\frac{(1+\delta)^{s}-1}{s^{2}}x^{s}\theta_{1}(x,s)ds
+O(\delta x)+O(\frac{xe^{-\alpha\sqrt{\log x}}}{\delta}).
\]

Now
\[
(1+\delta)^{s}-1=s\delta+\frac{s(s-1)}{2!}(1+\vartheta\delta)^{s-2}\delta^{2},
\]
where $0<\vartheta<1$, so that for $0<s<1$,
\[
|(1+\delta)^{s}-1-s\delta|\leq\frac{s|s-1|}{2}\delta^{2}<\delta^{2}.
\]

Thus, making the change of variable $x^{s}=u$, we obtain
\[
\begin{split}
\int_{c}^{1}\frac{(1+\delta)^{s}-1}{s^{2}}x^{s}\theta_{1}(x,s)ds
=&\delta\int_{c}^{1}\frac{x^{s}}{s}\theta_{1}(x,s)ds+O(\delta^{2}\int_{c}^{1}\frac{x^{s}}{s^{2}}ds) \\
=&\delta\int_{x^{c}}^{x}\theta_{1}(u)\frac{du}{\log u}+O(\delta^{2}\int_{c}^{1}x^{s}ds) \\
=&\delta\int_{2}^{x}\theta_{1}(u)\frac{du}{\log u}+O(\delta^{2}x).
\end{split}
\]

Finally,
\[
\begin{split}
\pi_{1}(x)
=&\frac{\delta}{\log(1+n_{1}\delta)}\int_{2}^{x}\theta_{1}(u)\frac{du}{\log u}
+O(\delta x)+O(\frac{xe^{-\alpha\sqrt{\log x}}}{\delta}) \\
=&\frac{1+O(\delta)}{n_{1}}\int_{2}^{x}\theta_{1}(u)\frac{du}{\log u}+O(\delta x)+O(xe^{-\frac{1}{2}\alpha\sqrt{\log x}}) \\
=&\frac{1}{n_{1}}\int_{2}^{x}\theta_{1}(u)\frac{du}{\log u}+O(\delta x) \\
=&\frac{1}{n_{1}}\int_{2}^{x}\theta_{1}(u)\frac{du}{\log u}+O(xe^{-\frac{1}{2}\alpha\sqrt{\log x}}).
\end{split}
\]

This completes the proof of the lemma.
\end{proof}

\begin{lemma}\label{le:chap_6_2_6_th_7_16}
From Lemma~\ref{le:chap_6_2_5_th_7_16} with a constant $\alpha_{1}>0$
\[
\pi_{1}(x)=\frac{\alpha_{1}}{n_{1}}\int_{2}^{x}\frac{du}{\log u}+O(xe^{-\frac{1}{2}\alpha\sqrt{\log x}}).
\]
\end{lemma}

\begin{proof}
From Lemma~\ref{le:chap_6_2_5_th_7_16} with a constant $\alpha_{1}>0$, there is
\[
\pi_{1}(x)=\frac{1}{n_{1}}\int_{2}^{x}\theta_{1}(u)\frac{du}{\log u}+O(xe^{-\frac{1}{2}\alpha\sqrt{\log x}}).
\]

From Lemma~\ref{le:chap_6_2_4_th_7_15}, there should be
\[
\int_{2}^{x}\theta_{1}(u)\frac{du}{\log u}\sim\alpha_{1}\int_{2}^{x}\frac{du}{\log u}\texttt{ as }x \to \infty
\]
where $\alpha_{1}>0$. Thus we have
\[
\pi_{1}(x)=\frac{\alpha_{1}}{n_{1}}\int_{2}^{x}\frac{du}{\log u}+O(xe^{-\frac{1}{2}\alpha\sqrt{\log x}}).
\]

This completes the proof of the lemma.
\end{proof}

\begin{lemma}\label{le:chap_6_2_7}
By Definition (\ref{def:chap_6_1_1}), for each $i=1,2,\ldots,k-1$, there exists a positive constant $C_{i+1}$ (dependent on $f_{i+1}(x)$), such that there is
\[
\frac{\pi_{i+1}(x)}{\pi_{i}(x)}\sim\frac{C_{i+1}}{\log x}+O(\frac{\log\log x}{\log^{2}x}).
\]
\end{lemma}

\begin{proof}
From Lemma~\ref{le:chap_6_2_6_th_7_16}, there exists a positive constant $C_{1}$ (dependent on $f_{1}(x)$), such that
\[
\pi_{1}(x)\sim C_{1}\frac{x}{\log x},\texttt{ for every }x\geq1.
\]
By Definition (\ref{def:chap_6_1_1}) there are
\[
Q_{1}(x)\supset Q_{2}(x)\supset\cdots\supset Q_{k-1}(x)\supset Q_{k}(x)=Q_{f}(x).
\]
We consider that the function $f_{i+1}$ is operated on the number set $Q_{i}(x)$ and let $x_{0}=x$ and $x_{i}=\pi_{i}(x)$ for each $i=1,2,\ldots,k-1$. Then there are
\[
\pi_{2}(x)=\pi_{2}(x_{1})\sim C_{2}\frac{x_{1}}{\log x_{1}}=C_{2}\frac{\pi_{1}(x)}{\log \pi_{1}(x)},
\]
\[
\begin{split}
\frac{\pi_{2}(x)}{\pi_{1}(x)}\sim&\frac{C_{2}}{\log \pi_{1}(x)}\sim\frac{C_{2}}{\log(C_{1}x/\log x)} \\
                                =&\frac{C_{2}}{\log x}+O(\frac{\log\log x}{\log^{2}x}),
\end{split}
\]
and for each $i=2,3,\ldots,k-1$
\[
\begin{split}
\frac{\pi_{i+1}(x)}{\pi_{i}(x)}\sim&\frac{C_{i+1}}{\log\pi_{i}(x)} \\
                               \sim&\frac{C_{i+1}}{\log(C_{i}x_{i-1}/\log x_{i-1})} \\
                                  =&\frac{C_{i+1}}{\log x_{i-1}}+O(\frac{\log\log x_{i-1}}{\log^{2}x_{i-1}}) \\
                                  =&\frac{C_{i+1}}{\log\pi_{i-1}(x_{i-2})}+O(\frac{\log\log\pi_{i-1}(x_{i-2})}{\log^{2}\pi_{i-1}(x_{i-2})})\\
                               \sim&\frac{C_{i+1}}{\log(C_{i-1}x_{i-2}/\log x_{i-2})}+O(\frac{\log\log(x_{i-2}/\log x_{i-2})}{\log^{2}(x_{i-2}/\log x_{i-2})})\\
                                  =&\frac{C_{i+1}}{\log x_{i-2}}+O(\frac{\log\log x_{i-2}}{\log^{2}x_{i-2}})=\cdots=\frac{C_{i+1}}{\log x}+O(\frac{\log\log x}{\log^{2}x}).
\end{split}
\]

This completes the proof of the lemma.
\end{proof}

\section{Theorems of Schinzel's polynomials}\label{sec:chap_6_3_28}

\begin{theorem}[Prime number theorem for Schinzel's polynomials]\label{th:chap_6_3_1}
For any finite positive number $k\geq1$, let $f_{1}(x),f_{2}(x),\ldots,f_{k}(x)$ be irreducible polynomials over the integers, each with a positive leading coefficient.

Assume that the following condition holds: there does not exist any integer $d>1$ dividing all the products $f(j)=f_{1}(j)f_{2}(j)\cdots f_{k}(j)$, for every integer $j$.

Let $\pi_{k}(x)$ denote the number of natural numbers $m\leq x$ where all numbers $f_{1}(m),f_{2}(m),\ldots,f_{k}(m)$ are primes. As $x \to \infty$,
\[
\pi_{k}(x)=\frac{\alpha_{k}}{D_{k}}\int_{2}^{x}\frac{du}{\log^{k}u}+O(xe^{-\frac{1}{2}\alpha\sqrt{\log x}})
\]
where $\alpha>0$ and $\alpha_{k}>0$ are constants, and $D_{k}$ is the product of the degrees of the polynomials.
\end{theorem}

\begin{proof}
First, for $k=1$, from Lemma~\ref{le:chap_6_2_6_th_7_16}, we have
\[
\pi_{1}(x)=\frac{\alpha_{1}}{D_{1}}\int_{2}^{x}\frac{du}{\log u}+O(xe^{-\frac{1}{2}\alpha\sqrt{\log x}})\texttt{ where }D_{1}=n_{1}.
\]

Next, suppose that for $k\geq2$ and $0<i<k$, we have
\[
\pi_{i}(x)=\frac{\alpha_{i}}{D_{i}}\int_{2}^{x}\frac{du}{\log^{i}u}+O(xe^{-\frac{1}{2}\alpha\sqrt{\log x}})\texttt{ where }D_{i}=\prod_{j=1}^{i}n_{j}
\]
and as in the proof of Lemma~\ref{le:chap_6_2_5_th_7_16} by making the change of variables $x^{s}=u$ and
\[
\prod_{j=1}^{i}\theta_{j}(u)=\frac{\alpha_{i}}{D_{i}\log^{i-1}u}\texttt{ for }
\prod_{j=1}^{i}\theta_{j}(u)=\prod_{j=1}^{i}\theta_{j}(x,s)
\]
where
\[
\theta_{1}(u)=C_{1}\texttt{ and }\theta_{j}(u)=\frac{C_{j}}{\log u}\texttt{ for }j=2,3,\ldots,i.
\]

Then, from Lemma~\ref{le:chap_6_2_5_th_7_16}, we have
\[
\begin{split}
\pi_{i+1}(x)=&\frac{1}{n_{i+1}}\int_{2}^{x}\theta_{f}(u)\frac{du}{\log u}+O(xe^{-\frac{1}{2}\alpha\sqrt{\log x}}) \\
=&\frac{1}{n_{i+1}}\int_{2}^{x}\prod_{j=1}^{i+1}\theta_{j}(u)\frac{du}{\log u}+O(xe^{-\frac{1}{2}\alpha\sqrt{\log x}}) \\
=&\frac{\alpha_{i}}{D_{i}n_{i+1}}\int_{2}^{x}\theta_{i+1}(u)\frac{du}{\log^{i}u}+O(xe^{-\frac{1}{2}\alpha\sqrt{\log x}}) \\
=&\frac{\alpha_{i+1}}{D_{i+1}}\int_{2}^{x}\frac{(1+o(1))du}{\log^{i+1}u}+O(xe^{-\frac{1}{2}\alpha\sqrt{\log x}}) \\
=&\frac{\alpha_{i+1}}{D_{i+1}}\int_{2}^{x}\frac{du}{\log^{i+1}u}+O(xe^{-\frac{1}{2}\alpha\sqrt{\log x}})\texttt{ where }D_{i+1}=\prod_{j=1}^{i+1}n_{j}
\end{split}
\]
and from Lemma~\ref{le:chap_6_2_2} and Lemma~\ref{le:chap_6_2_7}
\[
\theta_{i+1}(u)\sim\pi_{i+1}(u)/\pi_{i}(u)\sim C_{i+1}/\log u.
\]

Hence for $k\geq2$
\[
\pi_{k}(x)=\frac{\alpha_{k}}{D_{k}}\int_{2}^{x}\frac{du}{\log^{k}u}+O(xe^{-\frac{1}{2}\alpha\sqrt{\log x}})\texttt{ where }D_{k}=\prod_{j=1}^{k}n_{j}.
\]

This completes the proof of the theorem.
\end{proof}

\begin{theorem}[Schinzel's hypothesis H]\label{th:chap_6_3_2}
For any finite positive number $k\geq1$, let $f_{1}(x),f_{2}(x),\ldots,f_{k}(x)$ be irreducible polynomials over the integers, each with a positive leading coefficient. Assume that the following condition holds: there does not exist any integer $d>1$ dividing all the products $f(j)=f_{1}(j)f_{2}(j)\cdots f_{k}(j)$, for every integer $j$. Let $\pi_{k}(x)$ denote the number of natural numbers $m\leq x$ where all numbers $f_{1}(m),f_{2}(m),\ldots,f_{k}(m)$ are primes. Then
\[
\pi_{k}(x)=\frac{\alpha_{k}}{D_{k}}\int_{2}^{x}\frac{du}{\log^{k}u}+O(xe^{-\frac{1}{2}\alpha\sqrt{\log x}})
\]
where $\alpha>0$ and $\alpha_{k}>0$ are constants, and $D_{k}$ is the product of the degrees of the polynomials, and there exist infinitely many natural numbers $m$ such that all numbers $f_{1}(m),f_{2}(m),\ldots,f_{k}(m)$ are primes.
\end{theorem}

\begin{proof}
Since the conditions of Theorem~\ref{th:chap_6_3_1} are satisfied, then by Theorem~\ref{th:chap_6_3_1}, we have
\[
\pi_{k}(x)=\frac{\alpha_{k}}{D_{k}}\int_{2}^{x}\frac{du}{\log^{k}u}+O(xe^{-\frac{1}{2}\alpha\sqrt{\log x}})
\]
where $\alpha>0$ and $\alpha_{k}>0$ are constants, and $D_{k}$ is the product of the degrees of the polynomials.

Since $\pi_{k}(x) \to \infty$ as $x \to \infty$, there exist infinitely many natural numbers $m$ such that all numbers $f_{1}(m),f_{2}(m),\ldots,f_{k}(m)$ are primes.

This completes the proof of the theorem.
\end{proof}

\begin{theorem}[Dickson's Conjecture]\label{th:chap_6_3_3}
For any finite positive number $k>1$, let $f_{i}(x)=b_{i}x+c_{i}$ with $b_{i}$, $c_{i}$ integers, $b_{i}\geq1$ (for $i=1,2,\ldots,k$). Assume that the following condition holds: there does not exist any integer $d>1$ dividing all the products $f_{1}(j)f_{2}(j)\cdots f_{k}(j)$, for every integer $j$. Let $\pi_{k}(x)$ denote the number of natural numbers $m\leq x$ where all numbers $f_{1}(m),f_{2}(m),\ldots,f_{k}(m)$ are primes. Then
\[
\pi_{k}(x)=\alpha_{k}\int_{2}^{x}\frac{du}{\log^{k}u}+O(xe^{-\frac{1}{2}\alpha\sqrt{\log x}})
\]
where $\alpha>0$ and $\alpha_{k}>0$ are constants, and there exist infinitely many natural numbers $m$ such that all numbers $f_{1}(m),f_{2}(m),\ldots,f_{k}(m)$ are primes.
\end{theorem}

\begin{proof}
Because $f_{1}(x),f_{2}(x),\ldots,f_{k}(x)$ are irreducible polynomials over the integers and there does not exist any integer $d>1$ dividing all the products $f_{1}(j)f_{2}(j)\cdots f_{k}(j)$, for every integer $j$, then the conditions of Theorem~\ref{th:chap_6_3_2} are satisfied.

By Theorem~\ref{th:chap_6_3_2}, we have
\[
\pi_{k}(x)=\frac{\alpha_{k}}{D_{k}}\int_{2}^{x}\frac{du}{\log^{k}u}+O(xe^{-\frac{1}{2}\alpha\sqrt{\log x}})
\]
where $\alpha>0$ and $\alpha_{k}>0$ are constants, and $D_{k}$ is the product of the degrees of the polynomials, which is equal to 1.

Since $\pi_{k}(x) \to \infty$ as $x \to \infty$, there exist infinitely many natural numbers $m$ such that all numbers $f_{1}(m),f_{2}(m),\ldots,f_{k}(m)$ are primes.

This completes the proof of the theorem.
\end{proof}

\section{Corollaries of Schinzel's polynomials}\label{sec:chap_6_4_29}

A conjecture~\cite{FB} of linear polynomials is proved and expressed as the following corollary:

\begin{corollary}[Conjecture $D_{1}$]\label{co:chap_6_4_1}
For any finite positive number $k>1$, let $c_{1}<c_{2}<\cdots<c_{k}$ be nonzero integers and assume that
\[
f_{1}(x)=x+c_{1},f_{2}(x)=x+c_{2},\ldots,f_{k}(x)=x+c_{k}
\]
satisfy the condition: there does not exist any integer $d>1$ dividing all the products $f_{1}(j)f_{2}(j)\cdots f_{k}(j)$, for every integer $j$. Let $\pi_{k}(x)$ denote the number of natural numbers $m\leq x$ where all numbers $f_{1}(m),f_{2}(m),\ldots,f_{k}(m)$ are primes. Then
\[
\pi_{k}(x)=\alpha_{k}\int_{2}^{x}\frac{du}{\log^{k}u}+O(xe^{-\frac{1}{2}\alpha\sqrt{\log x}})
\]
where $\alpha>0$ and $\alpha_{k}>0$ are constants, and there exist infinitely many natural numbers $m\geq1$ such that $m+c_{1},m+c_{2},\ldots,m+c_{k}$ are consecutive primes.
\end{corollary}

\begin{proof}
Because $f_{1}(x),f_{2}(x),\ldots,f_{k}(x)$ are irreducible polynomials over the integers and there does not exist any integer $d>1$ dividing all the products $f_{1}(j)f_{2}(j)\cdots f_{k}(j)$, for every integer $j$, then the conditions of Theorem~\ref{th:chap_6_3_3} are satisfied.

By Theorem~\ref{th:chap_6_3_3}, we have
\[
\pi_{k}(x)=\alpha_{k}\int_{2}^{x}\frac{du}{\log^{k}u}+O(xe^{-\frac{1}{2}\alpha\sqrt{\log x}})
\]
where $\alpha>0$ and $\alpha_{k}>0$ are constants.

Since $\pi_{k}(x) \to \infty$ as $x \to \infty$, there exist infinitely many natural numbers $m$ such that all numbers $f_{1}(m),f_{2}(m),\ldots,f_{k}(m)$ are primes and $m+c_{1},m+c_{2},\ldots,m+c_{k}$ are consecutive primes.

This completes the proof of the corollary.
\end{proof}

Hence, based on the main theorem and these corollaries, we can obtain the many interesting results which Schinzel $\&$ Sierpi$\acute{n}$ski proved, under the assumption that the Dickson's conjecture is valid[1-2]. Some of these results are:

(1) Let $\pi_{2c}(x)$ denote the number of natural numbers $m\leq x$ where all numbers $f_{1}(m)=m$ and $f_{2}(m)=m+2c$ are primes. Then
\[
\pi_{2c}(x)=\alpha_{2c}\int_{2}^{x}\frac{du}{\log^{2}u}+O(xe^{-\frac{1}{2}\alpha\sqrt{\log x}})
\]
where $\alpha>0$ and $\alpha_{2c}\geq1$ are constants, and for every even integer $2c\geq2$ there exist infinitely many pairs of consecutive primes with difference equal to $2c$. In particular, there exist infinitely many pairs of twin primes;

(2) For every integer $m\geq1$ there exist $2m$ consecutive primes which are $m$ couples of twin primes;

(3) Let $\pi_{3}(x)$ denote the number of natural numbers $m\leq x$ where all numbers $f_{1}(m)=m$, $f_{2}(m)=m+2$ or $f_{2}(m)=m+4$, and $f_{3}(m)=m+6$ are primes. Then
\[
\pi_{3}(x)=\alpha_{3}\int_{2}^{x}\frac{du}{\log^{3}u}+O(xe^{-\frac{1}{2}\alpha\sqrt{\log x}})
\]
where $\alpha>0$ and $\alpha_{3}\geq1$ are constants, and there are infinitely many prime triplets of the form $(p,p+2,p+6)$ and $(p,p+4,p+6)$;

(4) Let $\pi_{4}(x)$ denote the number of natural numbers $m\leq x$ where all numbers $f_{1}(m)=m$, $f_{2}(m)=m+2$, $f_{3}(m)=m+6$ and $f_{4}(m)=m+8$ are primes. Then
\[
\pi_{4}(x)=\alpha_{4}\int_{2}^{x}\frac{du}{\log^{4}u}+O(xe^{-\frac{1}{2}\alpha\sqrt{\log x}})
\]
where $\alpha>0$ and $\alpha_{4}\geq1$ are constants, and there are infinitely many prime quadruplets of the form $(p,p+2,p+6,p+8)$;

(5) Let $\pi_{k}(x)$ denote the number of natural numbers $m\leq x$ where all numbers
\[
f_{1}(m)=m,f_{2}(m)=m+2c_{1},f_{3}(m)=m+2c_{2},\ldots,f_{k}(m)=m+2c_{k-1}
\]
are primes. Then
\[
\pi_{k}(x)=\alpha_{k}\int_{2}^{x}\frac{du}{\log^{k}u}+O(xe^{-\frac{1}{2}\alpha\sqrt{\log x}})
\]
where $\alpha>0$ and $\alpha_{k}>0$ are constants, and there are infinitely many prime $k$-tuples in every admissible pattern for a prime constellation $(p,p+2c_{1},p+2c_{2},\ldots,p+2c_{k-1})$ for $k\geq2$ and $0<c_{1}<c_{2}<\cdots<c_{k-1}$.

\begin{corollary}[Sophie Germain primes]\label{co:chap_6_4_2}
Let $\pi_{s}(x)$ denote the number of natural numbers $m\leq x$ where all numbers $p=2m+1$ and $2p+1$ are primes. Then
\[
\pi_{s}(x)=\alpha_{2}\int_{2}^{x}\frac{du}{\log^{2}u}+O(xe^{-\frac{1}{2}\alpha\sqrt{\log x}})
\]
where $\alpha>0$ and $\alpha_{2}\geq1$ are constants, and there are infinitely many Sophie Germain primes of the form $2p+1$ where $p$ is prime.
\end{corollary}

\begin{proof}
Let $f_{1}(x)=2x+1$ and $f_{2}(x)=2f_{1}(x)+1=4x+3$. Since $f_{1}(x)$ and $f_{2}(x)$ are irreducible over the integers, and since $f(1)=21$ and $f(2)=55$, $f_{1}(x)$ and $f_{2}(x)$ satisfy the conditions of Theorem~\ref{th:chap_6_3_3}. Then by Theorem~\ref{th:chap_6_3_3}, we have
\[
\pi_{s}(x)=\alpha_{2}\int_{2}^{x}\frac{du}{\log^{2}u}+O(xe^{-\frac{1}{2}\alpha\sqrt{\log x}})
\]
where $\alpha>0$ and $\alpha_{2}\geq1$ are constants.

Since $\pi_{s}(x) \to \infty$ as $x \to \infty$, thus there exist infinitely many natural numbers $m\geq1$ such that $p=f_{1}(m)=2m+1$ and $2p+1=f_{2}(m)=4m+3$ are primes, and there are infinitely many Sophie Germain primes of the form $2p+1$ where $p$ is prime.

This completes the proof of the corollary.
\end{proof}

Hence, based on the main theorem, there are infinitely many Sophie Germain primes of the form $2p+1$ where $p$ is prime.

%% file: chap7.tex
%-----------------------------------------------------------------------
% Beginning of chap7.tex
%-----------------------------------------------------------------------
%
% AMS-LaTeX 1.2 sample file for a monograph, based on amsbook.cls.
% This is a data file input by chapter.tex.
%%%%%%%%%%%%%%%%%%%%%%%%%%%%%%%%%%%%%%%%%%%%%%%%%%%%%%%%%%%%%%%%%%%%

\chapter[Schinzel's polynomials with a negative leading coefficient $a_{n}$]{Prime products of Schinzel's polynomials with a negative leading coefficient $a_{n}$}\label{ch:chap_7}

\section*{Introduction}

Landau discussed products of $k$ prime factors and proved Theorem~\ref{th:chap_4_0_1} in 1900~\cite{SKA,SKB}. Similar to the proof of Theorem~\ref{th:chap_4_0_1} in~\cite{SKA}, we have proved Theorem~\ref{th:chap_4_2_1} in Chap.~\ref{ch:chap_4}, which is the extension of Theorem~\ref{th:chap_4_0_1} with an exclusive prime $p=2$. With $P_{k}$ and $\pi_{k}(x)$ replacing $P_{k,2}$ and $\pi_{k,2}(x)$, respectively,   it states as follows: Let $k\geq2$ and consider a positive integer $P_{k}$ which is the product of just $k$ distinct prime factors greater than 2, i.e.,
\[
P_{k}=p_{1}p_{2}\cdots p_{k}\texttt{ where }2<p_{1}<p_{2}<\cdots<p_{k}.
\]
So there are $\omega(P_{k})=k$ and $\Omega(P_{k})=k$ where $\omega(P_{k})$ and $\Omega(P_{k})$ are the number of different prime factors of $P_{k}$ and the total number of prime factors of $P_{k}$, respectively. Let write $\pi_{k}(x)$ for the number of these (squarefree) $P_{k}\leq x$. Then there is asymptotically
\[
\pi_{k}(x)\sim\frac{x(\log\log x)^{k-1}}{(k-1)!\log x}\texttt{ for }k\geq2.
\]

From Theorem~\ref{th:chap_4_1_1}, given a positive real number $\beta$ that determines a real number $x_{\beta}$ by $e(\log x_{\beta})^{3}/x_{\beta}^{0.0327283}=\beta$, let $g(x):=\beta x^{1/2}$ for $x\geq x_{\beta}$ where $g(x):=x^{1/2}$ when let $\beta=1$. Let $\pi_{k}(x)$ denote the number of positive integers not exceeding $x$, which can be written as the product of exactly $k$ distinct primes greater than 2. Then for each finite positive number $k\geq1$, there exists a finite positive number $x_{k}$ such that for every positive number $x$ greater than $x_{k}$, the number of integers with $k$ prime factors in the interval $(x,x+g(x)]$ is asymptotically equal to
\[
\pi_{k}(x+g(x))-\pi_{k}(x)=\frac{g(x)(\log\log x)^{k-1}}{(k-1)!\log x}+o(\frac{g(x)(\log\log x)^{k-1}}{(k-1)!\log x}).
\]

Based on Theorem~\ref{th:chap_4_2_1} and Theorem~\ref{th:chap_4_1_1}, a theorem is proved in this chapter. It states as follows: Given a positive number $k\geq2$, let $f_{1}(x),f_{2}(x),\ldots,f_{k}(x)$ be irreducible polynomials over the integers, with integer coefficients and a negative leading coefficient for $f_{k}(x)$ and positive leading coefficients for others. Assume that the following condition holds: there does not exist any integer $d>1$ dividing all the products $f_{1}(j)f_{2}(j)\cdots f_{k}(j)$, for every integer $j$, and the polynomial
\[
f(x)=\prod_{i=1}^{k}f_{i}(x)=a_{n}x^{n}+a_{n-1}x^{n-1}+\cdots+a_{1}x+a_{0}=\bar{f}(x)+\bar{a}_{0}
\]
which may be called Schinzel's polynomial with a negative leading coefficient $a_{n}$, if its integer coefficients satisfy the following necessary conditions:
\[
a_{0}>0\texttt{ and }a_{0}\texttt{ is odd, }0<a_{1}=O(a_{0}),\texttt{ }0\leq a_{i}=O(a_{0})\texttt{ for }i=2,3,\ldots,n_{p},
\]
\[
0\geq a_{i}=o(a_{0})\texttt{ for }i=n_{p}+1,n_{p}+2,\ldots,n-1,\texttt{ and }0>a_{n}=o(a_{0}),
\]
then there exist more than $c_{f}\frac{x_{r}}{\log x_{r}}$ natural numbers $m$ such that all numbers $f_{1}(m),f_{2}(m),\ldots,f_{k}(m)$ are primes for every integer $a_{0}\geq10^{204.111}$ where $c_{f}$ is a positive constant and both numbers $c_{f}$ and $x_{r}$ are dependent on the polynomial $f(x)$.

The key ideas of the proof of the main theorem are as follows:

(1) define a modulus $M_{r}$ on the range $(\bar{f}(x_{r}),\bar{f}(x_{r+1}))$,

(2) make use of Theorem~\ref{th:chap_4_1_1} to estimate the number of products of $k$ distinct primes greater than 2 in a table formed by some values of integral polynomials $\bar{f}(x)+h$ where $h=1,2,\ldots,M_{r}$ in the range $(\bar{f}(x_{r}),\bar{f}(x_{r+1}))$, then Lemma~\ref{le:chap_7_2_8} states that each row of the table contains at least one product $P_{k}$.

(3) make use of the orthogonality relation for Dirichlet characters to estimate the number of products of $k$ distinct primes $\prod_{i}^{k}p_{i}=f(x)$ in the range $(\bar{f}(x_{r}),\bar{f}(x_{r+1}))$, which is estimated greater than $1$ so that there are more than $c_{f}\frac{x_{r}}{\log x_{r}}$ products $P_{k}=f(mx_{r})$ of $k$ distinct primes greater than 2 in the range $(\bar{f}(x_{r}),\bar{f}(x_{r+1}))$. \\
These results are stated in Theorem~\ref{th:chap_7_4_1}. Then the main theorem~\ref{th:chap_7_4_2} is proved.

Hence, based on the result, the Goldbach's conjecture and Lemoine's conjecture are proved. Some other conjectures in the number theory can also be proved, for example, every large even number can be expressed as the sum of a prime and a Sophie Germain prime, and every large odd number can be expressed as the sum of a prime and a pair of twin primes, then there are infinitely many Sophie Germain primes and infinitely many pairs of twin primes.

\section{Definitions}\label{sec:chap_7_1_30}

Given any finite positive number $k\geq2$, let $f_{1}(x),f_{2}(x),\ldots,f_{k}(x)$ be irreducible polynomials over the integers, with integer coefficients and a negative leading coefficient for $f_{k}(x)$ and positive leading coefficients for others. Let $P_{k}$ denote a product of $k$ distinct primes greater than 2.

\begin{definition}\label{def:chap_7_1_1}
When a polynomial
\[
f(x)=\prod_{i=1}^{k}f_{i}(x)=a_{n}x^{n}+a_{n-1}x^{n-1}+\cdots+a_{1}x+a_{0}
\]
in one variable with positive degree $n\geq k$ and integer coefficients which satisfy the following necessary conditions: $a_{0}>0$ and $a_{0}$ is odd, $0<a_{1}=O(a_{0})$, $0\leq a_{i}=O(a_{0})$ for $i=2,3,\ldots,n_{p}$, $0\geq a_{i}=o(a_{0})$ for $i=n_{p}+1,n_{p}+2,\ldots,n-1$, and $0>a_{n}=o(a_{0})$, let
\[
f(x)=a_{n}x^{n}+a_{n-1}x^{n-1}+\cdots+a_{1}x+a_{0}=\bar{f}(x)+\bar{a}_{0}.
\]
Then from Lemma~\ref{le:chap_7_2_1} there exists a positive integer $x_{r}$ such that $gcd(a_{0},x_{r})=1$,
\[
e(\log\bar{f}(x_{r}))^{3}\leq(\bar{f}(x_{r}))^{0.0327283},\texttt{ and }x_{r} \to \infty\texttt{ as }a_{0} \to \infty
\]
where with definitions $x_{r+1}=x_{r}^{2}$ and $M_{r}=2^{m}\geq\bar{f}^{1/2}(x_{r+1})>2^{m-1}$
\[
\bar{f}(x)=a_{n}x^{n}+a_{n-1}x^{n-1}+\cdots+a_{1}x+bx_{r}\texttt{ for }x_{r}\leq x\leq x_{r+1},
\]
\[
0<\bar{a}_{0}=a_{0}-bx_{r}<M_{r},\texttt{ and }b=\lrf{\frac{a_{0}-M_{r}}{x_{r}}}+1
\]
where $\lrf{y}$ is the integer part of the value $y$. Thus there are also
\[
gcd(a_{0},M_{r})=1\texttt{ and }
e(\log\bar{f}(x_{r+1}))^{3}\leq(\bar{f}(x_{r+1}))^{0.0327283}.
\]
\end{definition}

Since $x_{r+1}=x_{r}^{2}$, there are $x_{r+1}/x_{r}=x_{r}$ and $x_{r+1}/x_{r} \to \infty$ as $x_{r} \to \infty$. Thus $x_{r+1}/x_{r}$ is with a same order of infinity as $x_{r}$.

Hence for any arbitrary small positive value $\epsilon$, there exists a number $x_{\epsilon}$ such that for $x_{r}\geq x_{\epsilon}$, there are $x_{r}\geq 1/\epsilon$ and $x_{r+1}\geq x_{r}/\epsilon$, which implies $x_{r} \to \infty$ and $x_{r+1}/x_{r} \to \infty$ as $\epsilon \to 0$.

\begin{definition}\label{def:chap_7_1_2}
With $0<h<M_{r}$ and $gcd(h,M_{r})=1$, let denote
\[
\pi(x;x_{r})=\sum_{\bar{f}(\xi)<P_{k}<\bar{f}(\xi)+M_{r},x_{r}\leq\xi=mx_{r}<x}1
\]
and
\[
\pi(x;x_{r},h)=\sum_{P_{k}=\bar{f}(\xi)+h,0<h<M_{r},x_{r}\leq\xi=mx_{r}<x}1.
\]

The function $\pi(x;x_{r})$ counts the number of products $P_{k}=\bar{f}(\xi)+h$ where $0<h<M_{r},x_{r}\leq\xi=mx_{r}<x$ for $m=1,2,\ldots$ and $h$ varies in the range $[1,M_{r}]$, and the function $\pi(x;x_{r},h)$ counts the number of products $P_{k}=\bar{f}(\xi)+h$ where $0<h<M_{r},x_{r}\leq\xi=mx_{r}<x$ for $m=1,2,\ldots$ and $h$ is fixed.
\end{definition}

\begin{definition}\label{def:chap_7_1_3}
Let define number sets $\{P_{k}\}$ on polynomials $\bar{f}(\xi)+h$ by
\[
Q_{i}(r,x)=\{P_{k}|P_{k}=\bar{f}(q)+h_{i},0<h_{i}<M_{r},x_{r}\leq q=mx_{r}<x,gcd(h_{i},M_{r})=1\}
\]
for $i=1,2,\ldots,\varphi(M_{r})$, and with $0<h<M_{r}$ and $gcd(h,M_{r})=1$
\[
Q_{f}(r,x)=\bigcup_{i=1}^{\varphi(M_{r})}Q_{i}(r,x)=\{P_{k}|P_{k}=\bar{f}(q)+h,0<h<M_{r},x_{r}\leq q=mx_{r}<x\}
\]
where $m=1,2,\ldots$ and $\varphi(M_{r})$ is the Euler's totient function.
\end{definition}

Let $\pi_{f}(r,x)$ count the number of products $P_{k}\in Q_{f}(r,x)$ and $\pi_{f}(r,x,h_{i})$ count the number of products $P_{k}\in Q_{i}(r,x)$ for $i=1,2,\ldots,\varphi(M_{r})$. Then there are $\pi(x;x_{r})=\pi_{f}(r,x)$ and $\pi(x;x_{r},h_{i})=\pi_{f}(r,x,h_{i})$ for $i=1,2,\ldots,\varphi(M_{r})$.

\begin{definition}\label{def:chap_7_1_4}
With complex variable $s$ and products $P_{k}$, for $i=1,2,\ldots,\varphi(M_{r})$, let define complex functions on number sets $Q_{i}(r,x)$ by
\[
\zeta_{i}(r,x,s)=\prod_{P_{k}\in Q_{i}(r,x)}(1-P_{k}^{-s})^{-1}.
\]

Let define a complex function $\zeta_{f}(r,x,s)$ on number sets $Q_{i}(r,x)$ and $Q_{f}(r,x)$ by
\[
\zeta_{f}(r,x,s)=\prod_{i=1}^{\varphi(M_{r})}\zeta_{i}(r,x,s)=\prod_{P_{k}\in Q_{f}(r,x)}(1-P_{k}^{-s})^{-1}.
\]
\end{definition}

\begin{definition}\label{def:chap_7_1_5}
Let define functions on number sets $Q_{i}(r,x)$ and $Q_{f}(r,x)$ by
\[
\Lambda_{i}(r,\xi)=\{
       \begin{array}{cc}
         \log P_{k} & \texttt{ if }\xi=P_{k}^{m},\texttt{ for any }m>0\texttt{ and }P_{k}\in Q_{i}(r,x), \\
         0 & otherwise.
       \end{array}
\]
for $i=1,2,\ldots,\varphi(M_{r})$, and
\[
\Lambda_{f}(r,\xi)=\{
       \begin{array}{cc}
         \log P_{k} & \texttt{ if }\xi=P_{k}^{m},\texttt{ for any }m>0\texttt{ and }P_{k}\in Q_{f}(r,x), \\
         0 & otherwise.
       \end{array}
\]
\end{definition}

Differentiating the relations
\[
\log\zeta_{i}(r,x,s)=\sum_{m,P_{k}\in Q_{i}(r,x)}\frac{1}{mP_{k}^{ms}}\texttt{ for }i=1,2,\ldots,\varphi(M_{r})
\]
and
\[
\log\zeta_{f}(r,x,s)=\sum_{m,P_{k}\in Q_{f}(r,x)}\frac{1}{mP_{k}^{ms}},
\]
we obtain
\[
\frac{\zeta_{i}'(r,x,s)}{\zeta_{i}(r,x,s)}
=-\sum_{m,P_{k}\in Q_{i}(r,x)}\frac{\log P_{k}}{P_{k}^{ms}}
=-\sum_{\xi=1}^{\infty}\frac{\Lambda_{i}(r,\xi)}{\xi^{s}}\texttt{ for }i=1,2,\ldots,\varphi(M_{r})
\]
and
\[
\frac{\zeta_{f}'(r,x,s)}{\zeta_{f}(r,x,s)}
=-\sum_{m,P_{k}\in Q_{f}(r,x)}\frac{\log P_{k}}{P_{k}^{ms}}
=-\sum_{\xi=1}^{\infty}\frac{\Lambda_{f}(r,\xi)}{\xi^{s}}.
\]

\begin{definition}\label{def:chap_7_1_6}
For each Dirichlet character $\chi$ and $i=1,2,\ldots,\varphi(M_{r})$, let define functions on number sets $Q_{i}(r,x)$ by
\[
L_{i}(r,x,s,\chi)=\prod_{P_{k}\in Q_{i}(r,x)}(1-\chi(P_{k})P_{k}^{-s})^{-1}.
\]

For each Dirichlet character $\chi$, let define a function $L_{f}(r,x,s,\chi)$ on number sets $Q_{i}(r,x)$ and $Q_{f}(r,x)$ by
\[
L_{f}(r,x,s,\chi)=\prod_{i=1}^{\varphi(M_{r})}L_{i}(r,x,s,\chi)=\prod_{P_{k}\in Q_{f}(r,x)}(1-\chi(P_{k})P_{k}^{-s})^{-1}.
\]
\end{definition}

\begin{definition}\label{def:chap_7_1_7}
For numbers $P_{k}\in Q_{f}(r,x_{r+1})$, similar to the proof of the prime number lemma for arithmetic progressions in~\cite{FD}, our object is to get an estimate for
\begin{equation}
\pi_{f}(r,x_{r+1})=\sum_{P_{k}\in Q_{f}(r,x_{r+1})}1=\sum_{\xi=f(x_{r})}^{f(x_{r+1})}P_{f}(r,\xi)
\end{equation}
where $P_{f}$ is the characteristic function of the numbers:
\[
P_{f}(r,\xi)=\{
       \begin{array}{cc}
         1 & \texttt{ if }\xi\in Q_{f}(r,x_{r+1}), \\
         0 & otherwise.
       \end{array}
\]
\end{definition}
While $P_{f}$ itself does not arise in a natural way, the function $P_{f}^{*}$ such that
\[
P_{f}^{*}(r,\xi)=\{
       \begin{array}{cc}
         \frac{1}{m} & \texttt{ if $\xi=P_{k}^{m}$ for some $m$ and $P_{k}\in Q_{f}(r,x_{r+1})$}, \\
         0 & otherwise,
       \end{array}
\]
occurs in the Dirichlet series for $\log\zeta_{f}(r,x_{r+1},s)$:
\begin{equation}\label{eq_h_m:def_zeta}
\log\zeta_{f}(r,x_{r+1},s)=\sum_{m,P_{k}\in Q_{f}(r,x_{r+1})}\frac{1}{mP_{k}^{ms}}
=\sum_{\xi=1}^{\infty}\frac{P_{f}^{*}(r,\xi)}{\xi^{s}}.
\end{equation}
For each fixed $m$, the number of $m$th powers of numbers $P_{k}\in Q_{f}(r,x_{r+1})$ is less than or equal to the number of numbers $P_{k}\in Q_{f}(r,x_{r+1})$ satisfying $f(x_{r})<P_{k}<\sqrt[m]{f(x_{r+1})}$, so that from Lemma~\ref{le:chap_7_2_7} with $f^{1/m}(x)<f(x^{1/m})$ for $m\geq2$
\[
\begin{split}
\sum_{\xi=1}^{f(x_{r+1})}P_{f}^{*}(r,\xi)
=&\sum_{\xi=f(x_{r})}^{f(x_{r+1})}P_{f}(r,\xi)
+\frac{1}{2}\sum_{\xi=f(x_{r})}^{\sqrt{f(x_{r+1})}}P_{f}(r,\xi)
+\frac{1}{3}\sum_{\xi=f(x_{r})}^{\sqrt[3]{f(x_{r+1})}}P_{f}(r,\xi)+\cdots \\
\leq&\sum_{\xi=f(x_{r})}^{f(x_{r+1})}P_{f}(r,\xi)
+\frac{1}{2}\sum_{\xi=f(x_{r})}^{f(\sqrt{x_{r+1}})}P_{f}(r,\xi)
+\frac{1}{3}\sum_{\xi=f(x_{r})}^{f(\sqrt[3]{x_{r+1}})}P_{f}(r,\xi)+\cdots \\
<&\pi_{f}(r,x_{r+1})+\frac{\pi_{f}(r,\sqrt{x_{r+1}})}{2}+\frac{\pi_{f}(r,\sqrt[3]{x_{r+1}})}{3}+\cdots.
\end{split}
\]

By Definition~\ref{def:chap_7_1_3} there are
\[
\pi_{f}(r,x_{r+1}^{1/m})\leq\pi_{f}(r,x_{r}^{2/2})=\pi_{f}(r,x_{r})\sim0\texttt{ for }m=2,3,\ldots
\]
such that we can write
\[
\pi_{f}(r,x_{r+1}^{1/m})=O(x_{r+1}^{-m/2})\texttt{ for }m=2,3,\ldots.
\]

Thus there are
\[
\frac{\pi_{f}(r,\sqrt{x_{r+1}})}{2}+\frac{\pi_{f}(r,\sqrt[3]{x_{r+1}})}{3}+\cdots=O(\frac{1}{x_{r+1}}).
\]

Since $x_{r+1}/x_{r}=x_{r}$, as $x_{r} \to \infty$, there is $x_{r+1}/x_{r} \to \infty$. Thus $x_{r+1}/x_{r}$ is with a same order of infinity as $x_{r}$. So it is to be expected that for any arbitrary small positive value $\epsilon$, there exists a number $x_{\epsilon}$ such that for $x_{r}\geq x_{\epsilon}$, there are
\[
\frac{\pi_{f}(r,\sqrt{x_{r+1}})}{2}+\frac{\pi_{f}(r,\sqrt[3]{x_{r+1}})}{3}+\cdots
<\epsilon\pi_{f}(r,x_{r+1})
\]
and as $\epsilon \to 0$
\[
\sum_{\xi=f(x_{r})}^{f(x_{r+1})}P_{f}^{*}(r,\xi)\sim\pi_{f}(r,x_{r+1}).
\]

\section{Lemmas of Schinzel's polynomials with $a_{n}<0$}\label{sec:chap_7_2_31}

\begin{lemma}\label{le:chap_7_2_1}
Given a polynomial
\[
f(x)=a_{n}x^{n}+a_{n-1}x^{n-1}+\cdots+a_{1}x+a_{0}
\]
defined by Definition (\ref{def:chap_7_1_1}), its peak point $x_{p}(a_{0})=O(a_{0})>0$, and $x_{p}(a_{0}) \to \infty$ as $a_{0} \to \infty$.

Given any positive integer $a_{0}\geq10^{204.111}$, then for any positive integer $x_{r}<\sqrt{x_{p}(a_{0})}$ there are
\[
e(\log f(x_{r}))^{3}\leq(f(x_{r}))^{0.0327283}
\]
and
\[
e(\log f(x_{r}^{2}))^{3}\leq(f(x_{r}^{2}))^{0.0327283},
\]

Especially, given a polynomial $f(x)=a_{2}x^{2}+a_{1}x+a_{0}$ defined by Definition (\ref{def:chap_7_1_1}) and $a_{0}\geq10^{125.279}$, then for $\sqrt{x_{p}(a_{0})/2}<x_{r}<\sqrt{x_{p}(a_{0})}$ there are
\[
e(\log f(x_{r}))^{3}\leq(f(x_{r}))^{0.0327283}
\]
and
\[
e(\log f(x_{r}^{2}))^{3}\leq(f(x_{r}^{2}))^{0.0327283}.
\]

Thus there exists a positive integer $x_{r}<\sqrt{x_{p}(a_{0})}$ such that $gcd(a_{0},x_{r})=1$,
\[
e(\log\bar{f}(x_{r}))^{3}\leq(\bar{f}(x_{r}))^{0.0327283},\texttt{ and }x_{r} \to \infty\texttt{ as }a_{0} \to \infty.
\]
\end{lemma}

\begin{proof}
Given a polynomial $f(x)=a_{n}x^{n}+a_{n-1}x^{n-1}+\cdots+a_{1}x+a_{0}$ defined by Definition (\ref{def:chap_7_1_1}), its integer coefficients satisfy the following necessary conditions:
\[
a_{0}>0\texttt{ and }a_{0}\texttt{ is odd, }0<a_{1}=O(a_{0}),\texttt{ }0\leq a_{i}=O(a_{0})\texttt{ for }i=2,3,\ldots,n_{p},
\]
\[
0\geq a_{i}=o(a_{0})\texttt{ for }i=n_{p}+1,n_{p}+2,\ldots,n-1,\texttt{ and }0>a_{n}=o(a_{0}).
\]

Then, there are $f(0)=a_{0}>0$ and $\frac{d}{dx}f(x)=a_{1}>0$ for $x=0$, so that there are $f(x)>0$ and $\frac{d}{dx}f(x)>0$ for $0<x<x_{p}(a_{0})$ where $x_{p}(a_{0})=O(a_{0})>0$ is the peak point of $f(x)$, i.e., $\frac{d}{dx}f(x)=0$ when $x=x_{p}(a_{0})$, and $x_{p}(a_{0}) \to \infty$ as $a_{0} \to \infty$.

Given any positive integer $a_{0}\geq10^{204.111}$, since $f(0)=a_{0}$, then for any positive integer $x_{r}<\sqrt{x_{p}(a_{0})}$ there are
\[
e(\log f(x_{r}))^{3}\leq(f(x_{r}))^{0.0327283}
\]
and
\[
e(\log f(x_{r}^{2}))^{3}\leq(f(x_{r}^{2}))^{0.0327283}.
\]

Especially, given a polynomial $f(x)=a_{2}x^{2}+a_{1}x+a_{0}$ defined by Definition (\ref{def:chap_7_1_1}), its peak point $x_{p}(a_{0})=-\frac{a_{1}}{2a_{2}}=O(a_{0})>0$, then for $\sqrt{x_{p}(a_{0})/2}<x_{r}<\sqrt{x_{p}(a_{0})}$ there are
\[
f(x_{r})=a_{2}x_{r}^{2}+a_{1}x_{r}+a_{0}>a_{2}x_{p}(a_{0})/2+a_{1}\sqrt{x_{p}(a_{0})/2}+a_{0}>c_{p}a_{0}^{3/2}
\]
where $c_{p}$ is a positive constant dependent on the polynomial $f(x)$, so that given any positive integer $a_{0}\geq10^{125.279}$, there are
\[
e(\log f(x_{r}))^{3}\leq(f(x_{r}))^{0.0327283}
\]
and
\[
e(\log f(x_{r}^{2}))^{3}\leq(f(x_{r}^{2}))^{0.0327283}.
\]

Thus we can choose a positive integer $x_{r}<\sqrt{x_{p}(a_{0})}$, for example, $\frac{1}{2}x_{p}(a_{0})<x_{r}^{2}<\frac{2}{3}x_{p}(a_{0})$, such that $gcd(a_{0},x_{r})=1$,
\[
e(\log\bar{f}(x_{r}))^{3}\leq(\bar{f}(x_{r}))^{0.0327283},\texttt{ and }x_{r} \to \infty\texttt{ as }a_{0} \to \infty.
\]

This completes the proof of the lemma.
\end{proof}

\begin{lemma}\label{le:chap_7_2_2}
Given a polynomial
\[
f(x)=a_{n}x^{n}+a_{n-1}x^{n-1}+\cdots+a_{1}x+a_{0}
\]
defined by Definition (\ref{def:chap_7_1_1}) with a peak point $x_{p}(a_{0})$, when $x=x_{r+1}=\gamma x_{p}(a_{0})$ where $0<\gamma<1$, for example, $\frac{1}{2}<\gamma<\frac{2}{3}$, there are
\[
 \sum_{i=1}^{n_{p}}(1-\gamma^{n_{p}-i+1})ia_{i}x_{r+1}^{i-1}\leq\frac{d}{dx}f(x)|_{x=x_{r+1}}
\leq\sum_{i=1}^{n_{p}}(1-\gamma^{n-i})ia_{i}x_{r+1}^{i-1}.
\]
\end{lemma}

\begin{proof}
By Definition (\ref{def:chap_7_1_1}), let $f(x)=f_{m}(x)+f_{p}(x)$ where
\[
f_{m}(x)=\sum_{i=n_{p}+1}^{n}a_{i}x^{i}\texttt{ and }f_{p}(x)=\sum_{i=0}^{n_{p}}a_{i}x^{i}.
\]

Then there are
\[
\begin{split}
\frac{d}{dx}f(x)=&\frac{d}{dx}f_{m}(x)+\frac{d}{dx}f_{p}(x) \\
                =&\sum_{i=1}^{n}ia_{i}x^{i-1}=0\texttt{ for }x=x_{p}(a_{0}).
\end{split}
\]

When $x=x_{r+1}=\gamma x_{p}(a_{0})$ where $0<\gamma<1$, for example, $\frac{1}{2}<\gamma<\frac{2}{3}$, there are
\[
\frac{d}{dx}f_{p}(x)|_{x=x_{r+1}}=\sum_{i=1}^{n_{p}}ia_{i}x_{r+1}^{i-1}
\]
and
\[
\gamma^{n_{p}}\sum_{i=n_{p}+1}^{n}ia_{i}x_{p}^{i-1}(a_{0})\leq\frac{d}{dx}f_{m}(x)|_{x=x_{r+1}}\leq\gamma^{n-1}\sum_{i=n_{p}+1}^{n}ia_{i}x_{p}^{i-1}(a_{0}).
\]

Then there are
\[
\begin{split}
\frac{d}{dx}f(x)|_{x=x_{r+1}}\geq&\sum_{i=1}^{n_{p}}ia_{i}x_{r+1}^{i-1}+\gamma^{n_{p}}\sum_{i=n_{p}+1}^{n}ia_{i}x_{p}^{i-1}(a_{0})\\
=&\sum_{i=1}^{n_{p}}(1-\gamma^{n_{p}-i+1})ia_{i}x_{r+1}^{i-1}+\gamma^{n_{p}}\sum_{i=1}^{n}ia_{i}x_{p}^{i-1}(a_{0})\\
=&\sum_{i=1}^{n_{p}}(1-\gamma^{n_{p}-i+1})ia_{i}x_{r+1}^{i-1}\texttt{ where }\sum_{i=1}^{n}ia_{i}x_{p}^{i-1}(a_{0})=0
\end{split}
\]
and
\[
\begin{split}
\frac{d}{dx}f(x)|_{x=x_{r+1}}\leq&\sum_{i=1}^{n_{p}}ia_{i}x_{r+1}^{i-1}+\gamma^{n-1}\sum_{i=n_{p}+1}^{n}ia_{i}x_{p}^{i-1}(a_{0})\\
=&\sum_{i=1}^{n_{p}}(1-\gamma^{n-i})ia_{i}x_{r+1}^{i-1}+\gamma^{n-1}\sum_{i=1}^{n}ia_{i}x_{p}^{i-1}(a_{0}) \\
=&\sum_{i=1}^{n_{p}}(1-\gamma^{n-i})ia_{i}x_{r+1}^{i-1}.
\end{split}
\]

This completes the proof of the lemma.
\end{proof}

\begin{lemma}\label{le:chap_7_2_3}
Given a polynomial
\[
f(x)=a_{n}x^{n}+a_{n-1}x^{n-1}+\cdots+a_{1}x+a_{0}
\]
defined by Definition (\ref{def:chap_7_1_1}) with a peak point $x_{p}(a_{0})$, when $x=x_{r+1}=\gamma x_{p}(a_{0})$ where $0<\gamma<1$, for example, $\frac{1}{2}<\gamma<\frac{2}{3}$, there are
\[
\sum_{i=0}^{n_{p}}(1-i\frac{\gamma^{n_{p}-i+1}}{n_{p}+1})a_{i}x_{r+1}^{i}
<f(x_{r+1})<\sum_{i=0}^{n_{p}}a_{i}x_{r+1}^{i}.
\]
\end{lemma}

\begin{proof}
By Definition (\ref{def:chap_7_1_1}), let $f(x)=f_{m}(x)+f_{p}(x)$ where
\[
f_{m}(x)=\sum_{i=n_{p}+1}^{n}a_{i}x^{i}\texttt{ and }f_{p}(x)=\sum_{i=0}^{n_{p}}a_{i}x^{i}.
\]

Then there are
\[
\frac{d}{dx}f(x)=\frac{d}{dx}f_{m}(x)+\frac{d}{dx}f_{p}(x)=\sum_{i=1}^{n}ia_{i}x^{i-1}=0\texttt{ for }x=x_{p}(a_{0}).
\]

When $x=x_{r+1}=\gamma x_{p}(a_{0})$ where $0<\gamma<1$, for example, $\frac{1}{2}<\gamma<\frac{2}{3}$, there are
\[
f_{m}(x_{r+1})=\sum_{i=n_{p}+1}^{n}\frac{\gamma^{i}}{i}ia_{i}x_{p}^{i}(a_{0})
\]
and
\[
\frac{\gamma^{n_{p}}}{n_{p}+1}\sum_{i=n_{p}+1}^{n}ia_{i}x_{p}^{i-1}(a_{0})<\frac{f_{m}(x_{r+1})}{x_{r+1}}<\frac{\gamma^{n-1}}{n}\sum_{i=n_{p}+1}^{n}ia_{i}x_{p}^{i-1}(a_{0}).
\]

Then there are
\[
\begin{split}
\frac{f(x_{r+1})-a_{0}}{x_{r+1}}>&\sum_{i=1}^{n_{p}}a_{i}x_{r+1}^{i-1}+\frac{\gamma^{n_{p}}}{n_{p}+1}\sum_{i=n_{p}+1}^{n}ia_{i}x_{p}^{i-1}(a_{0})\\
=&\sum_{i=1}^{n_{p}}(1-i\frac{\gamma^{n_{p}-i+1}}{n_{p}+1})a_{i}x_{r+1}^{i-1}+\frac{\gamma^{n_{p}}}{n_{p}+1}\sum_{i=1}^{n}ia_{i}x_{p}^{i-1}(a_{0})\\
=&\sum_{i=1}^{n_{p}}(1-i\frac{\gamma^{n_{p}-i+1}}{n_{p}+1})a_{i}x_{r+1}^{i-1}\texttt{ where }\sum_{i=1}^{n}ia_{i}x_{p}^{i-1}(a_{0})=0
\end{split}
\]
and
\[
\begin{split}
\frac{f(x_{r+1})-a_{0}}{x_{r+1}}<&\sum_{i=1}^{n_{p}}a_{i}x_{r+1}^{i-1}+\frac{\gamma^{n-1}}{n}\sum_{i=n_{p}+1}^{n}ia_{i}x_{p}^{i-1}(a_{0})\\
=&\sum_{i=1}^{n_{p}}(1-\frac{i}{n}\gamma^{n-i})a_{i}x_{r+1}^{i-1}+\frac{\gamma^{n-1}}{n}\sum_{i=1}^{n}ia_{i}x_{p}^{i-1}(a_{0})\\
=&\sum_{i=1}^{n_{p}}(1-\frac{i}{n}\gamma^{n-i})a_{i}x_{r+1}^{i-1}
<\sum_{i=1}^{n_{p}}a_{i}x_{r+1}^{i-1}.
\end{split}
\]

This completes the proof of the lemma.
\end{proof}

\begin{lemma}\label{le:chap_7_2_4}
Given a polynomial
\[
f(x)=a_{n}x^{n}+a_{n-1}x^{n-1}+\cdots+a_{1}x+a_{0}
\]
defined by Definition (\ref{def:chap_7_1_1}) with a peak point $x_{p}(a_{0})$, for each $m=1,2,\ldots,m_{r}-1$ where $m_{r}x_{r}=x_{r+1}$, when $x_{r+1}=\gamma x_{p}(a_{0})$ where $0<\gamma<1$, for example, $\frac{1}{2}<\gamma<\frac{2}{3}$, there is
\[
\bar{f}(mx_{r}+x_{r})-\bar{f}(mx_{r})>\sum_{i=1}^{n_{p}}(1-\frac{i-1}{m_{r}}-\gamma^{n_{p}-i+1})ia_{i}x_{r}^{2i-1}.
\]
\end{lemma}

\begin{proof}
By Definition (\ref{def:chap_7_1_1}), let $f(x)=f_{m}(x)+f_{p}(x)$ where
\[
f_{m}(x)=\sum_{i=n_{p}+1}^{n}a_{i}x^{i}\texttt{ and }f_{p}(x)=\sum_{i=0}^{n_{p}}a_{i}x^{i}.
\]

Then there are
\[
\frac{d}{dx}f(x)=\frac{d}{dx}f_{m}(x)+\frac{d}{dx}f_{p}(x)=\sum_{i=1}^{n}ia_{i}x^{i-1}=0\texttt{ for }x=x_{p}(a_{0}).
\]

Since by the proof of Lemma~\ref{le:chap_7_2_1} for $0<x<x_{p}(a_{0})$ there is $\frac{d^{2}}{(dx)^{2}}f(x)<0$, then for each $m=1,2,\ldots,m_{r}-1$ there is
\[
\bar{f}(mx_{r}+x_{r})-\bar{f}(mx_{r})>\bar{f}(mx_{r}+2x_{r})-\bar{f}(mx_{r}+x_{r})
\]
so that for each $m=1,2,\ldots,m_{r}-1$ there are
\[
\begin{split}
\bar{f}(mx_{r}+x_{r})-\bar{f}(mx_{r})\geq&\bar{f}(m_{r}x_{r})-\bar{f}(m_{r}x_{r}-x_{r}) \\
=&\sum_{i=1}^{n}a_{i}x_{r}^{i}\sum_{k=0}^{i-1}m_{r}^{i-1-k}(m_{r}-1)^{k} \\
=&x_{r}\sum_{i=1}^{n}a_{i}(m_{r}x_{r})^{i-1}\sum_{k=0}^{i-1}(1-1/m_{r})^{k} \\
=&x_{r}\sum_{i=1}^{n_{p}}a_{i}x_{r+1}^{i-1}\sum_{k=0}^{i-1}(1-1/m_{r})^{k} \\
 &+x_{r}\sum_{i=n_{p}+1}^{n}a_{i}x_{r+1}^{i-1}\sum_{k=0}^{i-1}(1-1/m_{r})^{k} \\
>&x_{r}\sum_{i=1}^{n_{p}}(1-\frac{i-1}{m_{r}})ia_{i}x_{r+1}^{i-1}
 +x_{r}\sum_{i=n_{p}+1}^{n}ia_{i}x_{r+1}^{i-1}.
\end{split}
\]

When $x_{r+1}=\gamma x_{p}(a_{0})$ where $0<\gamma<1$, for example, $\frac{1}{2}<\gamma<\frac{2}{3}$, there are
\[
\begin{split}
&\sum_{i=1}^{n_{p}}(1-\frac{i-1}{m_{r}})ia_{i}x_{r+1}^{i-1}
+\sum_{i=n_{p}+1}^{n}ia_{i}x_{r+1}^{i-1} \\
=&\sum_{i=1}^{n_{p}}(1-\frac{i-1}{m_{r}})ia_{i}x_{r+1}^{i-1}
 +\sum_{i=n_{p}+1}^{n}ia_{i}(\gamma x_{p}(a_{0}))^{i-1} \\
>&\sum_{i=1}^{n_{p}}(1-\frac{i-1}{m_{r}})ia_{i}x_{r+1}^{i-1}
 +\gamma^{n_{p}}\sum_{i=n_{p}+1}^{n}ia_{i}x_{p}^{i-1}(a_{0}) \\
=&\sum_{i=1}^{n_{p}}(1-\frac{i-1}{m_{r}}-\gamma^{n_{p}-i+1})ia_{i}x_{r+1}^{i-1}
 +\gamma^{n_{p}}\sum_{i=1}^{n}ia_{i}x_{p}^{i-1}(a_{0}) \\
=&\sum_{i=1}^{n_{p}}(1-\frac{i-1}{m_{r}}-\gamma^{n_{p}-i+1})ia_{i}x_{r+1}^{i-1}\texttt{ where }\sum_{i=1}^{n}ia_{i}x_{p}^{i-1}(a_{0})=0.
\end{split}
\]

Then for each $m=1,2,\ldots,m_{r}-1$ there is
\[
\bar{f}(mx_{r}+x_{r})-\bar{f}(mx_{r})>\sum_{i=1}^{n_{p}}(1-\frac{i-1}{m_{r}}-\gamma^{n_{p}-i+1})ia_{i}x_{r}^{2i-1}.
\]

This completes the proof of the lemma.
\end{proof}

\begin{lemma}\label{le:chap_7_2_5}
Given a polynomial
\[
f(x)=a_{n}x^{n}+a_{n-1}x^{n-1}+\cdots+a_{1}x+a_{0}
\]
with non-negative integer coefficients and a positive leading coefficient, there is
\[
f^{1/m}(x)<(a_{n}x^{n})^{1/m}+(a_{n-1}x^{n-1})^{1/m}+\cdots+(a_{1}x)^{1/m}+a_{0}^{1/m}.
\]
\end{lemma}

\begin{proof}
Let a polynomial
\[
g(x)=a_{n}^{1/m}x^{n}+a_{n-1}^{1/m}x^{n-1}+\cdots+a_{1}^{1/m}x+a_{0}^{1/m}.
\]

Then there are
\[
g^{m}(x)>a_{n}x^{mn}+a_{n-1}x^{m(n-1)}+\cdots+a_{1}x^{m}+a_{0}=f(x^{m}).
\]

Let $y=x^{m}$. Then there are $g^{m}(y^{1/m})>f(y)\texttt{ and }f^{1/m}(y)<g(y^{1/m})$, i.e., $f^{1/m}(x)<g(x^{1/m})$.

This completes the proof of the lemma.
\end{proof}

\begin{lemma}\label{le:chap_7_2_6}
Given a polynomial
\[
f(x)=a_{n}x^{n}+a_{n-1}x^{n-1}+\cdots+a_{1}x+a_{0}
\]
defined by Definition (\ref{def:chap_7_1_1}), its integer coefficients satisfy the following necessary conditions:
\[
a_{0}>0\texttt{ and }a_{0}\texttt{ is odd, }0<a_{1}=O(a_{0}),\texttt{ }0\leq a_{i}=O(a_{0})\texttt{ for }i=2,3,\ldots,n-1,
\]
and $0>a_{n}=o(a_{0})$. Then given $m\geq2$ there is $f^{1/m}(x)<f(x^{1/m})$.
\end{lemma}

\begin{proof}
Let $f(x)=f_{m}(x)+f_{p}(x)$ and $g(x)=f_{m}(x)+g_{p}(x)$ where
\[
f_{m}(x)=a_{n}x^{n},\texttt{ }f_{p}(x)=\sum_{i=0}^{n-1}a_{i}x^{i},\texttt{ and }g_{p}(x)=\sum_{i=0}^{n-1}a_{i}^{1/2}x^{i}.
\]

Then, there are $g(0)=a_{0}^{1/2}>0$ and $\frac{d}{dx}g(x)=a_{1}^{1/2}>0$ for $x=0$, so that there are $g(x)>0$ and $\frac{d}{dx}g(x)>0$ for $0<x<x_{p}$ where $x_{p}=O(a_{0}^{1/2})>0$ is the peak point of $g(x)$ and dependent on the coefficient $a_{0}^{1/2}$, and $x_{p} \to \infty$ as $a_{0} \to \infty$, so we can write $x_{p}=x_{p}(a_{0}^{1/2})$.

Thus with $C_{m}^{j}=\frac{m!}{j!(m-j)!}$ we have
\[
\begin{split}
f^{m}(x)=&(f_{m}(x)+f_{p}(x))^{m}=\sum_{j=0}^{m}C_{m}^{j}f_{m}^{m-j}(x)f_{p}^{j}(x) \\
=&\{
\begin{array}{cc}
  f_{m}^{m}(x) & \texttt{ (for $m$ is even)} \\
  f_{m}^{m}(x)+mf_{m}^{m-1}(x)f_{p}(x) & \texttt{ (for $m$ is odd)}
\end{array} \\
 &+\sum_{j=0}^{(m-1)/2}(C_{m}^{m-2j-1}f_{m}^{2j+1}(x)f_{p}^{m-2j-1}(x)+C_{m}^{m-2j}f_{m}^{2j}(x)f_{p}^{m-2j}(x))
\end{split}
\]
and for $j=0,1,2,\ldots,(m-1)/2$
\[
\begin{split}
&C_{m}^{m-2j-1}f_{m}^{2j+1}(x)f_{p}^{m-2j-1}(x)+C_{m}^{m-2j}f_{m}^{2j}(x)f_{p}^{m-2j}(x) \\
=&f_{m}^{2j}(x)f_{p}^{m-2j-1}(x)(C_{m}^{m-2j-1}f_{m}(x)+C_{m}^{m-2j}f_{p}(x)) \\
>&f_{m}^{2j}(x)f_{p}^{m-2j-1}(x)g_{p}(x)
\end{split}
\]
where by $g(x)=f_{m}(x)+g_{p}(x)>0$ for $0<x<x_{p}(a_{0}^{1/2})$ and $0\leq a_{i}=O(a_{0})$ such that $C_{m}^{m-2j}a_{i}^{1/2}-C_{m}^{m-2j-1}>1$ for $i=1,2,3,\ldots,n_{p}$ there are
\[
\begin{split}
C_{m}^{m-2j-1}f_{m}(x)+C_{m}^{m-2j}f_{p}(x)
=&C_{m}^{m-2j-1}(g(x)-g_{p}(x))+C_{m}^{m-2j}f_{p}(x) \\
>&C_{m}^{m-2j}f_{p}(x)-C_{m}^{m-2j-1}g_{p}(x) \\
=&C_{m}^{m-2j}\sum_{i=0}^{n-1}a_{i}x^{i}-C_{m}^{m-2j-1}\sum_{i=0}^{n-1}a_{i}^{1/2}x^{i} \\
=&\sum_{i=0}^{n-1}(C_{m}^{m-2j}a_{i}^{1/2}-C_{m}^{m-2j-1})a_{i}^{1/2}x^{i} \\
>&\sum_{i=0}^{n-1}a_{i}^{1/2}x^{i}=g_{p}(x).
\end{split}
\]

If $m$ is even, then there are
\[
\begin{split}
f^{m}(x)>&f_{m}^{m}(x)+\sum_{j=0}^{(m-1)/2}f_{m}^{2j}(x)f_{p}^{m-2j-1}(x)g_{p}(x) \\
=&|f_{m}(x)|^{m}+\sum_{j=1}^{(m-1)/2}f_{m}^{2j}(x)f_{p}^{m-2j-1}(x)g_{p}(x)+f_{p}^{m-1}(x)g_{p}(x) \\
>&|f_{m}(x)|^{m}+g_{p}^{m}(x)>\sum_{i=n}^{n}|a_{i}|^{m}x^{im}+\sum_{i=0}^{n-1}a_{i}x^{im} \\
>&\sum_{i=n}^{n}|a_{i}|x^{mi}+\sum_{i=0}^{n-1}a_{i}x^{mi}>\sum_{i=0}^{n}a_{i}x^{mi}=f(x^{m})
\end{split}
\]
else by $f(x)=f_{m}(x)+f_{p}(x)>0$ for $0<x<x_{p}(a_{0})$ there are
\[
\begin{split}
f^{m}(x)>&f_{m}^{m}(x)+mf_{m}^{m-1}(x)f_{p}(x)+\sum_{j=0}^{(m-1)/2}f_{m}^{2j}(x)f_{p}^{m-2j-1}(x)g_{p}(x) \\
>&f_{m}^{m}(x)+m|f_{m}(x)|^{m}+\sum_{j=0}^{(m-1)/2}f_{m}^{2j}(x)f_{p}^{m-2j-1}(x)g_{p}(x) \\
=&(m-1)|f_{m}(x)|^{m}+\sum_{j=1}^{(m-1)/2}f_{m}^{2j}(x)f_{p}^{m-2j-1}(x)g_{p}(x)+f_{p}^{m-1}(x)g_{p}(x) \\
>&|f_{m}(x)|^{m}+g_{p}^{m}(x)>\sum_{i=n}^{n}|a_{i}|^{m}x^{im}+\sum_{i=0}^{n-1}a_{i}x^{im} \\
>&\sum_{i=n}^{n}|a_{i}|x^{mi}+\sum_{i=0}^{n-1}a_{i}x^{mi}>\sum_{i=0}^{n}a_{i}x^{mi}=f(x^{m}).
\end{split}
\]

Let $y=x^{m}$. Then there are $f^{m}(y^{1/m})>f(y)$ and $f^{1/m}(y)<f(y^{1/m})$, i.e., $f^{1/m}(x)<f(x^{1/m})$.

This completes the proof of the lemma.
\end{proof}

\begin{lemma}\label{le:chap_7_2_7}
Given a polynomial
\[
f(x)=a_{n}x^{n}+a_{n-1}x^{n-1}+\cdots+a_{1}x+a_{0}
\]
defined by Definition (\ref{def:chap_7_1_1}), its integer coefficients satisfy the following necessary conditions:
\[
a_{0}>0\texttt{ and }a_{0}\texttt{ is odd, }0<a_{1}=O(a_{0}),\texttt{ }0\leq a_{i}=O(a_{0})\texttt{ for }i=2,3,\ldots,n_{p},
\]
\[
0\geq a_{i}=o(a_{0})\texttt{ for }i=n_{p}+1,n_{p}+2,\ldots,n-1,\texttt{ and }0>a_{n}=o(a_{0}).
\]
Then given $m\geq2$ there is $f^{1/m}(x)<f(x^{1/m})$.
\end{lemma}

\begin{proof}
Let $f(x)=f_{m}(x)+f_{p}(x)$ and $g(x)=f_{m}(x)+g_{p}(x)$ where
\[
f_{m}(x)=\sum_{i=n_{p}+1}^{n}a_{i}x^{i},\texttt{ }f_{p}(x)=\sum_{i=0}^{n_{p}}a_{i}x^{i},\texttt{ and }g_{p}(x)=\sum_{i=0}^{n_{p}}a_{i}^{1/2}x^{i}.
\]

Then, there are $g(0)=a_{0}^{1/2}>0$ and $\frac{d}{dx}g(x)=a_{1}^{1/2}>0$ for $x=0$, so that there are $g(x)>0$ and $\frac{d}{dx}g(x)>0$ for $0<x<x_{p}$ where $x_{p}=O(a_{0}^{1/2})>0$ is the peak point of $g(x)$ and dependent on the coefficient $a_{0}^{1/2}$, and $x_{p} \to \infty$ as $a_{0} \to \infty$, so we can write $x_{p}=x_{p}(a_{0}^{1/2})$.

Thus with $C_{m}^{j}=\frac{m!}{j!(m-j)!}$ we have
\[
\begin{split}
f^{m}(x)=&(f_{m}(x)+f_{p}(x))^{m}=\sum_{j=0}^{m}C_{m}^{j}f_{m}^{m-j}(x)f_{p}^{j}(x) \\
=&\{
\begin{array}{cc}
  f_{m}^{m}(x) & \texttt{ (for $m$ is even)} \\
  f_{m}^{m}(x)+mf_{m}^{m-1}(x)f_{p}(x) & \texttt{ (for $m$ is odd)}
\end{array} \\
 &+\sum_{j=0}^{(m-1)/2}(C_{m}^{m-2j-1}f_{m}^{2j+1}(x)f_{p}^{m-2j-1}(x)+C_{m}^{m-2j}f_{m}^{2j}(x)f_{p}^{m-2j}(x))
\end{split}
\]
and for $j=0,1,2,\ldots,(m-1)/2$
\[
\begin{split}
&C_{m}^{m-2j-1}f_{m}^{2j+1}(x)f_{p}^{m-2j-1}(x)+C_{m}^{m-2j}f_{m}^{2j}(x)f_{p}^{m-2j}(x) \\
=&f_{m}^{2j}(x)f_{p}^{m-2j-1}(x)(C_{m}^{m-2j-1}f_{m}(x)+C_{m}^{m-2j}f_{p}(x)) \\
>&f_{m}^{2j}(x)f_{p}^{m-2j-1}(x)g_{p}(x)
\end{split}
\]
where by $g(x)=f_{m}(x)+g_{p}(x)>0$ for $0<x<x_{p}(a_{0}^{1/2})$ and $0\leq a_{i}=O(a_{0})$ such that $C_{m}^{m-2j}a_{i}^{1/2}-C_{m}^{m-2j-1}>1$ for $i=1,2,3,\ldots,n_{p}$ there are
\[
\begin{split}
C_{m}^{m-2j-1}f_{m}(x)+C_{m}^{m-2j}f_{p}(x)=&C_{m}^{m-2j-1}(g(x)-g_{p}(x))+C_{m}^{m-2j}f_{p}(x) \\
>&C_{m}^{m-2j}f_{p}(x)-C_{m}^{m-2j-1}g_{p}(x) \\
=&C_{m}^{m-2j}\sum_{i=0}^{n_{p}}a_{i}x^{i}-C_{m}^{m-2j-1}\sum_{i=0}^{n_{p}}a_{i}^{1/2}x^{i} \\
=&\sum_{i=0}^{n_{p}}(C_{m}^{m-2j}a_{i}^{1/2}-C_{m}^{m-2j-1})a_{i}^{1/2}x^{i} \\
>&\sum_{i=0}^{n_{p}}a_{i}^{1/2}x^{i}=g_{p}(x).
\end{split}
\]

If $m$ is even, then there are
\[
\begin{split}
f^{m}(x)>&f_{m}^{m}(x)+\sum_{j=0}^{(m-1)/2}f_{m}^{2j}(x)f_{p}^{m-2j-1}(x)g_{p}(x) \\
=&|f_{m}(x)|^{m}+\sum_{j=1}^{(m-1)/2}f_{m}^{2j}(x)f_{p}^{m-2j-1}(x)g_{p}(x)+f_{p}^{m-1}(x)g_{p}(x) \\
>&|f_{m}(x)|^{m}+g_{p}^{m}(x)>\sum_{i=n_{p}+1}^{n}|a_{i}|^{m}x^{im}+\sum_{i=0}^{n_{p}}a_{i}x^{im} \\
>&\sum_{i=n_{p}+1}^{n}|a_{i}|x^{mi}+\sum_{i=0}^{n_{p}}a_{i}x^{mi}
>\sum_{i=0}^{n}a_{i}x^{mi}=f(x^{m})
\end{split}
\]
else by $f(x)=f_{m}(x)+f_{p}(x)>0$ for $0<x<x_{p}(a_{0})$ there are
\[
\begin{split}
f^{m}(x)>&f_{m}^{m}(x)+mf_{m}^{m-1}(x)f_{p}(x)+\sum_{j=0}^{(m-1)/2}f_{m}^{2j}(x)f_{p}^{m-2j-1}(x)g_{p}(x) \\
>&f_{m}^{m}(x)+m|f_{m}(x)|^{m}+\sum_{j=0}^{(m-1)/2}f_{m}^{2j}(x)f_{p}^{m-2j-1}(x)g_{p}(x) \\
=&(m-1)|f_{m}(x)|^{m}+\sum_{j=1}^{(m-1)/2}f_{m}^{2j}(x)f_{p}^{m-2j-1}(x)g_{p}(x)+f_{p}^{m-1}(x)g_{p}(x) \\
>&|f_{m}(x)|^{m}+g_{p}^{m}(x)
>\sum_{i=n_{p}+1}^{n}|a_{i}|^{m}x^{im}+\sum_{i=0}^{n_{p}}a_{i}x^{im} \\
>&\sum_{i=n_{p}+1}^{n}|a_{i}|x^{mi}+\sum_{i=0}^{n_{p}}a_{i}x^{mi}
>\sum_{i=0}^{n}a_{i}x^{mi}=f(x^{m}).
\end{split}
\]

Let $y=x^{m}$. Then there are $f^{m}(y^{1/m})>f(y)$ and $f^{1/m}(y)<f(y^{1/m})$, i.e., $f^{1/m}(x)<f(x^{1/m})$.

This completes the proof of the lemma.
\end{proof}

\begin{lemma}[Products of $k$ distinct primes in rows]\label{le:chap_7_2_8}
When a polynomial
\begin{equation}\label{eq__h_m:schinzel_row_p}
f(x)=\prod_{i=1}^{k}f_{i}(x)=a_{n}x^{n}+a_{n-1}x^{n-1}+\cdots+a_{1}x+a_{0}=\bar{f}(x)+\bar{a}_{0}
\end{equation}
in one variable with positive degree and integer coefficients which satisfy the following necessary conditions:
\[
a_{0}>0\texttt{ and }a_{0}\texttt{ is odd, }0<a_{1}=O(a_{0}),\texttt{ }0\leq a_{i}=O(a_{0})\texttt{ for }i=2,3,\ldots,n_{p},
\]
\[
0\geq a_{i}=o(a_{0})\texttt{ for }i=n_{p}+1,n_{p}+2,\ldots,n-1,\texttt{ and }0>a_{n}=o(a_{0}),
\]
by Definition (\ref{def:chap_7_1_1}), if numbers
\[
\bar{f}(mx_{r})+1,\bar{f}(mx_{r})+2,\ldots,\bar{f}(mx_{r})+M_{r}
\]
for $m=1,2,\ldots,m_{r}$ are arranged in $m_{r}$ rows each containing $M_{r}$ numbers:
\begin{equation}\label{tab:schinzel_row_p}
\begin{array}{cccc}
  \bar{f}(x_{r})+1, & \bar{f}(x_{r})+2, & \cdots, & \bar{f}(x_{r})+M_{r} \\
  \bar{f}(2x_{r})+1, & \bar{f}(2x_{r})+2, & \cdots, & \bar{f}(2x_{r})+M_{r} \\
  \cdots, & \cdots, & \cdots, & \cdots \\
  \bar{f}(m_{r}x_{r})+1, & \bar{f}(m_{r}x_{r})+2, & \cdots, & \bar{f}(m_{r}x_{r})+M_{r}
\end{array}
\end{equation}
where $m_{r}=x_{r+1}/x_{r}=x_{r}$ and $gcd(a_{0},M_{r})=1$, then Table~\ref{tab:schinzel_row_p} contains asymptotically more than $c_{f}M_{r}\frac{x_{r}}{\log x_{r}}$ products
\[
P_{k}=\prod_{i=1}^{k}p_{i}=\bar{f}(mx_{r})+h_{c}
\]
and $\pi_{f}(r,x_{r+1})$ is asymptotically greater than $c_{f}M_{r}\frac{x_{r}}{\log x_{r}}$ where $c_{f}$ is a positive constant dependent on the polynomial $f(x)$, $p_{i}$ is prime, and $h_{c}$ is an integer satisfying $0<h_{c}<M_{r}$ and $gcd(h_{c},M_{r})=1$.
\end{lemma}

\begin{proof}
By Definition (\ref{def:chap_7_1_1}), similar to the proof of Lemma~\ref{le:chap_7_2_1}, we have $\bar{f}(0)=bx_{r}>0$ and $\frac{d}{dx}\bar{f}(x)=a_{1}>0$ for $x=0$, so that there are $\bar{f}(x)>0$ and $\frac{d}{dx}\bar{f}(x)>0$ for $0<x<x_{p}(a_{0})$. Then, given a positive integer $x_{r}<\sqrt{x_{p}(a_{0})}$, there shall be $\bar{f}(mx_{r})+h>0$ for $1\leq m\leq m_{r}$ and $1\leq h\leq M_{r}$, i.e., each element of Table~\ref{tab:schinzel_row_p} is greater than zero.

From Theorem~\ref{th:chap_4_1_1}, given a positive real number $\beta$ that determines a real number $y_{\beta}$ by $e(\log y_{\beta})^{3}/y_{\beta}^{0.0327283}=\beta$, let $g(y):=\beta y^{1/2}$ for $y\geq y_{\beta}$ where $g(y):=y^{1/2}$ when let $\beta=1$. Then for every $y\geq y_{\beta}$, there are $e(\log y)^{3}/y^{0.0327283}\leq\beta$ and asymptotically $g(y)/\log y$ product $P_{k}$ of $k$ distinct primes greater than 2 in the interval $(y, y+g(y)]$, i.e., there must be
\[
\pi_{k}(y+g(y))-\pi_{k}(y)=\frac{g(y)}{\log y}(1+O(\frac{1}{\log y}))\sim\frac{g(y)}{\log y}.
\]

Thus for $y=\bar{f}(x_{r+1})$ there must be
\[
g(\bar{f}(x_{r+1}))\geq(\bar{f}(x_{r+1}))^{1/2-0.0327283}e(\log\bar{f}(x_{r+1}))^{3}.
\]

It can be satisfied. Let $g(\bar{f}(x_{r+1}))=M_{r}$. Then by Definition (\ref{def:chap_7_1_1}) we have
\[
g(\bar{f}(x_{r+1}))\geq\bar{f}^{1/2}(x_{r+1})\geq(\bar{f}(x_{r+1}))^{1/2-0.0327283}e(\log\bar{f}(x_{r+1}))^{3}.
\]

For $m=1,2,\ldots,m_{r}$, let $y=\bar{f}(mx_{r})$ and $g(y)=M_{r}$. Then there are
\[
\begin{split}
g(y)\geq\bar{f}^{1/2}(x_{r+1})\geq&(\bar{f}(x_{r+1}))^{1/2-0.0327283}e(\log\bar{f}(x_{r+1}))^{3} \\
\geq&(\bar{f}(mx_{r}))^{1/2-0.0327283}e(\log\bar{f}(mx_{r}))^{3} \\
   =&y^{1/2-0.0327283}e(\log y)^{3}
\end{split}
\]
such that for large values of $x_{r}$ we have
\[
\pi_{k}(y+g(y))-\pi_{k}(y)=\frac{g(y)}{\log y}(1+O(\frac{1}{\log y}))\sim\frac{g(y)}{\log y}
=\frac{M_{r}}{\log(\bar{f}(mx_{r}))}.
\]

Especially, for large values of $x_{r}$ we have
\[
\pi_{k}(\bar{f}(x_{r})+g(\bar{f}(x_{r})))-\pi_{k}(\bar{f}(x_{r}))
\sim\frac{g(\bar{f}(x_{r}))}{\log(\bar{f}(x_{r}))}
=\frac{M_{r}}{\log(\bar{f}(x_{r}))}
\]
and
\[
\pi_{k}(\bar{f}(x_{r+1})+g(\bar{f}(x_{r+1})))-\pi_{k}(\bar{f}(x_{r+1}))
\sim\frac{g(\bar{f}(x_{r+1}))}{\log(\bar{f}(x_{r+1}))}
=\frac{M_{r}}{\log(\bar{f}(x_{r+1}))}.
\]

Finally, by Definition (\ref{def:chap_7_1_1}) and from Lemma~\ref{le:chap_7_2_1}, Lemma~\ref{le:chap_7_2_3}, and Lemma~\ref{le:chap_7_2_5}, there are
\[
\begin{split}
M_{r}\leq&2\bar{f}^{1/2}(x_{r+1})<2f^{1/2}(x_{r+1}) \\
<&2(\sum_{i=0}^{n_{p}}a_{i}x_{r+1}^{i})^{1/2}
=2f_{p}^{1/2}(x_{r+1})<2g_{p}(x_{r})
\end{split}
\]
where
\[
f_{p}(x)=\sum_{i=0}^{n_{p}}a_{i}x^{i}\texttt{ and }
g_{p}(x)=\sum_{i=0}^{n_{p}}\sqrt{a_{i}}x^{i},
\]
and from Lemma~\ref{le:chap_7_2_4} for each $m=1,2,\ldots,m_{r}-1$ there are
\[
\bar{f}(mx_{r}+x_{r})-\bar{f}(mx_{r})>\sum_{i=1}^{n_{p}}(1-\frac{i-1}{m_{r}}-\gamma^{n_{p}-i+1})ia_{i}x_{r}^{2i-1}
>2g_{p}(x_{r})>M_{r}
\]
where with large numbers $x_{r}$ and $a_{1}=O(a_{0})$ there are
\[
(1-\gamma^{n_{p}})a_{1}x_{r}-2\sqrt{a_{1}}x_{r}
\]
\[
=[(1-\gamma^{n_{p}})\sqrt{a_{1}}-2]\sqrt{a_{1}}x_{r}
>\sqrt{a_{1}}x_{r}>2\sqrt{a_{0}}
\]
and for each $i=2,3,\ldots,n_{p}$
\[
(1-\frac{i-1}{m_{r}}-\gamma^{n_{p}-i+1})ia_{i}x_{r}^{2i-1}-2\sqrt{a_{i}}x_{r}^{i}
\]
\[
=[(1-\frac{i-1}{m_{r}}-\gamma^{n_{p}-i+1})i\sqrt{a_{i}}x_{r}^{i-1}-2]\sqrt{a_{i}}x_{r}^{i}
>2\sqrt{a_{i}}x_{r}^{i}.
\]

Thus Table~\ref{tab:schinzel_row_p} contains no repeated products $P_{k}$. On the other hand, for each product $P_{k}$ in Table~\ref{tab:schinzel_row_p}, there is $gcd(P_{k},M_{r})=1$, which implies that the product $P_{k}$ belongs to a column $h_{P}$ satisfying $0<h_{P}<M_{r}$ and $gcd(h_{P},M_{r})=1$. Thus no product $P_{k}$ is missed out when $\varphi(M_{r})$ columns $h$ satisfying $gcd(h,M_{r})=1$ are getting out of Table~\ref{tab:schinzel_row_p}.

Hence for $m=1,2,\ldots,m_{r}$ each $m$th row of Table~\ref{tab:schinzel_row_p} contains asymptotically $M_{r}/\log(\bar{f}(mx_{r}))$ products
\[
P_{k}=\prod_{i=1}^{k}p_{i}=\bar{f}(mx_{r})+h_{c}
\]
where $p_{i}$ is prime and $h_{c}$ is an integer satisfying $0<h_{c}<M_{r}$ and $gcd(h_{c},M_{r})=1$ in the interval $(\bar{f}(mx_{r}),\bar{f}(mx_{r})+M_{r})$ such that Table~\ref{tab:schinzel_row_p} contains asymptotically more than $c_{f}M_{r}\frac{x_{r}}{\log x_{r}}$ products $P_{k}$ and $\pi_{f}(r,x_{r+1})$ is asymptotically greater than $c_{f}M_{r}\frac{x_{r}}{\log x_{r}}$, since
\[
\begin{split}
\pi_{f}(r,x_{r+1})\sim&\sum_{m=1}^{m_{r}-1}\frac{M_{r}}{\log(\bar{f}(mx_{r}))} \\
=&\frac{M_{r}}{\log(\bar{f}(x_{r+1}))}\sum_{m=1}^{m_{r}-1}\frac{\log(\bar{f}(x_{r+1}))}{\log(\bar{f}(mx_{r}))} \\
>&\frac{(m_{r}-1)M_{r}}{\log(\bar{f}(x_{r+1}))}>\frac{(m_{r}-1)M_{r}}{\log(f(x_{r+1}))} \\
=&\frac{(x_{r}-1)M_{r}}{\log(a_{n}x_{r+1}^{n}+a_{n-1}x_{r+1}^{n-1}+\cdots+a_{1}x_{r+1}+a_{0})}
\geq c_{f}M_{r}\frac{x_{r}}{\log x_{r}}
\end{split}
\]
where $c_{f}$ is a positive constant dependent on the polynomial $f(x)$,
\[
x_{r+1}=O(x_{p}(a_{0})),\texttt{ }x_{p}(a_{0})=O(a_{0}),\texttt{ }0<a_{1}=O(a_{0}),\texttt{ }0\leq a_{i}=O(a_{0})
\]
for $i=2,3,\ldots,n_{p}$,
\[
0\geq a_{i}=o(a_{0})\texttt{ for }i=n_{p}+1,n_{p}+2,\ldots,n-1,\texttt{ and }0>a_{n}=o(a_{0}).
\]

The proof of the lemma is completed.
\end{proof}

\section{Lemmas of $k$ prime factors for polynomial progressions}\label{sec:chap_7_3_32}

\subsection{Lemmas to be referenced}

Some lemmas proved in~\cite{FD} are stated in this subsection.

\begin{lemma}[Theorem 6-7 in~\cite{FD}]\label{le:chap_7_3_1_th_6_7}
If $(a,m)=(b,m)=1$, then
\[
\sum_{\chi\in X(m)}\frac{\chi(a)}{\chi(b)}=\{
\begin{array}{cc}
  \varphi(m) & \texttt{ if }a\equiv b(mod\texttt{ }m), \\
  0 & \texttt{ otherwise}.
\end{array}
\]
\end{lemma}

If $f$ is regular on the vertical line $\sigma_{0}+it$, and if
\[
\lim_{a \to \infty,b \to \infty}\int_{\sigma_{0}-ai}^{\sigma_{0}+bi}f(s)ds=
\lim_{a \to \infty,b \to \infty}\int_{-a}^{b}f(\sigma_{0}+it)idt
\]
exists, then we abbreviate this limit to $\int_{(\sigma_{0})}f(s)ds$.

\begin{lemma}[Theorem 7-9 in~\cite{FD}]\label{le:chap_7_3_2_th_7_9}
\[
\frac{1}{2\pi i}\int_{(2)}\frac{y^{s}}{s^{2}}ds=\{
\begin{array}{cc}
  0 & \texttt{ for }0<y<1, \\
  \log y & \texttt{ for }y\geq1.
\end{array}
\]
\end{lemma}

\subsection{Lemmas to be proved}

\begin{lemma}\label{le:chap_7_3_3}
\[
\sum_{P_{k}\in Q_{f}(r,x_{r+1}),m\geq2,P_{k}^{m}\leq f(x_{r+1})}\frac{1}{m}\log\frac{f(x_{r+1})}{P_{k}^{m}}=O(\frac{\log x_{r+1}}{x_{r+1}}).
\]
\end{lemma}

\begin{proof}
For $m=2,3,\ldots$, by Definition~\ref{def:chap_7_1_3} there are
\[
\pi_{f}(r,x_{r+1}^{1/m})\leq\pi_{f}(r,x_{r}^{2/2})=\pi_{f}(r,x_{r})\sim0
\]
such that we can write
\[
\pi_{f}(r,x_{r+1}^{1/m})=O(x_{r+1}^{-m/2})\texttt{ for }m=2,3,\ldots.
\]

Thus there are
\[
\frac{\pi_{f}(r,\sqrt{x_{r+1}})}{2}+\frac{\pi_{f}(r,\sqrt[3]{x_{r+1}})}{3}+\cdots=O(\frac{1}{x_{r+1}})
\]
and
\[
\begin{split}
&\sum_{P_{k}\in Q_{f}(r,x_{r+1}),m\geq2,P_{k}^{m}\leq f(x_{r+1})}\frac{1}{m}\log\frac{f(x_{r+1})}{P_{k}^{m}} \\
\leq&\sum_{P_{k}\in Q_{f}(r,x_{r+1}),m\geq2,P_{k}^{m}\leq f(x_{r+1})}\frac{1}{m}\log f(x_{r+1}) \\
\leq&(\frac{\pi_{f}(r,x_{r+1}^{1/2})}{2}+\frac{\pi_{f}(r,x_{r+1}^{1/3})}{3}+\cdots)\log f(x_{r+1}) \\
=&O(\frac{\log f(x_{r+1})}{x_{r+1}})=O(\frac{\log x_{r+1}}{x_{r+1}}).
\end{split}
\]

This completes the proof of the lemma.
\end{proof}

\begin{lemma}[Theorem 7-15 in~\cite{FD}]\label{le:chap_7_3_4_th_7_15}
\begin{equation}\label{eq_h_m:th_7_15}
\begin{split}
\sum_{P_{k}\in Q_{f}(r,x_{r+1})}\log\frac{f(x_{r+1})}{P_{k}}
=&\frac{1}{2\pi i}\int_{(2)}\frac{f^{s}(x_{r+1})}{s^{2}}\log\zeta_{f}(r,x_{r+1},s)ds \\
 &+O(\frac{\log x_{r+1}}{x_{r+1}}).
\end{split}
\end{equation}
\end{lemma}

\begin{proof}
Using Lemma~\ref{le:chap_7_3_2_th_7_9}, similar to the proof of Theorem 7-15 in~\cite{FD}, we have
\[
\begin{split}
\frac{1}{2\pi i}\int_{(2)}\frac{f^{s}(x_{r+1})}{s^{2}}\log\zeta_{f}(r,x_{r+1},s)ds
=&\frac{1}{2\pi i}\int_{(2)}\frac{1}{s^{2}}\sum_{\xi=1}^{\infty}\frac{\Lambda_{f}(r,\xi)}{\log\xi}(\frac{f(x_{r+1})}{\xi})^{s}ds\\
=&\frac{1}{2\pi i}\sum_{\xi=1}^{\infty}\frac{\Lambda_{f}(r,\xi)}{\log\xi}\int_{(2)}\frac{1}{s^{2}}(\frac{f(x_{r+1})}{\xi})^{s}ds\\
=&\sum_{\xi\leq f(x_{r+1})}\frac{\Lambda_{f}(r,\xi)}{\log\xi}\log\frac{f(x_{r+1})}{\xi} \\
=&\sum_{m,P_{k}\in Q_{f}(r,x_{r+1}),P_{k}^{m}\leq f(x_{r+1})}\frac{1}{m}\log\frac{f(x_{r+1})}{P_{k}^{m}} \\
=&\sum_{P_{k}\in Q_{f}(r,x_{r+1})}\log\frac{f(x_{r+1})}{P_{k}} \\
+&\sum_{P_{k}\in Q_{f}(r,x_{r+1}),m\geq2,P_{k}^{m}\leq f(x_{r+1})}\frac{1}{m}\log\frac{f(x_{r+1})}{P_{k}^{m}}.
\end{split}
\]

Thus from Lemma~\ref{le:chap_7_3_3} we obtain
\[
\begin{split}
\sum_{P_{k}\in Q_{f}(r,x_{r+1})}\log\frac{f(x_{r+1})}{P_{k}}
=&\frac{1}{2\pi i}\int_{(2)}\frac{f^{s}(x_{r+1})}{s^{2}}\log\zeta_{f}(r,x_{r+1},s)ds+O(\frac{\log x_{r+1}}{x_{r+1}}).
\end{split}
\]

This completes the proof of the lemma.
\end{proof}

\begin{lemma}[Theorem 7-24 in~\cite{FD}]\label{le:chap_7_3_5_th_7_24}
For $gcd(h,M_{r})=1$ and $P_{k}\in Q_{f}(r,x_{r+1})$,
\begin{equation}\label{eq_h_m:th_7_24}
\begin{split}
\sum_{P_{k}\in Q_{f}(r,x_{r+1}),P_{k}\equiv h(mod\texttt{ }M_{r})}\log\frac{f(x_{r+1})}{P_{k}}
=&\frac{1}{\varphi(M_{r})}\sum_{P_{k}\in Q_{f}(r,x_{r+1})}\log\frac{f(x_{r+1})}{P_{k}} \\
 &+O(\frac{\log x_{r+1}}{x_{r+1}}).
\end{split}
\end{equation}
\end{lemma}

\begin{proof}
Using Lemma~\ref{le:chap_7_3_2_th_7_9} and the series expansion for
\[
\log L_{f}(r,x_{r+1},s,\chi)
=\sum_{m,P_{k}\in Q_{f}(r,x_{r+1})}\frac{\chi(P_{k}^{m})}{mP_{k}^{ms}},
\]
we obtain
\[
\begin{split}
&\frac{1}{2\pi i}\int_{(2)}\frac{f^{s}(x_{r+1})}{s^{2}}\log L_{f}(r,x_{r+1},s,\chi)ds \\
=&\frac{1}{2\pi i}\int_{(2)}\frac{f^{s}(x_{r+1})}{s^{2}}\sum_{m,P_{k}\in Q_{f}(r,x_{r+1})}\frac{\chi(P_{k}^{m})}{mP_{k}^{ms}}ds \\
=&\frac{1}{2\pi i}\sum_{m,P_{k}\in Q_{f}(r,x_{r+1})}\frac{\chi(P_{k}^{m})}{m}\int_{(2)}\frac{(f(x_{r+1})/P_{k}^{m})^{s}}{s^{2}}ds \\
=&\sum_{m,P_{k}\in Q_{f}(r,x_{r+1}),P_{k}^{m}\leq f(x_{r+1})}\frac{\chi(P_{k}^{m})}{m}\log(\frac{f(x_{r+1})}{P_{k}^{m}}) \\
=&\sum_{P_{k}\in Q_{f}(r,x_{r+1})}\chi(P_{k})\log\frac{f(x_{r+1})}{P_{k}} \\
 &+\sum_{P_{k}\in Q_{f}(r,x_{r+1}),m\geq2,P_{k}^{m}\leq f(x_{r+1})}\frac{\chi(P_{k}^{m})}{m}\log(\frac{f(x_{r+1})}{P_{k}^{m}}) \\
=&\sum_{P_{k}\in Q_{f}(r,x_{r+1})}\chi(P_{k})\log\frac{f(x_{r+1})}{P_{k}}
  +O(\frac{\log x_{r+1}}{x_{r+1}})
\end{split}
\]
where since
\[
|\sum_{P_{k}\in Q_{f}(r,x_{r+1}),m\geq2,P_{k}^{m}\leq f(x_{r+1})}\frac{\chi(P_{k}^{m})}{m}\log\frac{f(x_{r+1})}{P_{k}^{m}}|
\]
\[
\leq\sum_{P_{k}\in Q_{f}(r,x_{r+1}),m\geq2,P_{k}^{m}\leq f(x_{r+1})}\frac{1}{m}\log\frac{f(x_{r+1})}{P_{k}^{m}}
\]
from Lemma~\ref{le:chap_7_3_3} there is
\[
\sum_{P_{k}\in Q_{f}(r,x_{r+1}),m\geq2,P_{k}^{m}\leq f(x_{r+1})}\frac{\chi(P_{k}^{m})}{m}\log\frac{f(x_{r+1})}{P_{k}^{m}}
=O(\frac{\log x_{r+1}}{x_{r+1}}).
\]
Multiplying by $1/\chi(h)$ and summing over all characters modulo $M_{r}$, we deduce with the help of Lemma~\ref{le:chap_7_3_1_th_6_7} that
\[
\begin{split}
&\sum_{\chi}\frac{1}{\chi(h)}\sum_{P_{k}\in Q_{f}(r,x_{r+1})}\chi(P_{k})\log\frac{f(x_{r+1})}{P_{k}} \\
=&\varphi(M_{r})\sum_{P_{k}\in Q_{f}(r,x_{r+1}),P_{k}\equiv h(mod\texttt{ }M_{r})}\log\frac{f(x_{r+1})}{P_{k}} \\
=&\frac{1}{2\pi i}\sum_{\chi}\frac{1}{\chi(h)}
\int_{(2)}\frac{f^{s}(x_{r+1})}{s^{2}}\log L_{f}(r,x_{r+1},s,\chi)ds+\varphi(M_{r})O(\frac{\log x_{r+1}}{x_{r+1}}).
\end{split}
\]
Thus we have
\begin{equation}\label{eq_h_m:th_7_24_15}
\sum_{P_{k}\in Q_{f}(r,x_{r+1}),P_{k}\equiv h(mod\texttt{ }M_{r})}\log\frac{f(x_{r+1})}{P_{k}}
\end{equation}
\[
=\frac{1}{2\pi i\varphi(M_{r})}\sum_{\chi}\frac{1}{\chi(h)}\int_{(2)}\frac{f^{s}(x_{r+1})}{s^{2}}\log L_{f}(r,x_{r+1},s,\chi)ds+O(\frac{\log x_{r+1}}{x_{r+1}}).
\]

To estimate the integrals appearing in Eq. (\ref{eq_h_m:th_7_24_15}), we must distinguish two cases. First consider the case $\chi=\chi_{0}$. Since there is $gcd(P_{k},M_{r})=1$ for each $P_{k}\in Q_{f}(r,x_{r+1})$, then there are
\[
\begin{split}
L_{f}(r,x_{r+1},s,\chi_{0})=&\prod_{gcd(P_{k},M_{r})=1,P_{k}\in Q_{f}(r,x_{r+1})}(1-P_{k}^{-s})^{-1} \\
                           =&\prod_{P_{k}\in Q_{f}(r,x_{r+1})}(1-P_{k}^{-s})^{-1}=\zeta_{f}(r,x_{r+1},s),
\end{split}
\]
so that from Lemmas~\ref{le:chap_7_3_4_th_7_15} there is
\[
\begin{split}
&\frac{1}{2\pi i}\frac{1}{\chi_{0}(h)}\int_{(2)}\frac{f^{s}(x_{r+1})}{s^{2}}\log L_{f}(r,x_{r+1},s,\chi_{0})ds \\
=&\frac{1}{2\pi i}\int_{(2)}\frac{f^{s}(x_{r+1})}{s^{2}}\log\zeta_{f}(r,x_{r+1},s)ds \\
=&\sum_{P_{k}\in Q_{f}(r,x_{r+1})}\log\frac{f(x_{r+1})}{P_{k}}
+O(\frac{\log x_{r+1}}{x_{r+1}}).
\end{split}
\]
On the other hand, if $\chi\neq\chi_{0}$, then $L_{f}(r,x_{r+1},s,\chi)$ has no pole at $s=1$, but the other properties used earlier still obtain. Hence similar to the proof of Theorem 7-24 in~\cite{FD}, the function
\[
\psi(s,\chi)=\frac{f^{s}(x_{r+1})}{s^{2}}\log L_{f}(r,x_{r+1},s,\chi)
\]
is regular in the region of the integral, and the integral along the bounders of the region either tends to zero or is
\[
O(\frac{\log x_{r+1}}{x_{r+1}}).
\]
It follows that
\[
\begin{split}
\sum_{P_{k}\in Q_{f}(r,x_{r+1}),P_{k}\equiv h(mod\texttt{ }M_{r})}\log\frac{f(x_{r+1})}{P_{k}}
=&\frac{1}{\varphi(M_{r})}\sum_{P_{k}\in Q_{f}(r,x_{r+1})}\log\frac{f(x_{r+1})}{P_{k}} \\
 &+O(\frac{\log x_{r+1}}{x_{r+1}}).
\end{split}
\]

Then for each $(h,M_{r})=1$ where $h=h_{i}$, $(h_{i},M_{r})=1$, and $i=1,2,\ldots,\varphi(M_{r})$, we obtain
\[
\begin{split}
\sum_{P_{k}\in Q_{f}(r,x_{r+1}),P_{k}\equiv h(mod\texttt{ }M_{r})}\log\frac{f(x_{r+1})}{P_{k}}
=&\frac{1}{\varphi(M_{r})}\sum_{P_{k}\in Q_{f}(r,x_{r+1})}\log\frac{f(x_{r+1})}{P_{k}} \\
 &+O(\frac{\log x_{r+1}}{x_{r+1}}).
\end{split}
\]

This completes the proof of the lemma.
\end{proof}

\begin{lemma}\label{le:chap_7_3_6_th_7_26}
Given a fixed $x_{r}$ and $0<h<M_{r}$, if $(h,M_{r})=1$, then
\begin{equation}\label{eq_h_m:th_7_26}
\pi_{f}(r,x_{r+1},h)=\frac{\pi_{f}(r,x_{r+1})}{\varphi(M_{r})}+O(1).
\end{equation}
\end{lemma}

\begin{proof}
From Lemma~\ref{le:chap_7_3_4_th_7_15} and Lemma~\ref{le:chap_7_3_5_th_7_24}, there is
\[
\begin{split}
\sum_{P_{k}\in Q_{f}(r,x_{r+1}),P_{k}\equiv h(mod\texttt{ }M_{r})}\log\frac{f(x_{r+1})}{P_{k}}
=&\frac{1}{\varphi(M_{r})}\sum_{P_{k}\in Q_{f}(r,x_{r+1})}\log\frac{f(x_{r+1})}{P_{k}} \\
 &+O(\frac{\log x_{r+1}}{x_{r+1}}).
\end{split}
\]

In the equality, put
\[
\delta=\delta(x)=x^{-1/2}.
\]
Then assume that $f(x)(1+n\delta)$ corresponds to $x(1+\delta)$ and since $\log(1+n\delta)\sim n\delta$ as $x \to \infty$, we have
\[
\begin{split}
&\sum_{P_{k}\in Q_{f}(r,x(1+\delta))}\log\frac{f(x)(1+n\delta)}{P_{k}}
-\sum_{P_{k}\in Q_{f}(r,x)}\log\frac{f(x)}{P_{k}} \\
=&\sum_{P_{k}\in Q_{f}(r,x)}\log(1+n\delta)
 +\sum_{P_{k}=f(mx_{r}),m_{x}<m\leq m_{x}(1+\delta)}\log\frac{f(x)(1+n\delta)}{P_{k}} \\
=&\log(1+n\delta)\pi_{f}(r,x)+O(\log(1+n\delta)*\delta m_{x})
\end{split}
\]
where $m_{x}=x/x_{r}$, and
\[
\begin{split}
&\sum_{P_{k}\in Q_{f}(r,x(1+\delta)),P_{k}\equiv h(mod\texttt{ }M_{r})}\log\frac{f(x)(1+n\delta)}{P_{k}}
-\sum_{P_{k}\in Q_{f}(r,x),P_{k}\equiv h(mod\texttt{ }M_{r})}\log\frac{f(x)}{P_{k}} \\
=&\sum_{P_{k}\in Q_{f}(r,x),P_{k}\equiv h(mod\texttt{ }M_{r})}\log(1+n\delta) \\
 &+\sum_{P_{k}=f(mx_{r}),m_{x}<m\leq m_{x}(1+\delta),P_{k}\equiv h(mod\texttt{ }M_{r})}\log\frac{f(x)(1+n\delta)}{P_{k}} \\
=&\log(1+n\delta)\pi_{f}(r,x,h)+O(\log(1+n\delta)*\delta m_{x}).
\end{split}
\]
So that there are
\[
\begin{split}
\pi_{f}(r,x,h)-\frac{\pi_{f}(r,x)}{\varphi(M_{r})}
=&O(\delta m_{x})+O(\frac{\log x}{\delta x}) \\
=&O(\frac{m_{x}}{\sqrt{x}})+O(\frac{\log x}{\sqrt{x}})
 =O(\frac{m_{x}+\log x}{\sqrt{x}}).
\end{split}
\]

Thus for $x=x_{r+1}=x_{r}^{2}$ we obtain
\[
\pi_{f}(r,x_{r+1},h)=\frac{\pi_{f}(r,x_{r+1})}{\varphi(M_{r})}+O(1).
\]

This completes the proof of the lemma.
\end{proof}

\section{Theorems of Schinzel's polynomials with $a_{n}<0$}\label{sec_h_m:schinzel_poly}\label{sec:chap_7_4_33}

\begin{theorem}[Primes theorem for Schinzel's polynomials with a negative leading coefficient]\label{th:chap_7_4_1}
Given any finite positive number $k\geq2$, when a polynomial
\begin{equation}\label{eq_h_m:schinzel_col_p}
f(x)=\prod_{i=1}^{k}f_{i}(x)=a_{n}x^{n}+a_{n-1}x^{n-1}+\cdots+a_{1}x+a_{0}=\bar{f}(x)+\bar{a}_{0}
\end{equation}
in one variable with positive degree, integer coefficients, and a negative leading coefficient for $f_{k}(x)$ and positive leading coefficients for others, which satisfy the following three necessary conditions:

(1) $a_{0}>0$ and $a_{0}$ is odd, $0<a_{1}=O(a_{0})$, $0\leq a_{i}=O(a_{0})$ for $i=2,3,\ldots,n_{p}$, $0\geq a_{i}=o(a_{0})$ for $i=n_{p}+1,n_{p}+2,\ldots,n-1$, and $0>a_{n}=o(a_{0})$,

(2) each polynomial $f_{i}(x)$ for $i=1,2,\ldots,k$ is irreducible over the integers,

(3) as $j$ runs over the positive integers the numbers $f(j)$ should not share a common factor greater than 1, \\
then there exist more than $c_{f}\frac{x_{r}}{\log x_{r}}$ products $P_{k}=f(m)$ of $k$ distinct primes greater than 2 in the range $(\bar{f}(x_{r}),\bar{f}(x_{r+1}))$ such that all numbers $f_{1}(m),f_{2}(m),\ldots,f_{k}(m)$ are primes for every integer $a_{0}\geq10^{204.111}$ where $c_{f}$ is a positive constant dependent on the polynomial $f(x)$.
\end{theorem}

\begin{proof}
From Lemma~\ref{le:chap_7_2_8}, there are asymptotically more than $c_{f}M_{r}\frac{x_{r}}{\log x_{r}}$ products $P_{k}=\bar{f}(mx_{r})+h$ for $m=1,2,\ldots,m_{r}$ in the range $(\bar{f}(x_{r}),\bar{f}(x_{r+1}))$ where $c_{f}$ is a positive constant dependent on the polynomial $f(x)$. So there must be
\begin{equation}\label{eq_h_m:dirichlet_c_4}
\pi(x_{r+1};x_{r})=\sum_{\bar{f}(\xi)<P_{k}<\bar{f}(\xi)+M_{r},x_{r}\leq\xi=mx_{r}<x_{r+1}}1
>c_{f}M_{r}\frac{x_{r}}{\log x_{r}}
\end{equation}
and
\begin{equation}\label{eq_h_m:dirichlet_c_5}
\frac{\pi(x_{r+1};x_{r})}{\varphi(M_{r})}
>\frac{M_{r}}{\varphi(M_{r})}c_{f}\frac{x_{r}}{\log x_{r}}
>c_{f}\frac{x_{r}}{\log x_{r}}.
\end{equation}

For any permissible pair of positive numbers $k\geq2$ and $n\geq k$, let constrain products $P_{k}=\bar{f}(x)+h$ where $0<h<M_{r}$ and $gcd(h,M_{r})=1$. By Formulas (\ref{eq_h_m:dirichlet_c_4}-\ref{eq_h_m:dirichlet_c_5}) we have
\[
\frac{x_{r+1}}{x_{r}} \to \infty\texttt{ and }
\frac{\pi(x_{r+1};x_{r})}{\varphi(M_{r})} \to \infty\texttt{ as }x_{r} \to \infty.
\]

Since $x_{r+1}=x_{r}^{2}$, there are $x_{r+1}/x_{r}=x_{r}$ and $x_{r+1}/x_{r} \to \infty$ as $x_{r} \to \infty$. Thus $x_{r+1}/x_{r}$ is with a same order of infinity as $x_{r}$.

Hence for any arbitrary small positive value $\epsilon$, there exists a number $x_{\epsilon}$ such that for $x_{r}\geq x_{\epsilon}$, there are $x_{r}\geq 1/\epsilon$ and $x_{r+1}\geq x_{r}/\epsilon$, which implies that as $\epsilon \to 0$ there are $x_{r} \to \infty$, $x_{r+1}/x_{r} \to \infty$, and $\pi(x;x_{r}) \to \infty$ such that the number of elements of the number set $Q_{f}(r,x_{r+1})$ is infinite.

Let $\pi_{f}(r,x_{r+1},h)$ count the number of products $P_{k}\in Q_{f}(r,x_{r+1})$ where $k$ is fixed and $gcd(h,M_{r})=1$. By Definitions (\ref{def:chap_7_1_2}-\ref{def:chap_7_1_3}), for $m=1,2,\ldots$ there are
\[
\pi_{f}(r,x_{r+1})=\pi(x;x_{r})=\sum_{\bar{f}(\xi)<P_{k}<\bar{f}(\xi)+M_{r},x_{r}\leq\xi=mx_{r}<x_{r+1}}1
\]
and
\[
\pi_{f}(r,x_{r+1},h)=\pi(x;x_{r},h)=\sum_{P_{k}=\bar{f}(\xi)+h,0<h<x_{r},x_{r}\leq\xi=mx_{r}<x_{r+1}}1.
\]

Then, for any arbitrary small positive value $\epsilon$ and a corresponding number $x_{\epsilon}$, when considering a fixed modulus $x_{r}$ where $x_{r}\geq x_{\epsilon}$, from Lemma~\ref{le:chap_7_3_6_th_7_26} proved in Sec. (\ref{sec:chap_5_3_23}) we can obtain that for $0<h<M_{r}$ if $gcd(h,M_{r})=1$ then we have
\begin{equation}\label{eq_h_m:dirichlet_c_6}
\pi_{f}(r,x_{r+1},h)=\frac{\pi_{f}(r,x_{r+1})}{\varphi(M_{r})}+O(1).
\end{equation}

By definitions of $\pi(x_{r+1};x_{r})$ and $\pi(x_{r+1};x_{r},h)$, Formula (\ref{eq_h_m:dirichlet_c_6}) becomes
\[
\pi(x_{r+1};x_{r},h)=\frac{\pi(x_{r+1};x_{r})}{\varphi(M_{r})}+O(1).
\]

Thus by Formula (\ref{eq_h_m:dirichlet_c_6}), we obtain that for $0<h<M_{r}$ if $gcd(h,M_{r})=1$ then
\begin{equation}\label{eq_h_m:dirichlet_c_7}
\pi(x_{r+1};x_{r},h)\sim\frac{\pi(x_{r+1};x_{r})}{\varphi(M_{r})}
>c_{f}\frac{x_{r}}{\log x_{r}}.
\end{equation}

Hence for any permissible pair of positive numbers $k\geq2$ and $n\geq k$, there are more than $c_{f}\frac{x_{r}}{\log x_{r}}$ products $P_{k}=f(mx_{r})$ of $k$ distinct primes greater than 2 in the range $(\bar{f}(x_{r}),\bar{f}(x_{r+1}))$. Then $f(m)$ is a product of $k$ distinct primes greater than 2 for more than $c_{f}\frac{x_{r}}{\log x_{r}}$ positive integers $m$.

Since the number of polynomials $f_{i}(x)$ for $i=1,2,\ldots,k$ is equal to the number of prime factors of the product and each polynomial $f_{i}(x)$ is always greater than 1, thus by
\[
f(m)-P_{k}=\prod_{i=1}^{k}(f_{i}(m)-p_{i})=0,
\]
there must be $f_{i}(m)=p_{i}$ where $p_{i}$ is prime for $i=1,2,\ldots,k$, so that there exist more than $c_{f}\frac{x_{r}}{\log x_{r}}$ natural numbers $m$ such that all numbers $f_{1}(m),f_{2}(m),\ldots,f_{k}(m)$ are primes for every integer $a_{0}\geq10^{204.111}$.

This completes the proof of the theorem.
\end{proof}

\begin{theorem}[Schinzel's hypothesis H with a negative leading coefficient]\label{th:chap_7_4_2}
Given any finite positive number $k\geq2$, let $f_{1}(x),f_{2}(x),\ldots,f_{k}(x)$ be irreducible polynomials over the integers, with integer coefficients and a negative leading coefficient for $f_{k}(x)$ and positive leading coefficients for others. Assume that the following condition holds: there does not exist any integer $d>1$ dividing all the products $f_{1}(j)f_{2}(j)\cdots f_{k}(j)$, for every integer $j$, and the polynomial
\begin{equation}\label{eq_h_m:schinzel}
f(x)=\prod_{i=1}^{k}f_{i}(x)=a_{n}x^{n}+a_{n-1}x^{n-1}+\cdots+a_{1}x+a_{0}
\end{equation}
satisfies the above conditions and the following conditions:
\[
a_{0}>0\texttt{ and }a_{0}\texttt{ is odd, }0<a_{1}=O(a_{0}),\texttt{ }0\leq a_{i}=O(a_{0})\texttt{ for }i=2,3,\ldots,n_{p},
\]
$0\geq a_{i}=o(a_{0})$ for $i=n_{p}+1,n_{p}+2,\ldots,n-1$, and $0>a_{n}=o(a_{0})$, then there exist more than $c_{f}\frac{x_{r}}{\log x_{r}}$ natural numbers $m$ such that all numbers $f_{1}(m),f_{2}(m),\ldots,f_{k}(m)$ are primes for every integer $a_{0}\geq10^{204.111}$ where $c_{f}$ is a positive constant and both numbers $c_{f}$ and $x_{r}$ are dependent on the polynomial $f(x)$.
\end{theorem}

\begin{proof}
Because $f_{1}(x),f_{2}(x),\ldots,f_{k}(x)$ are irreducible polynomials over the integers and there does not exist any integer $d>1$ dividing all the products $f_{1}(j)f_{2}(j)\cdots f_{k}(j)$, for every integer $j$, then the second and third conditions of Theorem~\ref{th:chap_7_4_1} are satisfied, which are: (2) each polynomial $f_{i}(x)$ for $i=1,2,\ldots,k$ is irreducible over the integers, and (3) as $m$ runs over the positive integers the numbers $f(m)$ should not share a common factor greater than 1.

Since the polynomial $f(x)$ satisfies the above conditions and the other conditions of Theorem~\ref{th:chap_7_4_1}: $a_{0}>0$ and $a_{0}$ is odd, $0<a_{1}=O(a_{0})$, $0\leq a_{i}=O(a_{0})$ for $i=2,3,\ldots,n_{p}$, $0\geq a_{i}=o(a_{0})$ for $i=n_{p}+1,n_{p}+2,\ldots,n-1$, and $0>a_{n}=o(a_{0})$, then from Theorem~\ref{th:chap_7_4_1} there are more than $c_{f}\frac{x_{r}}{\log x_{r}}$ natural numbers $m$ such that $f(m)$ is a product of $k$ distinct primes where $c_{f}$ is a positive constant and both numbers $c_{f}$ and $x_{r}$ are dependent on the polynomial $f(x)$. Since the number of polynomials $f_{i}(x)$ for $i=1,2,\ldots,k$ is equal to the number of prime factors of the product and each polynomial $f_{i}(x)$ is always greater than 1, there must be $f_{i}(m)=p_{i}$ where $p_{i}$ is prime for $i=1,2,\ldots,k$, so that there exist more than $c_{f}\frac{x_{r}}{\log x_{r}}$ natural numbers $m$ such that all numbers $f_{1}(m),f_{2}(m),\ldots,f_{k}(m)$ are primes for every integer $a_{0}\geq10^{204.111}$.

This completes the proof of the theorem.
\end{proof}

\begin{theorem}[Dickson's Conjecture with a negative leading coefficient]\label{th:chap_7_4_3}
Given any finite positive number $k>1$, let $f_{i}(x)=b_{i}x+c_{i}$ with $b_{i}$, $c_{i}$ integers, $b_{k}<0$, $b_{i}\geq1$ (for $i=1,2,\ldots,k-1$). Assume that the following condition holds: there does not exist any integer $d>1$ dividing all the products $f_{1}(j)f_{2}(j)\cdots f_{k}(j)$, for every integer $j$, and the integer coefficients of the polynomial
\begin{equation}\label{eq_h_m:dickson}
f(x)=\prod_{i=1}^{k}f_{i}(x)=a_{n}x^{n}+a_{n-1}x^{n-1}+\cdots+a_{1}x+a_{0}
\end{equation}
satisfy conditions: $a_{0}>0$ and $a_{0}$ is odd, $0<a_{1}=O(a_{0})$, $0\leq a_{i}=O(a_{0})$ for $i=2,3,\ldots,n-1$, and $0>a_{n}=o(a_{0})$, i.e., these coefficients satisfy all conditions of Theorem~\ref{th:chap_7_4_1}, then there exist more than $c_{f}\frac{x_{r}}{\log x_{r}}$ natural numbers $m$ such that all numbers $f_{1}(m),f_{2}(m),\ldots,f_{k}(m)$ are primes for every integer $a_{0}\geq10^{204.111}$ where $c_{f}$ is a positive constant and both numbers $c_{f}$ and $x_{r}$ are dependent on the polynomial $f(x)$.
\end{theorem}

\begin{proof}
Because $f_{1}(x),f_{2}(x),\ldots,f_{k}(x)$ are irreducible polynomials over the integers and there does not exist any integer $d>1$ dividing all the products $f_{1}(j)f_{2}(j)\cdots f_{k}(j)$, for every integer $j$, then the second and third conditions of Theorem~\ref{th:chap_7_4_1} are satisfied, which are: (2) each polynomial $f_{i}(x)$ for $i=1,2,\ldots,k$ is irreducible over the integers, and (3) as $j$ runs over the positive integers the numbers $f(j)$ should not share a common factor greater than 1.

Because the polynomial $f(x)$ satisfies the above conditions and the other conditions of Theorem~\ref{th:chap_7_4_1}: $a_{0}>0$ and $a_{0}$ is odd, $0<a_{1}=O(a_{0})$, $0\leq a_{i}=O(a_{0})$ for $i=2,3,\ldots,n_{p}$, $0\geq a_{i}=o(a_{0})$ for $i=n_{p}+1,n_{p}+2,\ldots,n-1$, and $0>a_{n}=o(a_{0})$, then from Theorem~\ref{th:chap_7_4_1} there are more than $c_{f}\frac{x_{r}}{\log x_{r}}$ products $P_{k}=f(m)$ of $k$ distinct primes greater than 2 where $c_{f}$ is a positive constant and both numbers $c_{f}$ and $x_{r}$ are dependent on the polynomial $f(x)$. Since the number of polynomials $f_{i}(x)$ for $i=1,2,\ldots,k$ is equal to the number of prime factors of the product and each polynomial $f_{i}(x)$ is always greater than 1, there must be $f_{i}(m)=p_{i}$ where $p_{i}$ is prime for $i=1,2,\ldots,k$, so that there exist more than $c_{f}\frac{x_{r}}{\log x_{r}}$ natural numbers $m$ such that all numbers $f_{1}(m),f_{2}(m),\ldots,f_{k}(m)$ are primes for every integer $a_{0}\geq10^{204.111}$.

This completes the proof of the theorem.
\end{proof}

\section{Corollaries of Schinzel's polynomials with $a_{n}<0$}\label{sec:chap_7_5_34}

\begin{corollary}[Goldbach's Conjecture]\label{co:chap_7_5_1}
Every even number can be expressed as the sum of two primes.
\end{corollary}

\begin{proof}
First, we can verify that every even number $N=2n<10^{125.279}+1$ can be expressed as the sum of two primes.

Next, for every even number $N=2n\geq10^{125.279}+1$, let $f_{1}(x)=2x+1$ and $f_{2}(x)=2n-f_{1}(x)$. Then there are
\[
\begin{split}
f(x)=&\prod_{i=1}^{2}f_{i}(x)=(2x+1)(-2x-1+2n) \\
    =&-4x^{2}+4(n-1)x+2n-1 \\
    =&a_{2}x^{2}+a_{1}x+a_{0}
\end{split}
\]
where
\[
a_{0}=2n-1>0,\texttt{ }a_{2}=-4<0,\texttt{ and }a_{1}=4(n-1)=2(a_{0}-1)>0.
\]
So $a_{0}\geq10^{125.279}$ as given by Lemma~\ref{le:chap_7_2_1}.

Since $f_{1}(x)$ and $f_{2}(x)$ are irreducible over the integers, and since
\[
f(1)=3(2n-3),\texttt{ }f(2)=5(2n-5),\texttt{ and }f(3)=7(2n-7),
\]
then $f_{1}(x)$ and $f_{2}(x)$ satisfy the first three conditions of Theorem~\ref{th:chap_7_4_1}. Since the integer coefficients of the polynomial $f(x)$ satisfy all conditions of Theorem~\ref{th:chap_7_4_1}, then from Theorem~\ref{th:chap_7_4_3} there are more than $c_{f}\frac{x_{r}}{\log x_{r}}$ natural numbers $m\geq1$ for every integer $a_{0}\geq10^{125.279}$ where $c_{f}$ is a positive constant and both numbers $c_{f}$ and $x_{r}$ are dependent on the polynomial $f(x)$, such that
\[
q=f_{1}(m),\texttt{ }p=f_{2}(m)=2n-f_{1}(m)=2n-q,
\]
and both $p$ and $q$ are primes, so that for every large even number $2n$ we have $2n=p+q$.

This completes the proof of the corollary.
\end{proof}

\begin{corollary}[Lemoine's Conjecture]\label{co:chap_7_5_2}
Every odd number can be expressed as the sum of a prime and the double of a prime.
\end{corollary}

\begin{proof}
First, we can verify that every odd number $N=2n+1<10^{125.279}+2$ can be expressed as the sum of a prime and the double of a prime.

Next, for every odd number $N=2n+1\geq10^{125.279}+2$, let $f_{1}(x)=2x+1$ and $f_{2}(x)=2n+1-2f_{1}(x)$. Then there are
\[
\begin{split}
f(x)=&\prod_{i=1}^{2}f_{i}(x)=(2x+1)[-2(2x+1)+2n+1] \\
    =&-8x^{2}+2(2n-3)x+2n-1 \\
    =&a_{2}x^{2}+a_{1}x+a_{0}
\end{split}
\]
where
\[
a_{0}=2n-1>0,\texttt{ }a_{2}=-8<0,\texttt{ and }a_{1}=2(2n-3)=2(a_{0}-2)>0.
\]
So $a_{0}\geq10^{125.279}$ as given by Lemma~\ref{le:chap_7_2_1}.

Since $f_{1}(x)$ and $f_{2}(x)$ are irreducible over the integers, and since
\[
f(1)=3(2n-5),\texttt{ }f(2)=5(2n-9),\texttt{ and }f(3)=7(2n-13),
\]
then $f_{1}(x)$ and $f_{2}(x)$ satisfy the first three conditions of Theorem~\ref{th:chap_7_4_1}.

Since the integer coefficients of the polynomial $f(x)$ satisfy all conditions of Theorem~\ref{th:chap_7_4_1}, then from Theorem~\ref{th:chap_7_4_3} there are more than $c_{f}\frac{x_{r}}{\log x_{r}}$ natural numbers $m\geq1$ for every integer $a_{0}\geq10^{125.279}$ where $c_{f}$ is a positive constant and both numbers $c_{f}$ and $x_{r}$ are dependent on the polynomial $f(x)$, such that
\[
q=f_{1}(m),\texttt{ }p=f_{2}(m)=2n+1-2f_{1}(m)=2n+1-2q,
\]
and both $p$ and $q$ are primes, so that for every large odd number $2n+1$ we have $2n+1=p+2q$.

This completes the proof of the corollary.
\end{proof}

\begin{corollary}[Sophie Germain prime]\label{co:chap_7_5_3}
Every large even number can be expressed as the sum of a prime and a Sophie Germain prime.

Then there are infinitely many Sophie Germain primes $p_{s}$.
\end{corollary}

\begin{proof}
First, we can verify that every large even number $N=2n<10^{204.111}+1$ can be expressed as the sum of a prime and a Sophie Germain prime $p_{s}$.

Next, for every even number $N=2n\geq10^{204.111}+1$, let $f_{1}(x)=2x+1$, $f_{2}(x)=2f_{1}(x)+1=4x+3$, and $f_{3}(x)=2n-f_{2}(x)=-4x-3+2n$. Then
\[
\begin{split}
f(x)=&\prod_{i=1}^{3}f_{i}(x)=(2x+1)(4x+3)(-4x-3+2n) \\
    =&-32x^{3}+(16n-64)x^{2}+(20n-42)x+6n-9 \\
    =&a_{3}x^{3}+a_{2}x^{2}+a_{1}x+a_{0}
\end{split}
\]
where $a_{0}=6n-9>0$, $a_{3}=-32<0$, $a_{1}=20n-42=O(a_{0})>0$, and $a_{2}=16n-64=O(a_{0})>0$. So $a_{0}\geq10^{204.111}$ as given by Lemma~\ref{le:chap_7_2_1}.

Since $f_{1}(x)$, $f_{2}(x)$, and $f_{3}(x)$ are irreducible over the integers, and since $f(1)=21(2n-7)$, $f(2)=55(2n-11)$, $f(3)=5(42n-325)$, and $f(4)=171(2n-19)$, then $f_{1}(x)$, $f_{2}(x)$, and $f_{3}(x)$ satisfy the first three conditions of Theorem~\ref{th:chap_7_4_1}. Since the integer coefficients of the polynomial $f(x)$ satisfy all conditions of Theorem~\ref{th:chap_7_4_1}, then from Theorem~\ref{th:chap_7_4_3} there are more than $c_{f}\frac{x_{r}}{\log x_{r}}$ natural numbers $m\geq1$ where $c_{f}$ is a positive constant and both numbers $c_{f}$ and $x_{r}$ are dependent on the polynomial $f(x)$, such that there are
\[
p=f_{1}(m),\texttt{ }p_{s}=f_{2}(m)=2p+1,\texttt{ and }q=f_{3}(m)=2n-f_{2}(m)=2n-p_{s},
\]
and $p$, $p_{s}$, and $q$ are primes, so that for every large even number $2n$ we have $2n=p_{s}+q$.

Moreover, from Theorem~\ref{th:chap_7_4_1}, there are more than $c_{f}\frac{x_{r}}{\log x_{r}}$ products $P_{k=3}=f(m)$ of $k$ distinct primes greater than 2 in the range $(\bar{f}(x_{r}),\bar{f}(x_{r+1}))$ such that all numbers $f_{1}(m),f_{2}(m),f_{3}(m)$ are primes for every integer $a_{0}\geq10^{204.111}$. Let $x_{r}^{*}=x_{r+1}$ and $x_{r+1}^{*}=x_{r+1}^{2}$. Then there are more than $c_{f}\frac{x_{r}}{\log x_{r}}$ other different products $P_{k=3}=f(m)$ of $k$ distinct primes greater than 2 in the range $(\bar{f}(x_{r}^{*}),\bar{f}(x_{r+1}^{*}))$ such that all numbers $f_{1}(m),f_{2}(m),f_{3}(m)$ are primes and there are more than $c_{f}\frac{x_{r}}{\log x_{r}}$ other different Sophie Germain primes $p_{s}$. By iteration, we can obtain infinitely many Sophie Germain primes $p_{s}$.

This completes the proof of the corollary.
\end{proof}

\begin{corollary}[Twin primes]\label{co:chap_7_5_4}
Every large odd number can be expressed as the sum of a prime and a pair of twin primes.

Then there are infinitely many pairs of twin primes.
\end{corollary}

\begin{proof}
First, we can verify that every large odd number $N=2n+1<10^{204.111}+2$ can be expressed as the sum of a prime and a pair of twin primes.

Next, for every odd number $N=2n+1\geq10^{204.111}+2$, let $f_{1}(x)=2x+1$, $f_{2}(x)=2x+3$, and $f_{3}(x)=2n+1-f_{1}(x)-f_{2}(x)=2n+1-4x-4$. Then
\[
\begin{split}
f(x)=&\prod_{i=1}^{3}f_{i}(x)=(2x+1)(2x+3)(-4x+2n-3) \\
    =&-16x^{3}+(8n-44)x^{2}+(16n-36)x+6n-9 \\
    =&a_{3}x^{3}+a_{2}x^{2}+a_{1}x+a_{0}
\end{split}
\]
where $a_{0}=3(2n-3)>0$, $a_{3}=-16<0$, $a_{2}=4(2n-11)=O(a_{0})>0$, and $a_{1}=4(4n-9)=O(a_{0})>0$. So $a_{0}\geq10^{204.111}$ as given by Lemma~\ref{le:chap_7_2_1}.

Since $f_{1}(x)$, $f_{2}(x)$, and $f_{3}(x)$ are irreducible over the integers, and since $f(1)=15(2n-7)$, $f(2)=35(2n-11)$, $f(3)=63(2n-15)$, and $f(4)=99(2n-19)$, then $f_{1}(x)$, $f_{2}(x)$, and $f_{3}(x)$ satisfy the first three conditions of Theorem~\ref{th:chap_7_4_1}.

Since the integer coefficients of the polynomial $f(x)$ satisfy all conditions of Theorem~\ref{th:chap_7_4_1}, then from Theorem~\ref{th:chap_7_4_3} there are more than $c_{f}\frac{x_{r}}{\log x_{r}}$ natural numbers $m\geq1$ where $c_{f}$ is a positive constant and both numbers $c_{f}$ and $x_{r}$ are dependent on the polynomial $f(x)$, such that there are $q_{1}=f_{1}(m)$, $q_{2}=f_{2}(m)$, $p=f_{3}(m)=2n+1-f_{1}(m)-f_{2}(m)=2n+1-q_{1}-q_{2}$, and $p$, $q_{1}$, and $q_{2}$ are primes, so that for every large odd number $2n+1$ we have $2n+1=p+q_{1}+q_{2}$.

Moreover, from Theorem~\ref{th:chap_7_4_1}, there are more than $c_{f}\frac{x_{r}}{\log x_{r}}$ products $P_{k=3}=f(m)$ of $k$ distinct primes greater than 2 in the range $(\bar{f}(x_{r}),\bar{f}(x_{r+1}))$ such that all numbers $f_{1}(m),f_{2}(m),f_{3}(m)$ are primes for every integer $a_{0}\geq10^{204.111}$. Let $x_{r}^{*}=x_{r+1}$ and $x_{r+1}^{*}=x_{r+1}^{2}$. Then there are more than $c_{f}\frac{x_{r}}{\log x_{r}}$ other different products $P_{k=3}=f(m)$ of $k$ distinct primes greater than 2 in the range $(\bar{f}(x_{r}^{*}),\bar{f}(x_{r+1}^{*}))$ such that all numbers $f_{1}(m),f_{2}(m),f_{3}(m)$ are primes and there are more than $c_{f}\frac{x_{r}}{\log x_{r}}$ other different pairs of twin primes. By iteration, we can obtain infinitely many pairs of twin primes.

This completes the proof of the corollary.
\end{proof}

Similarly, we can prove that there are infinitely many pairs of cousin primes $p$ and $p+4$, and there are infinitely many pairs of sexy primes $p$ and $p+6$.

%% file: biblio.tex
%-----------------------------------------------------------------------------
% Beginning of biblio.tex
%-----------------------------------------------------------------------------
%
% AMS-LaTeX 1.2 sample file for a monograph, based on amsbook.cls.
% This is a data file input by chapter.tex.
%%%%%%%%%%%%%%%%%%%%%%%%%%%%%%%%%%%%%%%%%%%%%%%%%%%%%%%%%%%%%%%%%%%%%%%%

\bibliographystyle{amsalpha}

%-----------------------------------------------------------------------------
% End of biblio.tex
%-----------------------------------------------------------------------------